\documentclass[12pt]{article}

\usepackage{varwidth}
\usepackage{lipsum}
\usepackage{threeparttable}
\usepackage{makecell}
\usepackage{amsfonts}
\usepackage{amsmath}
\usepackage{mathrsfs}
\usepackage{fancyhdr}
\usepackage{epsfig,color,morefloats}
\usepackage{xcolor}
\usepackage[authoryear]{natbib}
\usepackage{booktabs,multirow,multicol,longtable}
\usepackage{amsthm}
\usepackage{amssymb}
\usepackage{color}
\usepackage{lscape}
\usepackage{rotating}
\usepackage{caption}
\usepackage{subfigure}
\usepackage{subcaption}
\usepackage{graphicx}
\usepackage[breaklinks,colorlinks,linkcolor=blue,citecolor=blue,urlcolor=blue]{hyperref}
\usepackage{lineno}
\usepackage{cancel}
\usepackage{multirow}
\usepackage[plain,noend]{algorithm2e}
\makeatletter
\renewcommand{\@algocf@capt@plain}{above}
\makeatother

\SetAlgoCaptionLayout{algocaptionleft}
\SetAlCapFnt{\footnotesize}
\SetAlCapNameFnt{\footnotesize}

\usepackage[breaklinks,colorlinks,linkcolor=blue,citecolor=blue,urlcolor=blue]{hyperref}

\usepackage{setspace}

\usepackage{amsmath, amssymb, amsthm, latexsym,color, natbib}
\usepackage{fancyhdr}
\usepackage{algpseudocode}
\usepackage{algpseudocode}
\usepackage{enumitem}

\usepackage{pdflscape}

\usepackage[resetlabels]{multibib}
\allowdisplaybreaks
\newcites{S}{References}

\allowdisplaybreaks

\usepackage{pdflscape}

\makeatletter
\newcommand{\fdsy@scale}{1.0}
\newcommand\fdsy@mweight@normal{Book}
\newcommand\fdsy@mweight@small{Book}
\newcommand\fdsy@bweight@normal{Medium}
\newcommand\fdsy@bweight@small{Medium}

\DeclareFontFamily{U}{FdSymbolC}{}

\DeclareFontShape{U}{FdSymbolC}{m}{n}{
	<-7.1> s * [\fdsy@scale] FdSymbolC-\fdsy@mweight@small
	<7.1-> s * [\fdsy@scale] FdSymbolC-\fdsy@mweight@normal
}{}
\DeclareFontShape{U}{FdSymbolC}{b}{n}{
	<-7.1> s * [\fdsy@scale] FdSymbolC-\fdsy@bweight@small
	<7.1-> s * [\fdsy@scale] FdSymbolC-\fdsy@bweight@normal
}{}

\DeclareSymbolFont{arrows}{U}{FdSymbolC}{m}{n}
\SetSymbolFont{arrows}{bold}{U}{FdSymbolC}{b}{n}

\DeclareMathSymbol{\upvDash}{\mathrel}{arrows}{233}

\DeclareMathSymbol{\upmodels}{\mathrel}{arrows}{237}
\makeatother

\DeclareMathSymbol{\upvDash}{\mathrel}{arrows}{233}

\DeclareMathSymbol{\upmodels}{\mathrel}{arrows}{237}
\makeatother

\theoremstyle{definition}

\newcommand{\x}{X}

\newcommand{\w}{W}

\newcommand{\A}{ A}
\newcommand{\B}{ B}
\newcommand{\C}{ C}
\newcommand{\D}{ D}

\newcommand{\bS}{S}

\newcommand{\I}{ I}

\newcommand{\y}{Y}
\newcommand{\z}{Z}

\newcommand{\br}{r}

\newcommand{\E}{E}
\newcommand{\Var}{{\rm var}}
 
\newcommand{\pr}{{\rm pr}}
\newcommand{\Cov}{{\rm cov}}

\newcommand{\bV}{V}

\newcommand{\tr}{\mathop{\text{\rm tr}}}

\newcommand{\bSigma}{\Sigma}

\newtheorem{lemma}{Lemma}
\newtheorem{theorem}{Theorem}
\newtheorem{proposition}{Proposition}
\newtheorem{assumption}{Assumption}

\newtheorem{definition}{Definition}
\newtheorem{remark}{Remark}

\usepackage{authblk}

\allowdisplaybreaks

\makeatletter
\DeclareFontFamily{U}{FdSymbolC}{}

\DeclareFontShape{U}{FdSymbolC}{m}{n}{
	<-7.1> s * [\fdsy@scale] FdSymbolC-\fdsy@mweight@small
	<7.1-> s * [\fdsy@scale] FdSymbolC-\fdsy@mweight@normal
}{}
\DeclareFontShape{U}{FdSymbolC}{b}{n}{
	<-7.1> s * [\fdsy@scale] FdSymbolC-\fdsy@bweight@small
	<7.1-> s * [\fdsy@scale] FdSymbolC-\fdsy@bweight@normal
}{}

\DeclareSymbolFont{arrows}{U}{FdSymbolC}{m}{n}
\SetSymbolFont{arrows}{bold}{U}{FdSymbolC}{b}{n}

\DeclareMathSymbol{\upvDash}{\mathrel}{arrows}{233}

\DeclareMathSymbol{\upmodels}{\mathrel}{arrows}{237}
\makeatother

\def\spacingset#1{\renewcommand{\baselinestretch}%
	{#1}\small\normalsize} \spacingset{1}

\newcommand{\blind}{1}

\def\singlespace{\def\baselinestretch{1}\@normalsize}

\begin{document}
	
	% \def\spacingset#1{\renewcommand{\baselinestretch}%
		% {#1}\small\normalsize} \spacingset{1}

	%\renewcommand{\baselinestretch}{1.0}
	
	%%%%%%%%%%%%%%%%%%%%%%%%%%%%%%%%%%%%%%%%%%%%%%%%%%%%%%%%%%%%%%%%%%%%%%%%%%%%%%
	
	\if1\blind
	{
        \spacingset{1.25}
		\title{
			An Adaptive $L_2$-type Test for High-dimensional White Noise
			}
		
		\author[a,b]{Jinyuan Chang}
		\author[a]{Jing He}
        \author[c]{Weiming Li}
        \author[a]{Chen Lin}
		\affil[a]{\it \small Joint Laboratory of Data Science and Business
			Intelligence, Institute of Statistical Interdisciplinary Research,  Southwestern University of Finance and Economics, Chengdu, China}
		\affil[b]{\it \small State Key Laboratory of Mathematical Sciences, Academy of Mathematics and Systems Science, Chinese Academy of Sciences, Beijing, China}
		\affil[c]{\it \small School of Statistics and Data Science, Shanghai University of Finance and Economics, Shanghai, China}		
        
        \setcounter{Maxaffil}{0}
		
		\renewcommand\Affilfont{\itshape\small}
		%\date{\today}
		\date{\vspace{-5ex}}
		
		\maketitle
	} \fi
	
	\if0\blind
	{
		\bigskip
		\bigskip
		\bigskip
		\begin{center}
			{\LARGE\bf An Adaptive $L_2$-type Test for High-dimensional White Noise }
		\end{center}
		\medskip
	} \fi

    \spacingset{1.3}
	\begin{abstract}
	We propose a new $L_2$-type test for white noise which allows the dimension $p$ of the time series to either (i) be a fixed constant, or (ii) diverge with the sample size $n$. The proposed test statistic exhibits an interesting phase transition, following two different regimes of behavior: $p$ is fixed, and $p\rightarrow\infty$.	Because identification of the operable regime is difficult,
if not impossible in practice, we devise a novel adaptive bootstrap method
to construct a unified testing procedure across different phases. Numerical experiments confirm the
good finite sample performance of the proposed adaptive $L_2$-type test in comparison to the existing methods in the literature. The proposed testing procedure has been implemented
in \textsf{R} package \texttt{HDTSA}.

	\end{abstract}
	
\bigskip 
	 

\spacingset{1.69}
\setlength{\abovedisplayskip}{0.2\baselineskip}
\setlength{\belowdisplayskip}{0.2\baselineskip}
\setlength{\abovedisplayshortskip}{0.2\baselineskip}
\setlength{\belowdisplayshortskip}{0.2\baselineskip}

%\newpage
\section{Introduction}

Testing for white noise is a fundamental problem in statistical inference, particularly in the context of diagnostic checking for linear time series modeling. The Box–Pierce portmanteau test is the most widely used method for univariate time series, which tests the joint significance of the first $K$ autocorrelations of a time series. When the dimension $p$ of the time series is small or moderate relative to the sample size $n$, it has also been extended for multivariate scenario. See, for example,  \cite{Hosking:1980} and \cite{Li+McLeod:1981}. Under the assumption that the observations are independent and identically distributed (i.i.d.), the null distributions of the test statistics involved in Box–Pierce portmanteau test and its multivariate extensions converge weakly to chi-square distributions.
There are also some other portmanteau tests for white noise, e.g., \cite{Durlauf:1991}, \cite{Deo:2000}, \cite{Lobato:2001}, and \cite{Shao:2011}. When $p$ is large, these omnibus tests often suffer from slow convergence to their asymptotic null distributions, which usually results in inaccurate size control.

With the increasing availability of large-scale data, high-dimensional time series are commonly observed in practice, which have a critical need to develop associated valid and powerful white noise tests. To do this, 
\cite{Chang+Zhou+Yao:2017} propose an $L_\infty$-type test based on the maximum absolute sample autocorrelations and cross-correlations, and apply the Gaussian approximation technique \citep{Chang:2024} to approximate the null distribution of the test statistic. 
\cite{Tsay:2020} consider an $L_\infty$-type test statistic based on the Spearman’s rank correlations, and derive its limiting distribution using the extreme value theory. 
 Although these two methods are theoretically valid when $p$ diverges with $n$, they are often conservative in practice when $p$ is large.
\cite{Wang+kong+Xia:2023} propose test statistics based on the average of the largest $s$ absolute values of the elements in autocorrelation matrices, and develop a bootstrap method to determine the critical values. 
Although these three tests have good power performance under sparse alternatives, they tend to have relatively low power against dense alternatives.
Meanwhile, the tests of \cite{Chang+Zhou+Yao:2017} and \cite{Wang+kong+Xia:2023} are computationally intensive when $p\gg n$. 
For dense alternatives, \cite{LZ19} introduce an $L_2$-type test to enhance power. To guarantee correct size control, it requires that:  (i) $p/n\rightarrow c$ for some constant $c\in(0,\infty)$, and (ii) the data follow the independent components model. Numerical results in Section \ref{sec:size} indicate that the $L_2$-type test of \cite{LZ19} fails to control size when these two requirements are not satisfied. Moreover, the first requirement  may be too restrictive, which cannot cover the scenarios $p \ll n$ and $p \gg n$.

To understand the size inflation phenomenon of the test proposed by \cite{LZ19}, we carefully study in this paper the asymptotic properties of $G_{n,\tau}$ defined as \eqref{eq:G_ntau} in Section \ref{sec:pre}, which is the key element used to construct the $L_2$-type test statistic of \cite{LZ19}. More specifically, we decompose $G_{n,\tau}$ into two terms $G_{n,\tau,1}$ and $G_{n,\tau,2}$ such that $G_{n,\tau}=G_{n,\tau,1}+G_{n,\tau,2}$. When $p \gg n$, or the underlying distribution of the data shifts away from the independent components model,  $G_{n,\tau,2}$ may diverge in probability which is the primary reason of size inflation. Instead, we suggest to construct the test statistic just based on  $G_{n,\tau,1}$, see \eqref{eq:ts} in Section \ref{sec:wn} for details. 
Our theoretical analysis shows that the asymptotic distribution of the newly proposed $L_2$-type test statistic exhibits an interesting phase transition, following two different regimes of behavior: $p$ is fixed, and $p\rightarrow\infty$. This makes it difficult to determine the critical value via the limiting distribution of the newly proposed test statistic, since identification of the operable regime is difficult if not impossible in practice due to the fact we only have one $p$ and one $n$ in the given data. To overcome this difficulty, we propose a unified adaptive bootstrap procedure in Section \ref{sec:wn} that yields valid inference across different phases. Our proposed test has been implemented in the \textsf{R} package \texttt{HDTSA} \citep{HDTSA}.

{\it Notation.} For any positive integer $m$, write $[m] = \{1, \ldots , m\}$. 
For any $x \in \mathbb{R}$, let $\lfloor x \rfloor$ denote the largest integer less than or equal to $x$. For a vector $a = (a_1,\ldots, a_k)^{\top} \in \mathbb{R}^k$, let $|a|_2=(\sum_{i=1}^{k}a_{i}^2)^{1/2}$ be its $L_2$-norm. For a matrix $A \in \mathbb{R}^{m \times q}$, let $\|A\|_2= \sup_{|x|_2=1} |Ax|_2$ denote its spectral norm. 
 
\section{Preliminary}\label{sec:pre}

Let $\{X_t\}_{t=1}^n$ be $n$ observations from a $p$-dimensional stationary linear process with mean zero. 
By Wold's Decomposition Theorem \citep[pp.~25]{lutkepohl:2005}, the stationary linear process $X_t$ can be formulated as a moving average representation:
\begin{align}\label{eq:VMA}
X_t=Z_{t}+\sum_{j=1}^{\infty} A_jZ_{t-j}\,, 
\end{align}
where each $A_j\in\mathbb R^{p\times p}$ represents the coefficient matrix and $\{Z_t\} _{t\in \mathbb Z}$ is a $p$-dimensional white-noise process with mean zero. Based on the representation \eqref{eq:VMA} of $X_t$, our goal is to test the null hypothesis 
\begin{align}\label{eq:WH_null}
    H_0: X_t=Z_t\,.
\end{align}

For any $\tau\geq0$, define
$
S_\tau=n^{-1}\sum_{t=1}^{n}X_tX_{t+\tau}^\top$, 
where, by convention, we write $X_{t+\tau}=X_{t+\tau-n}$ if $t+\tau>n$.
\cite{LZ19} propose the $L_2$-type test statistic $
    L_q = \sum_{\tau=1}^q \tr(S_\tau S_\tau^\top) - qn^{-1} {\tr}^2(S_0)$ 
 for the white noise hypothesis, where $q \ge 1$ is some fixed integer, and suggest to reject the null hypothesis at the  significance level $\alpha\in(0,1)$ if 
\begin{align}\label{eq:Lq_reject}
    L_q >  \frac{\sqrt{2q} }{n} \bigg\{{\tr}( S_0^2)-\frac{1}{n}{\tr}^2( S_0)\bigg\}  z_{1-\alpha}\,,
\end{align}
where $z_{1-\alpha}$ denotes the lower $(1-\alpha)$-quantile of $\mathcal{N}(0,1)$. To establish the validity of such a test, \cite{LZ19} impose two key requirements: (i) $p/n\rightarrow c$ for some constant $c\in(0,\infty)$, and (ii) $Z_t$ satisfies the independent components model defined in Definition \ref{def:ICmodel}.
\begin{definition}[Independent components  model]\label{def:ICmodel}
    Write 
    $Z_t=\Sigma^{1/2}W_t$ with $\Sigma=\Cov(Z_t)$ and $ W_t=(w_{1,t},\ldots,w_{p,t})^\top$, 
    where $\{W_t\}_{t\in\mathbb{Z}}$ is a sequence of i.i.d. $p$-dimensional random vectors with independent components $w_{1,t},\ldots,w_{p,t}$ satisfying $E(w_{i,t})=0$, $E(w_{i,t}^2)=1$ and $\max_{i\in [p]}E(w_{i,t}^4)$ is uniformly bounded away from infinity.
\end{definition}

In practice, these two requirements may be too restrictive. A natural question is whether the test \eqref{eq:Lq_reject} remains valid when these requirements are not satisfied. Simulation results in Section \ref{sec:size} show that the test \eqref{eq:Lq_reject} suffers from size inflation when $p \gg n$ and $Z_t$ does not satisfy the independent components model. To understand this size inflation phenomenon, we need to carefully study the asymptotic properties of the $L_2$-type test statistic $L_q$.  
Let
\begin{align}\label{eq:G_ntau}
    G_{n,\tau} = \frac{n}{p}\operatorname{tr}( S_\tau S_\tau^\top ) - \frac{1}{p} 
    \operatorname{tr}^2( S_0)\,, \quad \tau\in [q]\,.
\end{align}
Then 
$
    L_q = pn^{-1}\sum_{\tau=1}^q G_{n,\tau}$. 
We decompose $G_{n,\tau}$ as
$
  G_{n,\tau}  =   G_{n,\tau,1}+  G_{n,\tau,2}$, 
where
\begin{align}\label{eq:Gnt1}
  G_{n,\tau,1}
= \frac {1}{np}\sum_{t\neq s}X_t^\top X_s X_{t+\tau}^\top X_{s+\tau}\,,\quad
  G_{n,\tau,2}
=\frac 1{np}\sum_{t=1}^n|X_t|_2^2 |X_{t+\tau}|_2^2- \frac{1}{p}\bigg(\frac1n\sum_{t=1}^n|X_t|_2^2\bigg)^2\,.
\end{align}
The asymptotic behaviors of $G_{n,\tau,1}$ and $G_{n,\tau,2}$ are summarized in Proposition \ref{null-lemma} below which indicates that $G_{n,\tau, 1}$ and $G_{n,\tau,2}$ have different convergence rates. Specifically, the variance of $G_{n,\tau, 1}$ is of constant order, while the variance of $G_{n,\tau,2}$ will be  affected by both the underlying distribution of $Z_t$ and relative growth rate of $p$ and $n$. When $p \gg n$ or the underlying distribution of $Z_t$ shifts away from the independent components model,  $G_{n,\tau,2}$ may diverge in probability. As a result, the established null distribution of $L_q$ in \cite{LZ19} becomes invalid, which mainly causes the size inflation phenomenon of the test \eqref{eq:Lq_reject}. See more  discussion in Remark \ref{rek1} below.  

Write $Z_t=(z_{1,t},\ldots,z_{p,t})^\top$. The following assumptions are required for Proposition \ref{null-lemma}.

\begin{assumption}\label{as:A2} 
 $\{Z_t\}_{t\in \mathbb{Z}}$ is a sequence of i.i.d. random vectors such that $\E\{(Z_t^\top Z_s)^4\}=O(p^2)$ for $t\neq s$, and $\max_{i\in[p]}\E (z_{i,t}^{4})\leq C_1$ for some universal constant $C_1>0$.
\end{assumption}

\begin{assumption}\label{as:A4}
The largest eigenvalue of the covariance matrix $\Sigma$ is uniformly bounded away from infinity, and either of the two conditions holds: (i) if $p$ is fixed, then ${\tr}(\Sigma)>0$, and (ii) if $p \to \infty$, then 
$
\lim_{p\to\infty}  p^{-1}\tr (\Sigma^\ell) >0$ for $\ell=1, 2$.
\end{assumption}

\begin{assumption}\label{as:A3}
$\max_{i\in[p]}\E (z_{i,t}^{16})\leq C_2$ for some universal constant $C_2>0$, $\E[\{|Z_t|_2^2-\tr(\Sigma)\}^8]=O\{(\E[\{|Z_t|_2^2-\tr(\bSigma)\}^2])^4\}$ and $
 |\E (Z_t|Z_t|_2^2) |_2^2=O(\E[\{|Z_t|_2^2-\tr(\Sigma)\}^2])$.
\end{assumption}

Assumption \ref{as:A2} is less restrictive than the independent components model  required in \cite{LZ19}, which is satisfied automatically under the independent components model. Assumption \ref{as:A4} is also required in  \cite{LZ19}, which is a common assumption in the high-dimensional data analysis literature. Assumption \ref{as:A3} is used only for establishing the asymptotic normality of $G_{n,\tau,2}$ and the joint convergence of $G_{n,\tau,1}$ and $G_{n,\tau,2}$. If we just focus on the asymptotic behavior of $G_{n,\tau,1}$, Assumption \ref{as:A3} is not necessary. Under the independent components model, Assumption \ref{as:A3} can be weakened via  truncation techniques.

\begin{proposition}\label{null-lemma}
Let Assumptions {\rm\ref{as:A2}--\ref{as:A3}} hold with $\min\{p, n\}\to\infty$. Write $\sigma_{n1}^2=2p^{-2}{\tr}^2(\Sigma^2)$ and $
\sigma_{n2}^2=p^{-2}n^{-3}E[\{|Z_t|_2^2-\tr(\Sigma)\}^4]+(n-2-3n^{-1})\mu_2^2$ with $
\mu_2=-p^{-1}n^{-1}E[\{|Z_t|_2^2-\tr(\Sigma)\}^2]$. Under the null hypothesis $H_0$, if $\sigma_{n2}>0$ for all sufficiently large $n$, then 
$$
\bigg(
\frac{G_{n,1,1}}{\sigma_{n1}}, \ldots, \frac{G_{n,q,1}}{\sigma_{n1}}, 
\frac{G_{n,1,2}-\mu_2}{\sigma_{n2}}, \ldots, \frac{G_{n,q,2}-\mu_2}{\sigma_{n2}}\bigg) \stackrel{\mathrm{d}}{\to} \mathcal{N}( 0,I_{2q}) ~~\mbox{as}~~ n \to \infty\,.
$$  
\end{proposition}

\begin{remark}\label{rek1}
 Proposition \ref{null-lemma} indicates that $G_{n,\tau, 1}$ and $G_{n,\tau,2}$ exhibit different convergence rates in the high-dimensional setting ($p\rightarrow\infty$ as $n\rightarrow\infty$). Assumption \ref{as:A4} implies that $\sigma_{n1}$ is of constant order, while $\sigma_{n2}$ is influenced by both the underlying distribution of $Z_t$ and the relative growth rate of $p$ and $n$. If $Z_t$ follows the independent components model, it can be shown that 
 $
E[\{|Z_t|_2^2-\tr(\Sigma)\}^2]=O(p)$ and $ E[\{|Z_t|_2^2-\tr(\Sigma)\}^4]=O(p^2)$ under certain conditions, 
which implies $G_{n,\tau,2}=O_{\mathrm{p}}(n^{-1/2})$.
However, if the distribution of $\z_t$ shifts away from the independent components model, $\sigma_{n2}$ may diverge with $p$ and $n$. For instance, if $Z_t$ belongs to the family of elliptical distributions \citep{FZ90}, $
Z_t\stackrel{\mathrm{d}}{=} \xi\Sigma^{1/2}U$, 
where $\bSigma$ denotes the shape matrix of $Z_t$, $\xi>0$ is a scalar variable representing a random length, and $U$ is a $p$-dimensional vector representing a random direction independent of $\xi$ and uniformly distributed on the unit sphere of $\mathbb R^p$. 
For model identifiability, we regularize $\xi$ by $\E(\xi^2)=p$ so that ${\rm cov}(Z_t)=\Sigma$.
Hence, for any nondegenerate variable $\xi$ such that $\Var(p^{-1/2}\xi)$ does not vanish and $\E(\xi^4)=O(p^2)$, we have
\begin{align*}
\E[\{|\z_t|_2^2-\tr(\bSigma)\}^2 ]=
\frac{2 \E(\xi^4)}{p(p+2)} \tr(\bSigma^2)+\left\{\frac{\E(\xi^4)}{p(p+2)}-1\right\} {\tr}^2(\bSigma) = O(p^2)\,.
\end{align*} 
In this case, we have $G_{n,\tau,2}=O_{\mathrm{p}}(pn^{-1/2})$, which may diverge in probability when $p \gg n^{1/2}$.	
\end{remark}

\section{An adaptive $L_2$-type test for white noise}\label{sec:wn}

To address the variance inflation issue of $G_{n,\tau}$, which mainly causes the size inflation phenomenon of the test \eqref{eq:Lq_reject}, we consider a new $L_2$-type test statistic only based on $G_{n,\tau,1}$:
\begin{align}\label{eq:ts}
T_n =  \sum_{\tau=1}^qH_{n,\tau}~~\textrm{with}~~H_{n,\tau}=\frac{G_{n,\tau,1}}{\hat\sigma_{n1}}\,,
\end{align}
where $G_{n,\tau,1}$ is defined in \eqref{eq:Gnt1}, and 
$
 \hat\sigma_{n1}^2 =2p^{-2}\{n^{-1}(n-1)^{-1}\sum_{t\neq s} (\x_t^\top \x_s)^2\}^2$. Theorem \ref{null-Hn} establishes the asymptotic behaviors of $H_{n,1}, \ldots, H_{n,q}$.

\begin{theorem}\label{null-Hn}
Let Assumptions {\rm\ref{as:A2}} and {\rm\ref{as:A4}} hold. Denote by $\lambda_1\geq\cdots\geq\lambda_p$ the eigenvalues of $\bSigma$. Under the null hypothesis $H_0$, as $n\rightarrow\infty$, it holds that 
$
    (H_{n,1}, \ldots, H_{n,q}) \stackrel{\mathrm{d}}{\to} (H_{1}, \ldots, H_{q})$, 
where $H_{1}, \ldots, H_{q}$ are independent random variables such that 
\begin{align*}
 H_{\tau}	
\stackrel{\mathrm{d}}{=}
\begin{cases}
\displaystyle \frac{\sum_{i, j=1}^p \lambda_i \lambda_j (\varepsilon_{i,j,\tau}^2-1)}{\sqrt{2}\sum_{i=1}^p\lambda_i^2}\,,& \text{if $p$ is fixed},	\\
\displaystyle \mathcal{N}(0,1)\,,& \text{if $p \to \infty$},	
\end{cases}
\end{align*}
with  i.i.d. standard Gaussian random variables $\{\varepsilon_{i,j,\tau}\}_{i,j\in [p], \tau\in [q]}$.
\end{theorem}

Theorem \ref{null-Hn} gives a general limiting distribution of $H_{n,\tau}$ under the null hypothesis, yielding a phase transition for the null distribution of the proposed test statistic $T_n$, which causes two main obstacles for determining the associated critical value: (i) it is hard to identify the appropriate phase in practice, and (ii) even if we can identify the appropriate phase, there are still $p$ unknown parameters $\lambda_1,\ldots,\lambda_p$  that need to be estimated in practice when $p$ is fixed.    To overcome these two obstacles, we propose a novel adaptive bootstrap method in Algorithm \ref{alg:adaptive}, which provides a unified testing procedure across two different phases.

\begin{algorithm}[!ht]
  \setstretch{0.9}
  \caption{Adaptive bootstrap procedure for white noise test}\label{alg:adaptive}
  {\footnotesize
    \textbf{Input:} (i) observations $\{\x_t\}_{t=1}^n$, (ii) the number of bootstrap samples $B$, and (iii) the significance level $\alpha$. \\
    1. Calculate the test statistic $T_n$ specified in \eqref{eq:ts}. \\
    2. \For{$b\in[B]$}{
        Generate independent random variables $\{e_t^{(b)}\}_{t=1}^n$ such that $\pr\{e_t^{(b)}=1\}=1/2=\pr\{e_t^{(b)}=-1\}$;
        
        Calculate $\y_t^{(b)}=e_t^{(b)}\x_t$ for each $t\in [n]$;
        
        Calculate bootstrap statistic $T_{n}^{e,(b)}=\sum_{\tau=1}^qH_{n,\tau}^{e,(b)}$, where $H_{n,\tau}^{e,(b)}$ is calculated in the same manner as $H_{n,\tau}$ but with replacing $\{\x_t\}_{t=1}^n$ by $\{\y_t^{(b)}\}_{t=1}^n$;
    }
    3. Calculate the critical value $\hat{{\rm cv}}_{\alpha}$ as the $\lfloor B\alpha \rfloor$-th largest value among $\{T_{n}^{e,(b)}\}_{b=1}^B$. \\
   4.  Reject the null hypothesis if $T_n > \hat{{\rm cv}}_{\alpha}$.
  }
\end{algorithm}

Let $\{e_t\}_{t=1}^n$ be i.i.d. Rademacher random variables. Write $\y_t = e_t \x_t$ and $T_{n}^{e}=\sum_{\tau=1}^q H_{n,\tau}^{e}$, where $H_{n,\tau}^{e}$ is calculated in the same manner as $H_{n,\tau}$ but with replacing $\{\x_t\}_{t=1}^n$ by $\{\y_t\}_{t=1}^n$. To establish the theoretical guarantee of the adaptive testing procedure given in Algorithm \ref{alg:adaptive}, we define
$\tilde{{\rm cv}}_{\alpha} = \inf\{t \in \mathbb{R}: \pr(T_{n}^{e} \le t\,|\, \{\x_t\}_{t=1}^n) \ge 1-\alpha\}$. The critical value $\hat{{\rm cv}}_{\alpha}$ determined in Algorithm \ref{alg:adaptive} with sufficiently large $B$ will provide a good approximation to $\tilde{{\rm cv}}_\alpha$. Theorem \ref{tm:boot_size}  shows that our test has good size control. 

\begin{theorem}\label{tm:boot_size}
Let Assumptions {\rm\ref{as:A2}} and {\rm\ref{as:A4}} hold. Under the null hypothesis $H_0$, as $n\rightarrow\infty$, it holds that $\pr( T_n > \tilde{{\rm cv}}_{\alpha} )\to \alpha$. 
\end{theorem}

To study the power performance of our test, we impose the following assumption on \eqref{eq:VMA}, which is more general than those required in \cite{LZ19} that assume $\A_j = 0$ in \eqref{eq:VMA} for any $j \geq 2$.

\begin{assumption}\label{as:A5}
 $\z_t = \bSigma^{1/2}\w_t$ satisfies the independent components model specified in Definition \ref{def:ICmodel} with $\sup_{i \in [p]}\E(w_{i,t}^8)\leq\nu < \infty$ for some universal constant $\nu>0$.
The coefficient matrices in \eqref{eq:VMA} satisfy $\sum_{j=0}^\infty j\|\A_j\|_2\leq C_3$ for some universal constant $C_3>0$.
\end{assumption}

Let $\mu_n = \sum_{\tau=1}^q E(G_{n,\tau,1})$. Theorem \ref{boot} establishes the consistency of the proposed test under certain local alternatives.

\begin{theorem}\label{boot}
Let Assumptions {\rm\ref{as:A2}}, {\rm\ref{as:A4}} and {\rm\ref{as:A5}} hold. Under the alternative hypothesis, as $n\rightarrow\infty$, if $\mu_n/\max\{p/\sqrt{n}, \sqrt{n/p}\}\to\infty$, it holds that $\pr( T_n > \tilde{{\rm cv}}_{\alpha} )\to 1$.
\end{theorem}

\section{Numerical studies}\label{sec:num}

 In this section, we conduct numerical studies to evaluate the finite-sample performance of the newly proposed test. 
 To implement the proposed test, we set the number of bootstrap samples $B=1000$, and conduct numerical simulations for $q\in \{2, 4, 6, 8\}$. The results for $q=2$ are reported in this section, and the results for $q\in\{4, 6, 8\}$ are reported in the supplementary material. Simulation results indicate that the finite-sample performance of the proposed test is robust to the choice of $q$.  We also compare the proposed test with eight existing methods: (i) the $L_2$-type test (denoted as LLYY) in \cite{LZ19}, (ii) the $L_\infty$-type test (denoted as CYZ) in \cite{Chang+Zhou+Yao:2017}, (iii) the test (denoted as WKX) in \cite{Wang+kong+Xia:2023}, (iv) the $L_\infty$-type test (denoted as Tsay) in \cite{Tsay:2020}, (v) the Fisher's combination test (denoted as $\mathrm{FLM}_{\mathrm{FC}}$) in \cite{Feng2022}, and (vi)--(viii) the $L_\infty$-type tests in \cite{Chen2025} based on Hoeffding's $D$ statistic (denoted as $\mathrm{CSF}_{L_D}$), Blum-Kiefer-Rosenblatt's $R$ statistic (denoted as $\mathrm{CSF}_{L_R}$), and Bergsma-Dassios-Yanagimoto's $\tau^*$ statistic (denoted as $\mathrm{CSF}_{L_{\tau^*}}$), respectively. 
Note that these competing methods are also designed to test whether the first $K$ autocorrelations are zero. Thus, we set $K$ in these methods to be the same as $q$ in our proposed procedure. All simulation results are based on 5000 replications and at the significance level $0.05$.  The replication code for the numerical studies is available at the GitHub repository: \url{https://github.com/JinyuanChang-Lab/AdaptiveHDWhiteNoiseTest}.

\subsection{Empirical sizes}\label{sec:size}

To evaluate the empirical size of the proposed test, we let $\x_t = \z_t=(z_{1,t},\ldots,z_{p,t})^{\top}$, where $\{\z_t\}_{t=1}^n$ is a sequence of i.i.d. $p$-dimensional random vectors satisfying the following three models:

\begin{description} 
	\item[{\bf Model 1.}]  $z_{1,t},\ldots,z_{p,t}$ are i.i.d. random variables from $\rm{Gamma}(4,0.5)-2$.
    \item[{\bf Model 2.}] $\z_t = \sqrt{7}U_t/3$, where $U_t$ follows a multivariate Student's $t$-distribution $\mathcal{T}_9(0, I_p)$. 
    \item[{\bf Model 3.}] $\z_t$ follows a mixture distribution $0.5\mathcal{N}( 0,0.2\I_p)+0.5\mathcal{N}(0,1.8\I_p)$.
\end{description}

\begin{table}
\centering
\footnotesize
\caption{\label{table_size} {Empirical sizes $(\%)$ of the proposed  $L_2$-type test with $q=2$ and the competing methods. `NA' indicates that the results are omitted due to long computation time.}}
\begin{threeparttable}
\resizebox{1\textwidth}{!}{
\def\arraystretch{1.15}
\begin{tabular}{cccccccccccc}
\hline
& $p$ & $n$ & Proposed & LLYY & $\rm{CYZ}$ & $\rm{WKX}$ & $\rm{Tsay}$ & $\mathrm{FLM}_{\mathrm{FC}}$ & $\mathrm{CSF}_{L_D}$ & $\mathrm{CSF}_{L_R}$ & $\mathrm{CSF}_{L_{\tau^*}}$ \\
\hline

Model 1 & 5 & $50$ & 4.96 & 5.10 & 4.84 & 5.20 & 2.70 & 5.44 & 5.58 & 3.60 & 3.86 \\
        &   & $100$ & 4.76 & 5.20 & 6.54 & 4.36 & 3.08 & 5.96 & 4.98 & 4.00 & 4.12 \\
        &   & $200$ & 5.42 & 6.06 & 10.18 & 5.20 & 3.60 & 7.26 & 4.98 & 4.60 & 4.56 \\

        & $\lfloor n^{1/2} \rfloor$ & $50$ & 5.22 & 4.80 & 4.44 & 5.24 & 2.34 & 4.98 & 6.50 & 4.06 & 4.34 \\
        &                            & $100$ & 4.62 & 4.90 & 8.58 & 5.12 & 2.96 & 4.66 & 5.90 & 4.46 & 4.56 \\
        &                            & $200$ & 5.06 & 5.24 & 14.78 & 4.86 & 3.56 & 5.50 & 5.04 & 4.12 & 4.42 \\

        & $n$ & $50$ & 5.32 & 4.60 & 0.14 & 5.16 & 1.44 & 2.24 & 11.54 & 3.04 & 4.44 \\
        &     & $100$ & 4.68 & 4.70 & 0.00 & 5.72 & 2.06 & 3.00 & 7.68 & 3.38 & 4.00 \\
        &     & $200$ & 5.48 & 5.18 & 0.00 & 5.92 & 3.44 & 4.70 & 6.18 & 3.56 & 4.04 \\

        & $n^2$ & $50$ & 5.24 & 4.46 & NA & NA & 1.30 & 1.08 & 21.44 & NA & NA \\
        &       & $100$ & 4.92 & 4.76 & NA & NA & 2.22 & 1.08 & NA & NA & NA \\
        &       & $200$ & 5.14 & 5.14 & NA & NA & NA & NA & NA & NA & NA \\

Model 2 & 5 & $50$ & 5.22 & 5.00 & 4.48 & 4.78 & 2.12 & 5.12 & 6.42 & 4.26 & 4.68 \\
        &   & $100$ & 5.22 & 5.66 & 6.00 & 5.18 & 2.78 & 6.20 & 5.24 & 4.30 & 4.56 \\
        &   & $200$ & 4.72 & 5.30 & 8.82 & 4.94 & 2.94 & 6.98 & 4.64 & 3.96 & 4.04 \\

        & $\lfloor n^{1/2} \rfloor$ & $50$ & 5.36 & 4.76 & 3.64 & 4.86 & 2.10 & 4.48 & 6.68 & 4.24 & 4.56 \\
        &                            & $100$ & 5.34 & 5.44 & 7.04 & 4.96 & 3.36 & 4.50 & 5.80 & 4.24 & 4.74 \\
        &                            & $200$ & 4.74 & 5.20 & 12.40 & 4.90 & 3.22 & 5.12 & 4.88 & 3.66 & 3.96 \\

        & $n$ & $50$ & 5.80 & 9.24 & 0.24 & 4.32 & 1.52 & 1.94 & 11.14 & 2.58 & 3.88 \\
        &     & $100$ & 4.86 & 14.70 & 0.00 & 4.12 & 2.46 & 2.34 & 7.56 & 2.86 & 3.64 \\
        &     & $200$ & 5.24 & 22.62 & 0.00 & 4.48 & 2.96 & 3.38 & 6.34 & 3.70 & 4.18 \\

        & $n^2$ & $50$ & 5.32 & 24.78 & NA & NA & 1.34 & 0.64 & 19.80 & NA & NA \\
        &       & $100$ & 4.88 & 31.68 & NA & NA & 1.96 & 0.76 & NA & NA & NA \\
        &       & $200$ & 5.30 & 36.56 & NA & NA & NA & NA & NA & NA & NA \\

Model 3 & 5 & $50$ & 5.12 & 5.24 & 2.70 & 4.68 & 2.36 & 4.92 & 6.34 & 4.02 & 4.42 \\
        &   & $100$ & 5.72 & 5.96 & 4.54 & 4.84 & 3.22 & 6.34 & 5.42 & 4.58 & 4.84 \\
        &   & $200$ & 4.86 & 5.78 & 7.36 & 5.20 & 2.84 & 6.70 & 4.58 & 4.00 & 3.98 \\

        & $\lfloor n^{1/2} \rfloor$ & $50$ & 5.08 & 5.96 & 2.00 & 4.72 & 2.28 & 4.96 & 6.82 & 4.40 & 4.56 \\
        &                            & $100$ & 5.22 & 6.78 & 3.30 & 5.28 & 3.52 & 5.38 & 6.28 & 4.76 & 5.00 \\
        &                            & $200$ & 4.60 & 6.80 & 8.08 & 4.64 & 3.16 & 5.44 & 4.38 & 3.42 & 3.70 \\

        & $n$ & $50$ & 5.68 & 16.40 & 0.30 & 3.50 & 1.50 & 2.08 & 10.60 & 2.46 & 3.76 \\
        &     & $100$ & 4.92 & 22.32 & 0.00 & 4.06 & 2.06 & 3.24 & 7.78 & 3.14 & 3.92 \\
        &     & $200$ & 5.36 & 30.28 & 0.00 & 4.94 & 2.56 & 3.84 & 5.70 & 3.28 & 3.60 \\

        & $n^2$ & $50$ & 5.22 & 28.70 & NA & NA & 1.38 & 0.94 & 18.44 & NA & NA \\
        &       & $100$ & 4.86 & 34.08 & NA & NA & 2.40 & 1.04 & NA & NA & NA \\
        &       & $200$ & 4.50 & 40.10 & NA & NA & NA & NA & NA & NA & NA \\

\hline
\end{tabular}}
\end{threeparttable}
\end{table}

Model 1 is used in \cite{LZ19} which satisfies the independent components model requirement, while Models 2 and 3 do not satisfy such requirement. 
Table \ref{table_size} reports the empirical sizes of all  methods, which indicates that except for our proposed method and the LLYY test, the other tests are computationally intensive when $p$ is large.
The LLYY test exhibits good size control  under Model 1, while it shows obvious size inflation under Models 2 and 3 when $p=n$ and $p = n^2$. This provides evidence for the discussion in Section \ref{sec:pre} that the established null distribution of the test statistic in \cite{LZ19} may be invalid when $\z_t$ does not satisfy the independent components model. 
The proposed test, the WKX test, the $\mathrm{CSF}_{L_R}$ test, and the $\mathrm{CSF}_{L_{\tau^*}}$ test have good size control in all reported cases. 
The CYZ and Tsay tests are relatively more conservative when $p$ is large. The $\mathrm{FLM}_{\mathrm{FC}}$ test exhibits good size control when $p \le n$, but tends to be conservative when $p \gg n$. The $\mathrm{CSF}_{L_D}$ test tends to suffer from size inflation as $p$ increases.

\subsection{Empirical powers}

To evaluate the empirical power of the proposed test, we consider the following three models:

\begin{description} 
	\item[{\bf Model 4.}] $X_t = 0.1\x_{t-1}+\z_t$ with $\z_t$ generated from Model 1. 
    \item[{\bf Model 5.}]  $\x_t = \z_t - \Theta\z_{t-1}$ with $\z_t$ generated from Model 2, where the coefficient matrix $\Theta=(\theta_{i,j})_{p\times p}$ is a banded-matrix. Specifically, let $\theta_{i,j} = 0$ if $|i - j| > 1$. For $|i - j| \le 1$, $\theta_{i,j}$ is independently generated such that $\theta_{i,j} = 0$ with probability $0.7$, and follows uniform distribution $\mathcal{U}(-0.95, 0.95)$ with probability $0.3$. 
    \item[{\bf Model 6.}] $\x_t= \z_t+0.07 \z_{t-1}$ with $\z_t$ generated from Model 3.
\end{description}

Model 4 also used in \cite{LZ19} is the first-order vector autoregressive model, where $\x_t$ can be formulated as \eqref{eq:VMA} with $\A_j = 0.1^j \I_p$ for $j \geq 1$. 
Models 5 and 6 represent the first-order vector moving average models, which follow  \eqref{eq:VMA} with $\A_j = 0$ for any $j \geq 2$. Model 5 is similar to the second setting for power evaluation in \cite{Tsay:2020}, but we use multivariate Student’s $t$-distributed innovations. 
Table \ref{table_power} reports the empirical powers of all methods. The results show that the proposed test outperforms the  CYZ test, the WKX test and the Tsay test across most of the reported scenarios, except that the empirical power of the proposed test is very close to that of the CYZ test when $(n,p)=(200, 5)$ under Model 5. 
Under Model 4, our proposed test achieves comparable powers to the LLYY test, while under Model 6, our proposed test performs comparably to the LLYY test when $p\leq n$, and outperforms the LLYY test when $p\gg n$. Moreover, under Models 4 and 6, our proposed test outperforms the $\mathrm{CSF}_{L_D}$, $\mathrm{CSF}_{L_R}$, and $\mathrm{CSF}_{L_{\tau^*}}$ tests with the advantage being more pronounced when $p \ge n$, and performs comparably to the $\mathrm{FLM}_{\mathrm{FC}}$ test. Under Model 5, our proposed test has lower empirical power than the LLYY and $\mathrm{FLM}_{\mathrm{FC}}$ tests when $n$ is small ($n=50$). However, when $n \ge 100$, the empirical powers of our proposed test become comparable to those of these two tests. Compared with the $\mathrm{CSF}_{L_D}$, $\mathrm{CSF}_{L_R}$, and $\mathrm{CSF}_{L_{\tau^*}}$ tests, our proposed test generally achieves comparable or even higher empirical powers under Model 5.

\begin{table}
\centering
\footnotesize
\caption{\label{table_power} {Empirical powers $(\%)$ of the proposed $L_2$-type test with $q=2$ and the competing methods. ‘NA’ indicates that the results are omitted due to long computation time. }}
\begin{threeparttable}
\resizebox{1\textwidth}{!}{
\def\arraystretch{1.15}
\begin{tabular}{cccccccccccc}
\hline
& $p$ & $n$ & Proposed & LLYY & $\rm{CYZ}$ & $\rm{WKX}$ & $\rm{Tsay}$ & $\mathrm{FLM}_{\mathrm{FC}}$ & $\mathrm{CSF}_{L_D}$ & $\mathrm{CSF}_{L_R}$ & $\mathrm{CSF}_{L_{\tau^*}}$ \\
\hline

Model 4 & 5 & $50$ & 11.08 & 11.54 & 5.24 & 6.44 & 3.16 & 11.30 & 7.18 & 4.92 & 5.34 \\
        &   & $100$ & 15.40 & 16.50 & 10.16 & 10.18 & 5.10 & 17.64 & 7.18 & 6.20 & 6.44 \\
        &   & $200$ & 27.38 & 29.30 & 23.74 & 18.64 & 12.36 & 32.44 & 11.62 & 11.42 & 11.46 \\

        & $\lfloor n^{1/2} \rfloor$ & $50$ & 11.28 & 11.42 & 4.80 & 6.96 & 2.92 & 10.62 & 7.82 & 5.06 & 5.28 \\
        &                            & $100$ & 17.76 & 18.22 & 12.46 & 10.02 & 4.42 & 17.90 & 8.20 & 6.42 & 6.84 \\
        &                            & $200$ & 34.64 & 34.78 & 33.08 & 19.18 & 10.34 & 35.28 & 11.56 & 9.88 & 10.10 \\

        & $n$ & $50$ & 33.72 & 30.66 & 0.14 & 25.52 & 1.88 & 20.58 & 12.28 & 3.00 & 4.64 \\
        &     & $100$ & 75.22 & 72.98 & 0.02 & 65.90 & 5.30 & 62.46 & 9.06 & 3.86 & 4.70 \\
        &     & $200$ & 99.80 & 99.84 & 0.00 & 99.62 & 35.50 & 99.58 & 10.04 & 6.28 & 7.00 \\

        & $n^2$ & $50$ & 100.00 & 100.00 & NA & NA & 2.60 & 100.00 & 26.26 & NA & NA \\
        &       & $100$ & 100.00 & 100.00 & NA & NA & 12.46 & 100.00 & NA & NA & NA \\
        &       & $200$ & 100.00 & 100.00 & NA & NA & NA & NA & NA & NA & NA \\

Model 5 & 5 & $50$ & 71.24 & 73.80 & 37.22 & 32.36 & 35.50 & 82.52 & 69.76 & 65.54 & 66.92 \\
        &   & $100$ & 90.62 & 91.00 & 83.90 & 58.38 & 78.84 & 93.78 & 91.36 & 90.84 & 91.04 \\
        &   & $200$ & 94.94 & 95.18 & 95.00 & 79.52 & 92.58 & 96.38 & 95.92 & 95.72 & 95.80 \\

        & $\lfloor n^{1/2} \rfloor$ & $50$ & 76.94 & 79.90 & 24.14 & 27.70 & 29.84 & 87.88 & 76.84 & 71.34 & 73.20 \\
        &                            & $100$ & 97.86 & 98.24 & 70.26 & 52.72 & 73.28 & 99.50 & 98.92 & 98.62 & 98.76 \\
        &                            & $200$ & 99.96 & 99.96 & 98.64 & 81.28 & 98.62 & 99.98 & 99.98 & 99.98 & 99.98 \\

        & $n$ & $50$ & 92.78 & 91.42 & 0.30 & 41.10 & 2.98 & 99.44 & 98.02 & 89.60 & 93.94 \\
        &     & $100$ & 99.98 & 100.00 & 0.04 & 74.74 & 7.74 & 100.00 & 100.00 & 100.00 & 100.00 \\
        &     & $200$ & 100.00 & 100.00 & 0.00 & 96.48 & 18.06 & 100.00 & 100.00 & 100.00 & 100.00 \\

        & $n^2$ & $50$ & 95.52 & 71.12 & NA & NA & 2.04 & 97.46 & 100.00 & NA & NA \\
        &       & $100$ & 100.00 & 89.98 & NA & NA & 1.90 & 100.00 & NA & NA & NA \\
        &       & $200$ & 100.00 & 98.16 & NA & NA & NA & NA & NA & NA & NA \\

Model 6 & 5 & $50$ & 7.52 & 7.66 & 2.48 & 5.16 & 2.44 & 7.08 & 6.56 & 4.34 & 4.74 \\
        &   & $100$ & 9.76 & 10.26 & 4.92 & 7.44 & 4.38 & 11.00 & 6.58 & 5.48 & 5.64 \\
        &   & $200$ & 13.58 & 15.18 & 10.42 & 11.14 & 6.74 & 16.68 & 8.40 & 7.94 & 7.86 \\

        & $\lfloor n^{1/2} \rfloor$ & $50$ & 8.02 & 8.56 & 2.24 & 5.78 & 2.56 & 7.28 & 6.84 & 4.62 & 4.90 \\
        &                            & $100$ & 10.36 & 11.94 & 4.14 & 7.46 & 3.54 & 9.94 & 7.34 & 5.54 & 5.70 \\
        &                            & $200$ & 15.78 & 16.96 & 11.32 & 9.80 & 6.66 & 15.50 & 7.58 & 6.80 & 6.92 \\

        & $n$ & $50$ & 15.22 & 20.86 & 0.32 & 9.40 & 1.52 & 8.00 & 10.64 & 2.64 & 3.92 \\
        &     & $100$ & 32.06 & 31.32 & 0.08 & 24.78 & 2.22 & 19.82 & 8.92 & 3.56 & 4.56 \\
        &     & $200$ & 74.58 & 45.28 & 0.00 & 65.02 & 6.46 & 60.44 & 7.14 & 4.00 & 4.46 \\

        & $n^2$ & $50$ & 100.00 & 34.26 & NA & NA & 1.58 & 100.00 & 20.76 & NA & NA \\
        &       & $100$ & 100.00 & 41.68 & NA & NA & 3.32 & 100.00 & NA & NA & NA \\
        &       & $200$ & 100.00 & 52.28 & NA & NA & NA & NA & NA & NA & NA \\

\hline
\end{tabular}}
\end{threeparttable}
\end{table}

\bibliographystyle{apalike}
\bibliography{ref.bib}

%\begin{spacing}{1.0}
%\printbibliography{} %   \bibliography{ref.bib}
%\end{spacing}
\clearpage
%\newpage
\onehalfspacing
\normalsize
\numberwithin{equation}{section}

\begin{center}
	{\bf  \Large
		Supplementary material for ``An Adaptive $L_2$-type Test for High-dimensional White Noise" by Jinyuan Chang, Jing He, Weiming Li and Chen Lin}  \\
\end{center}

\setcounter{page}{1}
\renewcommand{\thepage}{S\arabic{page}}

\bigskip

\setcounter{equation}{0}
\setcounter{section}{0}
\setcounter{table}{0}
\setcounter{figure}{0}
\renewcommand{\thetable}{T\arabic{table}}
\renewcommand{\thefigure}{F\arabic{figure}}
\renewcommand{\thesection}{\Alph{section}}

\newtheorem{lem}{Lemma} % 定义lem环境
\numberwithin{lem}{section} % 按 section 计数，并且每个 section 重新计数
\setcounter{lem}{0}
\renewcommand{\thelem}{\thesection\arabic{lem}} % 自定义

We first introduce some notation that will be used throughout the supplementary material. Let $K$ denote a generic positive constant that does not depend on $n$ and $p$,  which may be different in different places. 
The notation $I(\cdot)$ denotes the indicator function.  
In the sequel, we write $\alpha_{\ell p} = p^{-1}\tr (\Sigma^\ell)$ and denote $\lim_{p\to\infty}  \alpha_{\ell p} = \alpha_{\ell} >0$ for $\ell = 1, 2$. For convenience, we assume $n > Cq$ for some sufficiently large constant $C>0$.

\section{Comparison of our proposed test and the test of \cite{LZ19}}
Our method and that of \cite{LZ19} have the same starting point, namely,
\[
G_{n,\tau}
=\frac{n}{p}\operatorname{tr}(S_{\tau}S_{\tau}^{\top})
-\frac{1}{p}\operatorname{tr}^{2}(S_{0})\,,
\]
where $S_{\tau}$ is the sample autocovariance matrix at lag $\tau$. This statistic can be decomposed as
\[
G_{n,\tau}=G_{n,\tau,1}+G_{n,\tau,2}\,,
\]
where
\begin{align*}
G_{n,\tau,1}
&=\frac{1}{np}\sum_{t\ne s}
X_{t}^{\top}X_{s}
X_{t+\tau}^{\top}X_{s+\tau}\,,\\
G_{n,\tau,2}
&=\frac{1}{np}\sum_{t=1}^{n}
|X_{t}|_{2}^{2}|X_{t+\tau}|_{2}^{2}
-\frac{1}{p}\bigg(\frac{1}{n}\sum_{t=1}^{n}|X_{t}|_{2}^{2}\bigg)^{2}\,.
\end{align*}
\cite{LZ19} construct their test based on the whole statistic $G_{n,\tau}$, whereas our test just uses the first component $G_{n,\tau,1}$. Different asymptotic behaviors of $G_{n,\tau,1}$ and $G_{n,\tau,2}$ motivate our such choice. The first component $G_{n,\tau,1}$ is asymptotically normal under the null hypothesis in the high-dimensional setting, with variance $\sigma_{n1}^{2}=2p^{-2}\operatorname{tr}^{2}(\Sigma^{2}) \asymp 1$ under the mild Assumptions 1 and 2.
This result does not require the assumption of independent components model or a restriction on the relative growth rate of the dimension $p$ and the sample size $n$. However, the second component $G_{n,\tau,2}$ is more sensitive to the distribution of the innovations and the relative growth rate of $p$ and $n$: 
\begin{itemize}
\item Under the independent components model, we have $G_{n,\tau,2}=O_{\mathrm{p}}(n^{-1/2})$, 
so it is asymptotically negligible in comparison to the first component $G_{n,\tau,1}$.

\item If the assumption of independent components model is violated, $G_{n,\tau,2}$ may no longer be negligible in comparison to the first component $G_{n,\tau,1}$. For example, for elliptically distributed innovations with a nonvanishing variance of the normalized radial component, $G_{n,\tau,2}$ can have the order $O_{\mathrm{p}}(pn^{-1/2})$, which is not necessarily negligible if $p \gtrsim n^{1/2}$.
\end{itemize}

\noindent\begin{minipage}{\linewidth}
\centering
\captionsetup{type=figure}
\includegraphics[width=0.48\linewidth]{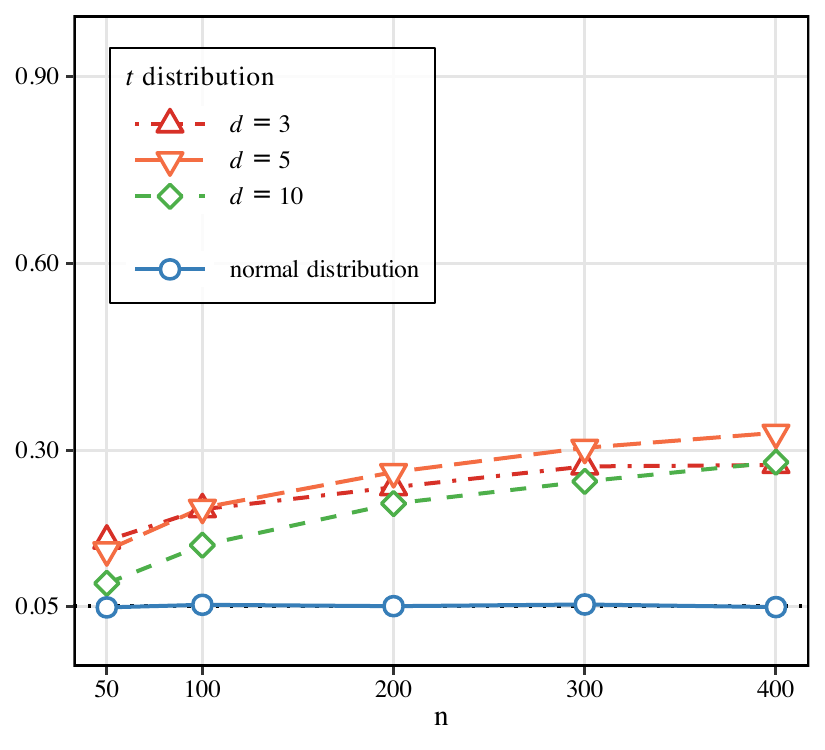}\hfill
\includegraphics[width=0.48\linewidth]{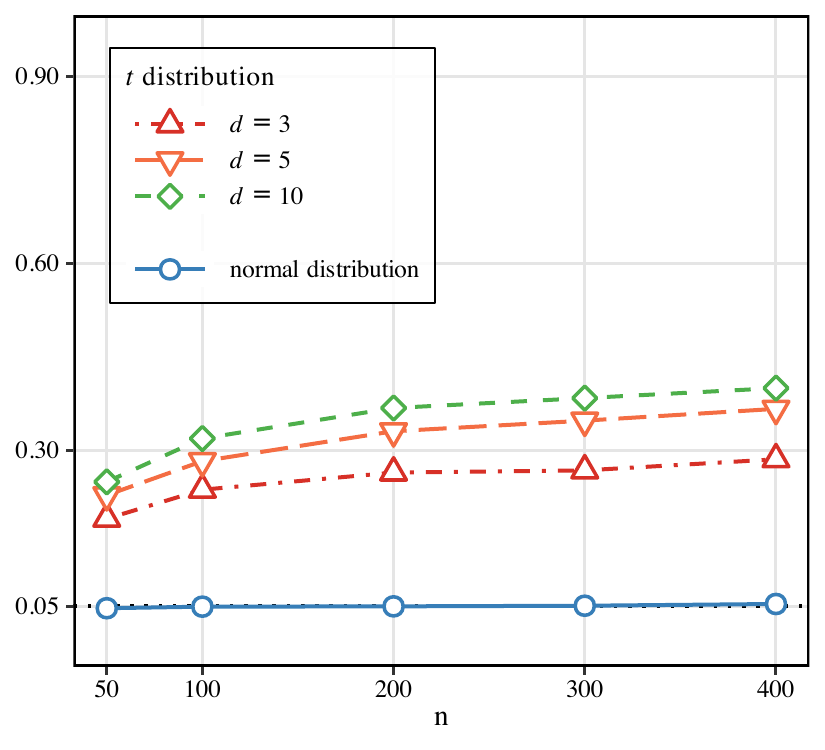}
\caption{Empirical sizes of the test proposed by \cite{LZ19}  for $p=n$ (left) and $p=n^2$ (right), $q=2$ and $n\in\{50,100,200,300,400\}$, based on 10000 independent replications. The data is  generated independently from (i) the multivariate Student's $t$-distribution $\mathcal{T}_d(0,I_p)$ with $d\in\{3,5,10\}$, and (ii) the multivariate normal distribution $\mathcal{N}(0,I_p)$. The nominal significance level is 5\%. The results show that the  test proposed by \cite{LZ19} has good size control under the multivariate normal distribution but exhibits size inflation under the multivariate Student's $t$-distributions.}
\label{fig:llyy-size-df}
\end{minipage}

\vspace{0.5cm}
The asymptotic behavior of $G_{n,\tau,2}$ can cause size inflation of the test proposed by \cite{LZ19}. See Figure \ref{fig:llyy-size-df} for details. In contrast, our proposed test just using $G_{n,\tau,1}$ can avoid this source of size distortion. More detailed comparison between our proposed test and the test of \cite{LZ19} is summarized as follows:

\begin{itemize}
\item \textbf{Size control.} 
The proposed test provides valid size control in a broader range of settings than the test of \cite{LZ19}. Theorem 2  shows that our proposed test achieves valid asymptotic size control without requiring either $p/n\to c \in (0,\infty)$ or the independent components model, which can accommodate more general cases. The simulation results in Table 1 further confirm the good finite-sample performance of our proposed test in terms of size control, showing that its empirical sizes remain close to the nominal level across cases with different relationships between $p$ and $n$, and even when the independent components model is violated.

\item \textbf{Power.} 
A potential limitation of our approach is that just using $G_{n,\tau,1}$ may discard useful signal contained in $G_{n,\tau,2}$ and thus reduce power under some alternatives. Nevertheless, this choice is necessary for establishing valid size control under more general settings considered in our framework, since the asymptotic behavior of $G_{n,\tau,2}$ generally requires stronger assumptions to characterize. Despite this potential loss of signal, our proposed test generally achieves empirical power comparable to that of \cite{LZ19} and even outperforms their test when $p$ is large. Specifically, under Model 4, where both tests have proper size control, our proposed test has slightly lower empirical power than the test of \cite{LZ19}, but these two tests remain comparable in terms of empirical power with the difference diminishing as $n$ increases. Under Models 5 and 6, a direct power comparison may not be fair because the test of \cite{LZ19} suffers from size inflation, which may lead to higher empirical power. Even if we make such a comparison, our proposed test performs comparably when $p\le n$ and achieves higher empirical power when $p\gg n$.

\end{itemize}

In general, by Wold's Decomposition Theorem \citep[pp.~25]{Lutkepohl2005}, a stationary linear process $X_t$ can be formulated as a moving average representation:
\begin{align}\label{eq:VMA1}
X_t=Z_{t}+\sum_{j=1}^{\infty} A_jZ_{t-j}\,,
\end{align}
where each $A_j\in\mathbb R^{p\times p}$ represents the coefficient matrix and $\{Z_t\} _{t\in \mathbb Z}$ is a $p$-dimensional white-noise process with mean zero. 
To establish the theoretical validity of our proposed test, we impose the independent and identically distributed assumption on $Z_t$ to simplify our technical arguments, which is weaker than the assumptions used in \cite{LZ19}. 
Regarding robustness to model misspecification, it is interesting to investigate whether our proposed method is still valid when $Z_t$'s in \eqref{eq:VMA1} are not independent and identically distributed.  
To answer this question, we first investigate the finite-sample performance of our proposed method in the following two examples, where the independent and identically distributed assumption does not hold:
\begin{description}

\item[\textbf{Model 7.}]
$X_t=Z_t = (z_{1,t},\ldots,z_{p,t})^{\top}$, where $Z_t$ consists of $p$ independent generalized autoregressive conditional heteroskedasticity processes, i.e. $z_{j,t}=\sigma_{j,t}u_{j,t}$, where $u_{j,t}$, $t \ge 0$, are i.i.d. random variables from $\mathcal{N}(0,1)$ and
$\sigma_{j,t}^2=\gamma_{0,j}+\gamma_{1,j}z_{j,t-1}^2+\beta_{1,j}\sigma_{j,t-1}^2$ with $\gamma_{0,j}$, $\gamma_{1,j}$, and $\beta_{1,j}$ independently generated from
$\mathcal{U}(0.25,0.50)$, $\mathcal{U}(0,0.25)$, and $\mathcal{U}(0,0.25)$, respectively. Set $\sigma_{j,0}^2=\gamma_{0,j}/(1-\gamma_{1,j}-\beta_{1,j})$.

\item[\textbf{Model 8.}]
$X_t=AX_{t-1}+Z_t$ with $Z_t$ generated from Model 7, where the coefficient matrix $A=\mathrm{diag}(a_1,\ldots,a_p)$ with $a_j=0.4$ for $j\leq k_0$ and $a_j=0$ for $j>k_0$, where $k_0=\max\{1,\lfloor p/4\rfloor\}$.

\end{description}

Models 7 and 8 are used to evaluate the empirical size and  empirical power, respectively. In both models, $\{z_{j,t}\}$ is generated from a GARCH(1,1) process with the initial value $z_{j,0}$ shifting away from the stationary distribution $\pi(\cdot)$, where $\pi(\cdot)$ is the distribution of 
\begin{align*}
u_{j,t}
\sqrt{
\gamma_{0,j}
\bigg[
1+
\sum_{k=1}^{\infty}
\prod_{\ell=1}^{k}
\bigg(
\gamma_{1,j}u_{j,t-\ell}^{2}
+\beta_{1,j}
\bigg)
\bigg]
}\,.
\end{align*}
Therefore, the marginal distribution of $z_{j,t}$ will vary with time $t$. 
As shown in Table~\ref{table_model78_2}, the LLYY test exhibits some size inflation under Model 7, while the proposed test maintains good size control under Model 7 and achieves comparable power to that of the LLYY test under Model 8. This indicates that the proposed test continues to perform well even when the independent and identically distributed assumption is violated. How to establish the theoretical guarantee of our proposed method without the independent and identically distributed assumption is nontrivial.
Even under the independent and identically distributed assumption, the theoretical analysis requires lengthy and involved proofs. Extending the theory to weaker assumptions would require substantial additional technical work and considerably complicated proofs. We therefore leave such a theoretical extension for future research.

\begin{table}
\centering
\footnotesize
\caption{\label{table_model78_2}
Empirical sizes $(\%)$ and powers $(\%)$ of the proposed $L_2$-type test and
the LLYY test under Models 7 and 8 for different choices of $q$.}

\begin{threeparttable}

\resizebox{\textwidth}{!}{
\def\arraystretch{1.15}

\begin{tabular}{ccccccccccccccccccc}
\hline

& &
\multicolumn{8}{c}{Model 7}
&
\multicolumn{8}{c}{Model 8}
\\
\cline{3-10}\cline{12-19}

& &
\multicolumn{2}{c}{$q=2$} &
\multicolumn{2}{c}{$q=4$} &
\multicolumn{2}{c}{$q=6$} &
\multicolumn{2}{c}{$q=8$} & &
\multicolumn{2}{c}{$q=2$} &
\multicolumn{2}{c}{$q=4$} &
\multicolumn{2}{c}{$q=6$} &
\multicolumn{2}{c}{$q=8$}
\\
\cmidrule(lr){3-4}
\cmidrule(lr){5-6}
\cmidrule(lr){7-8}
\cmidrule(lr){9-10}
\cmidrule(lr){12-13}
\cmidrule(lr){14-15}
\cmidrule(lr){16-17}
\cmidrule(lr){18-19}

$p$ & $n$
& Proposed & LLYY
& Proposed & LLYY
& Proposed & LLYY
& Proposed & LLYY &
& Proposed & LLYY
& Proposed & LLYY
& Proposed & LLYY
& Proposed & LLYY
\\
\hline

5 & $50$
& 4.46 & 7.24
& 4.40 & 6.16
& 4.94 & 6.36
& 4.84 & 5.44 &
& 25.96 & 32.54
& 21.82 & 25.70
& 19.20 & 22.38
& 17.68 & 19.90
\\

& $100$
& 4.62 & 8.08
& 4.74 & 7.70
& 4.56 & 6.82
& 4.92 & 6.98 &
& 50.78 & 58.76
& 41.50 & 47.20
& 36.36 & 40.74
& 32.86 & 36.72
\\

& $200$
& 5.20 & 8.96
& 5.62 & 8.32
& 5.06 & 7.30
& 5.28 & 7.26 &
& 82.04 & 87.20
& 71.96 & 77.58
& 65.08 & 70.12
& 59.70 & 64.36
\\

\cline{1-19}

$\lfloor n^{1/2}\rfloor$ & $50$
& 5.18 & 7.96
& 4.52 & 6.44
& 4.88 & 5.98
& 4.74 & 5.58 &
& 22.28 & 28.02
& 17.70 & 21.28
& 15.54 & 18.02
& 14.66 & 16.74
\\

& $100$
& 4.74 & 7.70
& 4.90 & 7.26
& 5.18 & 6.94
& 5.06 & 6.40 &
& 59.00 & 66.08
& 48.18 & 53.50
& 42.68 & 47.00
& 38.96 & 43.46
\\

& $200$
& 4.98 & 8.00
& 4.70 & 6.96
& 5.78 & 7.30
& 5.36 & 6.82 &
& 95.88 & 97.04
& 89.40 & 91.50
& 84.08 & 87.04
& 80.30 & 82.70
\\

\cline{1-19}

$n$ & $50$
& 5.22 & 6.96
& 4.96 & 6.22
& 4.94 & 5.80
& 5.20 & 5.42 &
& 69.26 & 74.22
& 70.42 & 73.30
& 74.12 & 75.40
& 77.30 & 77.84
\\

& $100$
& 5.18 & 7.54
& 5.00 & 6.88
& 5.70 & 6.62
& 5.32 & 6.14 &
& 99.52 & 99.74
& 99.58 & 99.72
& 99.78 & 99.84
& 99.92 & 99.94
\\

& $200$
& 4.86 & 7.30
& 4.54 & 6.22
& 4.50 & 5.82
& 4.64 & 5.76 &
& 100.00 & 100.00
& 100.00 & 100.00
& 100.00 & 100.00
& 100.00 & 100.00
\\

\cline{1-19}

$n^2$ & $50$
& 5.30 & 7.26
& 4.80 & 6.28
& 5.22 & 5.98
& 5.28 & 5.68 &
& 100.00 & 100.00
& 100.00 & 100.00
& 100.00 & 100.00
& 100.00 & 100.00
\\

& $100$
& 4.88 & 7.06
& 5.02 & 6.52
& 4.92 & 5.96
& 4.86 & 5.58 &
& 100.00 & 100.00
& 100.00 & 100.00
& 100.00 & 100.00
& 100.00 & 100.00
\\

& $200$
& 5.38 & 7.52
& 5.06 & 6.44
& 5.02 & 6.18
& 5.10 & 6.14 &
& 100.00 & 100.00
& 100.00 & 100.00
& 100.00 & 100.00
& 100.00 & 100.00
\\

\hline

\end{tabular}
}

\end{threeparttable}
\end{table}

\section{Additional simulation results}\label{sec:Additional simulations}
This section provides simulation results of the proposed adaptive $L_2$-type test with $q\in\{4,6,8\}$, for detailed simulation settings described in Section \ref{sec:num}. Tables \ref{suptable_size1}--\ref{suptable_size3} report the empirical sizes of the proposed test and the eight competing methods under Models 1--3, and Tables \ref{suptable_power1}--\ref{suptable_power3} report the empirical powers of the proposed test and the eight competing methods under Models 4--6.

\begin{table}[!htbp]
\centering
\footnotesize
\caption{\label{suptable_size1} {Empirical sizes $(\%)$ of the proposed  $L_2$-type test with $q=4$ and the competing methods. `NA' indicates that the results are omitted due to long computation time.}}
\begin{threeparttable}
\resizebox{1\textwidth}{!}{
\def\arraystretch{1.15}
\begin{tabular}{cccccccccccc}
\hline
& $p$ & $n$ & Proposed & LLYY & $\rm{CYZ}$ & $\rm{WKX}$ & $\rm{Tsay}$ & $\mathrm{FLM}_{\mathrm{FC}}$ & $\mathrm{CSF}_{L_D}$ & $\mathrm{CSF}_{L_R}$ & $\mathrm{CSF}_{L_{\tau^*}}$ \\
\hline

Model 1 & 5 & $50$ & 5.24 & 5.00 & 4.56 & 4.62 & 2.46 & 4.70 & 6.16 & 3.80 & 3.88 \\
        &   & $100$ & 4.84 & 5.42 & 6.54 & 4.88 & 3.04 & 5.08 & 5.08 & 4.18 & 4.34 \\
        &   & $200$ & 4.58 & 5.06 & 10.64 & 4.76 & 3.48 & 6.18 & 4.68 & 4.40 & 4.48 \\

        & $\lfloor n^{1/2} \rfloor$ & $50$ & 5.54 & 5.08 & 3.86 & 5.68 & 2.30 & 3.92 & 6.70 & 3.68 & 4.14 \\
        &                            & $100$ & 4.52 & 4.18 & 8.80 & 4.68 & 2.56 & 4.44 & 6.34 & 4.40 & 4.44 \\
        &                            & $200$ & 5.34 & 5.52 & 15.56 & 5.34 & 3.00 & 4.94 & 5.64 & 4.06 & 4.48 \\

        & $n$ & $50$ & 5.38 & 4.86 & 0.42 & 5.64 & 1.08 & 1.82 & 12.16 & 2.18 & 3.70 \\
        &     & $100$ & 4.88 & 4.52 & 0.02 & 5.38 & 2.50 & 2.92 & 8.22 & 2.90 & 3.86 \\
        &     & $200$ & 5.64 & 5.20 & 0.00 & 6.04 & 3.04 & 4.12 & 6.26 & 3.50 & 3.96 \\

        & $n^2$ & $50$ & 4.78 & 4.22 & NA & NA & 1.08 & 1.02 & 22.66 & NA & NA \\
        &       & $100$ & 5.02 & 4.26 & NA & NA & 1.76 & 1.28 & NA & NA & NA \\
        &       & $200$ & 5.06 & 4.80 & NA & NA & NA & NA & NA & NA & NA \\

Model 2 & 5 & $50$ & 5.20 & 5.04 & 3.86 & 5.16 & 2.32 & 4.42 & 6.92 & 4.42 & 4.68 \\
        &   & $100$ & 5.54 & 5.66 & 5.94 & 5.30 & 2.74 & 5.96 & 5.72 & 4.46 & 4.36 \\
        &   & $200$ & 5.24 & 5.68 & 9.02 & 5.34 & 3.06 & 6.32 & 4.92 & 4.00 & 4.22 \\

        & $\lfloor n^{1/2} \rfloor$ & $50$ & 5.22 & 4.78 & 3.00 & 5.18 & 1.90 & 3.96 & 7.16 & 3.40 & 3.90 \\
        &                            & $100$ & 5.32 & 5.40 & 7.02 & 5.26 & 2.90 & 4.46 & 6.34 & 4.06 & 4.56 \\
        &                            & $200$ & 5.00 & 5.16 & 12.56 & 5.22 & 3.22 & 4.88 & 5.14 & 4.16 & 4.32 \\

        & $n$ & $50$ & 5.58 & 7.76 & 0.56 & 4.20 & 0.96 & 1.58 & 11.10 & 2.26 & 3.68 \\
        &     & $100$ & 5.02 & 13.74 & 0.00 & 4.40 & 2.08 & 2.08 & 7.76 & 2.46 & 3.22 \\
        &     & $200$ & 4.50 & 20.74 & 0.00 & 4.12 & 2.62 & 2.50 & 5.90 & 3.40 & 3.92 \\

        & $n^2$ & $50$ & 5.18 & 20.30 & NA & NA & 1.10 & 0.62 & 21.26 & NA & NA \\
        &       & $100$ & 4.90 & 28.44 & NA & NA & 2.34 & 0.74 & NA & NA & NA \\
        &       & $200$ & 5.20 & 35.04 & NA & NA & NA & NA & NA & NA & NA \\

Model 3 & 5 & $50$ & 5.66 & 4.96 & 2.00 & 4.88 & 2.20 & 4.56 & 6.90 & 3.98 & 4.36 \\
        &   & $100$ & 5.28 & 5.76 & 3.34 & 5.36 & 3.08 & 5.60 & 5.86 & 4.32 & 4.60 \\
        &   & $200$ & 4.46 & 5.02 & 7.74 & 4.70 & 2.90 & 6.38 & 4.38 & 4.08 & 4.26 \\

        & $\lfloor n^{1/2} \rfloor$ & $50$ & 5.06 & 5.08 & 1.04 & 4.98 & 2.06 & 3.94 & 7.50 & 4.10 & 4.56 \\
        &                            & $100$ & 5.02 & 6.46 & 3.30 & 5.18 & 3.32 & 4.86 & 5.92 & 4.18 & 4.40 \\
        &                            & $200$ & 4.72 & 5.82 & 7.90 & 4.74 & 2.90 & 4.66 & 5.18 & 3.92 & 4.12 \\

        & $n$ & $50$ & 5.32 & 13.28 & 0.76 & 3.44 & 1.36 & 1.66 & 11.46 & 1.86 & 3.26 \\
        &     & $100$ & 5.18 & 20.74 & 0.06 & 3.86 & 1.98 & 2.50 & 7.96 & 2.76 & 3.70 \\
        &     & $200$ & 5.20 & 28.16 & 0.00 & 4.36 & 2.58 & 3.28 & 5.52 & 2.80 & 3.34 \\

        & $n^2$ & $50$ & 4.78 & 22.68 & NA & NA & 1.34 & 0.90 & 21.16 & NA & NA \\
        &       & $100$ & 4.96 & 28.74 & NA & NA & 1.86 & 0.82 & NA & NA & NA \\
        &       & $200$ & 4.74 & 36.54 & NA & NA & NA & NA & NA & NA & NA \\

\hline
\end{tabular}}
\end{threeparttable}
\end{table}

\begin{table}[!htbp]
\centering
\footnotesize
\caption{\label{suptable_size2} {Empirical sizes $(\%)$ of the proposed  $L_2$-type test with $q=6$ and the competing methods. `NA' indicates that the results are omitted due to long computation time.}}
\begin{threeparttable}
\resizebox{1\textwidth}{!}{
\def\arraystretch{1.15}
\begin{tabular}{cccccccccccc}
\hline
& $p$ & $n$ & Proposed & LLYY & $\rm{CYZ}$ & $\rm{WKX}$ & $\rm{Tsay}$ & $\mathrm{FLM}_{\mathrm{FC}}$ & $\mathrm{CSF}_{L_D}$ & $\mathrm{CSF}_{L_R}$ & $\mathrm{CSF}_{L_{\tau^*}}$ \\
\hline

Model 1 & 5 & $50$ & 4.72 & 4.52 & 3.18 & 5.08 & 2.20 & 3.62 & 6.82 & 4.04 & 4.14 \\
        &   & $100$ & 4.84 & 5.06 & 6.00 & 5.02 & 3.32 & 5.16 & 5.68 & 4.16 & 4.40 \\
        &   & $200$ & 4.62 & 5.00 & 9.46 & 4.80 & 3.38 & 5.80 & 5.04 & 4.42 & 4.54 \\

        & $\lfloor n^{1/2} \rfloor$ & $50$ & 4.90 & 4.48 & 3.18 & 5.36 & 2.30 & 3.36 & 7.50 & 4.04 & 4.62 \\
        &                            & $100$ & 4.76 & 4.40 & 7.38 & 4.66 & 2.90 & 3.94 & 6.44 & 4.30 & 4.84 \\
        &                            & $200$ & 5.08 & 5.16 & 15.04 & 4.92 & 3.24 & 4.32 & 5.06 & 3.90 & 4.10 \\

        & $n$ & $50$ & 5.10 & 4.38 & 0.58 & 4.96 & 0.94 & 1.58 & 12.54 & 1.70 & 3.20 \\
        &     & $100$ & 4.70 & 4.48 & 0.02 & 5.72 & 2.02 & 2.64 & 8.36 & 2.48 & 3.36 \\
        &     & $200$ & 5.46 & 5.16 & 0.00 & 6.74 & 3.06 & 4.06 & 6.80 & 3.72 & 3.96 \\

        & $n^2$ & $50$ & 4.78 & 4.42 & NA & NA & 0.98 & 0.90 & 23.72 & NA & NA \\
        &       & $100$ & 4.48 & 3.86 & NA & NA & 2.18 & 0.92 & NA & NA & NA \\
        &       & $200$ & 4.62 & 4.50 & NA & NA & NA & NA & NA & NA & NA \\

Model 2 & 5 & $50$ & 5.36 & 4.40 & 2.58 & 5.44 & 2.58 & 3.70 & 7.18 & 3.90 & 4.38 \\
        &   & $100$ & 4.86 & 5.18 & 4.90 & 4.64 & 2.56 & 5.16 & 5.74 & 4.44 & 4.56 \\
        &   & $200$ & 4.96 & 5.42 & 8.24 & 5.42 & 3.04 & 5.90 & 5.14 & 4.16 & 4.36 \\

        & $\lfloor n^{1/2} \rfloor$ & $50$ & 5.56 & 4.70 & 2.00 & 4.82 & 1.88 & 3.08 & 7.18 & 3.46 & 4.04 \\
        &                            & $100$ & 5.30 & 5.24 & 6.00 & 5.30 & 2.62 & 3.84 & 6.46 & 3.98 & 4.42 \\
        &                            & $200$ & 5.60 & 5.22 & 10.48 & 5.04 & 3.16 & 4.84 & 5.08 & 4.10 & 4.30 \\

        & $n$ & $50$ & 5.26 & 6.98 & 0.96 & 4.04 & 0.88 & 1.54 & 12.16 & 2.06 & 3.40 \\
        &     & $100$ & 4.94 & 11.76 & 0.02 & 4.62 & 2.10 & 2.18 & 8.08 & 2.36 & 3.30 \\
        &     & $200$ & 5.40 & 20.32 & 0.00 & 4.60 & 3.04 & 2.72 & 6.12 & 3.36 & 3.82 \\

        & $n^2$ & $50$ & 5.02 & 17.70 & NA & NA & 0.98 & 0.62 & 22.86 & NA & NA \\
        &       & $100$ & 5.06 & 26.24 & NA & NA & 1.46 & 0.72 & NA & NA & NA \\
        &       & $200$ & 5.20 & 32.82 & NA & NA & NA & NA & NA & NA & NA \\

Model 3 & 5 & $50$ & 5.58 & 4.44 & 1.78 & 4.94 & 2.04 & 3.88 & 6.82 & 3.82 & 4.34 \\
        &   & $100$ & 5.96 & 5.84 & 2.82 & 5.22 & 3.00 & 5.28 & 5.94 & 4.40 & 4.52 \\
        &   & $200$ & 4.94 & 5.46 & 6.42 & 5.00 & 2.84 & 5.96 & 4.76 & 3.80 & 3.98 \\

        & $\lfloor n^{1/2} \rfloor$ & $50$ & 4.90 & 4.34 & 1.02 & 4.80 & 1.86 & 3.46 & 7.84 & 3.68 & 4.46 \\
        &                            & $100$ & 5.26 & 5.64 & 2.38 & 5.54 & 2.90 & 4.34 & 6.06 & 4.18 & 4.46 \\
        &                            & $200$ & 4.66 & 5.78 & 6.56 & 5.14 & 3.18 & 4.30 & 5.00 & 3.90 & 3.94 \\

        & $n$ & $50$ & 5.52 & 11.56 & 1.86 & 3.30 & 1.22 & 1.62 & 11.56 & 1.54 & 3.02 \\
        &     & $100$ & 5.10 & 17.86 & 0.08 & 4.00 & 2.02 & 2.12 & 8.22 & 2.54 & 3.44 \\
        &     & $200$ & 4.86 & 25.78 & 0.00 & 4.24 & 2.62 & 2.94 & 5.90 & 2.74 & 3.38 \\

        & $n^2$ & $50$ & 5.02 & 19.82 & NA & NA & 1.06 & 0.58 & 22.08 & NA & NA \\
        &       & $100$ & 5.40 & 27.14 & NA & NA & 1.78 & 0.90 & NA & NA & NA \\
        &       & $200$ & 4.74 & 35.40 & NA & NA & NA & NA & NA & NA & NA \\

\hline
\end{tabular}}
\end{threeparttable}
\end{table}

\begin{table}[!htbp]
\centering
\footnotesize
\caption{\label{suptable_size3} {Empirical sizes $(\%)$ of the proposed  $L_2$-type test with $q=8$ and the competing methods. `NA' indicates that the results are omitted due to long computation time.}}
\begin{threeparttable}
\resizebox{1\textwidth}{!}{
\def\arraystretch{1.15}
\begin{tabular}{cccccccccccc}
\hline
& $p$ & $n$ & Proposed & LLYY & $\rm{CYZ}$ & $\rm{WKX}$ & $\rm{Tsay}$ & $\mathrm{FLM}_{\mathrm{FC}}$ & $\mathrm{CSF}_{L_D}$ & $\mathrm{CSF}_{L_R}$ & $\mathrm{CSF}_{L_{\tau^*}}$ \\
\hline

Model 1 & 5 & $50$ & 4.90 & 4.20 & 2.64 & 5.28 & 2.26 & 3.38 & 7.26 & 3.76 & 4.30 \\
        &   & $100$ & 5.06 & 5.30 & 4.64 & 4.96 & 3.00 & 4.86 & 5.94 & 4.16 & 4.44 \\
        &   & $200$ & 4.72 & 5.28 & 8.16 & 4.92 & 3.46 & 5.12 & 4.84 & 4.20 & 4.34 \\

        & $\lfloor n^{1/2} \rfloor$ & $50$ & 4.90 & 4.24 & 2.20 & 4.96 & 1.90 & 3.20 & 8.74 & 4.00 & 4.92 \\
        &                            & $100$ & 4.28 & 4.14 & 5.24 & 5.04 & 2.50 & 3.70 & 6.66 & 4.28 & 4.82 \\
        &                            & $200$ & 5.38 & 5.34 & 12.46 & 4.94 & 3.26 & 4.36 & 5.06 & 3.72 & 4.02 \\

        & $n$ & $50$ & 5.40 & 4.32 & 1.06 & 4.94 & 0.62 & 1.20 & 12.98 & 1.58 & 2.90 \\
        &     & $100$ & 4.98 & 4.66 & 0.00 & 5.38 & 1.66 & 2.28 & 8.58 & 2.44 & 3.44 \\
        &     & $200$ & 5.30 & 5.30 & 0.00 & 7.12 & 2.34 & 3.60 & 6.68 & 3.58 & 4.10 \\

        & $n^2$ & $50$ & 4.92 & 4.22 & NA & NA & 0.62 & 0.84 & 24.44 & NA & NA \\
        &       & $100$ & 4.60 & 3.84 & NA & NA & 1.84 & 0.76 & NA & NA & NA \\
        &       & $200$ & 4.58 & 4.54 & NA & NA & NA & NA & NA & NA & NA \\

Model 2 & 5 & $50$ & 4.96 & 4.42 & 2.18 & 5.02 & 2.76 & 3.08 & 7.66 & 3.94 & 4.46 \\
        &   & $100$ & 5.14 & 5.16 & 3.64 & 4.52 & 2.42 & 4.62 & 6.10 & 4.56 & 4.84 \\
        &   & $200$ & 5.50 & 5.42 & 6.68 & 5.62 & 3.08 & 5.50 & 5.08 & 3.90 & 4.24 \\

        & $\lfloor n^{1/2} \rfloor$ & $50$ & 4.92 & 4.28 & 1.48 & 4.20 & 1.80 & 2.22 & 7.52 & 3.42 & 4.20 \\
        &                            & $100$ & 5.36 & 4.88 & 4.48 & 5.12 & 2.66 & 3.32 & 6.74 & 4.18 & 4.82 \\
        &                            & $200$ & 5.00 & 5.32 & 8.24 & 5.36 & 3.24 & 4.38 & 5.12 & 3.98 & 4.20 \\

        & $n$ & $50$ & 5.18 & 5.54 & 1.20 & 3.60 & 0.76 & 1.08 & 12.88 & 1.76 & 3.32 \\
        &     & $100$ & 5.28 & 11.34 & 0.04 & 4.46 & 1.70 & 1.70 & 8.30 & 2.64 & 3.44 \\
        &     & $200$ & 4.62 & 18.44 & 0.00 & 4.54 & 2.44 & 2.38 & 6.42 & 3.52 & 4.04 \\

        & $n^2$ & $50$ & 4.54 & 15.56 & NA & NA & 0.88 & 0.60 & 23.96 & NA & NA \\
        &       & $100$ & 4.98 & 24.28 & NA & NA & 1.84 & 0.96 & NA & NA & NA \\
        &       & $200$ & 4.96 & 31.74 & NA & NA & NA & NA & NA & NA & NA \\

Model 3 & 5 & $50$ & 5.76 & 4.48 & 1.16 & 5.14 & 2.14 & 3.40 & 7.60 & 3.88 & 4.72 \\
        &   & $100$ & 5.02 & 5.14 & 2.14 & 5.04 & 3.06 & 4.86 & 6.16 & 4.52 & 4.82 \\
        &   & $200$ & 5.22 & 5.36 & 5.04 & 5.12 & 3.28 & 5.46 & 4.74 & 3.98 & 3.94 \\

        & $\lfloor n^{1/2} \rfloor$ & $50$ & 4.42 & 4.28 & 0.74 & 4.72 & 1.76 & 2.62 & 8.36 & 3.56 & 4.26 \\
        &                            & $100$ & 5.10 & 5.20 & 1.58 & 5.10 & 2.50 & 4.02 & 6.46 & 4.04 & 4.40 \\
        &                            & $200$ & 4.82 & 5.78 & 5.10 & 5.08 & 3.04 & 4.30 & 5.20 & 4.12 & 4.18 \\

        & $n$ & $50$ & 5.28 & 9.18 & 3.12 & 3.16 & 0.72 & 1.40 & 11.62 & 1.28 & 2.70 \\
        &     & $100$ & 5.08 & 16.30 & 0.08 & 4.08 & 1.72 & 2.06 & 8.22 & 2.76 & 3.38 \\
        &     & $200$ & 4.86 & 24.42 & 0.00 & 4.34 & 2.66 & 2.78 & 6.06 & 3.26 & 3.76 \\

        & $n^2$ & $50$ & 4.60 & 16.36 & NA & NA & 0.76 & 0.56 & 23.28 & NA & NA \\
        &       & $100$ & 5.22 & 25.36 & NA & NA & 1.54 & 0.94 & NA & NA & NA \\
        &       & $200$ & 4.86 & 34.66 & NA & NA & NA & NA & NA & NA & NA \\

\hline
\end{tabular}}
\end{threeparttable}
\end{table}

\begin{table}[!htbp]
\centering
\footnotesize
\caption{\label{suptable_power1} {Empirical powers $(\%)$ of the proposed  $L_2$-type test with $q=4$ and the competing methods. ‘NA’ indicates that the results are omitted due to long computation time. }}
\begin{threeparttable}
\resizebox{1\textwidth}{!}{
\def\arraystretch{1.15}
\begin{tabular}{cccccccccccc}
\hline
& $p$ & $n$ & Proposed & LLYY & $\rm{CYZ}$ & $\rm{WKX}$ & $\rm{Tsay}$ & $\mathrm{FLM}_{\mathrm{FC}}$ & $\mathrm{CSF}_{L_D}$ & $\mathrm{CSF}_{L_R}$ & $\mathrm{CSF}_{L_{\tau^*}}$ \\
\hline

Model 4 & 5 & $50$ & 10.36 & 9.60 & 4.82 & 7.38 & 2.76 & 9.02 & 8.30 & 4.80 & 5.42 \\
        &   & $100$ & 13.20 & 13.96 & 8.64 & 10.04 & 4.18 & 14.34 & 7.28 & 6.22 & 6.30 \\
        &   & $200$ & 21.46 & 22.96 & 19.30 & 15.10 & 9.44 & 25.02 & 9.76 & 9.12 & 9.06 \\

        & $\lfloor n^{1/2} \rfloor$ & $50$ & 11.26 & 10.80 & 3.58 & 7.52 & 2.34 & 9.24 & 8.70 & 4.68 & 5.28 \\
        &                            & $100$ & 16.70 & 16.80 & 10.48 & 10.84 & 4.32 & 15.02 & 7.60 & 5.28 & 5.56 \\
        &                            & $200$ & 29.68 & 29.52 & 30.04 & 18.02 & 7.74 & 28.56 & 9.56 & 7.64 & 8.18 \\

        & $n$ & $50$ & 43.66 & 39.84 & 0.42 & 36.20 & 1.40 & 27.24 & 13.04 & 2.38 & 4.08 \\
        &     & $100$ & 88.12 & 86.26 & 0.02 & 83.74 & 4.16 & 78.22 & 10.00 & 3.68 & 4.80 \\
        &     & $200$ & 100.00 & 100.00 & 0.00 & 99.98 & 25.24 & 99.98 & 9.24 & 5.70 & 6.16 \\

        & $n^2$ & $50$ & 100.00 & 100.00 & NA & NA & 3.62 & 100.00 & 27.30 & NA & NA \\
        &       & $100$ & 100.00 & 100.00 & NA & NA & 15.60 & 100.00 & NA & NA & NA \\
        &       & $200$ & 100.00 & 100.00 & NA & NA & NA & NA & NA & NA & NA \\

Model 5 & 5 & $50$ & 55.66 & 57.58 & 27.30 & 22.44 & 28.50 & 75.44 & 64.64 & 58.20 & 60.52 \\
        &   & $100$ & 85.58 & 86.08 & 78.66 & 40.30 & 74.56 & 92.34 & 89.62 & 88.78 & 89.24 \\
        &   & $200$ & 93.80 & 93.86 & 94.32 & 64.34 & 91.46 & 95.94 & 95.54 & 95.36 & 95.42 \\

        & $\lfloor n^{1/2} \rfloor$ & $50$ & 60.78 & 61.76 & 15.76 & 20.20 & 22.48 & 80.10 & 71.60 & 63.44 & 66.44 \\
        &                            & $100$ & 94.08 & 94.16 & 59.98 & 36.10 & 64.74 & 99.12 & 98.52 & 97.90 & 98.16 \\
        &                            & $200$ & 99.92 & 99.92 & 97.58 & 60.30 & 97.88 & 99.98 & 99.98 & 99.98 & 99.98 \\

        & $n$ & $50$ & 75.78 & 70.32 & 0.38 & 26.44 & 1.60 & 96.16 & 95.92 & 81.70 & 88.58 \\
        &     & $100$ & 99.86 & 97.90 & 0.00 & 49.50 & 5.98 & 100.00 & 100.00 & 100.00 & 100.00 \\
        &     & $200$ & 100.00 & 100.00 & 0.00 & 76.66 & 13.22 & 100.00 & 100.00 & 100.00 & 100.00 \\

        & $n^2$ & $50$ & 77.74 & 50.84 & NA & NA & 3.12 & 80.92 & 99.98 & NA & NA \\
        &       & $100$ & 99.92 & 70.64 & NA & NA & 3.30 & 100.00 & NA & NA & NA \\
        &       & $200$ & 100.00 & 87.20 & NA & NA & NA & NA & NA & NA & NA \\

Model 6 & 5 & $50$ & 7.72 & 6.70 & 2.14 & 5.60 & 2.36 & 6.04 & 7.40 & 4.50 & 4.90 \\
        &   & $100$ & 9.00 & 9.56 & 4.10 & 7.52 & 3.90 & 8.82 & 6.40 & 4.98 & 5.26 \\
        &   & $200$ & 11.52 & 12.04 & 9.30 & 9.02 & 5.52 & 13.24 & 7.22 & 6.52 & 6.52 \\

        & $\lfloor n^{1/2} \rfloor$ & $50$ & 7.88 & 7.14 & 1.52 & 5.98 & 2.36 & 6.14 & 8.30 & 4.52 & 4.94 \\
        &                            & $100$ & 9.68 & 10.66 & 3.44 & 7.52 & 3.32 & 8.52 & 7.12 & 4.84 & 5.38 \\
        &                            & $200$ & 13.44 & 14.58 & 10.10 & 9.96 & 5.60 & 12.30 & 7.28 & 6.12 & 6.38 \\  

        & $n$ & $50$ & 19.30 & 18.52 & 0.92 & 12.48 & 1.44 & 8.88 & 11.76 & 2.22 & 3.68 \\
        &     & $100$ & 42.56 & 31.52 & 0.02 & 34.58 & 1.96 & 26.92 & 9.32 & 3.10 & 3.98 \\
        &     & $200$ & 87.84 & 45.98 & 0.00 & 82.18 & 5.00 & 78.14 & 7.18 & 3.52 & 4.42 \\

        & $n^2$ & $50$ & 100.00 & 29.86 & NA & NA & 1.10 & 100.00 & 23.24 & NA & NA \\
        &       & $100$ & 100.00 & 41.58 & NA & NA & 3.42 & 100.00 & NA & NA & NA \\
        &       & $200$ & 100.00 & 55.56 & NA & NA & NA & NA & NA & NA & NA \\

\hline
\end{tabular}}
\end{threeparttable}
\end{table}

\begin{table}[!htbp]
\centering
\footnotesize
\caption{\label{suptable_power2} {Empirical powers $(\%)$ of the proposed  $L_2$-type test with $q=6$ and the competing methods. ‘NA’ indicates that the results are omitted due to long computation time. }}
\begin{threeparttable}
\resizebox{1\textwidth}{!}{
\def\arraystretch{1.15}
\begin{tabular}{cccccccccccc}
\hline
& $p$ & $n$ & Proposed & LLYY & $\rm{CYZ}$ & $\rm{WKX}$ & $\rm{Tsay}$ & $\mathrm{FLM}_{\mathrm{FC}}$ & $\mathrm{CSF}_{L_D}$ & $\mathrm{CSF}_{L_R}$ & $\mathrm{CSF}_{L_{\tau^*}}$ \\
\hline

Model 4 & 5 & $50$ & 9.76 & 9.48 & 3.90 & 7.44 & 2.66 & 7.82 & 8.60 & 4.32 & 5.04 \\
        &   & $100$ & 12.62 & 12.74 & 7.30 & 9.32 & 3.88 & 12.60 & 7.38 & 5.94 & 6.24 \\
        &   & $200$ & 19.50 & 20.58 & 16.14 & 13.86 & 8.14 & 21.94 & 8.60 & 8.08 & 8.02 \\

        & $\lfloor n^{1/2} \rfloor$ & $50$ & 11.78 & 11.16 & 2.76 & 8.24 & 2.28 & 8.40 & 8.62 & 4.44 & 5.24 \\
        &                            & $100$ & 16.66 & 15.92 & 8.36 & 11.04 & 3.72 & 13.94 & 7.06 & 4.78 & 5.14 \\
        &                            & $200$ & 28.46 & 28.14 & 25.84 & 18.74 & 6.78 & 25.90 & 8.80 & 7.12 & 7.26 \\

        & $n$ & $50$ & 51.60 & 47.10 & 0.64 & 44.10 & 0.98 & 33.28 & 14.12 & 2.02 & 3.86 \\
        &     & $100$ & 94.60 & 93.38 & 0.00 & 91.92 & 3.34 & 87.96 & 10.04 & 3.44 & 4.50 \\
        &     & $200$ & 100.00 & 100.00 & 0.00 & 100.00 & 21.48 & 100.00 & 8.50 & 5.12 & 5.82 \\

        & $n^2$ & $50$ & 100.00 & 100.00 & NA & NA & 3.70 & 100.00 & 28.90 & NA & NA \\
        &       & $100$ & 100.00 & 100.00 & NA & NA & 21.36 & 100.00 & NA & NA & NA \\
        &       & $200$ & 100.00 & 100.00 & NA & NA & NA & NA & NA & NA & NA \\

Model 5 & 5 & $50$ & 46.54 & 47.92 & 21.66 & 19.52 & 23.76 & 70.80 & 61.76 & 53.56 & 56.48 \\
        &   & $100$ & 81.38 & 81.90 & 74.60 & 32.64 & 71.78 & 91.36 & 88.62 & 87.50 & 87.94 \\
        &   & $200$ & 92.80 & 92.92 & 93.58 & 53.84 & 90.76 & 95.80 & 95.34 & 95.14 & 95.26 \\

        & $\lfloor n^{1/2} \rfloor$ & $50$ & 50.76 & 50.98 & 11.08 & 17.44 & 18.48 & 75.74 & 68.48 & 58.46 & 61.92 \\
        &                            & $100$ & 89.56 & 89.78 & 52.84 & 30.30 & 60.64 & 98.88 & 98.12 & 97.44 & 97.70 \\
        &                            & $200$ & 99.80 & 99.80 & 96.76 & 49.48 & 96.98 & 99.98 & 99.98 & 99.98 & 99.98 \\

        & $n$ & $50$ & 62.98 & 55.36 & 0.90 & 20.82 & 1.36 & 92.14 & 94.18 & 75.58 & 84.58 \\
        &     & $100$ & 97.90 & 92.14 & 0.04 & 37.58 & 5.36 & 100.00 & 100.00 & 100.00 & 100.00 \\
        &     & $200$ & 100.00 & 99.98 & 0.00 & 59.52 & 11.40 & 100.00 & 100.00 & 100.00 & 100.00 \\

        & $n^2$ & $50$ & 65.50 & 40.94 & NA & NA & 4.18 & 64.42 & 99.98 & NA & NA \\
        &       & $100$ & 98.40 & 59.90 & NA & NA & 8.40 & 100.00 & NA & NA & NA \\
        &       & $200$ & 100.00 & 78.18 & NA & NA & NA & NA & NA & NA & NA \\

Model 6 & 5 & $50$ & 7.88 & 6.20 & 1.30 & 5.90 & 2.24 & 5.50 & 7.68 & 4.36 & 4.74 \\
        &   & $100$ & 9.10 & 8.86 & 3.24 & 7.62 & 3.76 & 8.14 & 6.42 & 4.68 & 4.92 \\
        &   & $200$ & 11.14 & 11.44 & 7.66 & 9.00 & 5.48 & 11.98 & 7.04 & 6.02 & 6.34 \\

        & $\lfloor n^{1/2} \rfloor$ & $50$ & 7.68 & 6.48 & 1.06 & 6.10 & 1.82 & 5.44 & 8.28 & 3.88 & 4.60 \\
        &                            & $100$ & 10.38 & 9.96 & 2.22 & 7.80 & 3.44 & 8.24 & 7.32 & 4.60 & 4.90 \\
        &                            & $200$ & 13.86 & 14.56 & 7.98 & 10.28 & 4.98 & 11.38 & 6.76 & 5.46 & 5.62 \\

        & $n$ & $50$ & 23.24 & 16.52 & 1.88 & 14.80 & 1.10 & 9.88 & 12.68 & 1.72 & 3.50 \\
        &     & $100$ & 51.38 & 30.10 & 0.06 & 42.68 & 2.18 & 34.06 & 8.82 & 3.00 & 3.92 \\
        &     & $200$ & 94.58 & 47.06 & 0.00 & 91.54 & 4.00 & 88.38 & 7.42 & 4.00 & 4.40 \\

        & $n^2$ & $50$ & 100.00 & 27.38 & NA & NA & 1.24 & 100.00 & 24.02 & NA & NA \\
        &       & $100$ & 100.00 & 41.06 & NA & NA & 3.10 & 100.00 & NA & NA & NA \\
        &       & $200$ & 100.00 & 56.78 & NA & NA & NA & NA & NA & NA & NA \\

\hline
\end{tabular}}
\end{threeparttable}
\end{table}

\begin{table}[!htbp]
\centering
\footnotesize
\caption{\label{suptable_power3} {Empirical powers $(\%)$ of the proposed  $L_2$-type test with $q=8$ and the competing methods. ‘NA’ indicates that the results are omitted due to long computation time. }}
\begin{threeparttable}
\resizebox{1\textwidth}{!}{
\def\arraystretch{1.15}
\begin{tabular}{cccccccccccc}
\hline
& $p$ & $n$ & Proposed & LLYY & $\rm{CYZ}$ & $\rm{WKX}$ & $\rm{Tsay}$ & $\mathrm{FLM}_{\mathrm{FC}}$ & $\mathrm{CSF}_{L_D}$ & $\mathrm{CSF}_{L_R}$ & $\mathrm{CSF}_{L_{\tau^*}}$ \\
\hline

Model 4 & 5 & $50$ & 10.80 & 9.92 & 2.80 & 7.68 & 2.38 & 7.50 & 8.78 & 4.40 & 5.24 \\
        &   & $100$ & 12.04 & 12.48 & 5.94 & 9.30 & 3.68 & 11.70 & 7.22 & 5.54 & 6.06 \\
        &   & $200$ & 17.48 & 18.90 & 13.96 & 13.62 & 7.86 & 19.92 & 8.38 & 7.48 & 7.72 \\

        & $\lfloor n^{1/2} \rfloor$ & $50$ & 12.10 & 11.18 & 1.92 & 8.44 & 2.36 & 7.56 & 9.20 & 4.30 & 5.08 \\
        &                            & $100$ & 17.46 & 16.60 & 6.68 & 11.84 & 3.84 & 13.30 & 7.18 & 4.76 & 5.06 \\
        &                            & $200$ & 27.58 & 27.58 & 22.12 & 18.86 & 6.14 & 25.10 & 8.56 & 6.74 & 6.98 \\

        & $n$ & $50$ & 58.12 & 53.32 & 1.18 & 48.78 & 1.20 & 37.16 & 15.40 & 1.88 & 3.78 \\
        &     & $100$ & 97.20 & 96.74 & 0.02 & 95.76 & 3.28 & 92.64 & 10.38 & 3.20 & 4.56 \\
        &     & $200$ & 100.00 & 100.00 & 0.00 & 100.00 & 19.34 & 100.00 & 8.48 & 5.14 & 5.74 \\

        & $n^2$ & $50$ & 100.00 & 100.00 & NA & NA & 4.78 & 100.00 & 29.40 & NA & NA \\
        &       & $100$ & 100.00 & 100.00 & NA & NA & 27.68 & 100.00 & NA & NA & NA \\
        &       & $200$ & 100.00 & 100.00 & NA & NA & NA & NA & NA & NA & NA \\

Model 5 & 5 & $50$ & 41.34 & 42.26 & 17.12 & 17.54 & 20.94 & 67.00 & 59.28 & 50.70 & 53.20 \\
        &   & $100$ & 76.64 & 77.20 & 71.12 & 27.98 & 69.76 & 90.64 & 87.78 & 86.52 & 86.86 \\
        &   & $200$ & 91.60 & 91.80 & 93.20 & 47.16 & 90.46 & 95.62 & 95.26 & 94.96 & 95.14 \\

        & $\lfloor n^{1/2} \rfloor$ & $50$ & 44.50 & 44.00 & 8.34 & 16.52 & 16.04 & 71.78 & 65.90 & 54.74 & 59.04 \\
        &                            & $100$ & 85.50 & 85.62 & 47.34 & 26.56 & 57.18 & 98.70 & 97.94 & 97.06 & 97.28 \\
        &                            & $200$ & 99.68 & 99.54 & 95.76 & 42.06 & 96.28 & 99.98 & 99.98 & 99.98 & 99.98 \\

        & $n$ & $50$ & 55.42 & 45.66 & 1.44 & 18.58 & 1.00 & 87.78 & 92.72 & 71.50 & 81.28 \\
        &     & $100$ & 95.12 & 85.44 & 0.10 & 31.36 & 3.92 & 100.00 & 100.00 & 100.00 & 100.00 \\
        &     & $200$ & 100.00 & 99.44 & 0.00 & 47.84 & 9.78 & 100.00 & 100.00 & 100.00 & 100.00 \\

        & $n^2$ & $50$ & 57.40 & 32.98 & NA & NA & 5.90 & 53.32 & 99.98 & NA & NA \\
        &       & $100$ & 95.78 & 53.46 & NA & NA & 17.36 & 100.00 & NA & NA & NA \\
        &       & $200$ & 100.00 & 71.26 & NA & NA & NA & NA & NA & NA & NA \\

Model 6 & 5 & $50$ & 8.02 & 6.08 & 1.10 & 6.02 & 2.28 & 5.30 & 7.88 & 4.24 & 4.94 \\
        &   & $100$ & 8.38 & 8.10 & 2.30 & 7.54 & 3.78 & 7.38 & 6.70 & 5.06 & 5.30 \\
        &   & $200$ & 10.88 & 10.94 & 5.76 & 8.62 & 5.38 & 10.94 & 6.60 & 5.60 & 5.94 \\

        & $\lfloor n^{1/2} \rfloor$ & $50$ & 7.90 & 5.92 & 1.10 & 5.80 & 1.72 & 4.46 & 8.48 & 3.76 & 4.60 \\
        &                            & $100$ & 10.36 & 9.56 & 1.60 & 7.70 & 3.10 & 7.60 & 6.88 & 4.52 & 5.00 \\
        &                            & $200$ & 13.66 & 13.42 & 6.24 & 10.38 & 4.78 & 10.74 & 6.72 & 5.46 & 5.60 \\

        & $n$ & $50$ & 26.02 & 14.02 & 3.00 & 16.12 & 0.64 & 10.72 & 13.32 & 1.32 & 3.02 \\
        &     & $100$ & 59.18 & 28.92 & 0.00 & 49.10 & 2.04 & 40.28 & 8.70 & 2.98 & 3.88 \\
        &     & $200$ & 97.92 & 48.38 & 0.00 & 95.86 & 3.90 & 94.10 & 7.42 & 3.78 & 4.72 \\

        & $n^2$ & $50$ & 100.00 & 24.92 & NA & NA & 1.34 & 100.00 & 25.74 & NA & NA \\
        &       & $100$ & 100.00 & 41.68 & NA & NA & 3.74 & 100.00 & NA & NA & NA \\
        &       & $200$ & 100.00 & 58.86 & NA & NA & NA & NA & NA & NA & NA \\

\hline
\end{tabular}}
\end{threeparttable}
\end{table}
\clearpage

\section{Real data analysis}\label{sec:realdata}

In this section, we apply the proposed white noise test to investigate changes in short-horizon serial dependence in S\&P 500 stock returns before, during, and after the global financial crisis. Our focus is on whether the test can capture both the emergence and subsequent weakening of serial dependence. To this end, we conduct a rolling-window analysis over the period from January 1, 2005 to December 31, 2010, which covers the global financial crisis from 2007 to 2009 as well as the preceding and subsequent periods. The daily total return data (in percentages) are obtained from the Center for Research in Security Prices through Wharton Research Data Services (WRDS)\footnote{The website of WRDS: \url{https://wrds-www.wharton.upenn.edu/}}. We consider 349 stocks that are constituents of the S\&P 500 index over the entire sample period and have complete daily return records.

\begin{figure}[!h]
    \centering
    % first row
    % \setlength{\fboxsep}{5pt}
 
    \includegraphics[width=1\textwidth]{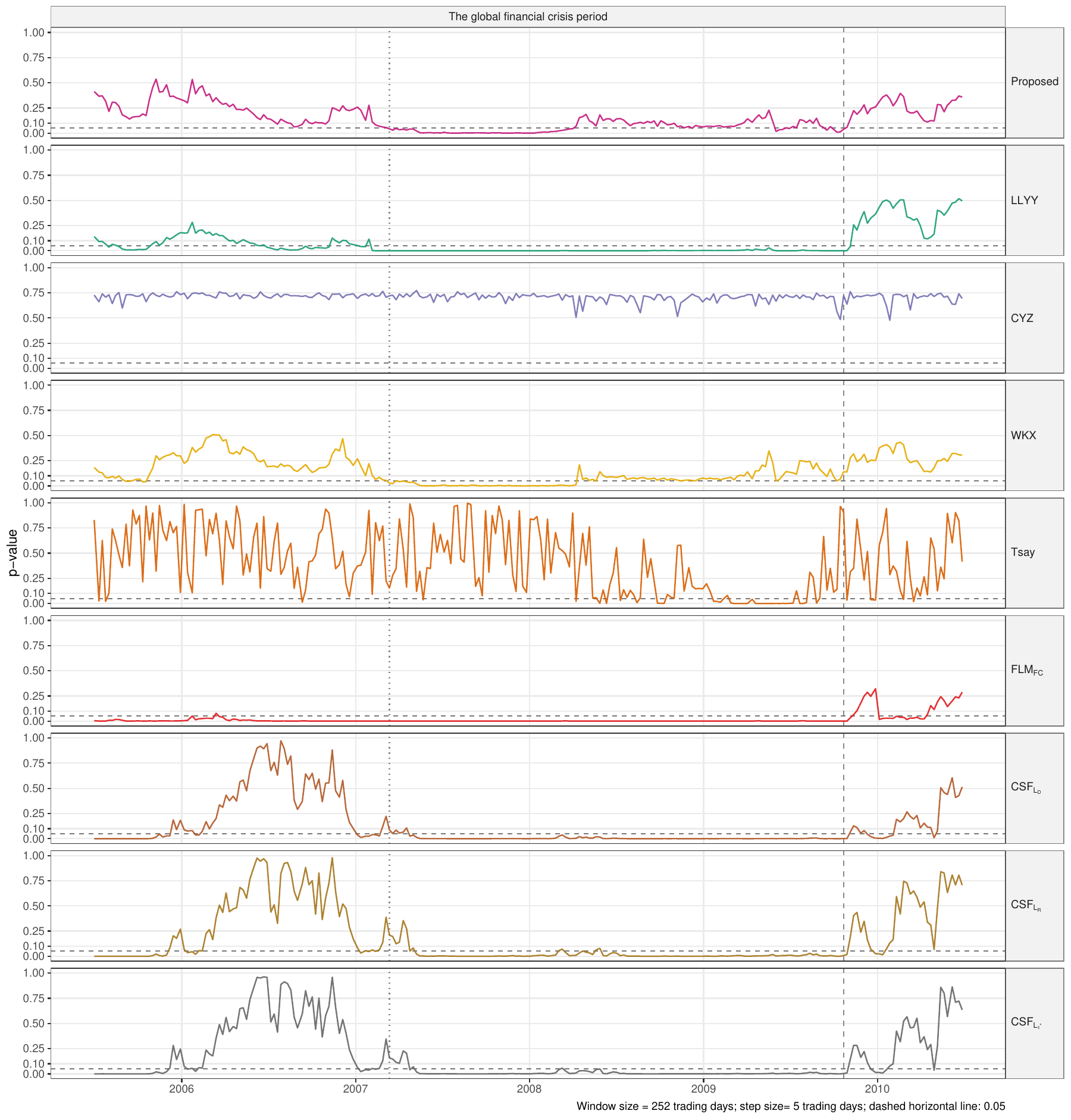}
    \caption{Rolling-window p-values of the white noise tests with $q=1$. Each rolling window contains 252 trading days and moves forward by five trading days at each step, with the resulting p-value plotted against the midpoint date of the window. Each row represents a different testing procedure. The dashed horizontal line indicates the 5\% significance level. The vertical dotted line marks the midpoint of the rolling window from September 11, 2006 to September 11, 2007, the first window in which our proposed test rejects the white noise hypothesis at the 5\% significance level. The vertical dashed line marks the midpoint of the rolling window from April 27, 2009 to April 26, 2010, the last window in which our proposed test rejects the white noise hypothesis at the 5\% significance level.
    }
    \label{fig:rolling}
\end{figure}

We apply all the tests considered in Section 4 to the demeaned stock returns within rolling windows of $n_w=252$ trading days, corresponding approximately to one trading year, and move the window forward by 5 trading days at each step. 
Figure~\ref{fig:rolling} presents the resulting p-values with $q=1$, plotted against the midpoint date of each window. From the perspective of market efficiency \citep{Fama1970,Fama1991}, stock returns are not predictable using past return information. However, during periods of financial crisis, disruptions to trading and information processing may give rise to strong short-horizon serial dependence. The white noise hypothesis is reasonable for stock returns in the periods before and after the financial crisis, while it is not reasonable for the period during the financial crisis. 
The proposed test exhibits a clear decline in p-values followed by a subsequent increase, suggesting that the evidence against the white noise hypothesis first strengthens and then weakens. More specifically, in the early part of the sample period, before the financial crisis, our proposed test does not reject the white noise hypothesis at the 5\% significance level. The first rejection occurs for the window spanning from September 11, 2006 to September 11, 2007, whose midpoint is indicated by the vertical dotted line in Figure~\ref{fig:rolling}. This window includes the onset of the financial crisis, which the Federal Reserve dates to August 2007. 
Toward the end of the sample period, the last rejection of the white noise hypothesis by our proposed test at the 5\% significance level occurs in the window spanning from April 27, 2009 to April 26, 2010, whose midpoint is indicated by the vertical dashed line in Figure~\ref{fig:rolling}.
Non-rejection persists in all subsequent windows, which cover a period of improving financial conditions and the beginning of the economic recovery: the Federal Reserve reported that financial markets began to show signs of improvement in March 2009, and the National Bureau of Economic Research dated the end of the U.S. recession to June 2009. These findings illustrate the ability of our proposed test to capture changes in the temporal dependence structure of stock returns across periods of financial crisis and recovery.
The results of the WKX test show a similar pattern to that of our proposed test. 
However, the LLYY, $\mathrm{FLM}_{\mathrm{FC}}$, $\mathrm{CSF}_{L_D}$, $\mathrm{CSF}_{L_R}$, and $\mathrm{CSF}_{L_{\tau^*}}$ tests reject the white noise hypothesis at the 5\% significance level in some rolling windows before the financial crisis, which is difficult to be interpreted based on the economic theory. 
The CYZ test produces consistently large p-values throughout the sample period, while the p-values of the Tsay test fluctuate substantially across adjacent windows. The results of the CYZ and Tsay tests provide less clear evidence of changes in serial dependence over time.

\section{Auxiliary lemmas}

In this section, we introduce some useful lemmas that will be used to prove the main results.

\begin{lemma}[Lemma B.26 of \cite{BSbook}] \label{ineq-main}
    Let $A$ be a $p\times p$ nonrandom matrix and $Z=(z_1,\ldots,z_p)^\top$ be a random vector of independent entries.
    Assume that $E(z_i)=0$, $E(z_i^2)=1$, and $E(|z_i|^\ell)\leq\nu_\ell$.
    Then, for any $k\geq 1$, we have
    \begin{align*}
    E\big\{|Z^\top A Z-\operatorname{tr}(A)|^{k} \big\}
    \leq
    C_k \big[ \nu_4^{k/2}\operatorname{tr}^{k/2}(AA^\top)+\nu_{2k}\operatorname{tr}\{(AA^\top)^{k/2}\}\big]\,,
    \end{align*}
    where $C_k$ is a constant depending only on $k$.

\end{lemma}

\begin{lemma}[Vector martingale central limit theorem]\label{lem:md-clt-vector}
Suppose for each $n\ge1$, $Y_{n,1},\ldots,Y_{n,r_n}$ is a $q$-dimensional real martingale difference sequence with respect to the increasing sequence of $\sigma$-fields $\{\mathcal F_{n,j}\}_{j=0}^{r_n}$, where $q\geq1$ is a fixed integer. Write $Y_{n,j}=(Y_{n,j,1},\ldots,Y_{n,j,q})^\top$. 
If, as $n\to\infty$,
\begin{align}
\sum_{j=1}^{r_n}
E(
Y_{n,j,a}Y_{n,j,b}\,|\, \mathcal F_{n,j-1})
\xrightarrow{{\rm p}}\sigma_{a, b}~\,
\mbox{for any}~\,a, b\in[q]\,, \label{eq:vector-cov}
\end{align}
where $\Sigma=(\sigma_{a, b})_{a,b\in[q]}$ is a positive definite matrix, and
\begin{align}
\sum_{j=1}^{r_n}E(Y_{n,j,a}^4)\to 0~\,
\mbox{for any}~\,a\in[q]\,, \label{eq:vector-4th}
\end{align}
then
\(
\sum_{j=1}^{r_n} Y_{n,j}\xrightarrow{{\rm d}}\mathcal N(0,\Sigma).
\)
\end{lemma}

\subsection{Proof of Lemma \ref{lem:md-clt-vector}}
Given a nonzero vector $t=(t_1,\ldots,t_q)^\top\in\mathbb R^q$, define
\[
Z_{n,j}=t^\top Y_{n,j}
=\sum_{a=1}^q t_aY_{n,j,a}\,,
\qquad j\in[r_n]\,.
\]
Then $Z_{n,1},\ldots,Z_{n,r_n}$ is a real martingale difference sequence with respect to the same increasing sequence of $\sigma$-fields $\{\mathcal F_{n,j}\}_{j=0}^{r_n}$.
By \eqref{eq:vector-cov}, it holds that
\begin{align*}
\sum_{j=1}^{r_n}E(Z_{n,j}^2\mid \mathcal F_{n,j-1})
&=
\sum_{a=1}^q\sum_{b=1}^q
t_at_b
\sum_{j=1}^{r_n}
E(
Y_{n,j,a}Y_{n,j,b}\,|\, \mathcal F_{n,j-1}) \\
&\xrightarrow{\rm p}
\sum_{a=1}^q\sum_{b=1}^q t_at_b\sigma_{a,b}
=t^\top\Sigma t>0\,.
\end{align*}
By Theorem 35.12 in \citet{bili95}, in order to show $ t^\top(\sum_{j=1}^{r_n} Y_{n,j})\xrightarrow{{\rm d}}\mathcal{N}(0,t^\top\Sigma t)$, it suffices to verify the following Lindeberg condition:
\begin{align}
\sum_{j=1}^{r_n}E\{
Z_{n,j}^2 I(|Z_{n,j}|\ge \varepsilon)
\}\to 0 \label{eq:scalar-lind}
\end{align}
 for any $\varepsilon>0$. Since
$
|Z_{n,j}|=|t^\top Y_{n,j}|
\le | t|_2 | Y_{n,j}|_2$, 
we have
\[
E\{
Z_{n,j}^2 I(|Z_{n,j}|\ge\varepsilon)
\}
\le
\varepsilon^{-2}E(Z_{n,j}^4)
\le
 \varepsilon^{-2}|t|_2^4\E (| Y_{n,j}|_2^4)
\le  \varepsilon^{-2}C\sum_{a=1}^q E (Y_{n,j,a}^4)\,,
\]
for some constant $C>0$. By \eqref{eq:vector-4th}, we have \eqref{eq:scalar-lind} holds.
Therefore, we have  
$
t^\top(\sum_{j=1}^{r_n} Y_{n,j})
\xrightarrow{\rm d}
\mathcal N(0, t^\top\Sigma t)
$ for any nonzero $ t\in\mathbb R^q$. This implies  
$
\sum_{j=1}^{r_n} Y_{n,j}\xrightarrow{\rm d}\mathcal N( 0,\Sigma)$. $\hfill\Box$

\section{Proof of Proposition \ref{null-lemma}}\label{sec:T1}

Define $r_t=p^{-1/2}X_t=p^{-1/2}Z_t$ for $t\in [n]$, and set $r_t=r_{t-n}$ for $t>n$. Recall $n > Cq$ for some sufficiently large constant $C>0$.
Due to $\max_{i\in[p]}E(z_{i,t}^{16})=O(1)$, for any matrix $A\in \mathbb R^{p\times p}$, we have
\begin{align}\label{moments}
E\big\{(r_t^\top A r_t)^k \big\}\leq \|A\|_2^k\cdot E\bigg\{\bigg(\frac1p\sum_{i=1}^pz_{i,t}^2\bigg)^k\bigg\}\leq \|A\|_2^k\cdot\max_{i\in [p]}E(z_{i,t}^{2k})=O\big(\|A\|_2^k \big)
\end{align}
for each $k\in[8]$.
The statistics $G_{n,\tau,1}$ and $G_{n,\tau,2}$ defined in \eqref{eq:Gnt1} can be written as
\begin{align*}
  G_{n,\tau,1}
=\frac {2p}{n}\sum_{t < s}r_t^\top r_s r_{t+\tau}^\top r_{s+\tau}~~\mbox{and}~~
  G_{n,\tau,2}
=\frac p{n}\sum_{t=1}^n r_t^\top r_t r_{t+\tau}^\top r_{t+\tau}- \frac{p}{n^2}\sum_{t,\ell=1}^n r_t^\top r_t r_\ell^\top r_\ell\,.
\end{align*}
Our proof includes the following three steps:

	{\bf Step 1.} As we will show in Section \ref{subsec:Gq1}, as $n \to \infty$, it holds that 
\begin{align}\label{eq:AsymNormGq1}
    \left(
\frac{G_{n,1,1}}{\sigma_{n1}}, \ldots, \frac{G_{n,q,1}}{\sigma_{n1}}\right) \stackrel{\mathrm{d}}{\to} \mathcal{N}( 0,I_{q})\,.
\end{align}

 {\bf Step 2.}  As we will show in Section \ref{subsec:Gq2}, as $n \to \infty$, it holds that 
\begin{align}\label{eq:AsymNormGq2}
    \left(\frac{G_{n,1,2}-\mu_2}{\sigma_{n2}}, \ldots, \frac{G_{n,q,2}-\mu_2}{\sigma_{n2}}\right) \stackrel{\mathrm{d}}{\to} \mathcal{N}( 0,I_{q})\,.
\end{align}

{\bf Step 3.} As we will show in Section \ref{subsec:Ind_Gq12}, as $n \to \infty$, it holds that 
\begin{align}\label{eq:JointAsymNorm}
\left(
\frac{G_{n,1,1}}{\sigma_{n1}}, \ldots, \frac{G_{n,q,1}}{\sigma_{n1}}, 
\frac{G_{n,1,2}-\mu_2}{\sigma_{n2}}, \ldots, \frac{G_{n,q,2}-\mu_2}{\sigma_{n2}}\right) \stackrel{\mathrm{d}}{\to} \mathcal{N}( 0, I_{2q})\,.
\end{align}
Then we have Proposition \ref{null-lemma}.  $\hfill\Box$

\subsection{Proof of Step 1}\label{subsec:Gq1}
We decompose $G_{n,\tau,1}$ into four terms such that
$
G_{n,\tau,1}=	G_{n,\tau,1,1}+  G_{n,\tau,1,2}+  G_{n,\tau,1,3} +  G_{n,\tau,1,4}$, 
where
\begin{align*}
&~~~~~G_{n,\tau,1,1}
=	\frac {2p}{n}\sum_{ t< s\neq t+\tau \atop t,s\in [n-\tau] }r_t^\top r_s r_{t+\tau}^\top r_{s+\tau}\,,
\quad~
G_{n,\tau,1,2}
=	\frac {2p}{n}\sum_{ t< s= t+\tau}r_t^\top r_s r_{t+\tau}^\top r_{s+\tau}\,,\\
&G_{n,\tau,1,3}
=	\frac {2p}{n}\sum_{ t< s\neq t+\tau \atop t\in [n-\tau],\, s\in\{n-\tau+1,\ldots,n\}}r_t^\top r_s r_{t+\tau}^\top r_{s+\tau}\,,
\quad
G_{n,\tau,1,4}
=	\frac {2p}{n}\sum_{ t< s\neq t+\tau \atop t, s\in\{n-\tau+1,\ldots,n\}}r_t^\top r_s r_{t+\tau}^\top r_{s+\tau}\,.
\end{align*}
Notice that, for $t\neq s\neq \ell\in [n]$, we have 
$E \{(r_t^\top r_{s} r_{s}^\top r_{\ell})^2\}=p^{-2}E\{ (r_{s}^\top \Sigma r_{s})^2\}$. 
Hence, under Assumptions \ref{as:A2} and \ref{as:A4}, by \eqref{moments} with $k=2$,  we have
\begin{align}\label{eq:Gnt12}
E (G_{n,\tau,1,2}^2)
= \frac{4p^2}{n^{2}}\sum_{t=1}^{n-\tau}E \{(r_t^\top r_{t+\tau} r_{t+\tau}^\top r_{t+2\tau})^2 \}=O(n^{-1})\,.
\end{align}
Meanwhile, for the third and fourth terms,
\begin{align}
E( G_{n,\tau,1,3}^2) &= \frac {4p^2}{n^2}\sum_{ t< s\neq t+\tau \atop t\in [n-\tau],\, s\in\{n-\tau+1,\ldots,n\}}E\{(r_t^\top r_s r_{t+\tau}^\top r_{s+\tau})^2\}\notag\\
&\leq 4\tau p^2n^{-1}E\{(r_1^\top r_2 r_{3}^\top r_{4})^2\} + 4\tau p^2 n^{-2} E\{(r_1^\top r_2 r_{2}^\top r_{3})^2\} \notag\\
& = 4\tau n^{-1}\alpha_{2p}^2 + 4\tau n^{-2}E\{(r_1^\top \Sigma r_{1})^2\} = O(n^{-1})\,,\label{eq:Gnt13}\\
E( G_{n,\tau,1,4}^2) &= \frac {4p^2}{n^2}\sum_{ t< s\neq t+\tau \atop t, s\in\{n-\tau+1,\ldots,n\}}E\{(r_t^\top r_s r_{t+\tau}^\top r_{s+\tau})^2\}\notag\\	
&\leq 2\tau^2p^2n^{-2}E\{(r_1^\top r_2 r_{3}^\top r_{4})^2 \} = 2\tau^2n^{-2}\alpha_{2p}^2 = O(n^{-2})\label{eq:Gnt14}\,,
\end{align}
which implies $G_{n,\tau,1}=G_{n,\tau,1,1}+O_{\rm p}(n^{-1/2})$. Elementary calculations can reveal that
\begin{align*}
    E (G_{n,\tau,1,1})=0 ~~\mbox{and}~~E (G_{n,\tau,1,1}^2)=2\alpha_{2p}^2\{1+O(n^{-1})\}\,.
\end{align*}

Let $E_0(\cdot)$ denote expectation, and $E_j(\cdot)$ denote conditional expectation with respect to the $\sigma$-field $\mathcal F_j=\sigma(X_1, \ldots,X_j)$ for $j \ge 1$. By convention, we set $r_\ell=0$ for $\ell\leq 0$. Applying the martingale decomposition, we obtain
\begin{align}
  G_{n,\tau,1,1}
=&~\frac {2p}{n}\sum_{j=1}^nE_j\Bigg( 
\sum_{t< s\neq t+\tau \atop t,s\in [n-\tau] }r_t^\top r_s r_{t+\tau}^\top r_{s+\tau}\Bigg) - \frac {2p}{n}\sum_{j=1}^n E_{j-1}\Bigg( \sum_{t< s\neq t+\tau \atop t,s\in [n-\tau] }r_t^\top r_s r_{t+\tau}^\top r_{s+\tau}\Bigg)	\nonumber\\
= &~\frac {2p}{n}\sum_{j=\tau+1}^n\sum_{j-\tau\neq t< j}r_{j-\tau}^\top r_{t-\tau}r_t^\top r_j =\sum_{j=\tau+1}^nD_{j,\tau}\,,\label{gnt11}
\end{align}
where $D_{j,\tau} = 2pn^{-1}\sum_{j-\tau\neq t< j}r_{j-\tau}^\top r_{t-\tau}r_t^\top r_j$.
Hence, $\{D_{j,\tau}\}$ forms a sequence of martingale differences with respect to $\{\mathcal F_j\}$. Recall $\sigma_{n1}^2=2\alpha_{2p}^2 \to 2\alpha_2^2$ as $p \to \infty$. 
As we will show in Sections \ref{subsubsec:c1}--\ref{subsubsec:c3}, 
\begin{align}
 \sum_{j=\tau+1}^n \frac{E_{j-1}(D_{j,\tau}^2)}{\sigma_{n1}^2} & \xrightarrow{\mathrm{p}}  1 ~\,\mbox{for any given}~\, \tau\in[q]\,, \label{mclt-conditiona1} \\ 
\sum_{j= \tau +1}^n \frac{E_{j-1}(D_{j,\tau}D_{j,\tau'})}{\sigma_{n1}^2} & \xrightarrow{\mathrm{p}} 0~\,\mbox{for any given}~\, (\tau,\tau')\in [q]^2 ~\mbox{satisfying}~\tau>\tau' \,, \label{mclt-conditiona2}\\
 \sum_{j=\tau+1}^n \frac{E (D_{j,\tau}^4)}{\sigma_{n1}^4} & \to0 ~\,\mbox{for any given}~\, \tau\in[q]\,.\label{mclt-conditiona3}
\end{align}
Notice that $D_{j,\tau}=0$ if $j \le \tau$. Hence, we have
\begin{align*}
    \left(
\frac{G_{n,1,1}}{\sigma_{n1}}, \ldots, \frac{G_{n,q,1}}{\sigma_{n1}}\right) = \sum_{j=2}^{n}\left(
\frac{D_{j,1}}{\sigma_{n1}}, \ldots, \frac{D_{j,q}}{\sigma_{n1}}\right) + O_{\mathrm{p}}(n^{-1/2})\,.
\end{align*}
Based on \eqref{mclt-conditiona1}--\eqref{mclt-conditiona3}, by Lemma \ref{lem:md-clt-vector}, we have \eqref{eq:AsymNormGq1} holds. $\hfill\Box$

\subsubsection{Proof of \eqref{mclt-conditiona1}}\label{subsubsec:c1}
Write 
\begin{align}
\sum_{j=\tau+1}^n E_{j-1}(D_{j,\tau}^2)
=&~\frac {4p^2}{n^2}\sum_{j=\tau+1}^n\sum_{j-\tau\neq t< j\atop j-\tau\neq s< j}r_{j-\tau}^\top r_{t-\tau}r_t^\top E_{j-1}	(r_j r_j^\top )r_s r_{s-\tau}^\top r_{j-\tau}\nonumber\\
=&~\frac {4p}{n^2}\sum_{j=\tau+1}^n	\sum_{j-\tau\neq t< j\atop j-\tau\neq s< j}r_{j-\tau}^\top r_{t-\tau}\cdot r_t^\top \Sigma r_s \cdot r_{s-\tau}^\top r_{j-\tau}\,.\label{lim-var}
\end{align}
Notice that the summands in \eqref{lim-var} have zero expectation for $t\neq s$. Then
\begin{align}
\frac{1}{\sigma_{n1}^2}E\left\{\sum_{j=\tau+1}^n E_{j-1}(D_{j,\tau}^2)\right\}
=&~\frac {4p}{n^2\sigma_{n1}^2}\sum_{j=\tau+1}^n	\sum_{j-\tau\neq t< j}E(r_t^\top \Sigma r_t \cdot r_{j-\tau}^\top r_{t-\tau}r_{t-\tau}^\top r_{j-\tau})\notag\\
=&~\frac {4p}{n^2\sigma_{n1}^2}\sum_{j=\tau+1}^n	\sum_{j-\tau\neq t< j}E(r_t^\top \Sigma r_t) \cdot E(r_{j-\tau}^\top r_{t-\tau}r_{t-\tau}^\top r_{j-\tau})\notag\\
=&~\frac {4p}{n^2\sigma_{n1}^2}\sum_{j=\tau+1}^n	\sum_{\tau<t< j \atop t\neq j-\tau}\tr\{E(r_t r_t^\top \Sigma)\} \cdot \tr\{E(r_{j-\tau}r_{j-\tau}^\top) E(r_{t-\tau}r_{t-\tau}^\top)\}\notag\\
=&~\frac {4}{n^2p^2\sigma_{n1}^2}\sum_{j=\tau+1}^n	\sum_{\tau<t< j \atop t\neq j-\tau}{\tr}^2(\Sigma^2)\notag\\
=&~\frac {4}{n^2\sigma_{n1}^2}\sum_{j=\tau+1}^n	\sum_{\tau<t< j \atop t\neq j-\tau}\alpha_{2p}^2\to 1\,. \label{eq:mean_conclusion1}
\end{align}
Moreover, it holds that
\begin{align}
\Var\left\{\sum_{j=\tau+1}^n E_{j-1}(D_{j,\tau}^2)\right\}
&\leq \frac {Kp^2}{n^4} E\Bigg[\bigg\{\sum_{j=\tau+1}^n	\sum_{\tau<t< j \atop t\neq j-\tau}\bigg(
r_{j-\tau}^\top r_{t-\tau}\cdot r_t^\top \Sigma r_t \cdot r_{t-\tau}^\top r_{j-\tau}-\frac1p\alpha_{2p}^2\bigg)\bigg\}^2\Bigg]\notag\\
&~~~+\frac {Kp^2}{n^4} E\Bigg\{\bigg(\sum_{j=\tau+1}^n	\sum_{{j-\tau\neq t< j\atop j-\tau\neq s< j}\atop t< s}r_{j-\tau}^\top r_{t-\tau}\cdot r_t^\top \Sigma r_s \cdot r_{s-\tau}^\top r_{j-\tau}\bigg)^2\Bigg\}\notag\\
&= : V_1+V_2\,.\label{eq:var_conclusion1}
\end{align}
As we will show in Section \ref{subsec:V1V2}, $V_1=o(1)$ and $V_2=o(1)$. Hence, based on \eqref{eq:mean_conclusion1} and \eqref{eq:var_conclusion1}, we have \eqref{mclt-conditiona1} holds.$\hfill\Box$

\subsubsection{Proof of $V_1=o(1)$ and $V_2=o(1)$}\label{subsec:V1V2}
For $V_1$, we have
\begin{align*}
V_1
= \frac {Kp^2}{n^4}\sum_{j_1, j_2=\tau+1}^n	\sum_{\tau<t_1< j_1 \atop t_1\neq j_1-\tau}\sum_{\tau<t_2< j_2 \atop t_2\neq j_2-\tau}\sigma(j_1,t_1,j_2,t_2)\,,
\end{align*}
where
\begin{align*}
\sigma(j_1,t_1,j_2,t_2)&= E\bigg\{\bigg(r_{j_1-\tau}^\top r_{t_1-\tau}r_{t_1}^\top \Sigma r_{t_1} r_{t_1-\tau}^\top r_{j_1-\tau}-\frac1p\alpha_{2p}^2\bigg)\\
&~~~~~~~~~~~~~~~\times\bigg(r_{j_2-\tau}^\top r_{t_2-\tau}r_{t_2}^\top \Sigma r_{t_2}r_{t_2-\tau}^\top r_{j_2-\tau}-\frac1p\alpha_{2p}^2\bigg) \bigg\}\,.
\end{align*}
If $\{j_1-\tau, t_1-\tau, t_1\} \cap \{j_2-\tau, t_2-\tau, t_2\}=\emptyset$, then $\sigma(j_1,t_1,j_2,t_2) = 0$. Same as \eqref{moments} with $k=2$, we have $\max_{t\in[n]}E\{(r_{t}^\top \Sigma r_{t})^2\}\leq K$ under Assumptions \ref{as:A2} and \ref{as:A4}. If $\{j_1-\tau, t_1-\tau, t_1\} \cap \{j_2-\tau, t_2-\tau, t_2\} \neq \emptyset$, by Assumption \ref{as:A2}, it holds that
\begin{align*}
    \sigma(j_1,t_1,j_2,t_2) \le &~ E\bigg\{\bigg(r_{3}^\top r_{2}\cdot r_{1}^\top \Sigma r_{1} \cdot r_{2}^\top r_{3}-\frac1p\alpha_{2p}^2\bigg)^2\bigg\}\\
    \leq &~ E\big\{(r_{3}^\top r_{2}\cdot r_{1}^\top \Sigma r_{1} \cdot r_{2}^\top r_{3})^2	\big\}
\leq K\cdot E\big\{(r_{3}^\top r_{2})^4\big\}=O(p^{-2})\,.
\end{align*}
Therefore, we obtain
\begin{align*}
V_1\leq 	\frac {K}{n^4}\sum_{j_1, j_2=\tau+1}^n	\sum_{\tau<t_1< j_1 \atop t_1\neq j_1-\tau}\sum_{\tau<t_2< j_2 \atop t_2\neq j_2-\tau}I\big(\{j_1-\tau, t_1-\tau, t_1\}\cap 	 \{j_2-\tau, t_2-\tau, t_2\}\neq\emptyset\big)=O(n^{-1})\,,
\end{align*}
which implies $V_1=o(1)$.

For $V_2$, write $\mathcal{T}_j = \{(t,s): j-\tau\neq t< j,~ j-\tau\neq s< j,~ t< s\}$ for $j \in \{\tau+1,\ldots,n\}$. By decomposing the summation over $(t,s)\in \mathcal{T}_j$ into two cases: (i)  $(t,s)\in \mathcal{T}_j$ with $s-\tau = t$, and (ii) $(t,s)\in \mathcal{T}_j$ with $s-\tau \neq t$, we can bound $V_2$  as follows:
\begin{align*}
   V_2 
   & \le \frac {Kp^2}{n^4}E\Bigg\{\bigg(\sum_{j=\tau+1}^n	\sum_{{j-\tau\neq t< j\atop j-\tau\neq s< j}\atop s-\tau = t}r_{j-\tau}^\top r_{t-\tau}\cdot r_t^\top \Sigma r_s \cdot r_{s-\tau}^\top r_{j-\tau}\bigg)^2\Bigg\} \\
   &~~~~~~+\frac {Kp^2}{n^4}E\Bigg\{\bigg(\sum_{j=\tau+1}^n	\sum_{{j-\tau\neq t< j\atop j-\tau\neq s< j}\atop t< s,\; s-\tau \neq t}r_{j-\tau}^\top r_{t-\tau}\cdot r_t^\top \Sigma r_s \cdot r_{s-\tau}^\top  r_{j-\tau}\bigg)^2\Bigg\} \\
   & =: V_{21}+V_{22}\,.
\end{align*}
Notice that $r_{j-\tau}^\top r_{t-\tau}r_t^\top \Sigma r_s r_{s-\tau}^\top r_{j-\tau}
=	r_{j-\tau}^\top r_{t-\tau}r_t^\top \Sigma r_{t+\tau} r_{t}^\top r_{j-\tau}$ for $(t,s)\in \mathcal{T}_j$ with $s-\tau=t$,
where the smallest index among all the indices of $r$'s is $t-\tau$, and $r_{t-\tau}$ is independent of all the others. This indicates that, if $t_1 \neq t_2$, then
\begin{align*}
E\{(r_{j_1-\tau}^\top r_{t_1-\tau}\cdot r_{t_1}^\top \Sigma r_{t_1+\tau} \cdot r_{t_1}^\top r_{j_1-\tau})\cdot 
(r_{j_2-\tau}^\top r_{t_2-\tau}\cdot r_{t_2}^\top \Sigma r_{t_2+\tau} \cdot r_{t_2}^\top r_{j_2-\tau})\}=0	
\end{align*}
for any $j_1, j_2 \in \{\tau+1,\ldots,n\}$ satisfying that $j_1-\tau \neq t_1 < j_1$ and $j_2-\tau \neq t_2 < j_2$.
Therefore, we obtain
\begin{align*}
V_{21}
= &~\frac {Kp^2}{n^4}E\Bigg\{\bigg(\sum_{j=\tau+1}^n	\sum_{ j-2\tau\neq t< j-\tau}r_{j-\tau}^\top r_{t-\tau}\cdot r_t^\top \Sigma r_{t+\tau} \cdot r_{t}^\top r_{j-\tau}\bigg)^2\Bigg\} \\
=&~\frac {Kp^2}{n^4}\sum_{j_1, j_2=\tau+1}^n	\sum_{\tau<t\atop{j_1-2\tau\neq t< j_1-\tau\atop j_2-2\tau\neq t< j_2-\tau}}E(r_{j_1-\tau}^\top r_{t-\tau}r_t^\top \Sigma r_{t+\tau} r_{t}^\top r_{j_1-\tau}r_{j_2-\tau}^\top r_{t-\tau}r_t^\top \Sigma r_{t+\tau} r_{t}^\top r_{j_2-\tau})\\
=&~\frac {K}{n^4}\sum_{j_1, j_2=\tau+1}^n	\sum_{\tau<t\atop{j_1-2\tau\neq t< j_1-\tau\atop j_2-2\tau\neq t< j_2-\tau}}
E(
r_{j_1-\tau}^\top \Sigma r_{j_2-\tau}\cdot r_t^\top \Sigma^3 r_{t} \cdot r_{t}^\top r_{j_1-\tau}
 \cdot r_{t}^\top r_{j_2-\tau})\\
 \leq&~\frac {K}{n^4}\sum_{j_1, j_2=\tau+1}^n	\sum_{\tau<t\atop{j_1-2\tau\neq t< j_1-\tau\atop j_2-2\tau\neq t< j_2-\tau}}
E\{
r_{j_1-\tau}^\top r_{j_1-\tau}\cdot(r_t^\top r_{t})^2\cdot r_{j_2-\tau}^\top r_{j_2-\tau}\}\\
 \leq &~ \frac {K}{n} E^2(r_{1}^\top r_{1})
E\{(r_{1}^\top r_{1})^2\} + \frac {K}{n^2}
E^2\{(r_{1}^\top r_{1})^2\}\to0 \,,
\end{align*}
where the first inequality follows from the Cauchy-Schwarz inequality and Assumption~\ref{as:A4}. Then $V_{21}=o(1)$.
On the other hand, notice that $\{t-\tau, s-\tau, j-\tau\}\cap\{t, s\}=\emptyset$ for $(t,s)\in \mathcal{T}_j$ with $s-\tau \neq t$. 
Moreover, given $j_1, j_2 \in \{\tau+1,\ldots,n\}$, for any $(t_1,s_1)\in \mathcal{T}_{j_1}$ and $(t_2,s_2)\in \mathcal{T}_{j_2}$ with $s_1-\tau \neq t_1$ and $s_2 -\tau \neq t_2$, if $t_1 \neq t_2$ or $s_1 \neq s_2$, it holds that
\begin{align*}
E\{
(r_{j_1-\tau}^\top r_{t_1-\tau}\cdot r_{t_1}^\top \Sigma r_{s_1} \cdot r_{s_1-\tau}^\top r_{j_1-\tau})\cdot 
(r_{j_2-\tau}^\top r_{t_2-\tau}\cdot r_{t_2}^\top \Sigma r_{s_2} \cdot r_{s_2-\tau}^\top r_{j_2-\tau})\}=0\,.
\end{align*}
Therefore, we obtain 
\begin{align*}
V_{22}
=&~\frac {Kp^2}{n^4}\sum_{j_1,j_2=\tau+1}^n	\sum_{{j_1-\tau\neq t< j_1,\;j_2-\tau\neq t< j_2\atop j_1-\tau\neq s< j_1,\;j_2-\tau\neq s< j_2}\atop t< s,\; s-\tau \neq t}E\big\{
r_{j_1-\tau}^\top r_{t-\tau}\cdot(r_t^\top \Sigma r_s)^2 \cdot r_{s-\tau}^\top r_{j_1-\tau}
r_{j_2-\tau}^\top r_{t-\tau}\cdot r_{s-\tau}^\top r_{j_2-\tau}\big\}
\\
=&~\frac {K}{n^4}\sum_{j_1,j_2=\tau+1}^n	\sum_{{j_1-\tau\neq t< j_1,\;j_2-\tau\neq t< j_2\atop j_1-\tau\neq s< j_1,\;j_2-\tau\neq s< j_2}\atop t< s,\; s-\tau \neq t}E\big\{
r_{j_1-\tau}^\top \Sigma r_{j_2-\tau} \cdot\left(r_t^\top \Sigma r_s\right)^2
\cdot r_{j_1-\tau}^\top \Sigma r_{j_2-\tau}\big\}
\\
= &~\frac {K}{n^4}\sum_{j_1,j_2=\tau+1}^n	\sum_{{j_1-\tau\neq t< j_1,\;j_2-\tau\neq t< j_2\atop j_1-\tau\neq s< j_1,\;j_2-\tau\neq s< j_2}\atop t< s,\; s-\tau \neq t}E\{(r_t^\top \Sigma r_s)^2\}E\{(
r_{j_1-\tau}^\top \Sigma r_{j_2-\tau})^2\} \\
= &~ \frac {K}{n^4}\sum_{j=\tau+1}^n	\sum_{{j-\tau\neq t< j \atop j-\tau\neq s< j}\atop t< s,\; s-\tau \neq t}E\{(r_t^\top \Sigma r_s)^2\}E\{(
r_{j-\tau}^\top \Sigma r_{j-\tau})^2\} \\
&~~~~~~~+ \frac {K}{n^4}\sum_{j_1\neq j_2}	\sum_{{j_1-\tau\neq t< j_1,\;j_2-\tau\neq t< j_2\atop j_1-\tau\neq s< j_1,\;j_2-\tau\neq s< j_2}\atop t< s,\; s-\tau \neq t}E\{(r_t^\top \Sigma  r_s)^2\}E\{(
r_{j_1-\tau}^\top \Sigma r_{j_2-\tau})^2\}\,.
\end{align*}
Under Assumptions \ref{as:A2} and \ref{as:A4}, by \eqref{moments} with $k=2$, $\max_{t\in[n]}\E\{(\br_t^\top \bSigma\br_t)^2\} =O(\|\bSigma\|_2^2) \le K$. Meanwhile, for any $t\neq s$, $\E\{(\br_t^\top \bSigma\br_s)^2\}=p^{-1}\E(\br_t^\top \bSigma^3\br_t) \le p^{-1}\|\bSigma\|_2^3 {\E}(\br_t^{\top} \br_t)\le Kp^{-1}$. Hence, it holds that $V_{22}
\le Kn^{-1}+Kp^{-2}=o(1)$. 
From the above results, we conclude that
$V_2\leq V_{21}+V_{22}=o(1)$. 
$\hfill\Box$

\subsubsection{Proof of \eqref{mclt-conditiona2}}\label{subsubsec:c2}
For any $\tau>\tau'$, we have 
\begin{align}\label{etts}
\sum_{j= \tau +1}^n\E_{j-1}(D_{j,\tau}D_{j,\tau'})
=&~\frac {4p^2}{n^2}\sum_{j= \tau +1}^n\E_{j-1}
\bigg(\sum_{j-\tau\neq t< j}\br_{j-\tau}^\top \br_{t-\tau}\br_t^\top \br_j \cdot
\sum_{j-\tau'\neq s< j}\br_{j-\tau'}^\top \br_{s-\tau'}\br_s^\top \br_j \bigg)
\nonumber\\
=&~\frac {4p}{n^2}\sum_{j= \tau +1}^n	\sum_{j-\tau\neq t< j\atop j-\tau'\neq s< j}
\br_{j-\tau}^\top \br_{t-\tau}\br_t^\top \bSigma\br_s\br_{s-\tau'}^\top \br_{j-\tau'}\,.
\end{align}
In addition, for any $j \ge \tau+1$, $t>\tau$ and $s>\tau'$, one can verify
\begin{align*}
\E(\br_{j-\tau}^\top \br_{t-\tau}\br_t^\top \bSigma\br_s\br_{s-\tau'}^\top \br_{j-\tau'})
=\begin{cases}
p^{-3}\tr(\bSigma^4)\,,& \mbox{if}~t=j-\tau', s=j-\tau\,,\\
0\,,&\text{otherwise}\,.
\end{cases}	
\end{align*}
Notice that $\br_{t} = 0$ for any $t \le 0$. Therefore, by Assumption \ref{as:A4}, the expectation of \eqref{etts} is $O(p^{-1}n^{-1})$.
The variance of  \eqref{etts} can be bounded as 
\begin{align*}
\Var\left\{\sum_{j= \tau +1}^n\E_{j-1}\left(D_{j,\tau}D_{j,\tau'}\right)\right\}
\leq &~\frac {Kp^2}{n^4}
\E\Bigg\{\bigg(\sum_{j= \tau +1}^n	\sum_{j-\tau\neq t< j\atop j-\tau'\neq t}
\br_{j-\tau}^\top \br_{t-\tau}\cdot\br_t^\top \bSigma\br_t\cdot\br_{t-\tau'}^\top \br_{j-\tau'}\bigg)^2\Bigg\}\\
&~~~+\frac {Kp^2}{n^4}
\E\Bigg\{\bigg(\sum_{j= \tau +1}^n	\sum_{t\neq s\atop{j-\tau\neq t< j\atop j-\tau'\neq s< j}}
\br_{j-\tau}^\top \br_{t-\tau}\br_t^\top \bSigma\br_s\br_{s-\tau'}^\top \br_{j-\tau'}\bigg)^2\Bigg\}\\
=&~ \tilde V_1+\tilde V_2\,.
\end{align*}
To show \eqref{mclt-conditiona2}, it suffices to show $\tilde{V}_1+\tilde{V}_2 = o(1)$.
We note that the order of $\tilde V_2$ can be derived by considering two cases: $t=s-\tau'$ and $t\neq s-\tau'$, respectively. Using the similar arguments for controlling $V_2$, we have $\tilde V_2 =o(1)$. Details are omitted.
To control $\tilde V_1$, notice that
\begin{align*}
\tilde V_1
=&~\frac {Kp^2}{n^4}\sum_{j_1,j_2= \tau +1}^n	\sum_{{j_1-\tau\neq t_1< j_1\atop j_1-\tau'\neq t_1}\atop {j_2-\tau\neq t_2< j_2\atop j_2-\tau'\neq t_2}}
\E
(\br_{j_1-\tau}^\top \br_{t_1-\tau}\br_{t_1}^\top \bSigma\br_{t_1}\br_{t_1-\tau'}^\top \br_{j_1-\tau'}
\br_{j_2-\tau}^\top \br_{t_2-\tau}\br_{t_2}^\top \bSigma\br_{t_2}\br_{t_2-\tau'}^\top \br_{j_2-\tau'})	\,.
\end{align*}
When $t_1 \neq t_2$, due to the facts $t_1-\tau=\min\{t_1-\tau, t_1-\tau', t_1, j_1-\tau, j_1-\tau'\}$ and $t_2-\tau=\min\{t_2-\tau, t_2-\tau', t_2, j_2-\tau, j_2-\tau'\}$, we have 
\[ \E
(\br_{j_1-\tau}^\top \br_{t_1-\tau}\br_{t_1}^\top \bSigma\br_{t_1}\br_{t_1-\tau'}^\top \br_{j_1-\tau'}\br_{j_2-\tau}^\top \br_{t_2-\tau}\br_{t_2}^\top \bSigma\br_{t_2}\br_{t_2-\tau'}^\top \br_{j_2-\tau'})=0 \]
for any $j_1,j_2 \ge \tau+1$ such that $j_1-\tau \neq t_1 < j_1$, $j_1-\tau' \neq t_1$, $j_2-\tau \neq t_2 < j_2$ and $j_2-\tau' \neq t_2$, which implies
\begin{align*}
\tilde V_1
=&~\frac {Kp^2}{n^4}\sum_{j_1,j_2= \tau +1}^n	\sum_{{j_1-\tau\neq t< j_1\atop j_1-\tau'\neq t}\atop {j_2-\tau\neq t< j_2\atop j_2-\tau'\neq t}}
\E
(\br_{j_1-\tau}^\top \br_{t-\tau}\br_{t}^\top \bSigma\br_{t}\br_{t-\tau'}^\top \br_{j_1-\tau'}
\br_{j_2-\tau}^\top \br_{t-\tau}\br_{t}^\top \bSigma\br_{t}\br_{t-\tau'}^\top \br_{j_2-\tau'}) \\
=&~\frac {Kp}{n^4}\sum_{j_1,j_2= \tau +1}^n	\sum_{{j_1-\tau\neq t< j_1\atop j_1-\tau'\neq t>\tau}\atop {j_2-\tau\neq t< j_2\atop j_2-\tau'\neq t}}
\E
\big\{\br_{j_1-\tau}^\top \bSigma\br_{j_2-\tau}(\br_{t}^\top \bSigma\br_{t})^2\br_{t-\tau'}^\top \br_{j_1-\tau'}
\br_{t-\tau'}^\top \br_{j_2-\tau'}\big\} \\
=&~\frac {Kp}{n^4}\sum_{j= \tau +1}^n	\sum_{j-\tau\neq t< j\atop j-\tau'\neq t>\tau}
\E
\big\{\br_{j-\tau}^\top \bSigma\br_{j-\tau}(\br_{t}^\top \bSigma\br_{t})^2\br_{t-\tau'}^\top \br_{j-\tau'}
\br_{t-\tau'}^\top \br_{j-\tau'}\big\} \\
& + \frac {2Kp}{n^4}\sum_{j_1>j_2}	\sum_{{j_1-\tau\neq t< j_1\atop j_1-\tau'\neq t>\tau}\atop {j_2-\tau\neq t< j_2\atop j_2-\tau'\neq t}}
\E
\big\{\br_{j_1-\tau}^\top \bSigma\br_{j_2-\tau}(\br_{t}^\top \bSigma\br_{t})^2\br_{t-\tau'}^\top \br_{j_1-\tau'}
\br_{t-\tau'}^\top \br_{j_2-\tau'}\big\} \\
=&~ \tilde{V}_{11}+\tilde{V}_{12}\,.
\end{align*}
For any $j \ge \tau+1$, if $j-\tau \neq t < j$ and $j-\tau' \neq t > \tau$, one can verify that
\begin{align*}
  &  \E
\big\{\br_{j-\tau}^\top \bSigma\br_{j-\tau}(\br_{t}^\top \bSigma\br_{t})^2\br_{t-\tau'}^\top \br_{j-\tau'}
\br_{t-\tau'}^\top \br_{j-\tau'}\big\} \\
&~~~~~~~~= \begin{cases} 
p^{-1}
\E\big\{\br_{j-\tau}^\top \bSigma\br_{j-\tau}(\br_{t}^\top \bSigma\br_{t})^2\br_{j-\tau'}^\top 
\bSigma\br_{j-\tau'}\big\}\,,& \mbox{if}~t-\tau' \neq j-\tau\,,\notag\\
p^{-1}\E
\big\{(\br_{j-\tau}^\top \bSigma\br_{j-\tau})^2(\br_{t}^\top \bSigma\br_{t})^2\big\}\,,&\mbox{if}~t-\tau'= j-\tau\,. \notag
\end{cases}
\end{align*}
By \eqref{moments}, under Assumptions \ref{as:A2} and \ref{as:A4}, for any $j \ge \tau+1$, $j-\tau \neq t < j$ and $j-\tau' \neq t > \tau$, we have
\begin{align*}
&\E\big\{\br_{j-\tau}^\top \bSigma\br_{j-\tau}(\br_{t}^\top \bSigma\br_{t})^2\br_{j-\tau'}^\top 
\bSigma\br_{j-\tau'}\big\}
=\E(\br_{j-\tau}^\top \bSigma\br_{j-\tau})\E(\br_{j-\tau'}^\top \bSigma\br_{j-\tau'})\E\big\{(\br_{t}^\top \bSigma\br_{t})^2\big\}\leq K\,, \\
&~~~~~~~~~~~~~\E\big\{(\br_{j-\tau}^\top \bSigma\br_{j-\tau})^2(\br_{t}^\top \bSigma\br_{t})^2\big\}
=\E\big\{(\br_{j-\tau}^\top \bSigma\br_{j-\tau})^2\big\}\E\big\{(\br_{t}^\top \bSigma\br_{t})^2\big\}\leq K\,,
\end{align*} 
which implies $\tilde{V}_{11} = O(n^{-2})$.
Additionally, one can verify that
\begin{align*}
    \E
\big\{\br_{j_1-\tau}^\top \bSigma\br_{j_2-\tau}(\br_{t}^\top \bSigma\br_{t})^2\br_{t-\tau'}^\top \br_{j_1-\tau'}
\br_{t-\tau'}^\top \br_{j_2-\tau'}\big\} = 0
\end{align*}
for any $j_1 > j_2$, $j_1-\tau \neq t < j_1$, $j_1-\tau' \neq t$, $j_2-\tau \neq t < j_2$ and $j_2-\tau' \neq t$, which implies $\tilde{V}_{12} = 0$.
Hence, $\tilde{V}_1 = O(n^{-2})$. 
Together with $\tilde V_2 =o(1)$, we have
$\tilde{V}_1+\tilde{V}_2 = o(1)$, which implies that
\eqref{mclt-conditiona2} holds. $\hfill\Box$

\subsubsection{Proof of \eqref{mclt-conditiona3}}\label{subsubsec:c3} 
Notice that
\begin{align}
\sum_{j=\tau+1}^n\E (D_{j,\tau}^4)
=&~\frac{16p^4}{n^4}\sum_{j=\tau+1}^n\E\left\{\bigg(\sum_{j-\tau\neq t< j}\br_{t-\tau}^\top \br_{j-\tau}\br_t^\top \br_j\bigg)^4\right\} \notag\\
\leq &~
\frac{Kp^4}{n^4}\sum_{j=\tau+1}^n\E\left\{\bigg(\sum_{t< j-\tau}\br_{t-\tau}^\top \br_{j-\tau}\br_t^\top \br_j\bigg)^4\right\} \notag\\
&~+\frac{Kp^4}{n^4}\sum_{j=\tau+1}^n\E\left\{\bigg(\sum_{j-\tau< t< j}\br_{t-\tau}^\top \br_{j-\tau}\br_t^\top \br_j\bigg)^4\right\}\notag\\
:= &~ R_1+R_2\,. \label{eq:4moment}
\end{align}
Since $\tau$ is a fixed constant, then
\begin{align*}
R_2
\leq &~\frac{Kp^4}{n^4}\sum_{j=\tau+1}^n \sum_{j-\tau< t< j}\E\big\{(\br_{t-\tau}^\top \br_{j-\tau}\br_t^\top \br_j)^4\big\} \\	=&~\frac{Kp^4}{n^4}\sum_{j=\tau+1}^n \sum_{j-\tau< t< j}{\E}\big\{(\br_{t-\tau}^\top \br_{j-\tau})^4\big\}{\E}\big\{(\br_t^\top \br_j)^4\big\} = O(n^{-3})\,,
\end{align*}
where the last step is based on Assumption \ref{as:A2}.
To bound $R_1$, we write
\begin{align*}
R_1
=&~\frac{Kp^4}{n^4}\sum_{j=\tau+1}^n\sum_{t_1,t_2,t_3,t_4< j-\tau}
\E \left(\prod_{i=1}^4   \br_{t_i-\tau}^\top \br_{j-\tau}\br_{t_i}^\top \br_j \right)\\	
\leq &~\frac{Kp^4}{n^4}\sum_{j=\tau+1}^n\sum_{t_1,t_2,t_3,t_4< j-\tau\atop t_1\leq t_2\leq t_3\leq t_4}
\left|\E \bigg( \prod_{i=1}^4   \br_{t_i-\tau}^\top \br_{j-\tau}\br_{t_i}^\top \br_j \bigg)\right|\,.
\end{align*}
To guarantee $\E ( \prod_{i=1}^4   \br_{t_i-\tau}^\top \br_{j-\tau}\br_{t_i}^\top \br_j) \neq 0$ with $t_1\leq t_2\leq t_3\leq t_4$, the indices $(t_1, t_2, t_3, t_4)$ must satisfy that $t_1 = t_2$ and $t_3 = t_4$, which implies
\begin{align*}
R_1
\leq &~\frac{Kp^4}{n^4}\sum_{j=\tau+1}^n\sum_{t_1,t_3< j-\tau\atop t_1\leq t_3}
\E \big\{  (\br_{t_1-\tau}^\top \br_{j-\tau}\br_{t_1}^\top \br_j)^2(\br_{t_3-\tau}^\top \br_{j-\tau}\br_{t_3}^\top \br_j)^2 \big\}\\
\leq&~ \frac{Kp^4}{n^4}\sum_{j=\tau+1}^n\left[\sum_{t< j-\tau}
{\E}^{1/2}  \big\{(\br_{t-\tau}^\top \br_{j-\tau}\br_{t}^\top \br_j)^4\big\} \right]^2\\
\leq&~ \frac{Kp^4}{n}
\E   \big\{(\br_{1}^\top \br_{3}\br_{2}^\top \br_4)^4\big\}
= \frac{Kp^4}{n}
   {\E}^2\big\{(\br_{1}^\top \br_{3})^4\big\}=O(n^{-1})\,.
\end{align*}
By \eqref{eq:4moment}, we have $\sum_{j=\tau+1}^n\E (D_{j,\tau}^4)
=O(n^{-1})$. Since $\sigma_{n1}^2 = 2p^{-2}{\tr}^2(\bSigma^2)$ is uniformly bounded away from zero under Assumption \ref{as:A4}, then
\eqref{mclt-conditiona3} holds. $\hfill\Box$

\subsection{Proof of Step 2}\label{subsec:Gq2}
Recall that $\alpha_{1p} = p^{-1}{\tr}(\bSigma)$ and
\begin{align*}
  G_{n,\tau,2}
=&~\frac p{n}\sum_{t=1}^n\br_t^\top \br_t \br_{t+\tau}^\top \br_{t+\tau}- \frac{p}{n^2}\sum_{t,\ell=1}^n \br_t^\top \br_t\br_\ell^\top \br_\ell\,.
\end{align*}
It then holds that 
\begin{align}\label{gn3-mean}
\E( G_{n,\tau,2})
=&~\frac p{n}\bigg[n\alpha_{1p}^2- \frac{1}{n}\sum_{t\neq \ell}\E(\br_t^\top \br_t\br_\ell^\top \br_\ell)-\frac{1}{n}\sum_{t=1}^n\E\big\{(\br_t^\top \br_t)^2\big\}\bigg]\nonumber\\
=&~\frac p{n}\left[\alpha_{1p}^2-\E\big\{(\br_1^\top \br_1)^2\big\}\right] = -\frac pn\Var(\br_1^\top \br_1)\,,
\end{align}
and
\begin{align}\label{eq:secondMoment_G_n2}
\E ( G_{n,\tau,2}^2)
=&~\frac{p^2}{n^2}\sum_{t,s=1}^n\E\bigg\{\bigg(\br_t^\top \br_t \br_{t+\tau}^\top \br_{t+\tau}- \br_t^\top \br_t \cdot \frac1n\sum_{\ell=1}^n\br_\ell^\top \br_\ell\bigg) \notag\\
&~~~~~~~~~~~~~~~~~~~\times\bigg(\br_s^\top \br_s \br_{s+\tau}^\top \br_{s+\tau}- \br_s^\top \br_s \cdot\frac1n\sum_{\ell=1}^n\br_\ell^\top \br_\ell\bigg) \bigg\}\notag\\
= &~\frac{p^2}{n^2}(E_1-2E_2+E_3)\,,
\end{align}
where 
\begin{align*}
E_1=&~\sum_{t,s=1}^n\E(\br_t^\top \br_t \br_{t+\tau}^\top \br_{t+\tau}\br_s^\top \br_s \br_{s+\tau}^\top \br_{s+\tau})\,,\\	
E_2=
&~\frac{1}{n}\sum_{t=1}^n\E \bigg\{(\br_t^\top \br_t \br_{t+\tau}^\top \br_{t+\tau}) \bigg(\sum_{\ell=1}^n\br_\ell^\top \br_\ell\bigg)^2 \bigg\}\,,\\
E_3=
&~\frac1{n^2}\E\bigg\{\bigg(\sum_{\ell=1}^n\br_\ell^\top \br_\ell\bigg)^4\bigg\}\,.	
\end{align*}
As we will show in Sections \ref{subsec:E1} and \ref{subsec:E2}, 
\begin{align}
E_1=&~ (n^2-3n){\E}^4(\br_1^\top \br_1)+2n\E\{(\br_1^\top \br_1)^2\}{\E}^2(\br_1^\top \br_1)+n{\E}^2\{(\br_1^\top \br_1)^2\}\,,\label{eq:E1}\\
E_2=&~
(n^2-5n+6){\E}^4(\br_1^\top \br_1)
+(5n-10)\E\{(\br_1^\top \br_1)^2\}{\E}^2(\br_1^\top \br_1)\notag\\
&+2{\E}^2\{(\br_1^\top \br_1)^2\}+2\E\{(\br_1^\top \br_1)^3\}\E(\br_1^\top \br_1)\,,\label{eq:E2}\\
E_3=&~\bigg(n^2- 6n+ 11 -\frac{6}{n}\bigg)\E^4(\br_1^\top \br_1)+\bigg( 6n- 18 +\frac{12}{n}\bigg)\E\{(\br_1^\top \br_1)^2\}{\E}^2(\br_1^\top \br_1) \notag\\
&+
\bigg(3-\frac{3}{n}\bigg){\E}^2\{(\br_1^\top \br_1)^2\} +\bigg(4-\frac{4}{n}\bigg) \E(\br_1^\top \br_1)\E\{(\br_1^\top \br_1)^3\}+\frac{1}{n}\E\{(\br_1^\top \br_1)^4\}\,.\label{eq:E3}
\end{align}
Collecting these results, by \eqref{eq:secondMoment_G_n2}, it holds that
\begin{align}\label{gn3-m2}
\E (G_{n,\tau,2}^2)
=&~
\frac{p^2}{n^2}\bigg[\bigg(n-1-\frac{6}{n}\bigg){\E}^4(\br_1^\top \br_1)+\bigg(2-2n+\frac{12}{n}\bigg)\E\{(\br_1^\top \br_1)^2\}{\E}^2(\br_1^\top \br_1)\\
&~~~~~~~+\bigg(n-1-\frac{3}{n}\bigg){\E}^2\{(\br_1^\top \br_1)^2\}
-\frac{4}{n}\E(\br_1^\top \br_1)\E\{(\br_1^\top \br_1)^3\}+\frac{1}{n}\E\{(\br_1^\top \br_1)^4\}\bigg]\,.	\nonumber
\end{align}
Since $\alpha_{1p} = \E(\br_1^\top \br_1)
$, then we have  
\begin{align*}
\E\{(\br_1^\top \br_1)^2\}
=&~\alpha_{1p}^2+\Var(\br_1^\top \br_1)\,,\\
\E\{(\br_1^\top \br_1)^3\}
=&~\E\{(\br_1^\top \br_1-\alpha_{1p})^3\}+3\alpha_{1p}\E\{(\br_1^\top \br_1-\alpha_{1p})^2\}+\alpha_{1p}^3\\
=&~\E\{(\br_1^\top \br_1-\alpha_{1p})^3\}+3\alpha_{1p}\Var(\br_1^\top \br_1)+\alpha_{1p}^3\,,\\
\E\{(\br_1^\top \br_1)^4\}
=&~\E\{(\br_1^\top \br_1-\alpha_{1p})^4\}+4\alpha_{1p}\E\{(\br_1^\top \br_1-\alpha_{1p})^3\} \\
&+6\alpha_{1p}^2\E\{(\br_1^\top \br_1-\alpha_{1p})^2\}+\alpha_{1p}^4\\	
=&~\E\{(\br_1^\top \br_1-\alpha_{1p})^4\}+4\alpha_{1p}\E\{(\br_1^\top \br_1-\alpha_{1p})^3\} \\
&+6\alpha_{1p}^2\Var(\br_1^\top \br_1)+\alpha_{1p}^4\,.
\end{align*}
By \eqref{gn3-mean} and \eqref{gn3-m2}, we have
\begin{align}\label{eq:Var_G2}
\Var(G_{n,\tau,2})=	\frac{p^2}{n^3}\E\{(\br_1^\top \br_1-\alpha_{1p})^4\}+\frac{p^2(n-2-3n^{-1})}{n^2}{\Var}^2(\br_1^\top \br_1)=\sigma_{n2}^2\,.
\end{align}

Let $\E_0(\cdot)$ denote expectation, and $\E_j(\cdot)$ denote conditional expectation with respect to the $\sigma$-field $\mathcal F_j=\sigma(\x_1, \ldots,\x_j)$ for $j \ge 1$. 
To establish \eqref{eq:AsymNormGq2}, we first apply the martingale decomposition to $ G_{n,\tau,2}$ as follows:
\begin{align}\label{eq:gn2}
  G_{n,\tau,2}-\E(G_{n,\tau,2})
=&~\sum_{j=1}^n\big\{\E_j(G_{n,\tau,2})-\E_{j-1}(G_{n,\tau,2})\big\} = \sum_{j=1}^nM_{j,\tau}\,,
\end{align}
where $\{M_{j,\tau}\}$ is a sequence of martingale differences with respect to $\{\mathcal{F}_j\}$. 
Notice that
\begin{align}\label{mjtl}
M_{j,\tau}
=&~ \frac {2p}{n} (\br_{j}^\top \br_{j}-\alpha_{1p})\alpha_{1p} I(1\leq j\leq \tau) \nonumber\\ 
&+ \frac {p}{n}(\br_{j}^\top \br_{j}-\alpha_{1p})( \br_{j-\tau}^\top \br_{j-\tau} +\alpha_{1p}) I(\tau+1\leq j\leq n-\tau) \nonumber\\
& + \frac {p}{n} (\br_{j}^\top \br_{j}-\alpha_{1p})( \br_{j-\tau}^\top \br_{j-\tau}+\br_{j+\tau}^\top \br_{j+\tau} ) I(n-\tau+1\leq j\leq n) \nonumber\\
&-\frac{p}{n}\bigg\{
(\br_j^\top \br_j-\alpha_{1p})\cdot\frac{2}{n}\sum_{t< j}\br_t^\top \br_t+
(\br_j^\top \br_j-\alpha_{1p})\cdot\frac{2}{n}\sum_{t> j}\alpha_{1p} \nonumber\\
&~~~~~~~~+
\frac{1}{n}\big[(\br_j^\top \br_j)^2-\E\{(\br_j^\top \br_j)^2\}\big]\bigg\} \nonumber\\
= &~ M_{j,\tau,1}+M_{j,\tau,2}+M_{j,\tau,3}\,,
\end{align}
where
\begin{align*}
    & M_{j,\tau,1} =  \frac pn \times \begin{cases}
 0\,,&	\mbox{if}~1\leq j\leq \tau\,,\\
 (\br_{j}^\top \br_{j}-\alpha_{1p})( \br_{j-\tau}^\top \br_{j-\tau} -\alpha_{1p})\,,& \mbox{if}~	\tau+1\leq j\leq n-\tau\,,\\
 (\br_{j}^\top \br_{j}-\alpha_{1p})( \br_{j-\tau}^\top \br_{j-\tau}+\br_{j+\tau}^\top \br_{j+\tau} -2\alpha_{1p})\,,&	\mbox{if}~n-\tau+1\leq j\leq n\,,
 \end{cases} \\
    & M_{j,\tau,2} = -\frac {2p}{n^2}(\br_{j}^\top \br_{j}-\alpha_{1p})\sum_{t< j}(\br_t^\top \br_t-\alpha_{1p}) \,, \\
    & M_{j,\tau,3} = - \frac{p}{n^2}\big\{(\br_j^\top \br_j-\alpha_{1p})^2-\Var(\br_j^\top \br_j)\big\}\,.
\end{align*}
As we will show in Sections \ref{subsubsec:b1}--\ref{subsubsec:b3}, 
\begin{align}
\sum_{j=1}^n\frac{\E_{j-1}(M_{j,\tau}^2)}{\Var(G_{n,\tau,2})} &\xrightarrow{\mathrm{p}} 1~\,\mbox{for any given}~\, \tau \in [q]\,, \label{mclt-conditionb1} \\
\sum_{j=1}^n\frac{\E_{j-1}(M_{j,\tau}M_{j,\tau'})}{\Var(G_{n,\tau,2})}&\xrightarrow{\mathrm{p}} 0 ~\, \mbox{for any given} ~\, (\tau,\tau')\in [q]^2 ~\mbox{satisfying}~ \tau> \tau'\,, \label{mclt-conditionb2}\\
\sum_{j=1}^n\frac{\E (M_{j,\tau}^4)}{\Var^2(G_{n,\tau,2})}&\to0 ~\,\mbox{for any given}~\, \tau \in [q]\,.\label{mclt-conditionb3}
\end{align} 
Combining \eqref{mclt-conditionb1}--\eqref{mclt-conditionb3}, by Lemma \ref{lem:md-clt-vector}, we have \eqref{eq:AsymNormGq2} holds. $\hfill\Box$

\subsubsection{Calculation of $E_1$}\label{subsec:E1}
To show \eqref{eq:E1}, we consider five cases for the relationship between $t$ and $s$. Based on the results of these five cases derived below, we know \eqref{eq:E1} holds by simple calculation.

\underline{\it Case 1: $t=s$.}
Write 
\begin{align*}
E_{11}
= \sum_{t=1}^n\E\{(\br_t^\top \br_t \br_{t+\tau}^\top \br_{t+\tau})^2\}\,.
\end{align*}
Based on Assumption \ref{as:A2}, we have $  E_{11}=n{\E}^2\{(\br_1^\top \br_1)^2\}$.

\underline{\it Case 2: 
$1\leq t\leq n-\tau$, $1\leq s\leq n-\tau$, and $t\neq s$.} Write
\begin{align*}
    E_{12}
= &~\sum_{1\le t\neq s\le n-\tau}
\E(\br_t^\top \br_t \br_{t+\tau}^\top \br_{t+\tau}\br_s^\top \br_s \br_{s+\tau}^\top \br_{s+\tau}) \\
=&~ 2\sum_{1\le t < s\le n-\tau}
\E(\br_t^\top \br_t \br_{t+\tau}^\top \br_{t+\tau}\br_s^\top \br_s \br_{s+\tau}^\top \br_{s+\tau})\,.
\end{align*}
Then, by Assumption \ref{as:A2}, we have
\begin{align*}
E_{12}
= &~2\sum_{t<s\neq t+\tau}
\E(\br_t^\top \br_t \br_{t+\tau}^\top \br_{t+\tau}\br_s^\top \br_s \br_{s+\tau}^\top \br_{s+\tau})+
2\sum_{t<s= t+\tau}\E(\br_t^\top \br_t \br_{t+\tau}^\top \br_{t+\tau}\br_s^\top \br_s \br_{s+\tau}^\top \br_{s+\tau})\\
= &~2\sum_{t<s\neq t+\tau}\E(\br_t^\top \br_t) \E(\br_{t+\tau}^\top \br_{t+\tau})\E(\br_s^\top \br_s) \E(\br_{s+\tau}^\top \br_{s+\tau})\\
&~~~+
2\sum_{t=1}^{n-2\tau}\E(\br_t^\top \br_t) \E\big\{(\br_{t+\tau}^\top \br_{t+\tau})^2\big\}\E(\br_{t+2\tau}^\top \br_{t+2\tau})\\
=&~\{(n-\tau)^2-(n-\tau)-2(n-2\tau)\}{\E}^4(\br_1^\top \br_1)+2(n-2\tau)\E\{(\br_1^\top \br_1)^2\}{\E}^2(\br_1^\top \br_1)\,.
\end{align*}

\underline{\it Case 3: 
$n-\tau+1\leq t\leq n$ and $1\leq s\leq n-\tau$.} Write
\begin{align*}
E_{13}
= &~\sum_{t=n-\tau+1}^{n}\sum_{s=1}^{n-\tau}\E(\br_t^\top \br_t \br_{t+\tau}^\top \br_{t+\tau}\br_s^\top \br_s \br_{s+\tau}^\top \br_{s+\tau})\,.
\end{align*}
Recall $\br_t=\br_{t-n}$ for $t>n$. By Assumption \ref{as:A2}, we have
\begin{align*}
    E_{13}=&~\sum_{t=1}^{\tau}\sum_{s=1}^{n-\tau}\E(\br_t^\top \br_t \br_{t+n-\tau}^\top \br_{t+n-\tau}\br_s^\top \br_s \br_{s+\tau}^\top \br_{s+\tau})\\
    =&~\bigg(\sum_{t=s=1}^{\tau}+\sum_{t,s=1\atop t \neq s}^{\tau}+\sum_{t=1}^{\tau}\sum_{s=\tau+1\atop s+\tau = t+n-\tau}^{n-\tau}+\sum_{t=1}^{\tau}\sum_{s=\tau+1\atop s+\tau \neq t+n-\tau}^{n-\tau}\bigg)\E(\br_t^\top \br_t \br_{t+n-\tau}^\top \br_{t+n-\tau}\br_s^\top \br_s \br_{s+\tau}^\top \br_{s+\tau})\\
    =&~\{\tau(n-\tau)-2\tau\}{\E}^4(\br_1^\top \br_1)+2\tau \E\{(\br_1^\top \br_1)^2\}{\E}^2(\br_1^\top \br_1)\,.
\end{align*}

\underline{\it Case 4: 
$1\leq t\leq n-\tau$ and $n-\tau+1\leq s\leq n$.} Write
\begin{align*}
E_{14}
= &~\sum_{t=1}^{n-\tau}\sum_{s=n-\tau+1}^{n}\E(\br_t^\top \br_t \br_{t+\tau}^\top \br_{t+\tau}\br_s^\top \br_s \br_{s+\tau}^\top \br_{s+\tau})\,.
\end{align*}
Then $E_{14}=E_{13}$.

\underline{\it Case 5: 
$n-\tau+1\leq t\leq n$, $n-\tau+1\leq s\leq n$, and $t\neq s$.} Write
\begin{align*}
E_{15}
=&~\sum_{t=n-\tau+1}^{n}\sum_{s=n-\tau+1\atop s \neq t}^{n}\E(\br_t^\top \br_t \br_{t+\tau}^\top \br_{t+\tau}\br_s^\top \br_s \br_{s+\tau}^\top \br_{s+\tau})\,.
\end{align*}
Recall $\br_t=\br_{t-n}$ for $t>n$. By Assumption \ref{as:A2}, we have
\begin{align*}
E_{15}
=&~\sum_{t=n-\tau+1}^{n}\sum_{s=n-\tau+1\atop s \neq t}^{n}\E(\br_t^\top \br_t \br_{t+\tau-n}^\top \br_{t+\tau-n}\br_s^\top \br_s \br_{s+\tau-n}^\top \br_{s+\tau-n})\\
=&~\sum_{t=1}^{\tau}\sum_{s=1\atop s \neq t}^{\tau}\E(\br_t^\top \br_t \br_{t+n-\tau}^\top \br_{t+n-\tau}\br_s^\top \br_s \br_{s+n-\tau}^\top \br_{s+n-\tau}) \\
=&~\sum_{t=1}^{\tau}\sum_{s=1\atop s \neq t}^{\tau}\E(\br_t^\top \br_t) \E(\br_{t+n-\tau}^\top \br_{t+n-\tau})\E(\br_s^\top \br_s) \E(\br_{s+n-\tau}^\top \br_{s+n-\tau})\,,
\end{align*}
which implies $E_{15}=(\tau^2-\tau){\E}^4(\br_1^\top \br_1)$. $\hfill\Box$

\subsubsection{Calculation of $E_2$ and $E_3$}\label{subsec:E2}
To calculate $E_2$, we write
\begin{align*}
E_2
=&~\frac1n\sum_{t=1}^{n}\sum_{\ell=1}^n\E\{\br_t^\top \br_t \br_{t+\tau}^\top \br_{t+\tau} (\br_\ell^\top \br_\ell)^2\}
+\frac1n\sum_{t=1}^{n}\sum_{\ell_1\neq \ell_2}\E(\br_t^\top \br_t \br_{t+\tau}^\top \br_{t+\tau} \br_{\ell_1}^\top \br_{\ell_1}\br_{\ell_2}^\top \br_{\ell_2})	\\
= &~E_{21}+E_{22}\,.
\end{align*}
It is straightforward to verify that 
\begin{align*}
E_{21}
=&~(n-2)\E\{(\br_1^\top \br_1)^2\}{\E}^2(\br_1^\top \br_1)+2\E\{(\br_1^\top \br_1)^3\}\E(\br_1^\top \br_1)\,,\\
E_{22}
=&~(n^2-5n+6){\E}^4(\br_1^\top \br_1)+2{\E}^2\{(\br_1^\top \br_1)^2
\}+(4n-8)\E\{(\br_1^\top \br_1)^2\}{\E}^2(\br_1^\top \br_1)\,,
\end{align*}
which implies
\begin{align*}
E_2=&~
(n^2-5n+6){\E}^4(\br_1^\top \br_1)
+(5n-10)\E\{(\br_1^\top \br_1)^2\}{\E}^2(\br_1^\top \br_1)\\
&+2{\E}^2\{(\br_1^\top \br_1)^2\}+2\E\{(\br_1^\top \br_1)^3\}\E(\br_1^\top \br_1)\,.	
\end{align*}
Hence, we have \eqref{eq:E2} holds.
For $E_3$, by Assumption \ref{as:A2}, we have
\begin{align*}
E_3
=&~\frac1{n^2}\E\bigg\{\bigg(\sum_{\ell=1}^n\br_\ell^\top \br_\ell\bigg)^4\bigg\} = \frac1{n^2}\sum_{\ell_1,\ell_2,\ell_3,\ell_4=1}^n \E(\br_{\ell_1}^\top \br_{\ell_1}\br_{\ell_2}^\top \br_{\ell_2}\br_{\ell_3}^\top \br_{\ell_3}\br_{\ell_4}^\top \br_{\ell_4})\\
=&~\frac1{n^2}\sum_{\ell=1}^n \E\{(\br_{\ell}^\top \br_{\ell})^4\} + \frac{4}{n^2}\sum_{\ell_1,\ell_2=1\atop \ell_1 \neq \ell_2}^n \E\{(\br_{\ell_1}^\top \br_{\ell_1})^3\}\E(\br_{\ell_2}^\top \br_{\ell_2})\\
&+\frac{3}{n^2}\sum_{\ell_1,\ell_2=1\atop \ell_1 \neq \ell_2}^n \E\{(\br_{\ell_1}^\top \br_{\ell_1})^2\}\E\{(\br_{\ell_2}^\top \br_{\ell_2})^2\}+\frac{6}{n^2}\sum_{\ell_1,\ell_2,\ell_3=1\atop \ell_1 \neq \ell_2 \neq \ell_3}^n \E\{(\br_{\ell_1}^\top \br_{\ell_1})^2\}\E(\br_{\ell_2}^\top \br_{\ell_2})\E(\br_{\ell_3}^\top \br_{\ell_3})\\
&+\frac1{n^2}\sum_{\ell_1,\ell_2,\ell_3,\ell_4=1 \atop \ell_1\neq \ell_2 \neq \ell_3 \neq \ell_4}^n \E(\br_{\ell_1}^\top \br_{\ell_1})\E(\br_{\ell_2}^\top \br_{\ell_2})\E(\br_{\ell_3}^\top \br_{\ell_3})\E(\br_{\ell_4}^\top \br_{\ell_4})\\
=&~\frac{1}{n}\E\{(\br_1^\top \br_1)^4\} +\bigg(4-\frac{4}{n}\bigg) \E(\br_1^\top \br_1)\E\{(\br_1^\top \br_1)^3\}+
\bigg(3-\frac{3}{n}\bigg){\E}^2\{(\br_1^\top \br_1)^2\} \\
&+\bigg( 6n- 18 +\frac{12}{n}\bigg)\E\{(\br_1^\top \br_1)^2\}{\E}^2(\br_1^\top \br_1)+\bigg(n^2- 6n+ 11 -\frac{6}{n}\bigg)\E^4(\br_1^\top \br_1)\,,
\end{align*}
which implies \eqref{eq:E3} holds. $\hfill\Box$

\subsubsection{Proof of \eqref{mclt-conditionb1}}\label{subsubsec:b1}
By the martingale decomposition of $G_{n,\tau,2}$, it holds naturally that
\begin{align*}
\E\Bigg\{\sum_{j=1}^n\frac{\E_{j-1}(M_{j,\tau}^2)}{\Var(G_{n,\tau,2})}\Bigg\}=1\,.
\end{align*}
Therefore, to show \eqref{mclt-conditionb1}, it suffices to show
\begin{align}\label{mjl123}
\Var\Bigg\{\sum_{j=1}^n\frac{\E_{j-1}(M_{j,\tau,\ell_1}M_{j,\tau,\ell_2})}{\Var(G_{n,\tau,2})}\Bigg\}\to0 ~\, \mbox{for any}~ \ell_1,\ell_2=1,2,3\,.	
\end{align}

For $\ell_1=\ell_2=1$, we have
\begin{align*}
\E_{j-1}(M_{j,\tau,1}^2)
	=&~\frac{p^2}{n^2} \times\begin{cases}
0\,,&	\mbox{if}~1\leq j\leq \tau\,,\\
\Var(\br_{j}^\top \br_{j})\cdot( \br_{j-\tau}^\top \br_{j-\tau} -\alpha_{1p})^2\,,& \mbox{if}~	\tau+1\leq j\leq n-\tau\,,\\
\Var(\br_{j}^\top \br_{j})\cdot( \br_{j-\tau}^\top \br_{j-\tau}+\br_{j+\tau}^\top \br_{j+\tau} -2\alpha_{1p})^2\,,&	\mbox{if}~ n-\tau+1\leq j\leq n\,.
\end{cases}
\end{align*}
Hence, it holds that
\begin{align*}
\Var\left\{\sum_{j=1}^n\frac{\E_{j-1}(M_{j,\tau,1}^2)}{\Var(G_{n,\tau,2})}\right\}	
=&~ 
\frac{p^4\Var^2(\br_{1}^\top \br_{1})}{n^{4}\Var^2(G_{n,\tau,2})}
\Bigg[\Var\Bigg\{\sum_{j=\tau+1}^{n-\tau}  ( \br_{j-\tau}^\top \br_{j-\tau} -\alpha_{1p})^2\\
&~~~~+\sum_{j=n-\tau+1}^{n}  ( \br_{j-\tau}^\top \br_{j-\tau}+\br_{j+\tau}^\top \br_{j+\tau} -2\alpha_{1p})^2\Bigg\}\Bigg]	\\
\leq &~ 
\frac{Kp^4\Var^2(\br_{1}^\top \br_{1})\Var\{(\br_{1}^\top \br_{1} -\alpha_{1p})^2\}}{n^{3}\Var^2(G_{n,\tau,2})} + \frac{Kp^4\Var^4(\br_{1}^\top \br_{1})}{n^{4}\Var^2(G_{n,\tau,2})}\\
\leq &~ 
\frac{K\E\{( \br_{1}^\top \br_{1} -\alpha_{1p})^4\}}{n\Var^2(\br_{1}^\top \br_{1})} + \frac{K}{n^2} \leq 
\frac{K}{n} \to0\,,
\end{align*}
where the last inequality is based on Assumption \ref{as:A3}. 

For $\ell_1=\ell_2=2$, we have
\begin{align*}
\E_{j-1}(M_{j,\tau,2}^2)
	=&~\frac{p^2}{n^2}\Var(\br_{j}^\top \br_{j})\Bigg\{\frac{2}{n}\sum_{t< j}(\br_t^\top \br_t-\alpha_{1p})\Bigg\}^2\,.
\end{align*}
Since $\Var(G_{n,\tau,2}) \gtrsim p^2n^{-1}\Var^2(\br_1^{\top}\br_1)$, it holds that
\begin{align*}
\Var\left\{\sum_{j=1}^n\frac{\E_{j-1}(M_{j,\tau,2}^2)}{\Var(G_{n,\tau,2})}\right\}
=&~ 
\frac{p^4\Var^2(\br_{1}^\top \br_{1})}{n^{4}\Var^2(G_{n,\tau,2})}
\Var\Bigg[\sum_{j=1}^n\Bigg\{\frac{2}{n}\sum_{t< j}(\br_t^\top \br_t-\alpha_{1p})\Bigg\}^2\Bigg]\\
\leq &~ 
\frac{p^4\Var^2(\br_{1}^\top \br_{1})}{n^{4}\Var^2(G_{n,\tau,2})} \times n^2
\max_{j\in[n]}\Var\Bigg[\Bigg\{\frac{2}{n}\sum_{t< j}(\br_t^\top \br_t-\alpha_{1p})\Bigg\}^2\Bigg]\\
\leq &~ 
\frac{K}{\Var^2(\br_{1}^\top \br_{1})}
\max_{j\in[n]}\Var\Bigg[\Bigg\{\frac{1}{n}\sum_{t< j}(\br_t^\top \br_t-\alpha_{1p})\Bigg\}^2\Bigg]\\
\leq &~ 
\frac{K}{n^4\Var^2(\br_{1}^\top \br_{1})}
\max_{j\in[n]}\sum_{t_1, t_2, t_3, t_4< j}\E\Bigg\{\prod_{i=1}^4(\br_{t_i}^\top \br_{t_i}-\alpha_{1p})\Bigg\}\\
\leq &~ 
\frac{K\E\{(\br_{t}^\top \br_{t}-\alpha_{1p})^4\}}{n^2\Var^2(\br_{1}^\top \br_{1})} \leq \frac{K}{n^2}\to0 \,,
\end{align*}
where the last inequality is based on Assumption \ref{as:A3}.

For $\ell_1=\ell_2=3$, notice that
\begin{align*}
    \E_{j-1}(M_{j,\tau,3}^2) =&~ \frac{p^2}{n^4}\E_{j-1}\big[\{(\br_j^\top \br_j-\alpha_{1p})^2-\Var(\br_j^\top \br_j)\}^2\big] \\
    =&~ \frac{p^2}{n^4}\E\big[\{(\br_j^\top \br_j-\alpha_{1p})^2-\Var(\br_j^\top \br_j)\}^2 \big]
\end{align*}
is a constant. Then 
\begin{align*}
    \Var\Bigg\{\sum_{j=1}^n\E_{j-1}(M_{j,\tau,3}^2)\Bigg\}=0\,.
\end{align*}

For $\ell_1=1$ and $\ell_2=2$, since $\Var(G_{n,\tau,2}) \gtrsim p^2n^{-1}\Var^2(\br_1^{\top}\br_1)$, we have
\begin{align*}
&\Var\left\{\sum_{j=1}^n\frac{\E_{j-1}(M_{j,\tau,1}M_{j,\tau,2})}{\Var(G_{n,\tau,2})}\right\}\\	
&~~~~~~~~=
\frac{p^4\Var^2(\br_{1}^\top \br_{1})}{n^{4}\Var^2(G_{n,\tau,2})}
\Bigg[\Var\bigg\{ \frac{2}{n}\sum_{j=\tau+1}^{n-\tau}  ( \br_{j-\tau}^\top \br_{j-\tau} -\alpha_{1p})\sum_{t< j}(\br_t^\top \br_t-\alpha_{1p})\\
&~~~~~~~~~~~~~~~~
+ \frac{2}{n}\sum_{j=n-\tau+1}^{n}  ( \br_{j-\tau}^\top \br_{j-\tau}+\br_{j+\tau}^\top \br_{j+\tau} -2\alpha_{1p})\sum_{t< j}(\br_t^\top \br_t-\alpha_{1p})\bigg\}\Bigg]\\	
&~~~~~~~~\leq
\frac{K}{n^{2}\Var^2(\br_1^\top\br_1)}
\Bigg[\Var\bigg\{ \frac{2}{n}\sum_{j=\tau+1}^{n}  ( \br_{j-\tau}^\top \br_{j-\tau} -\alpha_{1p})\sum_{t< j}(\br_t^\top \br_t-\alpha_{1p})\bigg\}\\
&~~~~~~~~~~~~~~~~
+ \Var\bigg\{\frac{2}{n}\sum_{j=n-\tau+1}^{n} (\br_{j+\tau}^\top \br_{j+\tau} -\alpha_{1p})\sum_{t< j}(\br_t^\top \br_t-\alpha_{1p})\bigg\}\Bigg]\\	
&~~~~~~~~\leq  \frac{K}{\Var^2(\br_{1}^\top \br_{1})}
\max_{j\in[n]}\Var\Bigg\{\frac{2}{n}( \br_{j-\tau}^\top \br_{j-\tau} -\alpha_{1p})\sum_{t< j}(\br_t^\top \br_t-\alpha_{1p})\Bigg\}\\
&~~~~~~~~~~~~~~~~
+  \frac{K\tau^2}{n^2\Var^2(\br_{1}^\top \br_{1})}
\max_{n-\tau+1 \le j\le n}\Var\Bigg\{\frac{2}{n}( \br_{j+\tau}^\top \br_{j+\tau} -\alpha_{1p})\sum_{t< j}(\br_t^\top \br_t-\alpha_{1p})\Bigg\}\\
&~~~~~~~~ \leq \frac{K\E\{(\br_{t}^\top \br_{t}-\alpha_{1p})^4\}}{n\Var^2(\br_{1}^\top \br_{1})}\leq \frac{K}{n}\to0 \,,
\end{align*}
where the last inequality is based on Assumption \ref{as:A3}.

For $\ell_1=1$ and $\ell_2=3$, since $\Var(G_{n,\tau,2}) \gtrsim p^2n^{-1}\Var^2(\br_1^{\top}\br_1)$, we have
\begin{align*}
&\Var\Bigg\{\sum_{j=1}^n\frac{\E_{j-1}(M_{j,\tau,1}M_{j,\tau,3})}{\Var(G_{n,\tau,2})}\Bigg\}	\\
&~~~~ = \frac{p^4{\E}^2\{(\br_{1}^\top \br_{1}-\alpha_{1p})^3\}}{n^{6}\Var^2(G_{n,\tau,2})}
\Var\Bigg\{\sum_{j=\tau+1}^{n}  ( \br_{j-\tau}^\top \br_{j-\tau} -\alpha_{1p})+\sum_{j=n-\tau+1}^{n}  ( \br_{j+\tau}^\top \br_{j+\tau} -\alpha_{1p})\Bigg\}\\
&~~~~\leq \frac{K{\E}^2\{(\br_{1}^\top \br_{1}-\alpha_{1p})^3\}}{n^{3}\Var^4(\br_1^\top \br_1)} 
\Var(\br_{1}^\top \br_{1} -\alpha_{1p}) =  
\frac{K{\E}^2\{(\br_{1}^\top \br_{1}-\alpha_{1p})^3\}\E\{(\br_{1}^\top \br_{1}-\alpha_{1p})^2\}}{n^3\Var^4(\br_1^\top \br_1)}\\
&~~~~\leq  
\frac{K{\E}\{(\br_{1}^\top \br_{1}-\alpha_{1p})^4\}}{n^3\Var^2(\br_1^\top \br_1)}\leq \frac{K}{n^3}\to0 \,,
\end{align*}
where the second inequality is based on the H\"older's inequality, and the last inequality is based on Assumption \ref{as:A3}.

For $\ell_1=2$ and $\ell_2=3$, since $\Var(G_{n,\tau,2}) \gtrsim p^2n^{-1}\Var^2(\br_1^{\top}\br_1)$, we have
\begin{align*}
\Var\Bigg\{\sum_{j=1}^n\frac{\E_{j-1}(M_{j,\tau,2}M_{j,\tau,3}) 
}{\Var(G_{n,\tau,2})}\Bigg\}
=&~ 
\frac{4p^4{\E}^2\{(\br_{1}^\top \br_{1}-\alpha_{1p})^3\}}{n^{8}\Var^2(G_{n,\tau,2})} \Var\Bigg\{\sum_{j=1}^{n} \sum_{t< j}(\br_t^\top \br_t-\alpha_{1p}) \Bigg\}\\
\leq &~ 
\frac{K{\E}^2\{(\br_{1}^\top \br_{1}-\alpha_{1p})^3\}}{n^4\Var^4(\br_1^{\top}\br_1)}
\max_{j\in [n]}\Var\Bigg\{\sum_{t< j}(\br_t^\top \br_t-\alpha_{1p}) \Bigg\}\\	
\leq &~ 
\frac{K{\E}^2\{(\br_{1}^\top \br_{1}-\alpha_{1p})^3\}\E\{(\br_{1}^\top \br_{1}-\alpha_{1p})^2\}}{n^3\Var^4(\br_1^{\top}\br_1)} \\
\leq&~ \frac{K{\E}\{(\br_{1}^\top \br_{1}-\alpha_{1p})^4\}}{n^3\Var^2(\br_1^{\top}\br_1)} \le \frac{K}{n^3}\to0 \,,
\end{align*}
where the third inequality is based on the H\"older's inequality, and the last inequality is based on Assumption \ref{as:A3}. $\hfill\Box$

\subsubsection{Proof of \eqref{mclt-conditionb2}}\label{subsubsec:b2}
Recall, as defined in \eqref{mjtl}, $M_{j,\tau}=M_{j,\tau,1}+M_{j,\tau,2}+M_{j,\tau,3}$. By definition, for any $\tau > \tau'$, we have
\begin{align*}
&~\E_{j-1}(M_{j,\tau,1}M_{j,\tau',1}) \\
=&~ \frac {p^2}{n^2}\Var(\br_{1}^\top \br_{1}) \cdot \big\{(\br_{j-\tau}^\top \br_{j-\tau} -\alpha_{1p})( \br_{j-\tau'}^\top \br_{j-\tau'} -\alpha_{1p}) I(\tau+1\leq j\leq n-\tau) \\
&~~~+ (\br_{j-\tau}^\top \br_{j-\tau}+\br_{j+\tau}^\top \br_{j+\tau} -2\alpha_{1p})( \br_{j-\tau'}^\top \br_{j-\tau'} -\alpha_{1p})I(n-\tau+1\leq j\leq n-\tau') \\
&~~~+ (\br_{j-\tau}^\top \br_{j-\tau}+\br_{j+\tau}^\top \br_{j+\tau} -2\alpha_{1p})( \br_{j-\tau'}^\top \br_{j-\tau'}+\br_{j+\tau'}^\top \br_{j+\tau'} -2\alpha_{1p}) I(n-\tau'+1\leq j\leq n)\big\} \,.
\end{align*}
Notice that $M_{j,\tau,2}$ and $M_{j,\tau,3}$ do not depend on $\tau$. Hence, for any $\tau > \tau'$, it holds that
\begin{align*}
    &~\E_{j-1}\{M_{j,\tau,1}(M_{j,\tau',2}+M_{j,\tau',3})\} = \E_{j-1}\{M_{j,\tau,1}(M_{j,\tau,2}+M_{j,\tau,3})\}\\ 
 = & -\frac {2p^2}{n^3}\Var(\br_{1}^\top \br_{1}) \cdot ( \br_{j-\tau}^\top \br_{j-\tau} -\alpha_{1p})\sum_{t< j}(\br_t^\top \br_t-\alpha_{1p})I(\tau+1\leq j\leq n-\tau) \\
& -\frac {2p^2}{n^3}\Var(\br_{1}^\top \br_{1}) \cdot( \br_{j-\tau}^\top \br_{j-\tau}+\br_{j+\tau}^\top \br_{j+\tau} -2\alpha_{1p})\sum_{t< j}(\br_t^\top \br_t-\alpha_{1p}) I(n-\tau+1\leq j\leq n) \\
& -\frac {p^2}{n^3}\E\{(\br_{1}^\top \br_{1}-\alpha_{1p})^3\} \cdot (\br_{j-\tau}^\top \br_{j-\tau} -\alpha_{1p}) I(\tau+1\leq j\leq n-\tau) \\
& -\frac {p^2}{n^3}\E\{(\br_{1}^\top \br_{1}-\alpha_{1p})^3\} \cdot (\br_{j-\tau}^\top \br_{j-\tau}+\br_{j+\tau}^\top \br_{j+\tau} -2\alpha_{1p}) I(n-\tau+1\leq j\leq n)\,.
\end{align*}
and
\begin{align*}
    &~\E_{j-1}\{M_{j,\tau',1}(M_{j,\tau,2}+M_{j,\tau,3})\} \\ 
 = & -\frac {2p^2}{n^3}\Var(\br_{1}^\top \br_{1}) \cdot ( \br_{j-\tau'}^\top \br_{j-\tau'} -\alpha_{1p})\sum_{t< j}(\br_t^\top \br_t-\alpha_{1p})I(\tau'+1\leq j\leq n-\tau') \\
& -\frac {2p^2}{n^3}\Var(\br_{1}^\top \br_{1}) \cdot( \br_{j-\tau'}^\top \br_{j-\tau'}+\br_{j+\tau'}^\top \br_{j+\tau'} -2\alpha_{1p})\sum_{t< j}(\br_t^\top \br_t-\alpha_{1p}) I(n-\tau'+1\leq j\leq n) \\
& -\frac {p^2}{n^3}\E\{(\br_{1}^\top \br_{1}-\alpha_{1p})^3\} \cdot (\br_{j-\tau'}^\top \br_{j-\tau'} -\alpha_{1p}) I(\tau'+1\leq j\leq n-\tau') \\
& -\frac {p^2}{n^3}\E\{(\br_{1}^\top \br_{1}-\alpha_{1p})^3\} \cdot (\br_{j-\tau'}^\top \br_{j-\tau'}+\br_{j+\tau'}^\top \br_{j+\tau'} -2\alpha_{1p}) I(n-\tau'+1\leq j\leq n)\,.
\end{align*}
Analogously, we have
\begin{align*}
&~\E_{j-1}\{(M_{j,\tau,2}+M_{j,\tau,3})(M_{j,\tau',2}+M_{j,\tau',3})\} = \E_{j-1}\{(M_{j,\tau,2}+M_{j,\tau,3})^2 \}\\ 
=&~\frac {4p^2}{n^4}\Var(\br_{1}^\top \br_{1}) \cdot\left\{\sum_{t< j}(\br_t^\top \br_t-\alpha_{1p})\right\}^2 +\frac {4p^2}{n^4}\E\{(\br_{1}^\top \br_{1}-\alpha_{1p})^3\} \cdot \sum_{t< j}(\br_t^\top \br_t-\alpha_{1p})\\
&~~+\frac{p^2}{n^4}\E\big[\{(\br_1^\top \br_1-\alpha_{1p})^2-\Var(\br_1^\top \br_1)\}^2 \big]\,.
\end{align*}
It follows that $\E( M_{j,\tau,1}M_{j,\tau',1})=0$ and
\begin{align*}
\E \{M_{j,\tau,1}(M_{j,\tau',2}+M_{j,\tau',3})\}
&=-\frac {2p^2}{n^3}\Var^2(\br_{1}^\top \br_{1}) \big\{I(\tau+1\leq j\leq n-\tau) \\
&~~~~~~~~~~~~~~~ + 2 I(n-\tau+1\leq j\leq n) \big\}\,,\\
\E \{M_{j,\tau',1}(M_{j,\tau,2}+M_{j,\tau,3})\}
&=-\frac {2p^2}{n^3}\Var^2(\br_{1}^\top \br_{1}) \big\{I(\tau'+1\leq j\leq n-\tau') \\
&~~~~~~~~~~~~~~~ + 2 I(n-\tau'+1\leq j\leq n) \big\}\,,\\
\E \{(M_{j,\tau,2}+M_{j,\tau,3})^2\}
&=\frac {4(j-1)p^2}{n^4}\Var^2(\br_{1}^\top \br_{1})
+\frac{p^2}{n^4}\E\big[\{(\br_1^\top \br_1-\alpha_{1p})^2-\Var(\br_1^\top \br_1)\}^2\big]\,.
\end{align*}
Based on these results, since $\Var(G_{n,\tau,2}) \gtrsim p^2n^{-1}\Var^2(\br_1^{\top}\br_1)$, we obtain
\begin{align*}
\bigg|\E\Bigg\{\sum_{j=1}^n\frac{\E_{j-1}(M_{j,\tau}M_{j,\tau'})}{\Var(G_{n,\tau,2})}\Bigg\}\bigg|
=&~\bigg|\frac{p^2 \E\{(\br_j^\top \br_j-\alpha_{1p})^4\}-p^2n(2+3n^{-1})\Var^2(\br_1^\top \br_1)}{n^{3}\Var(G_{n,\tau,2})}\bigg|\\
\leq &~ \frac{Kp^2 \E\{(\br_j^\top \br_j-\alpha_{1p})^4\}}{n^{3}\Var(G_{n,\tau,2})} + \frac{Kp^2 n(2+3n^{-1})\Var^2(\br_1^\top \br_1)}{n^{3}\Var(G_{n,\tau,2})}\\
\leq &~ \frac{Kp^2 \E\{(\br_j^\top \br_j-\alpha_{1p})^4\}}{n^{2}\Var(G_{n,\tau,2})} \leq \frac{ K\E\{(\br_j^\top \br_j-\alpha_{1p})^4\}}{n\Var^2(\br_1^\top \br_1)} \leq \frac{K}{n}\to0\,,
\end{align*}
where the last inequality is based on Assumption \ref{as:A3}. 
Moreover, following similar discussions as in the proof of \eqref{mjl123}, one can verify that 
\begin{align*}
\Var\Bigg\{\sum_{j=1}^n\frac{\E_{j-1}(M_{j,\tau}M_{j,\tau'})}{\Var(G_{n,\tau,2})}\Bigg\}
\leq &~ K\sum_{\ell_1,\ell_2=1}^3 \Var\Bigg\{\sum_{j=1}^n\frac{\E_{j-1}(M_{j,\tau,\ell_1}M_{j,\tau',\ell_2})}{\Var(G_{n,\tau,2})}\Bigg\}
\leq \frac{K}{n}\to0 \,.
\end{align*}
Hence, we have \eqref{mclt-conditionb2} holds. $\hfill\Box$

\subsubsection{Proof of \eqref{mclt-conditionb3}}\label{subsubsec:b3}

Since $\Var(G_{n,\tau,2}) \gtrsim p^2n^{-1}\Var^2(\br_1^{\top}\br_1)$, it holds that
\begin{align*}
\sum_{j=1}^n\frac{\E (M_{j,\tau}^4)}{\Var^2(G_{n,\tau,2})}
\leq &~
\frac{ K}{\Var^2(G_{n,\tau,2})}\sum_{j=1}^n
\E( M_{j,\tau,1}^4+ M_{j,\tau,2}^4+ M_{j,\tau,3}^4)\\
\leq &~\frac{ K}{\Var^2(G_{n,\tau,2})}\bigg(
\frac{p^4}{n^3}{\E}^2\{(\br_1^\top \br_1-\alpha_{1p})^4\}+\frac{p^4}{n^5}{\E}^2\{(\br_1^\top \br_1-\alpha_{1p})^4\}\\
&\qquad~~~~~~~~~~~~~~~~~~~~~~~+\frac{p^4}{n^7}\E\big[\{(\br_1^\top \br_1-\alpha_{1p})^2-\Var(\br_1^\top \br_1)\}^4\big]\bigg)\\
\leq &~\frac{ Kp^4\E\{(\br_1^\top \br_1-\alpha_{1p})^8\}}{n^{3}\Var^2(G_{n,\tau,2})} \le \frac{ K\E\{(\br_1^\top \br_1-\alpha_{1p})^8\}}{n\Var^4(\br_1^{\top}\br_1)} \leq \frac{K}{n}\to 0\,,
\end{align*}
where the third inequality is based on the H\"older's inequality, and the last inequality is based on Assumption \ref{as:A3}.
Hence, we have \eqref{mclt-conditionb3} holds. $\hfill\Box$

\subsection{Proof of Step 3}\label{subsec:Ind_Gq12}

Recall that $G_{n,\tau,1}=G_{n,\tau,1,1}+O_{\rm p}(n^{-1/2})$ and $\E(G_{n,\tau,1,1})=0$. By \eqref{gnt11}, we have
\begin{align*}
G_{n,\tau,1,1} = &\sum_{j=\tau+1}^nD_{j,\tau} ~\,\text{with}~\, D_{j,\tau}= \frac {2p}{n}\sum_{j-\tau\neq t< j}\br_{j-\tau}^\top \br_{t-\tau}\br_t^\top \br_j\,.
\end{align*}
Let 
 \begin{align*}
\tilde G_{n,\tau,1,1}
= &\sum_{j=\tau+1}^n \tilde D_{j,\tau} ~\,\text{with} ~\, \tilde D_{j,\tau}= \frac {2p}{n}\sum_{t< j-\tau}\br_{j-\tau}^\top \br_{t-\tau}\br_t^\top \br_j\,.
\end{align*}
Since $\E(\tilde G_{n,\tau,1,1})=0$, and
 \begin{align*}
\E\{( G_{n,\tau,1,1}-\tilde G_{n,\tau,1,1})^2\}
=&~\frac {4p^2}{n^2}\E\bigg\{\bigg(\sum_{j=\tau+1}^n \sum_{j-\tau<t<j}\br_{j-\tau}^\top \br_{t-\tau}\br_t^\top \br_j\bigg)^2\bigg\}\\
=&~\frac {4p^2}{n^2}\sum_{j=\tau+1}^n \sum_{j-\tau<t<j}\E\{(\br_{j-\tau}^\top \br_{t-\tau}\br_t^\top \br_j)^2\}\\
=&~\frac {4}{n^2}\sum_{j=\tau+1}^n \sum_{j-\tau<t<j} \alpha_{2p}^2\leq \frac{K}{n} \to 0\,,
 \end{align*}
it holds that $G_{n,\tau,1}=\tilde G_{n,\tau,1,1}+O_{\rm p}(n^{-1/2})$. Moreover, $\{\tilde D_{j,\tau}\}$ also forms a sequence of martingale differences with respect to $\{\mathcal F_j\}$, where $\mathcal F_j$ is the $\sigma$-field generated by $\{\x_1,\ldots,\x_j\}$. 
In addition, by \eqref{eq:gn2}, we have
\begin{align*}
    G_{n,\tau,2}-\mu_2 = \sum_{j=1}^nM_{j,\tau}\,,
\end{align*}
where $\{M_{j,\tau}\}$ is a sequence of martingale differences with respect to $\{\mathcal F_j\}$. Hence, to prove \eqref{eq:JointAsymNorm}, by Lemma \ref{lem:md-clt-vector}, it suffices to show 
\begin{align}
\sum_{j=\tau+1}^n\frac{\E_{j-1}(\tilde{D}_{j,\tau}^2)}{\sigma_{n1}^2} & \xrightarrow{\mathrm{p}} 1 ~\,\mbox{for any given}~\, \tau \in [q]\,, \label{mclt-jointcondition1} \\ 
\sum_{j= \tau +1}^n\frac{\E_{j-1}(\tilde{D}_{j,\tau}\tilde{D}_{j,\tau'})}{\sigma_{n1}^2} & \xrightarrow{\mathrm{p}} 0~\,\mbox{for any given}~\, (\tau,\tau')\in [q]^2 ~\mbox{satisfying}~ \tau > \tau'\,, \label{mclt-jointcondition2}\\
\sum_{j=\tau+1}^n\frac{\E (\tilde{D}_{j,\tau}^4)}{\sigma_{n1}^4} & \to0 ~\,\mbox{for any given}~\, \tau \in [q]\,.\label{mclt-jointcondition3} \\
\sum_{j=1}^n\frac{\E_{j-1}(M_{j,\tau}^2)}{\sigma_{n2}^2} & \xrightarrow{\mathrm{p}} 1~\,\mbox{for any given}~\, \tau \in [q]\,, \label{mclt-jointcondition4} \\
\sum_{j=1}^n\frac{\E_{j-1}(M_{j,\tau}M_{j,\tau'})}{\sigma_{n2}^2} & \xrightarrow{\mathrm{p}} 0 ~\, \mbox{for any given} ~\, (\tau,\tau')\in [q]^2 ~\mbox{satisfying}~ \tau> \tau'\,, \label{mclt-jointcondition5}\\
\sum_{j=1}^n\frac{\E (M_{j,\tau}^4)}{\sigma_{n2}^4} & \to0 ~\,\mbox{for any given}~\, \tau \in [q]\,.\label{mclt-jointcondition6} \\
\sum_{j=1}^n\frac{\E_{j-1} (\tilde D_{j,\tau}M_{j,\tau'})}{\sigma_{n1}\sigma_{n2}} & \xrightarrow{\mathrm{p}} 0 ~\, \mbox{for any given}~\, \tau, \tau' \in [q]\,, \label{dmj}
\end{align}
where $\tilde D_{j,\tau} = 0$ for $j \in [\tau]$.
Using the similar arguments for proving \eqref{mclt-conditiona1}--\eqref{mclt-conditiona3}, we can also show \eqref{mclt-jointcondition1}--\eqref{mclt-jointcondition3} hold. \eqref{mclt-jointcondition4}--\eqref{mclt-jointcondition6} are actually the same as \eqref{mclt-conditionb1}--\eqref{mclt-conditionb3}, which have been shown in Sections \ref{subsubsec:b1}--\ref{subsubsec:b3}. As we will show in Section \ref{subsubsec:c4}, \eqref{dmj} holds. This completes the proof of \eqref{eq:JointAsymNorm}. $\hfill\Box$

\subsubsection{Proof of \eqref{dmj}}\label{subsubsec:c4}
Recall, as defined in \eqref{mjtl}, $M_{j,\tau}=M_{j,\tau,1}+M_{j,\tau,2}+M_{j,\tau,3}$. By definition, for any $\tau, \tau' \in [q]$, we have
\begin{align*}
 \tilde D_{j,\tau}M_{j,\tau'}
= \tilde D_{j,\tau}M_{j,\tau',1}+ \tilde D_{j,\tau}M_{j,\tau',2}+ \tilde D_{j,\tau}M_{j,\tau',3}\,,
\end{align*}
where
\begin{align*}	
 \tilde D_{j,\tau}M_{j,\tau',1}=&~\frac {2p^2}{n^2}\sum_{t< j-\tau}\br_{j-\tau}^\top \br_{t-\tau}\br_t^\top \br_j \big\{ (\br_{j}^\top \br_{j}-\alpha_{1p})( \br_{j-\tau'}^\top \br_{j-\tau'} -\alpha_{1p}) I(\tau'+1\leq j\leq n-\tau')   \nonumber \\
 &~~~~~~~~~~~ + (\br_{j}^\top \br_{j}-\alpha_{1p})( \br_{j-\tau'}^\top \br_{j-\tau'}+\br_{j+\tau'}^\top \br_{j+\tau'} -2\alpha_{1p}) I(n-\tau'+1\leq j\leq n) \big\}\,, \nonumber\\
 \tilde D_{j,\tau}M_{j,\tau',2}=&-\frac {4p^2}{n^3}\sum_{t_1< j-\tau}\br_{j-\tau}^\top \br_{t_1-\tau}\br_{t_1}^\top \br_j(\br_{j}^\top \br_{j}-\alpha_{1p})\sum_{t_2< j}(\br_{t_2}^\top \br_{t_2}-\alpha_{1p})\,,\\
 \tilde D_{j,\tau}M_{j,\tau',3}=&- \frac {2p^2}{n^3}\sum_{t< j-\tau}\br_{j-\tau}^\top \br_{t-\tau}\br_t^\top \br_j \big\{(\br_j^\top \br_j-\alpha_{1p})^2-\Var(\br_j^\top \br_j)\big\}\,.
\end{align*}
It is easy to check that $\E(\tilde D_{j,\tau}M_{j,\tau',\ell})=0$ for $\ell=1,2,3$.
Recall $\sigma_{n1}^2=2\alpha_{2p}^2 \to 2\alpha_2^2$ as $p \to \infty$, and $\sigma_{n2}^2=\Var(G_{n,\tau,2})$. To show \eqref{dmj}, it suffices to show
\begin{align*}
 \Var\bigg\{ \sum_{j=1}^n\frac{\E_{j-1} (\tilde D_{j,\tau}M_{j,\tau'})}{\sqrt{\Var(G_{n,\tau,2})}}\bigg\} 
\leq K\sum_{\ell=1}^3\Var\bigg\{ \sum_{j=1}^n\frac{\E_{j-1} (\tilde D_{j,\tau}M_{j,\tau',\ell})}{\sqrt{\Var(G_{n,\tau,2})}}\bigg\}\to0\,.
\end{align*}
To do this, we first bound the variance of $\sum_{j=1}^n\E_{j-1} (\tilde D_{j,\tau}M_{j,\tau',\ell})$ for $\ell=1,2,3$. 

For $\ell=1$, we have
\begin{align*}
\E\bigg[\bigg\{ 
\sum_{j=1}^n\E_{j-1} (\tilde D_{j,\tau}M_{j,\tau',1} )\bigg\}^2\bigg]\leq 2(R_{11}+R_{12})\,,
\end{align*}
where
\begin{align*}
 R_{11}= &~\frac{4p^4}{n^{4}}
 \E \bigg(\bigg[\sum_{j=\tau'+1}^{n-\tau'}\sum_{t<j-\tau }\br_{j-\tau}^\top \br_{t-\tau}\br_t^\top \E\{\br_j (\br_{j}^\top \br_{j}-\alpha_{1p})\}( \br_{j-\tau'}^\top \br_{j-\tau'} -\alpha_{1p})\bigg]^2\bigg)\,,\\
 R_{12}= &~\frac{4p^4}{n^{4}}
\E\bigg(\bigg[\sum_{j=n-\tau'+1}^{n}\sum_{t<j-\tau}\br_{j-\tau}^\top \br_{t-\tau}\br_t^\top \E\{\br_j (\br_{j}^\top \br_{j}-\alpha_{1p})\}( \br_{j-\tau'}^\top \br_{j-\tau'}+\br_{j+\tau'}^\top \br_{j+\tau'} -2\alpha_{1p})
\bigg]^2\bigg)\,.
\end{align*}
As we will show in Sections \ref{subsec:R11} and \ref{subsec:R12}, for any $\tau, \tau' \in [q]$,
\begin{align}
     R_{11}& \leq  K(n^{-1}+p^2n^{-2})\Var^2(\br_1^\top \br_1)\,, \label{eq:R11_all}\\
    R_{12}&\leq  Kp^2n^{-2}\Var^2(\br_1^\top \br_1)\,,\label{eq:R12_all}
\end{align} 
which implies
\begin{align}\label{eq:DM1}
\E\bigg[\bigg\{ 
\sum_{j=1}^n\E_{j-1} (\tilde D_{j,\tau}M_{j,\tau',1} )\bigg\}^2\bigg]\leq K\bigg(\frac{1}{n}+\frac{p^2}{n^{2}}\bigg)\Var^2(\br_1^\top \br_1)
\end{align}
for any $\tau, \tau' \in [q]$.

For $\ell = 2$, we have
\begin{align*}
&\E\bigg[\bigg\{ 
\sum_{j=1}^n\E_{j-1}( \tilde D_{j,\tau}M_{j,\tau',2}) \bigg\}^2\bigg]\\
&~~~~~= \frac {16p^4}{n^6}\E\bigg(\bigg[ 
\sum_{j=\tau+1}^n\sum_{t_1<j-\tau \atop t_2< j}\br_{j-\tau}^\top \br_{t_1-\tau}\br_{t_1}^\top \E\{\br_j(\br_{j}^\top \br_{j}-\alpha_{1p})\}(\br_{t_2}^\top \br_{t_2}-\alpha_{1p})\bigg]^2\bigg) \\
&~~~~~\leq 2(R_{21}+R_{22})\,,
\end{align*}
where 
\begin{align*}
R_{21}
=&~\frac {16p^4}{n^6}\E\bigg(\bigg[
\sum_{j=\tau+1}^n\sum_{t_1<j-\tau \atop t_2< j-\tau}\br_{j-\tau}^\top \br_{t_1-\tau}\br_{t_1}^\top \E\{\br_j(\br_{j}^\top \br_{j}-\alpha_{1p})\}(\br_{t_2}^\top \br_{t_2}-\alpha_{1p})\bigg]^2\bigg),\\
R_{22}
=&~\frac {16p^4}{n^6}\E\bigg(\bigg[ 
\sum_{j=\tau+1}^n\sum_{t_1<j-\tau \atop j-\tau\leq t_2 <j}\br_{j-\tau}^\top \br_{t_1-\tau}\br_{t_1}^\top \E\{\br_j(\br_{j}^\top \br_{j}-\alpha_{1p})\}(\br_{t_2}^\top \br_{t_2}-\alpha_{1p})\bigg]^2\bigg)\,.	
\end{align*}
As we will show in Section \ref{subsec:R2}, 
\begin{align} 
R_{21}  
\leq Kp^2n^{-3}
\Var^2(\br_1^\top \br_1) ~~\mbox{and}~~
R_{22}  
\leq Kpn^{-2}
\Var^2(\br_1^\top \br_1) \label{eq:R2_all}
\end{align}
for any $\tau, \tau'\in [q]$, which implies
\begin{align}\label{eq:DM2}
&\E \bigg[\bigg\{ 
\sum_{j=1}^n\E_{j-1}( \tilde D_{j,\tau}M_{j,\tau',2}) \bigg\}^2 \bigg]\leq K\bigg(\frac{p^2}{n^{3}}+\frac{p}{n^{2}}\bigg)\Var^2(\br_1^\top \br_1)
\end{align}
for any $\tau, \tau'\in [q]$.

For $\ell=3$, by Assumption \ref{as:A4}, we have
\begin{align}\label{eq:DM3}
\E\bigg[\bigg\{
\sum_{j=1}^n\E_{j-1} (\tilde D_{j,\tau}M_{j,\tau',3})
\bigg\}^2\bigg] &=\frac {4p^4}{n^6}\E\bigg(\bigg[
\sum_{j=\tau+1}^n\sum_{t<j-\tau}\br_{j-\tau}^\top \br_{t-\tau}\br_t^\top
 \E\{\br_j (\br_j^\top \br_j-\alpha_{1p})^2\}
\bigg]^2\bigg) \notag\\
&=\frac {4p^4}{n^6}\sum_{j=\tau+1}^n\sum_{t<j-\tau}\E\big([ \br_{j-\tau}^\top \br_{t-\tau}\br_t^\top \E\{\br_j (\br_j^\top \br_j-\alpha_{1p})^2\}
]^2\big) \notag\\
&\leq \frac {Kp^4}{n^4}
\max_{\tau+1 \le j \le n \atop t<j-\tau}
\E\big([ \br_{j-\tau}^\top \br_{t-\tau}\br_t^\top \E\{\br_j (\br_j^\top \br_j-\alpha_{1p})^2\}
]^2\big) \notag\\
&\leq \frac {Kp^3}{n^4}
\max_{\tau+1 \le j \le n \atop t< j-\tau}\E\big([\br_t^\top \E\{\br_j (\br_j^\top \br_j-\alpha_{1p})^2\}
]^2\big) \notag\\
&\leq \frac {Kp^2}{n^4} \E\{(\br_1^\top \br_1-\alpha_{1p})^4\} \le \frac {Kp^2}{n^4}\Var^2(\br_1^\top \br_1)
\end{align}
for any $\tau, \tau'\in [q]$.

Combining \eqref{eq:DM1}, \eqref{eq:DM2}, and \eqref{eq:DM3}, since $\Var(G_{n,\tau,2}) \gtrsim p^2n^{-1}\Var^2(\br_1^{\top}\br_1)$, 
it yields
\begin{align*}
\Var\Bigg\{ \sum_{j=1}^n\frac{\E_{j-1} (\tilde D_{j,\tau}M_{j,\tau'})}{\sqrt{\Var(G_{n,\tau,2})}}\Bigg\} & \leq K\sum_{\ell=1}^3\Var\Bigg\{ \sum_{j=1}^n\frac{\E_{j-1} (\tilde D_{j,\tau}M_{j,\tau',\ell})}{\sqrt{\Var(G_{n,\tau,2})}}\Bigg\}\\
&\leq  \frac{K(n^{-1}+p^2n^{-2}) \Var^2(\br_1^\top \br_1)}{\Var(G_{n,\tau,2})}=O(n^{-1}+p^{-2})\to0
\end{align*}
as $\min\{p,n\} \to \infty$. 
Hence, we have \eqref{dmj} holds. $\hfill\Box$

\subsubsection{Proof of \eqref{eq:R11_all}}\label{subsec:R11}
For $R_{11}$ with $\tau=\tau'$, the indices satisfy $t-\tau<t<j-\tau=j-\tau'$, from which we obtain
\begin{align*}
 R_{11}
 = &~\frac{4p^4}{n^{4}}
 \E \Bigg[\sum_{j_1, j_2=\tau+1}^{n-\tau}\sum_{t<j_1-\tau\atop t<j_2-\tau }
 \br_{j_1-\tau}^\top ( \br_{j_1-\tau}^\top \br_{j_1-\tau} -\alpha_{1p})\br_{t-\tau}\br_t^\top \E\{\br_j (\br_{j}^\top \br_{j}-\alpha_{1p})\}\\
 &~~~~~~~~~~~~~~~~~~~~~~~~~~~~~~~\times  \br_{j_2-\tau}^\top ( \br_{j_2-\tau}^\top \br_{j_2-\tau} -\alpha_{1p})\br_{t-\tau}\br_t^\top \E\{\br_j (\br_{j}^\top \br_{j}-\alpha_{1p})\}\Bigg]\\
  = &~\frac{4p^2}{n^{4}}
\sum_{j_1, j_2=\tau+1}^{n-\tau}\sum_{t<j_1-\tau\atop t<j_2-\tau }
\big[ \E \{\br_{j_1-\tau}^\top ( \br_{j_1-\tau}^\top \br_{j_1-\tau} -\alpha_{1p})\bSigma
( \br_{j_2-\tau}^\top \br_{j_2-\tau} -\alpha_{1p}) \br_{j_2-\tau}\}\\
 &~~~~~~~~~~~~~~~~~~~~~~~~~~~~~~~\times  \E\{\br_j^\top  (\br_{j}^\top \br_{j}-\alpha_{1p})\} \bSigma \E\{\br_j (\br_{j}^\top \br_{j}-\alpha_{1p})\} \big]\,.
\end{align*}
Due to $\E(\br_1)=0$, by Assumption \ref{as:A4}, we have
\begin{align*}
    \E\{\br_1^\top  (\br_{1}^\top \br_{1}-\alpha_{1p})\} \bSigma \E\{\br_1 (\br_{1}^\top \br_{1}-\alpha_{1p})\} \le K\big|\E(\br_1\br_1^\top \br_1)\big|_2^2\,.
\end{align*}
By the Cauchy-Schwarz inequality and Assumption \ref{as:A3}, it holds that
\begin{align}
\E \{\br_1^\top \br_1( \br_1^\top \br_1 -\alpha_{1p})^2\}
\leq &~K{\E}^{1/2} \{(\br_1^\top \br_1)^2\}{\E}^{1/2}\{( \br_1^\top \br_1 -\alpha_{1p})^4\}
\leq K\Var(\br_1^\top \br_1)\,,\label{var111}
\end{align}
where the last inequality is based on the fact that ${\E}\{(\br_1^\top \br_1)^2\} \le K$ as shown in \eqref{moments}. Also, by Assumption \ref{as:A3},
\begin{align}
\big|\E(\br_1\br_1^\top \br_1)\big|_2^2	
=&~\frac1{p^3}\big|\E(\z_1\z_1^\top \z_1)\big|_2^2\leq \frac {K}{p^3}\Var(\z_1^\top \z_1)=\frac K{p}\Var(\br_1^\top \br_1)\,.
\label{var112}
\end{align}
Therefore, we have
\begin{align*}
 &\E \{\br_{j_1-\tau}^\top ( \br_{j_1-\tau}^\top \br_{j_1-\tau} -\alpha_{1p})\bSigma
( \br_{j_2-\tau}^\top \br_{j_2-\tau} -\alpha_{1p}) \br_{j_2-\tau}\}\\
&~~~~~~ \leq K \times
\begin{cases}
\E \{\br_1^\top \br_1( \br_1^\top \br_1 -\alpha_{1p})^2\}\,, &\mbox{if}~ j_1=j_2\,,\\
\big|\E(\br_1\br_1^\top \br_1) \big|_2^2\,, &\mbox{if}~ j_1\neq j_2\,,
\end{cases}\\
&~~~~~~ \leq K \times
\begin{cases}
\Var(\br_1^\top \br_1)\,,& \mbox{if}~j_1=j_2\,,\\
p^{-1} \Var(\br_1^\top \br_1)\,,&\mbox{if}~ j_1\neq j_2\,,
\end{cases}
\end{align*}
which implies that $R_{11}$ with $\tau=\tau'$ can be bounded as
\begin{align}\label{R111}
 R_{11}
 \leq&~ 
 \frac{Kp^2}{n^{4}}
 \big\{
 n^3p^{-2} \Var^2(\br_1^\top \br_1)+n^2p^{-1}\Var^2(\br_1^\top \br_1)
\big\}\nonumber\\
\leq&~ 
 K( n^{-1}+pn^{-2})\Var^2(\br_1^\top \br_1)\,.
\end{align}

For $R_{11}$ with $\tau>\tau'$, the indices satisfy $t-\tau<t<j-\tau<j-\tau'$, from which we have
\begin{align}\label{eq:R11_22}
 R_{11}
  = &~\frac{4p^4}{n^{4}}
\sum_{j=\tau'+1}^{n-\tau'}\sum_{t<j-\tau } \E \big([\br_{j-\tau}^\top \br_{t-\tau}\br_t^\top \E\{\br_j (\br_{j}^\top \br_{j}-\alpha_{1p})\}( \br_{j-\tau'}^\top \br_{j-\tau'} -\alpha_{1p})]^2 \big)\notag\\	
  \leq &~\frac{Kp^2}{n^{4}}
\sum_{j=\tau'+1}^{n-\tau'}\sum_{t<j-\tau } \big|\E\{\br_j (\br_{j}^\top \br_{j}-\alpha_{1p})\} \big|_2^2 \Var( \br_{j-\tau'}^\top \br_{j-\tau'} ) \notag\\	
  \leq &~\frac{Kp}{n^{2}}\Var^2( \br_{1}^\top \br_{1})\,.
\end{align}
Here the last inequality holds by \eqref{var112}.

For $R_{11}$ with $\tau < \tau'$, since $\E(\br_t)=0$, by Assumption \ref{as:A4}, we have
\begin{align*}
 R_{11}
  & = \frac{4p^4}{n^{4}}
\sum_{j=\tau'+1}^{n-\tau'}\sum_{t<j-\tau } \E \big([\br_{j-\tau}^\top \br_{t-\tau}\br_t^\top \E\{\br_j (\br_{j}^\top \br_{j}-\alpha_{1p})\}( \br_{j-\tau'}^\top \br_{j-\tau'} -\alpha_{1p})]^2 \big)\notag\\	
  & \leq \frac{Kp^4}{n^2}
\max_{\tau'+1 \le j \le n-\tau' \atop t<j-\tau} \E \big([\br_{j-\tau}^\top \br_{t-\tau}\br_t^\top \E\{\br_j (\br_{j}^\top \br_{j}-\alpha_{1p})\}( \br_{j-\tau'}^\top \br_{j-\tau'} -\alpha_{1p})]^2 \big) \notag\\
  & \leq \frac{Kp^3}{n^2}
\max_{\tau'+1 \le j \le n-\tau' \atop t<j-\tau} \E \big[\br_{t-\tau}^\top \br_{t-\tau}\{\br_t^\top \E(\br_j\br_{j}^\top \br_{j})\}^2( \br_{j-\tau'}^\top \br_{j-\tau'} -\alpha_{1p})^2 \big]\,. 
\end{align*}
Notice that
\begin{align*}
  &\E \big[\br_{t-\tau}^\top \br_{t-\tau}\{\br_t^\top \E(\br_j\br_{j}^\top \br_{j})\}^2( \br_{j-\tau'}^\top \br_{j-\tau'} -\alpha_{1p})^2 \big] \\
  &~~~~~ \le \E \{\br_{t-\tau}^\top \br_{t-\tau}( \br_{t-\tau}^\top \br_{t-\tau} -\alpha_{1p})^2 \} \E[\{\br_t^\top \E(\br_j\br_{j}^\top \br_{j})\}^2] \notag \\
 &~~~~~~~~~+ \E(\br_{t-\tau}^\top \br_{t-\tau})\E\big[\{\br_t^\top \E(\br_j\br_{j}^\top \br_{j})\}^2( \br_{j-\tau'}^\top \br_{j-\tau'} -\alpha_{1p})^2\big] \\
  &~~~~~ \le \frac{K}{p}\Var( \br_{1}^\top \br_{1})\big|\E(\br_1\br_1^\top \br_1) \big|_2^2 + K\E[\{\br_t^\top \E(\br_j\br_{j}^\top \br_{j})\}^2]\E\{( \br_{j-\tau'}^\top \br_{j-\tau'} -\alpha_{1p})^2\} \\
  &~~~~~~~~~+ K\E\big[\{\br_t^\top \E(\br_j\br_{j}^\top \br_{j})\}^2( \br_{t}^\top \br_{t} -\alpha_{1p})^2\big] \\
  &~~~~~ \le \frac{K}{p}\Var( \br_{1}^\top \br_{1})\big|\E(\br_1\br_1^\top \br_1) \big|_2^2 + K\E\{\br_t^\top \br_t( \br_{t}^\top \br_{t} -\alpha_{1p})^2\}\big|\E(\br_1\br_1^\top \br_1) \big|_2^2 \\
  &~~~~~ \le \frac{K}{p^2}\Var^2( \br_{1}^\top \br_{1}) + \frac{K}{p}\Var^2( \br_{1}^\top \br_{1}) \le \frac{K}{p}\Var^2( \br_{1}^\top \br_{1})\,,
\end{align*}
where the second inequality is based on the fact that ${\E}\{(\br_1^\top \br_1)^2\} \le K$ as shown in \eqref{moments}, and the third inequality is based on the Cauchy-Schwarz inequality. Hence, we have
\begin{align*}
 R_{11}\leq \frac{Kp^2}{n^{2}}\Var^2( \br_{1}^\top \br_{1})
\end{align*}
for $\tau < \tau'$.
Together with \eqref{R111} and \eqref{eq:R11_22}, we have \eqref{eq:R11_all} holds for any $\tau, \tau' \in [q]$. $\hfill\Box$

\subsubsection{Proof of \eqref{eq:R12_all}}\label{subsec:R12}
By the Cauchy-Schwarz inequality, we have
\begin{align*}
R_{12}
\leq &~\frac{Kp^4}{n^{2}}
\max_{n-\tau'+1\leq j\leq n \atop t<j-\tau}
\E\big([ \br_{j-\tau}^\top \br_{t-\tau}\br_t^\top \E\{\br_j (\br_{j}^\top \br_{j}-\alpha_{1p})\}( \br_{j-\tau'}^\top \br_{j-\tau'}+\br_{j+\tau'}^\top \br_{j+\tau'} -2\alpha_{1p})
]^2 \big)\\
= &~\frac{Kp^4}{n^{2}}
\max_{n-\tau'+1\leq j\leq n \atop t<j-\tau}
\E \big([ \br_{j-\tau}^\top \br_{t-\tau}\br_t^\top \E\{\br_j (\br_{j}^\top \br_{j}-\alpha_{1p})\}( \br_{j-\tau'}^\top \br_{j-\tau'}+\br_{j+\tau'-n}^\top \br_{j+\tau'-n} -2\alpha_{1p})
]^2 \big)\,,
\end{align*}
where the last step is based on the fact $\br_t = \br_{t-n}$ for $t > n$.

For $\tau = \tau'$, by Assumption \ref{as:A4}, due to $\E(\br_t)=0$, it holds that
\begin{align}\label{eq:R12_11}
&\E \big([ \br_{j-\tau}^\top \br_{t-\tau}\br_t^\top \E\{\br_j (\br_{j}^\top \br_{j}-\alpha_{1p})\}( \br_{j-\tau}^\top \br_{j-\tau}+\br_{j+\tau-n}^\top \br_{j+\tau-n} -2\alpha_{1p})
]^2 \big) \notag\\
&~~~~~~~=\E\big[\{ \br_{j-\tau}^\top \br_{t-\tau}( \br_{j-\tau}^\top \br_{j-\tau}+\br_{j+\tau-n}^\top \br_{j+\tau-n} -2\alpha_{1p})
\}^2 \big] \notag\\
&~~~~~~~~~~~~\times \E\big[\E\{\br_j^\top (\br_{j}^\top \br_{j}-\alpha_{1p})\}\br_t\br_t^\top \E\{\br_j (\br_{j}^\top \br_{j}-\alpha_{1p})\}\big] \notag\\
&~~~~~~~=\E\big[\{ \br_{j-\tau}^\top \br_{t-\tau}( \br_{j-\tau}^\top \br_{j-\tau}+\br_{j+\tau-n}^\top \br_{j+\tau-n} -2\alpha_{1p})
\}^2 \big] \notag\\
&~~~~~~~~~~~~\times \frac{1}{p}\E\{\br_j^\top (\br_{j}^\top \br_{j}-\alpha_{1p})\}\bSigma \E\{\br_j (\br_{j}^\top \br_{j}-\alpha_{1p})\} \notag\\
&~~~~~~~\le \frac{K}{p}\E\big[\{ \br_{j-\tau}^\top \br_{t-\tau}( \br_{j-\tau}^\top \br_{j-\tau}+\br_{j+\tau-n}^\top \br_{j+\tau-n} -2\alpha_{1p})
\}^2 \big] \big|\E(\br_1 \br_{1}^\top \br_{1})\big|_2^2\,,	
\end{align}
which implies
\begin{align}\label{eq:R12_1}
R_{12}
\leq &~\frac{Kp^3}{n^{2}}
\max_{n-\tau+1\leq j\leq n \atop t<j-\tau }\E\big[\{ \br_{j-\tau}^\top \br_{t-\tau}( \br_{j-\tau}^\top \br_{j-\tau}+\br_{j+\tau-n}^\top \br_{j+\tau-n} -2\alpha_{1p})
\}^2 \big] \big|\E(\br_1 \br_{1}^\top \br_{1})\big|_2^2 \notag\\	
\leq &~\frac{Kp^3}{n^{2}}
\E\big[\{ \br_{j-\tau}^\top \br_{t-\tau}( \br_{j-\tau}^\top \br_{j-\tau} -\alpha_{1p})
\}^2 \big] \big|\E(\br_1 \br_{1}^\top \br_{1})\big|_2^2 \notag\\
&~~~~~+\frac{Kp^3}{n^{2}}\E\big[\{ \br_{j-\tau}^\top \br_{t-\tau}( \br_{j+\tau-n}^\top \br_{j+\tau-n} -\alpha_{1p})
\}^2 \big] \big|\E(\br_1 \br_{1}^\top \br_{1})\big|_2^2 \notag\\	
\leq &~\frac{Kp^2}{n^{2}}
\E\{\br_1^\top \br_1( \br_{1}^\top \br_{1} -\alpha_{1p})^2\} \big|\E(\br_1 \br_{1}^\top \br_{1})\big|_2^2 + \frac{Kp^2}{n^{2}}
\Var(\br_1^\top \br_1) \big|\E(\br_1 \br_{1}^\top \br_{1})\big|_2^2 \notag\\	
\leq &~\frac{Kp}{n^{2}}\Var^2(\br_1^\top \br_1)
\end{align}
for $\tau=\tau'$. In above result, the last inequality is based on \eqref{var111} and \eqref{var112}.

For $\tau > \tau'$, analogous to \eqref{eq:R12_11}, we have
\begin{align*}
&\E \big([ \br_{j-\tau}^\top \br_{t-\tau}\br_t^\top \E\{\br_j (\br_{j}^\top \br_{j}-\alpha_{1p})\}( \br_{j-\tau'}^\top \br_{j-\tau'}+\br_{j+\tau'-n}^\top \br_{j+\tau'-n} -2\alpha_{1p})
]^2 \big) \\
&~~~~~~~=\E\big[\{ \br_{j-\tau}^\top \br_{t-\tau}( \br_{j-\tau'}^\top \br_{j-\tau'}+\br_{j+\tau'-n}^\top \br_{j+\tau'-n} -2\alpha_{1p})
\}^2 \big] \\
&~~~~~~~~~~~~\times \E\big[\E\{\br_j^\top (\br_{j}^\top \br_{j}-\alpha_{1p})\}\br_t\br_t^\top \E\{\br_j (\br_{j}^\top \br_{j}-\alpha_{1p})\}\big] \notag\\
&~~~~~~~=\E\big[\{ \br_{j-\tau}^\top \br_{t-\tau}( \br_{j-\tau'}^\top \br_{j-\tau'}+\br_{j+\tau'-n}^\top \br_{j+\tau'-n} -2\alpha_{1p})
\}^2 \big] \\
&~~~~~~~~~~~~\times \frac{1}{p}\E\{\br_j^\top (\br_{j}^\top \br_{j}-\alpha_{1p})\}\bSigma \E\{\br_j (\br_{j}^\top \br_{j}-\alpha_{1p})\} \\
&~~~~~~~\le \frac{K}{p}\E\big[\{ \br_{j-\tau}^\top \br_{t-\tau}( \br_{j-\tau'}^\top \br_{j-\tau'}+\br_{j+\tau'-n}^\top \br_{j+\tau'-n} -2\alpha_{1p})
\}^2 \big] \big|\E(\br_1 \br_{1}^\top \br_{1})\big|_2^2\,,	
\end{align*}
which implies
\begin{align}\label{eq:R12_2}
R_{12}
\leq &~\frac{Kp^3}{n^{2}}
\max_{t<j-\tau\atop n-\tau'+1\leq j\leq n}\E\big[\{ \br_{j-\tau}^\top \br_{t-\tau}( \br_{j-\tau'}^\top \br_{j-\tau'}+\br_{j+\tau'-n}^\top \br_{j+\tau'-n} -2\alpha_{1p})
\}^2 \big] \big|\E(\br_1 \br_{1}^\top \br_{1})\big|_2^2\notag\\	
\leq &~\frac{Kp^3}{n^{2}}
\E\big[\{ \br_{j-\tau}^\top \br_{t-\tau}( \br_{j-\tau'}^\top \br_{j-\tau'} -\alpha_{1p})
\}^2 \big] \big|\E(\br_1 \br_{1}^\top \br_{1})\big|_2^2 \notag\\
&~~~~~+\frac{Kp^3}{n^{2}}\E\big[\{ \br_{j-\tau}^\top \br_{t-\tau}( \br_{j+\tau'-n}^\top \br_{j+\tau'-n} -\alpha_{1p})
\}^2 \big] \big|\E(\br_1 \br_{1}^\top \br_{1})\big|_2^2 \notag\\	
\leq &~\frac{Kp^2}{n^{2}}
\Var( \br_{1}^\top \br_{1} ) \big|\E(\br_1 \br_{1}^\top \br_{1})\big|_2^2 + \frac{Kp^2}{n^{2}}\E \{\br_1^\top \br_1( \br_1^\top \br_1 -\alpha_{1p})^2\} \big|\E(\br_1 \br_{1}^\top \br_{1})\big|_2^2 \notag\\
\leq&~ \frac{Kp}{n^{2}}\Var^2(\br_1^\top \br_1)
\end{align}
for $\tau > \tau'$.

For $\tau < \tau'$, we have
\begin{align*}
&\E \big([ \br_{j-\tau}^\top \br_{t-\tau}\br_t^\top \E\{\br_j (\br_{j}^\top \br_{j}-\alpha_{1p})\}( \br_{j-\tau'}^\top \br_{j-\tau'}+\br_{j+\tau'-n}^\top \br_{j+\tau'-n} -2\alpha_{1p})
]^2 \big) \\
&~~~~~~~\le 2\E \big([ \br_{j-\tau}^\top \br_{t-\tau}\br_t^\top \E\{\br_j (\br_{j}^\top \br_{j}-\alpha_{1p})\}( \br_{j-\tau'}^\top \br_{j-\tau'} -\alpha_{1p})
]^2 \big) \\
&~~~~~~~~~~~~ + 2\E \big([ \br_{j-\tau}^\top \br_{t-\tau}\br_t^\top \E\{\br_j (\br_{j}^\top \br_{j}-\alpha_{1p})\}( \br_{j+\tau'-n}^\top \br_{j+\tau'-n} -\alpha_{1p}) ]^2 \big)\,.
\end{align*}
By the Cauchy-Schwarz inequality and Assumption \ref{as:A4}, we have
\begin{align*}
    &\E \big([ \br_{j-\tau}^\top \br_{t-\tau}\br_t^\top \E\{\br_j (\br_{j}^\top \br_{j}-\alpha_{1p})\}( \br_{j-\tau'}^\top \br_{j-\tau'} -\alpha_{1p})]^2 \big) \\
    &~~~~~~~\le \frac{K}{p} \E \big([ \br_t^\top \E\{\br_j (\br_{j}^\top \br_{j}-\alpha_{1p})\}( \br_{j-\tau'}^\top \br_{j-\tau'} -\alpha_{1p})]^2 \br_{t-\tau}^\top \br_{t-\tau} \big) \\
    &~~~~~~~\le \frac{K}{p} \E \big([ \br_t^\top \E\{\br_j (\br_{j}^\top \br_{j}-\alpha_{1p})\}]^2 \E\{( \br_{j-\tau'}^\top \br_{j-\tau'} -\alpha_{1p})^2\br_{t-\tau}^\top \br_{t-\tau}\} \big) \\
    &~~~~~~~~~~~~ + \frac{K}{p} \E \big([ ( \br_{t}^\top \br_{t} -\alpha_{1p})\br_t^\top \E\{\br_j (\br_{j}^\top \br_{j}-\alpha_{1p})\}]^2      \big)\E(\br_{t-\tau}^\top \br_{t-\tau}) \\
    &~~~~~~~\le \frac{K}{p^2}\big|\E(\br_1\br_1^\top \br_1) \big|_2^2 \big[ \E \{\br_1^\top \br_1( \br_1^\top \br_1 -\alpha_{1p})^2\} +  \Var(\br_1^\top \br_1) \E(\br_{t-\tau}^\top \br_{t-\tau}) \big] \\
    &~~~~~~~~~~~~ + \frac{K}{p}  \E^{1/2}\{( \br_t^\top \br_t -\alpha_{1p})^4\} \E^{1/2}\{ [\br_t^\top \E\{\br_j (\br_{j}^\top \br_{j}-\alpha_{1p})\}]^4\} \\
    &~~~~~~~\le \frac{K}{p^2}\big|\E(\br_1\br_1^\top \br_1) \big|_2^2 \big[ \E \{\br_1^\top \br_1( \br_1^\top \br_1 -\alpha_{1p})^2\} +  \Var(\br_1^\top \br_1) \E(\br_{t-\tau}^\top \br_{t-\tau}) \big] \\
    &~~~~~~~~~~~~ +  \frac{K}{p} \E^{1/2}\{( \br_t^\top \br_t -\alpha_{1p})^4\} \E^{1/2}\{(\br_t^\top\br_t)^2\} \big|\E(\br_1\br_1^\top \br_1) \big|_2^2 \\
    &~~~~~~~\le \frac{K}{p^3} \Var^2(\br_1^\top \br_1) + \frac{K}{p^2} \Var^2(\br_1^\top \br_1) \le \frac{K}{p^2} \Var^2(\br_1^\top \br_1)
\end{align*}
for $\tau < \tau'$.
Analogously, for $\tau < \tau'$, we can also show 
\begin{align*}
    \E \big([ \br_{j-\tau}^\top \br_{t-\tau}\br_t^\top \E\{\br_j (\br_{j}^\top \br_{j}-\alpha_{1p})\}( \br_{j+\tau'-n}^\top \br_{j+\tau'-n} -\alpha_{1p})]^2 \big) \le \frac{K}{p^2}\Var^2(\br_1^\top \br_1)\,.
\end{align*}
Hence, for $\tau < \tau'$,
\begin{align*}
    R_{12} & \le \frac{Kp^4}{n^{2}} \max_{n-\tau'+1\leq j\leq n \atop t<j-\tau}
\E \big([ \br_{j-\tau}^\top \br_{t-\tau}\br_t^\top \E\{\br_j (\br_{j}^\top \br_{j}-\alpha_{1p})\}( \br_{j-\tau'}^\top \br_{j-\tau'}+\br_{j+\tau'-n}^\top \br_{j+\tau'-n} -2\alpha_{1p})
]^2 \big)  \\
& \leq  \frac{Kp^2}{n^{2}}\Var^2(\br_1^\top \br_1)\,.
\end{align*}
Together with \eqref{eq:R12_1} and \eqref{eq:R12_2}, we have \eqref{eq:R12_all} holds for any $\tau, \tau' \in [q]$. $\hfill\Box$

\subsubsection{Proof of  \eqref{eq:R2_all}}\label{subsec:R2}

For $R_{21}$, we have 
\begin{align*}
R_{21}
=&~\frac {16p^4}{n^6} 
\sum_{j=\tau+1}^n \E\bigg(\bigg[ \sum_{t_1, t_2< j-\tau}\br_{j-\tau}^\top \br_{t_1-\tau}\br_{t_1}^\top \E\{\br_j(\br_{j}^\top \br_{j}-\alpha_{1p})\}(\br_{t_2}^\top \br_{t_2}-\alpha_{1p})\bigg]^2\bigg)\\
=&~\frac {16p^4}{n^6} 
\sum_{j=\tau+1}^n  \sum_{t_1, t_2, t_3, t_4 <j-\tau}
\E
\big[\br_{j-\tau}^\top \br_{t_1-\tau}\br_{t_1}^\top \E\{\br_j(\br_{j}^\top \br_{j}-\alpha_{1p})\}(\br_{t_2}^\top \br_{t_2}-\alpha_{1p})\\
&~~~\qquad \qquad \qquad \qquad \qquad  \times \br_{j-\tau}^\top \br_{t_3-\tau}\br_{t_3}^\top \E\{\br_j(\br_{j}^\top \br_{j}-\alpha_{1p})\}(\br_{t_4}^\top \br_{t_4}-\alpha_{1p})\big]
\end{align*}
for any $\tau, \tau' \in [q]$.
By the Cauchy-Schwarz inequality, we have
\begin{align*}
R_{21}\leq &~\frac {Kp^4}{n^3} 
\max_{\tau+1 \le j \le n \atop t_1, t_2<j-\tau}
\E\big([\br_{j-\tau}^\top \br_{t_1-\tau}\br_{t_1}^\top \E\{\br_j(\br_{j}^\top \br_{j}-\alpha_{1p})\}(\br_{t_2}^\top \br_{t_2}-\alpha_{1p})]^2\big) \notag\\
\leq &~\frac {Kp^3}{n^3} 
\max_{\tau+1 \le j \le n \atop t_1, t_2<j-\tau}
\E\big[\br_{t_1-\tau}^\top \br_{t_1-\tau}\br_{t_1}^\top \br_{t_1}
\big|\E\{\br_j(\br_{j}^\top \br_{j}-\alpha_{1p})\}\big|_2^2(\br_{t_2}^\top \br_{t_2}-\alpha_{1p})^2\big] \notag\\
\leq &~\frac {Kp^3}{n^3}
\Var(\br_1^\top \br_1)\big|\E(\br_1\br_1^\top \br_1)\big|_2^2 \leq \frac {Kp^2}{n^3}\Var^2(\br_1^\top \br_1)
\end{align*}
for any $\tau, \tau' \in [q]$, where the last inequality is based on \eqref{var112}. 
For $R_{22}$, 
by the Jensen's inequality and the Cauchy-Schwarz inequality, we have
\begin{align*}
R_{22}
\leq &~\frac {Kp^4}{n^2} 
\max_{\tau+1 \le j \le n \atop t_1<j-\tau,~ j-\tau\leq t_2 <j}
\E\big([\br_{j-\tau}^\top \br_{t_1-\tau}\br_{t_1}^\top \E\{\br_j(\br_{j}^\top \br_{j}-\alpha_{1p})\}(\br_{t_2}^\top \br_{t_2}-\alpha_{1p})]^2\big) \\
\leq &~\frac {Kp^2}{n^2}
\Var(\br_1^\top \br_1)\big|\E (\br_1\br_1^\top \br_1)\big|_2^2\
\leq \frac {Kp}{n^2}
\Var^2(\br_1^\top \br_1)
\end{align*}
for any $\tau, \tau' \in [q]$.
Hence, we have \eqref{eq:R2_all} holds for any $\tau,\tau' \in [q]$. $\hfill\Box$

\section{Proof of Theorem \ref{null-Hn}}

\subsection{Proof of Theorem  \ref{null-Hn} with $p\to \infty$}\label{subsec:null_pinf}

When $\min\{p,n\}\to\infty$,   
as shown in Section \ref{subsec:Gq1}, under Assumptions \ref{as:A2} and \ref{as:A4}, we have 
\begin{align*}
    \left(
\frac{G_{n,1,1}}{\sigma_{n1}}, \ldots, \frac{G_{n,q,1}}{\sigma_{n1}}\right) \stackrel{\mathrm{d}}{\to} \mathcal{N}( 0,I_{q})\,.
\end{align*}
It then suffices to show $\hat\sigma_{n1}- \sigma_{n1} \stackrel{\mathrm{p}}{\to} 0$, where $\sigma_{n1}=\sqrt{2}p^{-1}\tr(\bSigma^2)$.
Since
\begin{align*}
\E (\hat\sigma_{n1})
 =&\frac{\sqrt{2}}{p}\cdot\frac1{n(n-1)}\sum_{t\neq s} \E\{(\x_t^\top \x_s)^2\}=\sigma_{n1}\,,
\end{align*}
in order to prove $\hat\sigma_{n1}- \sigma_{n1} \stackrel{\mathrm{p}}{\to} 0$, we only need to show $\Var (\hat\sigma_{n1}) = o(1)$.
Notice that
 \begin{align*}
\Var (\hat\sigma_{n1})
=&~\E\bigg(\bigg[\frac{\sqrt{2}}{p}\cdot\frac2{n(n-1)}\sum_{t< s} \{(\x_t^\top \x_s)^2-\tr(\bSigma^2)\}\bigg]^2\bigg)\\
=&~\frac8{p^2n^2(n-1)^2}\sum_{t_1< s_1,~t_2< s_2} \E\big[\{(\x_{t_1}^\top \x_{s_1})^2-\tr(\bSigma^2)\}\{(\x_{t_2}^\top \x_{s_2})^2-\tr(\bSigma^2)\}\big]\\
\leq &~\frac K{p^2n^2(n-1)^2}\bigg(\sum_{t< s} \E\big[\{(\x_{t}^\top \x_{s})^2-\tr(\bSigma^2)\}^2\big]\\
&~\qquad\qquad\qquad\quad+\sum_{t_1, t_2< s}\big| \E\big[\{(\x_{t_1}^\top \x_{s})^2-\tr(\bSigma^2)\}\{(\x_{t_2}^\top \x_{s})^2-\tr(\bSigma^2)\}\big]\big| \\
&~\qquad\qquad\qquad\quad+\sum_{t < s_1, s_2}\big| \E\big[\{(\x_{t}^\top \x_{s_1})^2-\tr(\bSigma^2)\}\{(\x_{t}^\top \x_{s_2})^2-\tr(\bSigma^2)\}\big]\big| \\
&~\qquad\qquad\qquad\quad+\sum_{t_1< s_1< s_2}\big| \E\big[\{(\x_{t_1}^\top \x_{s_1})^2-\tr(\bSigma^2)\}\{(\x_{s_1}^\top \x_{s_2})^2-\tr(\bSigma^2)\}\big]\big| \bigg)\\
\leq &~\frac K{p^2n^2(n-1)^2}\bigg[\sum_{t< s} \E \{(\x_{t}^\top \x_{s})^4\}+\sum_{t_1< s_1< s_2} {\E}^{1/2}\{(\x_{t_1}^\top \x_{s_1})^4\}{\E}^{1/2}\{(\x_{s_1}^\top \x_{s_2})^4\} \\
&~~ +\sum_{t_1, t_2< s} {\E}^{1/2}\{(\x_{t_1}^\top \x_{s})^4\}{\E}^{1/2}\{(\x_{t_2}^\top \x_{s})^4\} + \sum_{t < s_1, s_2} {\E}^{1/2}\{(\x_{t}^\top \x_{s_1})^4\}{\E}^{1/2}\{(\x_{t}^\top \x_{s_2})^4\} \bigg]\\
\leq &~ \frac{K}{p^2n} \max_{t\neq s}{\E}\{(\x_{t}^\top \x_{s})^4\} = \frac{K}{p^2n} \max_{t\neq s}{\E}\{(\z_{t}^\top \z_{s})^4\} \le Kn^{-1}\to0\,,
\end{align*} 
where the last inequality is based on the condition $\E\{(\z_t^\top \z_s)^4\}=O(p^2)$ for any $t \ne s$ as imposed in Assumption \ref{as:A2}. 
We complete the proof of Theorem \ref{null-Hn} with $p \to \infty$. $\hfill\Box$

\subsection{Proof of Theorem  \ref{null-Hn} with fixed $p$}\label{subsec:null_pfix}
Notice that 
\begin{align}\label{eq:decomp}
    G_{n, \tau,1} = \frac{n}{p}{\tr}(\bS_\tau\bS_\tau^\top) - \frac {1}{np}\sum_{t=1}^n\x_t^\top \x_t\x_{t+\tau}^\top \x_{t+\tau} := \frac{n}{p}{\tr}(\bS_\tau\bS_\tau^\top) - M_{n,\tau}\,.
\end{align}
For any $\tau\geq0$, notice that
\begin{align*}
S_\tau=\frac{1}{n}\sum_{t=1}^{n-\tau}X_tX_{t+\tau}^{\top} + \frac{1}{n}\sum_{t=n-\tau+1}^{n}X_tX_{t+\tau-n}^{\top} := S_\tau^{(1)} + S_\tau^{(2)}\,.
\end{align*}
Under the null hypothesis and Assumption \ref{as:A2}, by Theorem 1 of \cite{Chitturi1976}, the sample autocovariance matrices $\bS_1^{(1)},\ldots,\bS_q^{(1)}$ are asymptotically uncorrelated and satisfy 
\begin{align*}
    \sqrt{n} \begin{pmatrix}
         \mathrm{vec}\{\bS_1^{(1)}\} \\
         \vdots \\
         \mathrm{ vec}\{\bS_q^{(1)}\}\end{pmatrix} \xrightarrow{\mathrm{d}} \mathcal{N}(0, I_{q}\otimes (\bSigma\otimes \bSigma))
\end{align*}
as $n \to \infty$. For each fixed $\tau\in[q]$, $S_\tau^{(2)}$ contains only $\tau$ summands. Since $p$ is fixed, we have $\sqrt{n}S_\tau^{(2)} = O_{\mathrm{p}}(n^{-1/2})$.
Hence, it holds that
\begin{align}\label{eq:Dist1}
    \sqrt{n} \begin{pmatrix}
         \mathrm{vec}(\bS_1) \\
         \vdots \\
         \mathrm{ vec}(\bS_q)\end{pmatrix} \xrightarrow{\mathrm{d}} \mathcal{N}(0, I_{q}\otimes (\bSigma\otimes \bSigma))
\end{align}
as $n \to \infty$. Let $ N_\tau=(\varepsilon_{i,j,\tau})\in\mathbb R^{p\times p}$, where  $\varepsilon_{i,j,\tau}$ are i.i.d.\ $\mathcal{N}(0,1)$ variables for any $i,j\in [p]$ and $\tau\in[q]$. Hence, we can reformulate \eqref{eq:Dist1} as follows:
\begin{align}\label{eq:SDist}
     \sqrt{n} \begin{pmatrix}
         \mathrm{vec}(\bS_1) \\
         \vdots \\
         \mathrm{ vec}(\bS_q)\end{pmatrix}\xrightarrow{\mathrm{d}} \begin{pmatrix}
         \bSigma^{1/2}\otimes \bSigma^{1/2} &  &  \\
         & \ddots & \\
         &  &  \bSigma^{1/2}\otimes \bSigma^{1/2}
         \end{pmatrix} \begin{pmatrix}
         \mathrm{vec}( N_1) \\
         \vdots \\
         \mathrm{ vec}( N_q)\end{pmatrix}
\end{align}
as $n \to \infty$. Elementary calculations combined with classical results on U-statistics yield that $ M_{n,\tau} \xrightarrow{\mathrm{p}}  p^{-1}{\tr}^2(\bSigma)$ and $\hat\sigma_{n1}\xrightarrow{\mathrm{p}}  \sqrt{2}p^{-1}\tr(\bSigma^2)$.
Hence, we have
\begin{align}\label{eq:jointHn}
    \begin{pmatrix}
        H_{n, 1} \\
        \vdots \\
        H_{n, q}
    \end{pmatrix} &=  \begin{pmatrix}
        \hat\sigma_{n1}^{-1} &   & \\
        &  \ddots & \\
        &  & \hat\sigma_{n1}^{-1}
\end{pmatrix}\begin{pmatrix}
    np^{-1}\operatorname{tr}(\bS_1\bS_1^\top ) - M_{n,1} \\
    \vdots \\
    np^{-1}\operatorname{tr}(\bS_q\bS_q^\top )-M_{n,q}
    \end{pmatrix} \\
    &\xrightarrow{\mathrm{d}} \frac{1}{\sqrt{2}} \begin{pmatrix}
    \tr^{-1}(\bSigma^2)\{\operatorname{tr}(\bSigma N_1\bSigma  N_1^\top )- {\tr}^2(\bSigma)\} \notag\\
    \vdots \\
    \tr^{-1}(\bSigma^2)\{\operatorname{tr}(\bSigma  N_q\bSigma  N_q^\top )- {\tr}^2(\bSigma)\}
    \end{pmatrix} 
\end{align}
as $n \to \infty$.
Write $H_{\tau} = \tr^{-1}(\bSigma^2)\{\operatorname{tr}(\bSigma  N_{\tau}\bSigma  N_{\tau}^\top )- {\tr}^2(\bSigma)\}/\sqrt{2}$. Since $ N_1, \ldots,  N_q$ are independent random matrices, we know $H_1, \ldots, H_q$ are  independent random variables.

Let $\lambda_1\geq\cdots\geq\lambda_p$ be the eigenvalues of $\bSigma$. Applying the eigenvalue decomposition of $\bSigma$, we have $\bSigma =  Q \Lambda  Q^{\top}$, where $\Lambda = {\rm diag}(\lambda_1,\ldots, \lambda_p)$ and $ Q$ is an orthogonal matrix. Notice that $( Q N_{1} Q^{\top}, \ldots ,  Q N_{q} Q^{\top}) \stackrel{\mathrm{d}}{=} (N_{1},\ldots, N_{q})$.
Hence, due to ${\tr}(\bSigma) = \sum_{i=1}^p \lambda_i$ and ${\tr}(\bSigma^2) = \sum_{i=1}^p \lambda_i^2$, we have
\begin{align*}
    \frac{\operatorname{tr}(\bSigma  N_\tau\bSigma  N_\tau^\top) - {\tr}^2(\bSigma)}{\sqrt{2}{\tr}(\bSigma^2)} &=  \frac{\operatorname{tr}(\Lambda Q^{\top}  N_\tau  Q \Lambda Q^{\top}  N_\tau^\top  Q) - \sum_{i=1}^{p} \sum_{j=1}^{p} \lambda_i \lambda_j}{\sqrt{2}\sum_{i=1}^{p} \lambda_i^2} \\
    &\stackrel{\mathrm{d}}{=} \frac{\operatorname{tr}(\Lambda  N_\tau \Lambda  N_\tau^\top ) - \sum_{i=1}^{p} \sum_{j=1}^{p} \lambda_i \lambda_j}{\sqrt{2}\sum_{i=1}^{p} \lambda_i^2} \\
    & = \frac{\sum_{i=1}^{p} \sum_{j=1}^{p} \lambda_i \lambda_j (\varepsilon_{i,j,\tau}^2-1)}{\sqrt{2}\sum_{i=1}^p\lambda_i^2}
\end{align*}
for any $\tau \in [q]$.
We complete the proof of Theorem \ref{null-Hn} with fixed $p$. $\hfill\Box$

\section{Proof of Theorem  \ref{tm:boot_size}}

Let $G_{n} = \sum_{\tau=1}^q G_{n,\tau,1}$. Recall 
\begin{align}\label{eq:Tn1}
T_n =  \sum_{\tau=1}^q H_{n,\tau} = \sum_{\tau=1}^q \frac{G_{n,\tau,1}}{\hat\sigma_{n1}} = \frac{G_{n}}{\hat\sigma_{n1}}\,.
\end{align}
Write 
\begin{align*}
    G_{n}^e = \sum_{\tau=1}^q G_{n,\tau,1}^e ~~ \mbox{with} ~~ G_{n,\tau,1}^e= \frac 1{np}\sum_{t\neq s}\y_t^\top \y_s\y_{t+\tau}^\top \y_{s+\tau}\,.
\end{align*}
Recall $T_{n}^{e}=\sum_{\tau=1}^q H_{n,\tau}^{e}$, where $H_{n,\tau}^{e}$ is calculated in the same manner as $H_{n,\tau}$ but with replacing $\{\x_t\}_{t=1}^n$ by $\{\y_t\}_{t=1}^n$.
Then
\begin{align*}
    T_{n}^{e} = \sum_{\tau=1}^q \frac{G_{n,\tau,1}^e}{\hat\sigma_{n1}^e} ~~\mbox{with}~~ \hat\sigma_{n1}^e = \frac{\sqrt{2}}{n(n-1)p}\sum_{t, s\in[n]:\,t\neq s} (\y_t^\top \y_s)^2\,.
\end{align*}
Since $\y_t = e_t \x_t$ and $e_t^2 = 1$, we have ${\hat\sigma_{n1}^e} = {\hat\sigma_{n1}}$,
which implies 
\begin{align}\label{eq:Tne}
    T_{n}^{e} = \sum_{\tau=1}^q \frac{G_{n,\tau,1}^e}{\hat\sigma_{n1}} = \frac{G_{n}^e}{\hat\sigma_{n1}}\,.
\end{align}

\subsection{Proof of Theorem  \ref{tm:boot_size} with $p\to \infty$}\label{subsec:size_pinf}

Write $\mathcal{X}_n = \{\x_1,\ldots,\x_n\}$. Recall $\tilde{{\rm cv}}_{\alpha} = \inf\{t \in \mathbb{R}: \pr(T_{n}^{e} \le t\,|\, \mathcal{X}_n) \ge 1-\alpha\}$. Define $\tilde{{\rm cv}}_{G,\alpha} = \inf\{t \in \mathbb{R}: \pr(G_{n}^{e} \le t\,|\, \mathcal{X}_n) \ge 1-\alpha\}$. Notice that $\hat\sigma_{n1}$ is a measurable function of $\mathcal{X}_n$. Then
\begin{align*}
\hat\sigma_{n1}\tilde{{\rm cv}}_{\alpha} & = \hat\sigma_{n1}\inf\{t \in \mathbb{R}: \pr(T_{n}^{e} \le t\,|\, \mathcal{X}_n) \ge 1-\alpha\} \\
    & = \hat\sigma_{n1}\inf\bigg\{\frac{t}{\hat\sigma_{n1}} \in \mathbb{R}: \pr\bigg(\frac{G_{n}^{e}}{\hat\sigma_{n1}} \le \frac{t}{\hat\sigma_{n1}}\,\bigg|\, \mathcal{X}_n\bigg) \ge 1-\alpha\bigg\} \\
    & = \inf\{t \in \mathbb{R}: \pr(G_{n}^{e} \le t\,|\, \mathcal{X}_n) \ge 1-\alpha\} = \tilde{{\rm cv}}_{G,\alpha}\,,
\end{align*}
which implies $\pr_{H_0}(T_n > \tilde{{\rm cv}}_{\alpha}) =\pr_{H_0}(G_{n}> \tilde{{\rm cv}}_{G,\alpha})$. Hence, to show $\pr_{H_0}(T_n > \tilde{{\rm cv}}_{\alpha}) \to \alpha$, it is equivalent to prove $\pr_{H_0}(G_{n}> \tilde{{\rm cv}}_{G,\alpha}) \to \alpha$ as $n \to \infty$.

Let $\E_0^e(\cdot) = \E^e(\cdot) = \E(\cdot \,|\, \mathcal{X}_n)$ denote the conditional expectation given $\mathcal{X}_n$, i.e. the expectation is taken with respect to $e_1, \ldots, e_n$. 
For $j \in [n]$, let $\E_j^e(\cdot) = \E(\cdot \,|\, \mathcal{X}_n, \mathcal F_j^e)$ denote the conditional expectation given both $\mathcal{X}_n$ and the $\sigma$-field $\mathcal F_j^e$ generated by $\{e_1,\ldots,e_j\}$.
Recall $\br_t=p^{-1/2}\x_t$ for $t\in [n]$. By convention, we let $e_{t} = e_{t-n}$ if $t > n$. Then
\begin{align*}
G_{n,\tau,1}^e = \frac pn\sum_{t\neq s}e_te_se_{t+\tau}e_{s+\tau}\br_t^\top \br_s\br_{t+\tau}^\top \br_{s+\tau} = \frac {2p}{n}\sum_{1 \le t < s \le n}e_te_se_{t+\tau}e_{s+\tau}\br_t^\top \br_s\br_{t+\tau}^\top \br_{s+\tau}\,.
\end{align*}
We decompose $G_{n,\tau,1}^e$ into four terms: 
\begin{align*}
    G_{n,\tau,1}^e = G_{n,\tau,1,1}^e + G_{n,\tau,1,2}^e + G_{n,\tau,1,3}^e + G_{n,\tau,1,4}^e,
\end{align*}
where 
\begin{align*}
    G_{n,\tau,1,1}^e & = \frac {2p}{n}\sum_{1 \le t < s \le n-\tau}e_te_se_{t+\tau}e_{s+\tau}\br_t^\top \br_s\br_{t+\tau}^\top \br_{s+\tau}\,,\\
    G_{n,\tau,1,2}^e & = \frac {2p}{n}\sum_{t\in [n-\tau] \atop { s\in \{n-\tau+1,\ldots,n\} \atop t < s \neq t+\tau}}e_te_se_{t+\tau}e_{s+\tau}\br_t^\top \br_s\br_{t+\tau}^\top \br_{s+\tau}\,,\\
    G_{n,\tau,1,3}^e & = \frac {2p}{n}\sum_{t=n-2\tau+1}^{n-\tau}e_te_{t+2\tau}\br_t^\top \br_{t+\tau}\br_{t+\tau}^\top \br_{t+2\tau}\,,\\
    G_{n,\tau,1,4}^e & = \frac {2p}{n}\sum_{n-\tau+1 \le t < s \le n}e_te_se_{t+\tau}e_{s+\tau}\br_t^\top \br_s\br_{t+\tau}^\top \br_{s+\tau}\,.
\end{align*}
We have
\begin{align*}
    G_{n,\tau,1,1}^e = &~ \frac {2p}{n} \sum_{j=1}^n \bigg\{\E_j^e \bigg(\sum_{1 \le t < s\le n-\tau}e_te_se_{t+\tau}e_{s+\tau}\br_t^\top \br_s\br_{t+\tau}^\top \br_{s+\tau}\bigg) \\
    &~~~~~~~~~~~~~~~~~~~- \E_{j-1}^e \bigg(\sum_{1 \le t < s\le n-\tau}e_te_se_{t+\tau}e_{s+\tau}\br_t^\top \br_s\br_{t+\tau}^\top \br_{s+\tau}\bigg) \bigg\}  \\
    = &~ \frac {2p}{n} \sum_{j=\tau+1}^n \sum_{t < j-\tau}e_{t-\tau}e_{j-\tau}e_{t}e_{j}\br_{t-\tau}^\top \br_{j-\tau}\br_{t}^\top \br_{j} \\
    &~~~~~~+ \frac {2p}{n} \sum_{j=\tau+1}^n \sum_{j-\tau \le t < j}e_{t-\tau}e_{j-\tau}e_{t}e_{j}\br_{t-\tau}^\top \br_{j-\tau}\br_{t}^\top \br_{j} \\
    = &~ J_{1,\tau} + J_{2,\tau}\,.
\end{align*}
Notice that $\E^e(J_{2,\tau}) = 0$ and
\begin{align}\label{eq:Q51}
	 \E^e (J_{2,\tau}^2) 
  &=\frac {4p^2}{n^2}\E^e\bigg\{ \bigg(\sum_{j=\tau+1}^n \sum_{j-\tau \le t < j}e_{t-\tau}e_{j-\tau}e_{t}e_{j}\br_{t-\tau}^\top \br_{j-\tau}\br_{t}^\top \br_{j} \bigg)^2 \bigg\}\notag\\
  &=\frac {4p^2}{n^2}\sum_{j=\tau+1}^n \sum_{j-\tau \le t < j}\E^e\{(e_{t-\tau}e_{j-\tau}e_{t}e_{j}\br_{t-\tau}^\top \br_{j-\tau}\br_{t}^\top \br_{j})^2 \} \notag\\
  & = \frac{4p^2}{n^2}\sum_{j=\tau+1}^n \sum_{j-\tau \le t < j}(\br_{t-\tau}^\top \br_{j-\tau}\br_{t}^\top \br_{j})^2 := Q_{n5,1}(\tau) \,.
\end{align}
As we will show in Section \ref{subsubsec:q345}, 
\begin{align}\label{q5}
    \E\{|Q_{n5,1}(\tau)|\} \le Kn^{-1}
\end{align}
for any $\tau \in [q]$. 
It is easy to verify that $\E^e(G_{n,\tau,1,2}^e)=0$, $\E^e(G_{n,\tau,1,3}^e)=0$ and $\E^e(G_{n,\tau,1,4}^e)=0$. Moreover, we have
\begin{align}
    \E^e\{(G_{n,\tau,1,2}^e)^2\} 
    & = \frac {4p^2}{n^2}\sum_{t\in [n-\tau] \atop { s\in \{n-\tau+1,\ldots,n\} \atop t < s \neq t+\tau}}\E^e\{(e_te_se_{t+\tau}e_{s+\tau}\br_t^\top \br_s\br_{t+\tau}^\top \br_{s+\tau})^2\}\notag\\
    & = \frac {4p^2}{n^2}\sum_{t\in [n-\tau] \atop { s\in \{n-\tau+1,\ldots,n\} \atop t < s \neq t+\tau}}(\br_t^\top \br_s\br_{t+\tau}^\top \br_{s+\tau})^2:=Q_{n5,2}(\tau)\,,\label{eq:Q52}\\
    \E^e\{(G_{n,\tau,1,3}^e)^2\} 
    & = \frac {4p^2}{n^2}\sum_{t=n-2\tau+1}^{n-\tau}\E^e\{(e_te_{t+2\tau}\br_t^\top \br_{t+\tau}\br_{t+\tau}^\top \br_{t+2\tau})^2\}\notag\\
    & = \frac {4p^2}{n^2}\sum_{t=n-2\tau+1}^{n-\tau}(\br_t^\top \br_{t+\tau}\br_{t+\tau}^\top \br_{t+2\tau})^2:=Q_{n5,3}(\tau)\,,\label{eq:Q53}\\
     \E^e\{(G_{n,\tau,1,4}^e)^2\} 
    & = \frac {4p^2}{n^2}\sum_{n-\tau+1 \le t < s \le n}\E^e\{(e_te_se_{t+\tau}e_{s+\tau}\br_t^\top \br_s\br_{t+\tau}^\top \br_{s+\tau})^2\}\notag\\
    & = \frac {4p^2}{n^2}\sum_{n-\tau+1 \le t < s \le n}(\br_t^\top \br_s\br_{t+\tau}^\top \br_{s+\tau})^2:=Q_{n5,4}(\tau)\,.\label{eq:Q54}
\end{align}
Using similar arguments for deriving \eqref{eq:Gnt12}--\eqref{eq:Gnt14}, we have
\begin{align}\label{eq:q5tt}
    \E\{|Q_{n5,2}(\tau)|\} \le Kn^{-1}\,, ~ \E\{|Q_{n5,3}(\tau)|\} \le Kn^{-1} ~~\mbox{and}~~ \E\{|Q_{n5,4}(\tau)|\} \le Kn^{-2}\,.
\end{align}
Moreover, write 
\begin{align*}
    J_{1,\tau} = \sum_{j=\tau+1}^n D_{j,\tau}^e ~~ \mbox{with} ~~
    D_{j,\tau}^e = \frac {2p}{n}\sum_{t<j-\tau} e_{j-\tau}e_{t-\tau}e_te_j\br_{j-\tau}^\top \br_{t-\tau}\br_t^\top \br_j\,.
\end{align*}
Given $\mathcal{X}_n$, it can be shown that $\{D_{j,\tau}^e\}$ forms a sequence of martingale differences with respect to $\{\mathcal F_j^e\}$. By convention, we set $\br_\ell=0$ for $\ell\leq 0$. Given $\mathcal{X}_n$, write
\begin{align*}
    \tilde{G}_n^e =  \sum_{j=1}^n \sum_{\tau=1}^q D_{j,\tau}^e\,,
\end{align*}
where $\{\sum_{\tau=1}^q D_{j,\tau}^e\}$ forms a sequence of martingale differences with respect to $\{\mathcal F_j^e\}$. 
In the sequel, we will show
\begin{align}\label{eq:distConverge2}
    \sup_{t \in \mathbb{R}}\bigg|\pr(\tilde{G}_{n}^{e} \le t\,|\, \mathcal{X}_n)-\Phi\bigg(\frac{t}{\sqrt{2q\alpha_2^2}}\bigg)\bigg| \stackrel{\mathrm{p}}{\to} 0
\end{align}
as $n \to \infty$. 

Write
\begin{align*}
    \frac{\tilde{G}_n^e}{\sqrt{2q\alpha_2^2}} = \sum_{j=1}^n Z_j^e ~~\mbox{with}~~
    Z_j^e = \sum_{\tau=1}^q \frac{D_{j,\tau}^e}{\sqrt{2q\alpha_2^2}}\,.
\end{align*}
By Theorem 3.9 of \cite{Hall1980} with $\delta=1$, for all $x \in \mathbb{R}$, there exists a universal constant $C>0$ that does not depend on $x$, such that whenever $L_n \le 1$,
\begin{align}\label{eq:BEbound}
    \bigg|\pr\bigg(\frac{\tilde{G}_{n}^{e}}{\sqrt{2q\alpha_2^2}} \le x\,\bigg|\, \mathcal{X}_n\bigg)-\Phi(x)\bigg| \le CL_n^{1/5}(1+|x|^{16/5})^{-1}\,,
\end{align}
 where
\begin{align}\label{eq:Ln}
    L_n = \E^e\{(V_n^2-1)^2\} + \sum_{j=1}^n \E^e\{(Z_j^e)^4\} ~~\mbox{with}~~ V_n^2 = \sum_{j=1}^n \E^e_{j-1}\{(Z_j^e)^2\}\,.
\end{align}
To show \eqref{eq:distConverge2}, it suffices to prove $L_n \stackrel{\mathrm{p}}{\to} 0$. 
Notice that
\begin{align}\label{eq:Vn2}
    V_n^2 & = \sum_{j=1}^n \E^e_{j-1}\bigg\{\bigg(\sum_{\tau=1}^q \frac{D_{j,\tau}^e}{\sqrt{2q\alpha_2^2}}\bigg)^2\bigg\} \notag\\
    & = \frac{1}{2q\alpha_2^2}\sum_{j=1}^n \sum_{\tau=1}^q \E^e_{j-1}\{(D_{j,\tau}^e)^2\} + \frac{2}{2q\alpha_2^2}\sum_{j=1}^n \sum_{1 \le \tau' < \tau\le q} \E^e_{j-1}(D_{j,\tau}^eD_{j,\tau'}^e)\,.
\end{align}
By the fact $e_\ell^2=1$ for $\ell\in [n]$, we have
\begin{align*}
 \sum_{j=1}^n\sum_{\tau=1}^q\E_{j-1}^e\{(D_{j,\tau}^e)^2\}
=&~ \frac {4p^2}{n^2}\sum_{j=1}^n\sum_{\tau=1}^q\sum_{t, s<j-\tau}e_{t-\tau}e_t e_{s-\tau}e_s
\br_{s-\tau}^\top \br_{j-\tau}\br_{j-\tau}^\top \br_{t-\tau}\br_s^\top \br_j \br_j^\top \br_t\\
=&~ \sum_{\tau=1}^q\frac {4p^2}{n^2}\sum_{j=1}^n\sum_{t<j-\tau }
(\br_{t-\tau}^\top \br_{j-\tau})^2 (\br_t^\top \br_j)^2\\
&~~~+ \sum_{\tau=1}^q \frac {8p^2}{n^2}\sum_{j=1}^n\sum_{t<s<j-\tau }e_{t-\tau}e_t e_{s-\tau}e_s
\br_{s-\tau}^\top \br_{j-\tau}\br_{j-\tau}^\top \br_{t-\tau}\br_s^\top \br_j \br_j^\top \br_t \\
:= &~ \sum_{\tau=1}^q Q_{n1}(\tau) + \sum_{\tau=1}^q \tilde{Q}_{n2}^e(\tau)\,,
\end{align*}
and 
\begin{align*}
&\sum_{j=1}^n\sum_{1 \le \tau' < \tau\le q}\E_{j-1}^e (D_{j,\tau}^eD_{j,\tau'}^e) \\
&~~~~~~~~= \sum_{1 \le \tau' < \tau\le q} \frac {4p^2}{n^2}\sum_{j=1}^n\sum_{t<j-\tau \atop s<j-\tau'}
e_{j-\tau}e_{t-\tau}  e_t\br_{j-\tau}^\top  \br_{t-\tau} \br_t^\top \br_j
e_{j-\tau'}e_{s-\tau'}e_s\br_{j-\tau'}^\top \br_{s-\tau'}\br_s^\top \br_j \\
&~~~~~~~~:= \sum_{1 \le \tau' < \tau\le q}\tilde{Q}_{n3}^e(\tau,\tau')\,.
\end{align*}
Hence, by \eqref{eq:Vn2}, we have
\begin{align*}
    V_n^2 & = \frac{1}{2q\alpha_2^2}\sum_{\tau=1}^q Q_{n1}(\tau) + \frac{1}{2q\alpha_2^2}\sum_{\tau=1}^q \tilde{Q}_{n2}^e(\tau) + \frac{1}{q\alpha_2^2}\sum_{1 \le \tau' < \tau\le q}\tilde{Q}_{n3}^e(\tau,\tau')\,.
\end{align*}
Notice that $\E^e\{\tilde{Q}_{n2}^e(\tau)\}=0$ and $\E^e\{\tilde{Q}_{n3}^e(\tau,\tau')\}=0$ for any $\tau, \tau'$ satisfying $\tau > \tau'$,
which implies that 
\begin{align*}
    \E^e(V_n^2) =  \frac{1}{2q\alpha_2^2}\sum_{\tau=1}^q Q_{n1}(\tau)\,.
\end{align*}
We then have
\begin{align}\label{eq:EVn2}
    &\E^e\{(V_n^2-1)^2\} = \E^e[\{V_n^2 -\E^e(V_n^2)\}^2] + \{\E^e(V_n^2) -1\}^2 \notag\\
    &~~~~~~ = \E^e\bigg[\bigg\{\frac{1}{2q\alpha_2^2}\sum_{\tau=1}^q \tilde{Q}_{n2}^e(\tau) + \frac{1}{q\alpha_2^2}\sum_{1 \le \tau' < \tau\le q}\tilde{Q}_{n3}^e(\tau,\tau')\bigg\}^2\bigg] + \bigg\{\frac{1}{2q\alpha_2^2}\sum_{\tau=1}^q Q_{n1}(\tau)-1\bigg\}^2\notag\\
    &~~~~~~ \le \frac{1}{2q^2\alpha_2^4}\E^e\bigg[\bigg\{\sum_{\tau=1}^q \tilde{Q}_{n2}^e(\tau)\bigg\}^2\bigg] + \frac{2}{q^2\alpha_2^4}\E^e\bigg[\bigg\{\sum_{1 \le \tau' < \tau\le q}\tilde{Q}_{n3}^e(\tau,\tau')\bigg\}^2\bigg] \notag\\
    &~~~~~~~~~~~~+ \bigg\{\frac{1}{2q\alpha_2^2}\sum_{\tau=1}^q Q_{n1}(\tau)-1\bigg\}^2 \notag \\
    &~~~~~~ \le \frac{1}{2\alpha_2^4}\max_{\tau \in [q]}\E^e[\{\tilde{Q}_{n2}^e(\tau)\}^2] + \frac{q^2}{2\alpha_2^4}\max_{1 \le \tau' < \tau\le q}\E^e[\{\tilde{Q}_{n3}^e(\tau,\tau')\}^2] \notag\\
    &~~~~~~~~~~~~+ \bigg\{\frac{1}{2q\alpha_2^2}\sum_{\tau=1}^q Q_{n1}(\tau)-1\bigg\}^2\,.
\end{align}
On the other hand, by the Cauchy-Schwarz inequality, we have
\begin{align*}
    \sum_{j=1}^n \E^e\{(Z_j^e)^4\} & = \frac{1}{4q^2\alpha_2^4}\sum_{j=1}^n \E^e\bigg\{\bigg(\sum_{\tau=1}^q D_{j,\tau}^e\bigg)^4\bigg\} \le \frac{1}{4\alpha_2^4} \sum_{j=1}^n \E^e\bigg[\bigg\{\sum_{\tau=1}^q (D_{j,\tau}^e)^2\bigg\}^2\bigg] \\
    &\le \frac{q}{4\alpha_2^4} \sum_{j=1}^n \sum_{\tau=1}^q \E^e\{(D_{j,\tau}^e)^4\} \le \frac{q^2}{4\alpha_2^4} \max_{\tau \in [q]} \tilde{Q}_{n4}(\tau)\,, 
\end{align*}
where $\tilde{Q}_{n4}(\tau) = \sum_{j=1}^n\E^e\{(D_{j,\tau}^e)^4\}$.
Hence, by \eqref{eq:Ln}, we have 
\begin{align}\label{eq:Ln2}
    L_n & \le \bigg\{\frac{1}{2q\alpha_2^2}\sum_{\tau=1}^q Q_{n1}(\tau)-1\bigg\}^2 + \frac{1}{2\alpha_2^4}\max_{\tau \in [q]}\E^e[\{\tilde{Q}_{n2}^e(\tau)\}^2] \notag \\
    &~~~~~~~~~+ \frac{q^2}{2\alpha_2^4}\max_{1 \le\tau'<\tau\le  q}\E^e[\{\tilde{Q}_{n3}^e(\tau,\tau')\}^2] + \frac{q^2}{4\alpha_2^4} \max_{\tau \in [q]} \tilde{Q}_{n4}(\tau)\,.
\end{align}
By the fact $e_\ell^2=1$ for $\ell\in [n]$, notice that
\begin{align*}
\E^e[\{\tilde{Q}_{n2}^e(\tau)\}^2] = &~ \E^e\bigg[\bigg\{\sum_{j=1}^n\frac {8p^2}{n^2}\sum_{t<s<j-\tau }e_{t-\tau}e_t e_{s-\tau}e_s
(\br_{s-\tau}^\top \br_{j-\tau}\br_{j-\tau}^\top \br_{t-\tau})  (\br_s^\top \br_j \br_j^\top \br_t)\bigg\}^2\bigg]\\
 \leq &~ 
\frac {Kp^4}{n^4}
\sum_{j_1, j_2=1}^n\sum_{t<s<\min\{j_1-\tau, j_2-\tau\}}\prod_{i=1,2}
(\br_{s-\tau}^\top \br_{j_i-\tau}\br_{j_i-\tau}^\top \br_{t-\tau}\br_s^\top \br_{j_i} \br_{j_i}^\top \br_t)\\
:=&~ KQ_{n2}(\tau)\,,\\
\E^e[\{\tilde{Q}_{n3}^e(\tau,\tau')\}^2] =&~ \frac {16p^4}{n^4} 
\sum_{j=\tau+1}^n\E^e\bigg\{\bigg(\sum_{t<j-\tau \atop s<j-\tau'}
e_{t-\tau}  e_t\br_{j-\tau}^\top  \br_{t-\tau} \br_t^\top \br_j
e_{s-\tau'}e_s\br_{j-\tau'}^\top \br_{s-\tau'}\br_s^\top \br_j\bigg)^2\bigg\}\\
=&~ \frac {16p^4}{n^4} \sum_{j=\tau+1}^n\sum_{t<j-\tau \atop s<j-\tau'}
(\br_{j-\tau}^\top  \br_{t-\tau} \br_t^\top \br_j \br_{j-\tau'}^\top \br_{s-\tau'}\br_s^\top \br_j)^2\\
&~~~~~~+\frac {64p^4}{n^4} \sum_{j=\tau+1}^n\sum_{t<j-\tau-\tau'}
\big\{\br_{j-\tau}^\top  \br_{t-\tau} \br_t^\top \br_j \br_{j-\tau'}^\top \br_{t}(\br_{t+\tau'}^\top \br_j)^2\\
&~~~~~~~~~~~~~~~~~~~~~~~~~~~~~~~~~~~~~\times\br_{j-\tau}^\top  \br_{t-\tau+\tau'}\br_{j-\tau'}^\top \br_{t-\tau}\br_{t-\tau+\tau'}^\top\br_{j}\big\}\\
:= &~ KQ_{n3}(\tau,\tau')\,,
\end{align*}
and 
\begin{align*}
\tilde{Q}_{n4}(\tau) =&~
\frac {16p^4}{n^4}\sum_{j=\tau+1}^n\sum_{t_1, t_2, t_3, t_4<j-\tau}
 \E^e \bigg\{\bigg(\prod_{i=1}^4 e_{t_i-\tau}e_{t_i}\bigg)\prod_{i=1}^4\br_{j-\tau}^\top \br_{t_i-\tau}\br_{t_i}^\top \br_j \bigg\}\\
\leq &~
\frac {Kp^4}{n^4}\sum_{j=\tau+1}^n\sum_{t_1\leq t_2\leq t_3\leq t_4<j-\tau}
 \Bigg|\E^e \bigg\{\bigg(\prod_{i=1}^4 e_{t_i-\tau}e_{t_i}\bigg)\prod_{i=1}^4\br_{j-\tau}^\top \br_{t_i-\tau}\br_{t_i}^\top \br_j \bigg\}\Bigg|\\ 
\leq &~
\frac {Kp^4}{n^4}\sum_{j=\tau+1}^n\sum_{t_1\leq t_3<j-\tau}
\prod_{i=1, 3}(\br_{j-\tau}^\top \br_{t_i-\tau}\br_{t_i}^\top \br_j)^2\\
:= &~KQ_{n4}(\tau)
\end{align*}
for any $\tau,\tau' \in [q]$ satisfying $\tau>\tau'$. Under the null hypothesis, as we will show in Sections \ref{subsubsec:q1}--\ref{subsubsec:q345}, 
\begin{align}
& \E\{Q_{n1}(\tau)\}
\to 2\alpha_2^2~~\mbox{and}~~
\Var\{Q_{n1}(\tau)\}
\leq Kn^{-1} \,,\label{q1}\\
&~~~~~~~~~~~\E\{|Q_{n2}(\tau)|\} \le K(p^{-2}+n^{-1/2}) \label{q2}\\
&~~~~~~~\,\E\{|Q_{n3}(\tau,\tau')|\}
\leq Kn^{-1}\,, \label{q3}\\
&~~~~~~~~~~~\E\{|Q_{n4}(\tau)|\}
\leq  Kn^{-1}\,, \label{q4}
\end{align}
for any $\tau,\tau' \in [q]$ satisfying $\tau>\tau'$. Together with \eqref{q5} and \eqref{eq:q5tt}, by Markov's inequality, there exists a deterministic sequence $\varepsilon_n \to 0$ such that the event
\begin{align*}
A_{\varepsilon_n} & = \bigg\{\{\br_t\}_{t=1}^n: \max_{\tau \in [q]}|Q_{n1}(\tau)- 2\alpha_2^2| +
\max_{\tau \in [q]}\{ |Q_{n2}(\tau)|+ |Q_{n4}(\tau)|\} \\
&~~~~~~~~~~~~~~~~~~~~~~~+ \max_{1 \le\tau'<\tau \leq q}|Q_{n3}(\tau, \tau')| + \sum_{\ell=1}^4\max_{\tau \in [q]}|Q_{n5,\ell}(\tau)| \le \varepsilon_n\bigg\}
\end{align*}
occurs with probability approaching one, i.e. $\pr(A_{\varepsilon_n}) \to 1$. Restricted on the event $A_{\varepsilon_n}$, by \eqref{eq:Ln2}, we have $L_n \le K \varepsilon_n$. Hence, \eqref{eq:BEbound} implies that, restricted on the event $A_{\varepsilon_n}$, 
\begin{align}\label{eq:BE2}
\sup_{t \in \mathbb{R}}\bigg|\pr(\tilde{G}_{n}^{e} \le t\,|\, \mathcal{X}_n)-\Phi\bigg(\frac{t}{\sqrt{2q\alpha_2^2}}\bigg)\bigg| =\sup_{x\in \mathbb{R}}\bigg|\pr\bigg(\frac{\tilde{G}_{n}^{e}}{\sqrt{2q\alpha_2^2}} \le x\,\bigg|\, \mathcal{X}_n\bigg)-\Phi(x)\bigg| \le \eta_n 
\end{align}
for any sufficiently large $n$, where $\eta_n = K\varepsilon_n^{1/5}$. 
Moreover, recall $G_n^e-\tilde{G}_n^e
=\sum_{\tau=1}^q
(G_{n,\tau,1,2}^e+G_{n,\tau,1,3}^e+G_{n,\tau,1,4}^e+J_{2,\tau})$
and $q$ is a fixed constant. Restricted on the event $A_{\varepsilon_n}$, by the Markov's inequality, \eqref{eq:Q51}, \eqref{eq:Q52}, \eqref{eq:Q53} and \eqref{eq:Q54}, we have
\begin{align*}
    \pr\big(|G_n^e-\tilde{G}_n^e|> \varepsilon_n^{1/4}\,|\,\mathcal{X}_n\big) &\le
    \sum_{\tau=1}^q\pr\bigg(
    |J_{2,\tau}|>\frac{\varepsilon_n^{1/4}}{4q}
    \,\bigg|\,\mathcal{X}_n\bigg)
    +\sum_{\tau=1}^q\pr\bigg(
    |G_{n,\tau,1,2}^e|>\frac{\varepsilon_n^{1/4}}{4q}
    \,\bigg|\,\mathcal{X}_n\bigg) \\
    &~~~~+
    \sum_{\tau=1}^q\pr\bigg(
    |G_{n,\tau,1,3}^e|>\frac{\varepsilon_n^{1/4}}{4q}
    \,\bigg|\,\mathcal{X}_n\bigg)
    +\sum_{\tau=1}^q\pr\bigg(
    |G_{n,\tau,1,4}^e|>\frac{\varepsilon_n^{1/4}}{4q}
    \,\bigg|\,\mathcal{X}_n\bigg) \\
    &\le
    \sum_{\tau=1}^q\frac{16q^2}{\varepsilon_n^{1/2}}
    \E^e(J_{2,\tau}^2)
    +\sum_{\tau=1}^q\frac{16q^2}{\varepsilon_n^{1/2}}
    \E^e\{(G_{n,\tau,1,2}^e)^2\} \\
    &~~~~+
    \sum_{\tau=1}^q\frac{16q^2}{\varepsilon_n^{1/2}}
    \E^e\{(G_{n,\tau,1,3}^e)^2\}
    +\sum_{\tau=1}^q\frac{16q^2}{\varepsilon_n^{1/2}}
    \E^e\{(G_{n,\tau,1,4}^e)^2\}\\
    &\le
    \frac{16q^2}{\varepsilon_n^{1/2}}\sum_{\tau=1}^q
    \{Q_{n5,1}(\tau)+Q_{n5,2}(\tau)+Q_{n5,3}(\tau)+Q_{n5,4}(\tau)\}
    \le K\varepsilon_n^{1/2}\,.
\end{align*}
Hence, restricted on the event $A_{\varepsilon_n}$, we have
\begin{align*}
    \pr\bigg(\frac{G_{n}^{e}}{\sqrt{2q\alpha_2^2}}  \le x\,\bigg|\, \mathcal{X}_n\bigg) & \le \pr\bigg(\frac{G_{n}^{e}}{\sqrt{2q\alpha_2^2}}  \le x,~|G_n^e-\tilde{G}_n^e|\le \varepsilon_n^{1/4} \,\bigg|\, \mathcal{X}_n\bigg) \\
    &~~~~+ \pr\big(|G_n^e-\tilde{G}_n^e|> \varepsilon_n^{1/4}\,|\,\mathcal{X}_n\big)\\
    & \le \pr\bigg(\frac{\tilde{G}_{n}^{e}}{\sqrt{2q\alpha_2^2}}  \le x + \frac{\varepsilon_n^{1/4}}{\sqrt{2q\alpha_2^2}}\,\bigg|\, \mathcal{X}_n\bigg) + K\varepsilon_n^{1/2}
\end{align*}
and
\begin{align*}
    \pr\bigg(\frac{G_{n}^{e}}{\sqrt{2q\alpha_2^2}}  \le x\,\bigg|\, \mathcal{X}_n\bigg) & \ge \pr\bigg(\frac{G_{n}^{e}}{\sqrt{2q\alpha_2^2}}  \le x,~|G_n^e-\tilde{G}_n^e|\le \varepsilon_n^{1/4} \,\bigg|\, \mathcal{X}_n\bigg) \\
    & \ge \pr\bigg(\frac{\tilde{G}_{n}^{e}}{\sqrt{2q\alpha_2^2}}  \le x - \frac{\varepsilon_n^{1/4}}{\sqrt{2q\alpha_2^2}}\,\bigg|\, \mathcal{X}_n\bigg) \\
    &~~~~-\pr\big(|G_n^e-\tilde{G}_n^e|> \varepsilon_n^{1/4}\,|\,\mathcal{X}_n\big)\\
    & \ge \pr\bigg(\frac{\tilde{G}_{n}^{e}}{\sqrt{2q\alpha_2^2}}  \le x - \frac{\varepsilon_n^{1/4}}{\sqrt{2q\alpha_2^2}}\,\bigg|\, \mathcal{X}_n\bigg) - K\varepsilon_n^{1/2}\,.
\end{align*}
Recall $\tilde{{\rm cv}}_{G,\alpha} = \inf\{t \in \mathbb{R}: \pr(G_{n}^{e} \le t\,|\, \mathcal{X}_n) \ge 1-\alpha\}$. 
Together with \eqref{eq:BE2}, restricted on the event $A_{\varepsilon_n}$, it implies that 
\begin{align*}
    \sqrt{2q\alpha_2^2}\Phi^{-1}(1-\alpha-\eta_n-K\varepsilon_n^{1/2})-\varepsilon_n^{1/4} \le  \tilde{{\rm cv}}_{G,\alpha} \le \sqrt{2q\alpha_2^2}\Phi^{-1}(1-\alpha+\eta_n+K\varepsilon_n^{1/2})+\varepsilon_n^{1/4}\,.
\end{align*}
Let $u_{1-\alpha}$ be the $(1-\alpha)$-quantile of $\mathcal{N}(0,1)$. For any fixed $\alpha \in (0,1)$, due to $\eta_n \to 0$ and $\varepsilon_n \to 0$ as $n \to \infty$, there exists a deterministic sequence $\epsilon_n \to 0$ such that $\Phi^{-1}(1-\alpha-\eta_n-K\varepsilon_n^{1/2}) \ge u_{1-\alpha}-\epsilon_n$ and $\Phi^{-1}(1-\alpha+\eta_n+K\varepsilon_n^{1/2}) \le u_{1-\alpha}+\epsilon_n$. Under the null hypothesis, due to $G_{n} \stackrel{\mathrm{d}}{\to} \mathcal{N}(0,2q\alpha_2^2)$ as $n \to \infty$, we have  
\begin{align*}
    \pr(G_{n}> \tilde{{\rm cv}}_{G,\alpha}) & \le \pr(G_{n}> \tilde{{\rm cv}}_{G,\alpha}, A_{\varepsilon_n}) + \pr(A_{\varepsilon_n}^{\mathrm{c}}) \\
    & \le \pr\big\{G_{n}> \alpha_2 \sqrt{2q}(u_{1-\alpha}-\epsilon_n)-\varepsilon_n^{1/4}, A_{\varepsilon_n}\big\} + o(1) \\
    & \le \pr\big\{G_{n}> \alpha_2 \sqrt{2q}(u_{1-\alpha}-\epsilon_n)-\varepsilon_n^{1/4}\big\} + o(1) \\
    & \le \alpha + o(1)\,.
\end{align*}
On the other hand, we also have
\begin{align*}
    \pr(G_{n}> \tilde{{\rm cv}}_{G,\alpha}) & \ge \pr(G_{n}> \tilde{{\rm cv}}_{G,\alpha}, A_{\varepsilon_n})  \\
    & \ge \pr\big\{G_{n}> \alpha_2\sqrt{2q}(u_{1-\alpha}+\epsilon_n)+\varepsilon_n^{1/4}, A_{\varepsilon_n}\big\} \\
    & \ge \pr\big\{G_{n}> \alpha_2\sqrt{2q}(u_{1-\alpha}+\epsilon_n)+\varepsilon_n^{1/4}\big\} - \pr(A_{\varepsilon_n}^{\mathrm{c}}) \\
    & \ge \alpha - o(1)\,.
\end{align*}
Hence, we have $\pr(G_{n}> \tilde{{\rm cv}}_{G,\alpha}) \to \alpha$ as $n \to \infty$. We complete the proof of Theorem \ref{tm:boot_size} with $p \to \infty$. $\hfill\Box$

\subsubsection{Proof of \eqref{q1}}\label{subsubsec:q1}
For $Q_{n1}(\tau)$, we have
\begin{align}
\E\{Q_{n1}(\tau)\}
=&~\frac {4p^2}{n^2}\sum_{j=\tau+1}^n\sum_{t<j-\tau }
\E\{(\br_{t-\tau}^\top \br_{j-\tau})^2\} \E\{(\br_t^\top \br_j)^2\}=\frac {4}{n^2}\sum_{j=\tau+1}^n\sum_{\tau <t<j-\tau }
\alpha_{2p}^2\to 2\alpha_2^2\,,\nonumber\\
\Var\{Q_{n1}(\tau)\}
=&~\E\bigg(\bigg[\frac {4p^2}{n^2}\sum_{j=\tau+1}^n\sum_{\tau <t<j-\tau }
\bigg\{(\br_{t-\tau}^\top \br_{j-\tau})^2(\br_t^\top \br_j)^2-\frac1{p^2}\alpha_{2p}^2\bigg\}\bigg]^2\bigg)\nonumber\\
= &~\frac {16p^4}{n^4}\sum_{j_1, j_2=\tau+1}^n\sum_{\tau <t_1<j_1-\tau\atop \tau < t_2<j_2-\tau}
\E\bigg[\prod_{i=1,2}\bigg\{(\br_{t_i-\tau}^\top \br_{j_i-\tau})^2(\br_{t_i}^\top \br_{j_i})^2-\frac1{p^2}\alpha_{2p}^2\bigg\}\bigg]\,.\label{vqn1}
\end{align}
The expectation in \eqref{vqn1} is zero if $\{t_1, j_1, t_1-\tau, j_1-\tau\}\cap \{t_2, j_2, t_2-\tau, j_2-\tau\}=\emptyset$. Otherwise, it can be controlled as  
\begin{align*}
\Bigg|\E \bigg[\prod_{i=1,2}\bigg\{(\br_{t_i-\tau}^\top \br_{j_i-\tau})^2(\br_{t_i}^\top \br_{j_i})^2-\frac1{p^2}\alpha_{2p}^2\bigg\} \bigg]\Bigg|
\leq&~ \E\bigg[\bigg\{(\br_{1}^\top \br_{2})^2(\br_{3}^\top \br_{4})^2-\frac1{p^2}\alpha_{2p}^2\bigg\}^2\bigg]\\	
\leq &~\E\{(\br_{1}^\top \br_{2})^4\}\E\{(\br_{3}^\top \br_{4})^4\} = O(p^{-4})\,,
\end{align*}
where the last step is based on Assumption \ref{as:A2}.
Thus, $\Var\{Q_{n1}(\tau)\}\leq Kn^{-1}$. 
$\hfill\Box$

\subsubsection{Proof of \eqref{q2}}\label{subsubsec:q2}
To demonstrate the convergence of $Q_{n2}(\tau)$, we split it into two parts as
\begin{align*}
Q_{n2}(\tau)=Q_{n2,1}(\tau)+2Q_{n2,2}(\tau)\,,	
\end{align*}
where 
\begin{align*}
Q_{n2,1}(\tau)=&~	
\frac {p^4}{n^4}\sum_{j=\tau+1}^n\sum_{t<s<j-\tau }
(\br_{s-\tau}^\top \br_{j-\tau}\br_{j-\tau}^\top \br_{t-\tau})^2 (\br_s^\top \br_{j} \br_{j}^\top \br_t)^2\,,\\
Q_{n2,2}(\tau)=&~
\frac {p^4}{n^4}\sum_{\tau+1 \le j_1<j_2 \le n}\sum_{t<s<j_1-\tau}
(\br_{s-\tau}^\top \br_{j_1-\tau}\br_{j_1-\tau}^\top \br_{t-\tau} \br_s^\top \br_{j_1} \br_{j_1}^\top \br_t)
(\br_{s-\tau}^\top \br_{j_2-\tau}\br_{j_2-\tau}^\top \br_{t-\tau} \br_s^\top \br_{j_2} \br_{j_2}^\top \br_t)\,.
\end{align*}
In the sequel, we will derive the bounds for $\E\{|Q_{n2,1}(\tau)|\}$ and $\E\{|Q_{n2,2}(\tau)|\}$.

For $Q_{n2,1}(\tau)$, we have
\begin{align*}
\E \{|Q_{n2,1}(\tau)|\}
\le&~\frac {p^2}{n^4}\sum_{j=\tau+1}^n\sum_{t<s<j-\tau }
|\E(\br_{s-\tau}^\top \br_{j-\tau}\br_{j-\tau}^\top \bSigma\br_{j-\tau}\br_{j-\tau}^\top \br_{s-\tau} \br_t^\top \br_{j} \br_{j}^\top \bSigma\br_{j} \br_{j}^\top \br_t)|\\
\leq&~\frac{Kp^2}{n^4}\sum_{j=\tau+1}^n\sum_{t<s<j-\tau }
{\E}^{1/2}\{(\br_{s-\tau}^\top \br_{j-\tau})^4(\br_{j}^\top \bSigma\br_{j})^2\}{\E}^{1/2}\{(\br_t^\top \br_{j})^4 (\br_{j-\tau}^\top \bSigma\br_{j-\tau})^2\} \\	
= &~ \frac{Kp^2}{n^4}\sum_{j=\tau+1}^n\sum_{t<s<j-\tau }
{\E}\{(\br_{1}^\top \br_{2})^4\}\E\{(\br_{3}^\top \bSigma\br_{3})^2\} \\	
\leq&~Kn^{-1}\,,
\end{align*}
where the first inequality is based on the Cauchy-Schwarz inequality, and the last inequality is based on \eqref{moments} and Assumption \ref{as:A2}. 
Write $\mathcal{G} = \{(j_1,j_2,t,s): \tau+1 \le j_1 < j_2 \le n, t<s<j_1-\tau\}$. Then
\begin{align*}
    Q_{n2,2}(\tau) & =
\frac {p^4}{n^4}\sum_{(j_1,j_2,t,s)\in \mathcal{G} \atop t < s-\tau, \, j_1<j_2-\tau}
(\br_{s-\tau}^\top \br_{j_1-\tau}\br_{j_1-\tau}^\top \br_{t-\tau} \br_s^\top \br_{j_1} \br_{j_1}^\top \br_t)
(\br_{s-\tau}^\top \br_{j_2-\tau}\br_{j_2-\tau}^\top \br_{t-\tau} \br_s^\top \br_{j_2} \br_{j_2}^\top \br_t) \\
&~~~~+ \frac {p^4}{n^4}\sum_{(j_1,j_2,t,s)\in \mathcal{G} \atop t \ge s-\tau \,{\rm or}\, j_1 \ge j_2-\tau}
(\br_{s-\tau}^\top \br_{j_1-\tau}\br_{j_1-\tau}^\top \br_{t-\tau} \br_s^\top \br_{j_1} \br_{j_1}^\top \br_t)
(\br_{s-\tau}^\top \br_{j_2-\tau}\br_{j_2-\tau}^\top \br_{t-\tau} \br_s^\top \br_{j_2} \br_{j_2}^\top \br_t) \\
& = Q_{n2,2,1}(\tau) + Q_{n2,2,2}(\tau)\,.
\end{align*}
To derive the bound of $\E\{|Q_{n2,2}(\tau) |\}$, it suffices to derive the bounds of $\E\{|Q_{n2,2,1}(\tau)|\}$ and $\E\{|Q_{n2,2,2}(\tau)|\}$.

For $Q_{n2,2,1}(\tau)$, it holds that
\begin{align*}
\E \{Q_{n2,2,1}^2(\tau)\}
=&~\frac {p^8}{n^8}\E\bigg[\bigg\{\sum_{(j_1,j_2,t,s)\in \mathcal{G} \atop t < s-\tau,\,j_1<j_2-\tau}
\prod_{i=1,2}(\br_{s-\tau}^\top \br_{j_i-\tau}\br_{j_i-\tau}^\top \br_{t-\tau} \br_s^\top \br_{j_i} \br_{j_i}^\top \br_t)\bigg\}^2\bigg]\\
= &~ \frac{p^8}{n^8}\sum_{(j_1,j_2,t,s)\in \mathcal{G} \atop t < s-\tau,\,j_1<j_2-\tau}\sum_{(j'_1,j'_2,t',s')\in \mathcal{G} \atop t' < s'-\tau,\,j'_1<j'_2-\tau}\E\bigg[\bigg\{
\prod_{i=1,2}(\br_{s-\tau}^\top \br_{j_i-\tau}\br_{j_i-\tau}^\top \br_{t-\tau} \br_s^\top \br_{j_i} \br_{j_i}^\top \br_t)\bigg\} \\
&~~~~~~~~~~~~~~~~~~~~~~~~~~~~~~~~~~~~~~~~~~~~\times \bigg\{
\prod_{i=1,2}(\br_{s'-\tau}^\top \br_{j'_i-\tau}\br_{j'_i-\tau}^\top \br_{t'-\tau} \br_{s'}^\top \br_{j'_i} \br_{j'_i}^\top \br_{t'})\bigg\}\bigg]\\
\leq &~\frac {Kp^8}{n^8}\cdot n^8\max_{(j_1,j_2,t,s)\in \mathcal{G} \atop t < s-\tau,\, j_1 < j_2-\tau}{\E}^2\bigg(\prod_{i=1,2}\br_{s-\tau}^\top \br_{j_i-\tau}\br_{j_i-\tau}^\top \br_{t-\tau} \br_s^\top \br_{j_i} \br_{j_i}^\top \br_t\bigg)\\
&~~~~~~+\frac {Kp^8}{n^8} \cdot n^7\max_{(j_1,j_2,t,s)\in \mathcal{G} \atop t < s-\tau,\, j_1 < j_2-\tau}\E\bigg[\bigg(\prod_{i=1,2}\br_{s-\tau}^\top \br_{j_i-\tau}\br_{j_i-\tau}^\top \br_{t-\tau} \br_s^\top \br_{j_i} \br_{j_i}^\top \br_t\bigg)^2\bigg]\\
\leq&~ K(p^{-4}+n^{-1})\,,	
\end{align*}
where the first inequality is based on the Cauchy-Schwarz inequality,
and the last inequality is based on  the following facts: (i)
\begin{align*}
&\E \bigg(\prod_{i=1,2}\br_{s-\tau}^\top \br_{j_i-\tau}\br_{j_i-\tau}^\top \br_{t-\tau} \br_s^\top \br_{j_i} \br_{j_i}^\top \br_t\bigg)\\
&~~~~~~=
\E (\br_{s-\tau}^\top \br_{j_1-\tau} \br_{j_1-\tau}^\top \br_{t-\tau}\br_{t-\tau}^\top \br_{j_2-\tau} \br_{j_2-\tau}^\top \br_{s-\tau})
\E (\br_s^\top \br_{j_1} \br_{j_1}^\top \br_t\br_t^\top \br_{j_2} \br_{j_2}^\top \br_s )
	\leq Kp^{-6}
\end{align*}
for any $(j_1,j_2,t,s) \in \mathcal{G}$ with $t < s-\tau$ and $j_1 < j_2-\tau$, and (ii)
\begin{align*}
& \E \bigg\{\bigg(\prod_{i=1,2}\br_{s-\tau}^\top \br_{j_i-\tau}\br_{j_i-\tau}^\top \br_{t-\tau} \br_s^\top \br_{j_i} \br_{j_i}^\top \br_t \bigg)^2\bigg\} \nonumber\\
&~~~~~~ = {\E}\{(\br_{s-\tau}^\top \br_{j_1-\tau} \br_{j_2-\tau}^\top \br_{t-\tau})^2(\br_{s-\tau}^\top \br_{j_2-\tau} \br_{j_1-\tau}^\top \br_{t-\tau})^2\}{\E}\{(\br_{s}^\top \br_{j_1} \br_{j_2}^\top \br_{t})^2(\br_{s}^\top \br_{j_2} \br_{j_1}^\top \br_{t})^2\}\nonumber\\
&~~~~~~ \leq {\E}^{1/2}\{(\br_{s-\tau}^\top \br_{j_1-\tau} \br_{j_2-\tau}^\top \br_{t-\tau})^4\}{\E}^{1/2}\{(\br_{s-\tau}^\top \br_{j_2-\tau} \br_{j_1-\tau}^\top \br_{t-\tau})^4\}\nonumber\\
&~~~~~~~~~~~~~~ \times {\E}^{1/2}\{(\br_{s}^\top \br_{j_1} \br_{j_2}^\top \br_{t})^4\}{\E}^{1/2}\{(\br_{s}^\top \br_{j_2} \br_{j_1}^\top \br_{t})^4\}\nonumber\\
&~~~~~~ \leq \E^2\{(\br_1^\top \br_{2}\br_3^\top \br_{4})^4\} =  {\E}^4\{(\br_1^\top \br_{2})^4\}
\leq Kp^{-8}
\end{align*}
for any $(j_1,j_2,t,s) \in \mathcal{G}$ with $t < s-\tau$ and $j_1 < j_2-\tau$. By the Cauchy-Schwarz inequality, we have $\E\{|Q_{n2,2,1}(\tau)|\} \le \E^{1/2}\{Q_{n2,2,1}^2(\tau)\} \le K(p^{-2}+n^{-1/2})$.

For any  $(j_1,j_2,t,s)\in \mathcal{G}$ with $t \ge s-\tau$ or $j_1 \ge j_2-\tau$, the indices $t$ and $j_1$ must satisfy $s-\tau \le t < s$ or $j_2 - \tau \le j_1 < j_2$. Hence, 
\begin{align*}
\E\{|Q_{n2,2,2}(\tau)|\}
\leq &~\frac{K p^4}{n^4} \cdot n^3
\max_{(j_1,j_2,t,s)\in \mathcal{G} \atop t \ge s-\tau \,{\rm or}\, j_1 \ge j_2-\tau} \E\bigg\{\bigg|\prod_{i=1,2}(\br_{s-\tau}^\top \br_{j_i-\tau}\br_{j_i-\tau}^\top \br_{t-\tau} \br_s^\top \br_{j_i} \br_{j_i}^\top \br_t)\bigg|\bigg\}\\
\leq &~\frac{K p^4}{n}
\max_{(j_1,j_2,t,s)\in \mathcal{G} \atop t \ge s-\tau \,{\rm or}\, j_1 \ge j_2-\tau} \E\bigg\{\bigg(\prod_{i=1,2}\br_s^\top \br_{j_i} \br_{j_i}^\top \br_t\bigg)^2\bigg\} \leq \frac{K}{n}\,,
\end{align*}
where the second inequality is based on the Cauchy-Schwarz inequality, and the last inequality is based on the fact
\begin{align*}
\E\bigg\{\bigg(\prod_{i=1,2}\br_s^\top \br_{j_i} \br_{j_i}^\top \br_t\bigg)^2\bigg\} &=  {\E}\{(\br_{s}^\top \br_{j_1} \br_{j_2}^\top \br_{t})^2(\br_{s}^\top \br_{j_2} \br_{j_1}^\top \br_{t})^2\}\nonumber\\
&\leq {\E}^{1/2}\{(\br_{s}^\top \br_{j_1} \br_{j_2}^\top \br_{t})^4\}{\E}^{1/2}\{(\br_{s}^\top \br_{j_2} \br_{j_1}^\top \br_{t})^4\}\nonumber\\
&= \E\{(\br_1^\top \br_{2}\br_3^\top \br_{4})^4\} =  {\E}^2\{(\br_1^\top \br_{2})^4\}
\leq Kp^{-4}
\end{align*}
for any $(j_1,j_2,t,s) \in \mathcal{G}$ with $t \ge s-\tau$ or $j_1 \ge j_2-\tau$.
Collecting all the results, we have \eqref{q2}.
$\hfill\Box$

\subsubsection{Proofs of \eqref{q5}, \eqref{q3} and \eqref{q4}}\label{subsubsec:q345}

For $Q_{n5,1}(\tau)$, notice that
\begin{align*}
    \E\{|Q_{n5,1}(\tau)|\} & = \frac{4p^2}{n^2}\sum_{j=\tau+1}^n \sum_{j-\tau \le t < j}\E\{(\br_{t-\tau}^\top \br_{j-\tau}\br_{t}^\top \br_{j})^2 \} \\ 
    & \le \frac{Kp^2}{n}\max_{\tau+1 \le j \le n \atop j-\tau \le t < j }\E\{(\br_{t-\tau}^\top \br_{j-\tau}\br_{t}^\top \br_{j})^2 \}\\ 
& \le \frac {Kp^2}{n}\E^{1/2} \{(\br_1^\top \br_2)^4\}\E^{1/2} \{(\br_3^\top \br_4)^4\} \leq \frac{K}{n}
\end{align*}
for any $\tau \in [q]$, where the second inequality is based on the Cauchy-Schwarz inequality, and the last inequality is based on Assumption \ref{as:A2}. Hence, we have \eqref{q5}.

For $Q_{n3}(\tau,\tau')$ with any $\tau > \tau'$, by the Cauchy-Schwarz inequality, we have
\begin{align*}
\E\{|Q_{n3}(\tau,\tau')|\}
\le &~\frac {Kp^4}{n^4} \sum_{j=\tau+1}^n\sum_{t<j-\tau \atop s<j-\tau'}
\E\{(\br_{j-\tau}^\top  \br_{t-\tau} \br_t^\top \br_j \br_{j-\tau'}^\top \br_{s-\tau'}\br_s^\top \br_j)^2\}\\
&~~~~~~+\frac {Kp^4}{n^4} \sum_{j=\tau+1}^n\sum_{t<j-\tau-\tau'}
\E\big\{\br_{j-\tau}^\top  \br_{t-\tau} \br_t^\top \br_j \br_{j-\tau'}^\top \br_{t}(\br_{t+\tau'}^\top \br_j)^2\\
&~~~~~~~~~~~~~~~~~~~~~~~~~~~~~~~~~~~~~\times\br_{j-\tau}^\top  \br_{t-\tau+\tau'}\br_{j-\tau'}^\top \br_{t-\tau}\br_{t-\tau+\tau'}^\top\br_{j}\big\}\\
\leq &~ \frac {Kp^4}{n} \max_{\tau+1 \le j \le n \atop t<j-\tau,\; s<j-\tau'}
\E^{1/2}\{(\br_{j-\tau}^\top  \br_{t-\tau} \br_t^\top \br_j)^4\}\E^{1/2}\{(\br_{j-\tau'}^\top  \br_{s-\tau'} \br_s^\top \br_j)^4\} \\
&~~~~~~+ \frac {Kp^4}{n^2} \max_{\tau+1 \le j \le n \atop t<j-\tau-\tau'}
\E^{1/4}\{(\br_{j-\tau'}^\top  \br_{t-\tau} \br_{t+\tau'}^\top \br_j)^4\}\E^{1/4}\{(\br_{j-\tau'}^\top  \br_{t} \br_{t+\tau'}^\top \br_j)^4\} \\
&~~~~~~~~~~~~~~~~~~~~~~~~~~~\times\E^{1/4}\{(\br_{j-\tau}^\top  \br_{t-\tau+\tau'} \br_t^\top \br_j)^4\}\E^{1/4}\{(\br_{j-\tau}^\top  \br_{t-\tau} \br_{t-\tau+\tau'}^\top \br_j)^4\} \\
\le &~ \frac {Kp^4}{n} \E\{(\br_{1}^\top  \br_{2} \br_3^\top \br_4)^4\} = \frac {Kp^4}{n}\E\{(\br_{1}^\top \br_{2})^4\}\E\{(\br_{3}^\top \br_{4})^4\}
\leq \frac{K}{n}\,,
\end{align*}
where the last inequality is based on Assumption \ref{as:A2}. 
Hence, we have \eqref{q3}.

For $Q_{n4}(\tau)$, by the Cauchy-Schwarz inequality, we have
\begin{align*}
\E\{|Q_{n4}(\tau)|\}
\le &~
\frac {p^4}{n^4}\sum_{j=\tau+1}^n\sum_{t_1\leq t_3<j-\tau}\E\bigg\{\bigg(\prod_{i=1, 3}\br_{j-\tau}^\top \br_{t_i-\tau}\br_{t_i}^\top \br_j\bigg)^2\bigg\} \\
\leq&~ \frac{Kp^4}{n}\E\{(\br_{1}^\top \br_{2})^4\}\E\{(\br_{3}^\top \br_{4})^4\} \leq  \frac{K}{n}\,,
\end{align*}
where the last inequality is based on Assumption \ref{as:A2}. 
Hence, we have \eqref{q4}.
 $\hfill\Box$

\subsection{Proof of Theorem  \ref{tm:boot_size} with fixed $p$}\label{subsec:size_pfix}
Recall $\mathcal{X}_n = \{\x_1,\ldots,\x_n\}$. To show $\pr_{H_0}(T_n > \tilde{{\rm cv}}_{\alpha}) \to \alpha$ with $\tilde{{\rm cv}}_{\alpha} = \inf\{t \in \mathbb{R}: \pr(T_{n}^{e} \le t\,|\, \mathcal{X}_n) \ge 1-\alpha\}$, it suffices to show that under the null hypothesis, the conditional distribution of $T_n^e$ given $\mathcal{X}_n$ converges to the limiting distribution of $T_n$. 
Recall 
\begin{align*}
    G_{n, \tau,1} = \frac {1}{np}\sum_{t \ne s}\x_t^\top \x_s\x_{t+\tau}^\top \x_{s+\tau} := \frac{n}{p}{\tr}(\bS_\tau\bS_\tau^\top) - M_{n,\tau}
\end{align*}
as defined in \eqref{eq:decomp} in Section \ref{subsec:null_pfix}. By convention, we let  $e_{t}=e_{t-n}$ if $t>n$, and
$$
\bS_\tau^e=\frac1{n}\sum_{t=1}^{n}e_te_{t+\tau}\x_t\x_{t+\tau}^\top\,.
$$
By the fact $e_t^2 = 1$ for $t \in [n]$, we have
\begin{align*}
    G_{n, \tau,1}^e = \frac 1{np}\sum_{t\neq s}e_t e_s e_{t+\tau} e_{s+\tau}\x_t^\top \x_s\x_{t+\tau}^\top \x_{s+\tau}  = \frac{n}{p}{\tr}\{\bS_\tau^e (\bS_\tau^e)^\top\} - M_{n,\tau}\,.
\end{align*}

Recall 
$$T_n = \sum_{\tau=1}^q H_{n,\tau}=\frac{1}{\hat\sigma_{n1}}\sum_{\tau=1}^q G_{n, \tau,1} = \frac{1}{\hat\sigma_{n1}}
    \sum_{\tau=1}^q
    \bigg[\frac np\tr(\bS_\tau\bS_\tau^\top)
    -M_{n,\tau}\bigg]\,.$$ 
By \eqref{eq:jointHn}, we have 
$T_n\xrightarrow{\mathrm d}T_\infty$ under the null hypothesis, where
\begin{align}\label{eq:Tinf_fixedp}
    T_\infty
    &=\frac{1}{\sqrt2\,\tr(\bSigma^2)}
    \sum_{\tau=1}^q
    \left\{\tr(\bSigma  N_\tau\bSigma N_\tau^\top)
    -{\tr}^2(\bSigma)\right\} \stackrel{\mathrm d}{=}
    \frac{\sum_{\tau=1}^q\sum_{i,j=1}^p
    \lambda_i\lambda_j(\varepsilon_{i,j,\tau}^2-1)}
    {\sqrt2\sum_{i=1}^p\lambda_i^2}\,.
\end{align}
Here $N_\tau=(\varepsilon_{i,j,\tau})\in\mathbb R^{p\times p}$, where $\varepsilon_{i,j,\tau}$ are i.i.d.\ $\mathcal{N}(0,1)$ variables for any $i,j\in [p]$ and $\tau\in[q]$, and $\lambda_1 \ge\cdots \ge\lambda_p$ are the eigenvalues of $\bSigma$.
The distribution function $F(\cdot)$ of $T_\infty$ is continuous and
strictly increasing on the interior of its support. For a fixed significance level $\alpha$, let $c_\alpha=F^{-1}(1-\alpha)$. To show $\pr_{H_0}(T_n > \tilde{{\rm cv}}_{\alpha}) \to \alpha$, it suffices to show
\begin{align}\label{eq:probConverge}
        \tilde{\mathrm{cv}}_{\alpha,n}:=\tilde{\mathrm{cv}}_\alpha
    \stackrel{\mathrm p}{\to}c_\alpha
    ~~\text{as }n\to\infty\,,
\end{align}
which is equivalent to show that, for any subsequence $\{n_k\}$,
there exists a further subsequence $\{n_{k_\ell}\}$ such that
\begin{align}\label{eq:as}
        \tilde{\mathrm{cv}}_{\alpha,n_{k_\ell}}
    \xrightarrow{\mathrm{a.s.}}c_\alpha ~~\text{as }\ell\to\infty\,.
\end{align}

Define
\begin{align}
    W_{n1} & = \frac{1}{n}\sum_{\tau=1}^q\sum_{t=n-\tau+1}^{n}\sum_{k,l=1}^p x_{k,t}^2 x_{l,t+\tau}^2\,,\label{eq:Wn1}\\
    W_{n2} & = \frac{1}{n^2}\sum_{\tau=1}^q \sum_{k,l,k',l'=1}^p\bigg\{\sum_{j=\tau+1}^n(x_{k,j-\tau} x_{k',j-\tau}x_{l,j}x_{l',j}-\sigma_{k,k'} \sigma_{l,l'})\bigg\}^2\,,\label{eq:Wn2}\\
    W_{n3}(\tau,\tau') & = \frac{1}{n^2}\sum_{j=1}^n\sum_{k,l,k',l'=1}^p x_{k,j-\tau}^2 x_{l,j}^2x_{k',j-\tau'}^2 x_{l',j}^2\,,\label{eq:Wn3}\\
    W_{n4} &= \frac{1}{n^2}\sum_{j=1}^n \sum_{\tau=1}^q \sum_{k,l=1}^p x_{k,j-\tau}^4 x_{l,j}^4\,.\label{eq:Wn4}
\end{align}
Under the null hypothesis, as we will show in Section \ref{subsubsec:q123_pfix}, 
\begin{align}
W_{n1} &\stackrel{\mathrm{p}}{\to }0 \,,\label{q1_pfix} \\
W_{n2} &\stackrel{\mathrm{p}}{\to} 0 \,,\label{eq:q2_pfix}\\
W_{n3}(\tau,\tau')&\stackrel{\mathrm{p}}{\to} 0 ~~\mbox{for any}~~ 1 \le \tau' < \tau \le q\,, \label{eq:q3_pfix}\\
W_{n4} &\stackrel{\mathrm{p}}{\to} 0 \label{eq:q4_pfix}\,.
\end{align}
Fix an arbitrary subsequence $\{n_k\}$. Based on \eqref{q1_pfix}--\eqref{eq:q4_pfix}, we can select a further subsequence $\{n_{k_\ell}\}$ such that
$\pr(\Omega_1)=1$, where
\begin{align*}
    \Omega_1 &=\bigg\{\omega: \lim_{\ell \to \infty}W_{n_{k_\ell},1}(\omega)=0,~
        \lim_{\ell \to \infty}W_{n_{k_\ell},2}(\omega)= 0,~
        \lim_{\ell \to \infty}W_{n_{k_\ell},4}(\omega)=0,\\
        &~~~~~~~~~~\lim_{\ell \to \infty}\max_{1\le\tau'<\tau\le q}
        |W_{n_{k_\ell},3}(\tau,\tau')(\omega)| = 0, ~\lim_{\ell \to \infty}\max_{\tau\in[q]}
        |M_{n_{k_\ell},\tau}(\omega)-p^{-1}{\tr}^2(\bSigma)| = 0,\\
        &~~~~~~~~~~\lim_{\ell \to \infty}|\hat\sigma_{n_{k_\ell}1}(\omega)
        -\sqrt{2}p^{-1}\tr(\bSigma^2)|=0 \bigg\}\,.
\end{align*}
Write 
\begin{align*}
  \bV_n^e = \sqrt{n}\begin{pmatrix}\mathrm{vec}(\bS_1^e)\\
  \vdots\\ \mathrm{vec}(\bS_q^e)
  \end{pmatrix} \,.
\end{align*}
In the sequel, we first show 
\begin{align}\label{eq:joint_cond}
        \bV_{n_{k_\ell}}^e\,|\,
    \mathcal X_{n_{k_\ell}}
    \xrightarrow{\mathrm d}
    \mathcal N(0,\I_q\otimes\bSigma\otimes\bSigma)~~\mbox{as}~\ell \to \infty
\end{align}
almost surely. Based on the Cram\'er--Wold theorem, it suffices to show that there exists $\Omega_0$, a subset  of the sample space $\Omega$, with $\pr(\Omega_0)=1$ such that 
\begin{align}\label{eq:joint_cond_ge0}
\beta^\top\bV_{n_{k_\ell}}^e\,|\,
    \mathcal X_{n_{k_\ell}}=\mathcal X_{n_{k_\ell}}(\omega)
    \xrightarrow{\mathrm d}
    \mathcal N(0,\beta^\top(\I_q\otimes\bSigma\otimes\bSigma)\beta) ~~\mbox{as}~\ell \to \infty
\end{align}
for any $\omega \in \Omega_0$ and nonzero $\beta$.
Write 
\begin{align*}
  \beta = \begin{pmatrix}\mathrm{vec}(\C_1)\\
  \vdots\\ \mathrm{vec}(\C_q)
  \end{pmatrix}\,,
\end{align*}
where $\C_\tau = (c_{\tau,k,l})_{p\times p}$ for each $\tau \in [q]$. 
Let $\bSigma = (\sigma_{i,j})_{p\times p}$, and $$\sigma_{\beta}^2 = \beta^\top (\I_q \otimes \bSigma \otimes \bSigma) \beta= \sum_{\tau=1}^q \sum_{u,v,k,l=1}^p c_{\tau,u,v}c_{\tau,k,l}\sigma_{u,k}\sigma_{v,l}\,.$$ 
To prove \eqref{eq:joint_cond_ge0}, we need to consider two scenarios: (i) the smallest eigenvalue of $\bSigma$ is positive, and (ii) the smallest eigenvalue of $\bSigma$ is zero. 

If the smallest eigenvalue of $\bSigma$ is positive, we have $\sigma_{\beta}^2 > 0$ for any nonzero $\beta$. As we will show in Section \ref{subsec:projection_pfix}, for any $\varepsilon>0$, we have
\begin{align}\label{eq:projection_fixedp}
    &\sup_{x\in \mathbb{R}}\bigg|\pr\bigg(\frac{\beta^\top\bV_n^e}{\sigma_{\beta}} \le x\,\bigg|\, \mathcal{X}_n\bigg)-\Phi(x)\bigg| \\
    &~~~~~~~~~\le \frac{\varepsilon}{\sigma_{\beta}\sqrt{2\pi}} + K^\dagger(\beta) \bigg\{n^{-2}+\max_{1 \le \tau' < \tau \le q}W_{n3}(\tau,\tau') + W_{n2} + W_{n4}\bigg\}^{1/5}+ K_1(\beta)\varepsilon^{-2}W_{n1}\,,\notag
\end{align}
where $K^\dagger(\beta)>0$ and $K_1(\beta)>0$ are some constants independent of $n$. Now fix any $\omega\in\Omega_1$. Applying \eqref{eq:projection_fixedp} along the selected subsequence $\{n_{k_{\ell}}\}$, first letting $\ell\to\infty$ and then $\varepsilon \to 0$, we have
\[
    \sup_{x\in\mathbb R}
    \bigg|\pr\big\{
    \beta^\top\bV_{n_{k_\ell}}^e\le x
    \,|\,
    \mathcal X_{n_{k_\ell}}=\mathcal X_{n_{k_\ell}}(\omega)
    \big\}-\Phi\bigg(\frac{x}{\sigma_{\beta}}\bigg)\bigg|
    \to 0\,,
\]
which implies \eqref{eq:joint_cond_ge0} holds 
with $\Omega_0=\Omega_1$.

If the smallest eigenvalue of $\bSigma$ is zero, we know the smallest eigenvalue of
$\I_q\otimes\bSigma\otimes\bSigma$ is also zero. Let $\ker(\I_q\otimes\bSigma\otimes\bSigma) = \{\alpha: (\I_q\otimes\bSigma\otimes\bSigma)\alpha=0\}$ and $r = \mathrm{dim}\{\ker(\I_q\otimes\bSigma\otimes\bSigma)\} \ge 1$. For any $\beta \notin \ker(\I_q\otimes\bSigma\otimes\bSigma)$, identical to the arguments for the scenario with positive definite $\bSigma$, we know \eqref{eq:joint_cond_ge0} holds with $\Omega_0=\Omega_1$. Notice that the normal distribution $\mathcal{N}(0,\beta^\top(\I_q\otimes\bSigma\otimes\bSigma)\beta)$ in \eqref{eq:joint_cond_ge0} 
degenerates to the point mass at zero when $\beta \in \ker(\I_q\otimes\bSigma\otimes\bSigma)$. 
Write  $\x_t=(x_{1,t},\ldots,x_{p,t})^\top$. Notice that the $(k,l)$-th entry of $\sqrt{n}\bS_\tau^e$, denoted by $s_{\tau,k,l}^e$, can be expressed as $s_{\tau,k,l}^e=n^{-1/2}\sum_{t=1}^{n}e_te_{t+\tau} x_{k,t} x_{l,t+\tau}$.
It is easy to verify that $\E^e(\beta^\top\bV_n^e)=0$.
By the fact $e_t^2=1$ for $t\in [n]$, we have 
\begin{align*}
    \E^e\{(\beta^\top\bV_n^e)^2\} & = \frac{1}{n}\sum_{t=1}^{n}\sum_{\tau=1}^q\sum_{k,l,k',l'=1}^p c_{\tau,k,l}c_{\tau,k',l'} x_{k,t} x_{k',t}x_{l,t+\tau}x_{l',t+\tau}\,,
\end{align*}
which implies
\begin{align*}
    \E[\E^e\{(\beta^\top\bV_n^e)^2\}]& = \frac{1}{n}\sum_{t=1}^{n}\sum_{\tau=1}^q\sum_{k,l,k',l'=1}^p c_{\tau,k,l}c_{\tau,k',l'} \E(x_{k,t} x_{k',t}) \E(x_{l,t+\tau}x_{l',t+\tau})\\
     & = \frac{1}{n}\sum_{t=1}^{n}\sum_{\tau=1}^q\sum_{k,l,k',l'=1}^p c_{\tau,k,l}c_{\tau,k',l'} \sigma_{k,k'} \sigma_{l,l'} =\beta^{\top}(\I_q \otimes \bSigma \otimes \bSigma)\beta\,.
\end{align*}
Hence, $\E\{{\Var}(\beta^\top\bV_n^e\,|\,\mathcal{X}_n)\} = \beta^{\top}(\I_q \otimes \bSigma \otimes \bSigma)\beta$, which implies $\E\{{\Var}(\beta^\top\bV_n^e\,|\,\mathcal{X}_n)\} =0$ if and only if $\beta \in \ker(\I_q\otimes\bSigma\otimes\bSigma)$. Let $\{\alpha_1,\ldots,\alpha_r\}$ be a given basis of $\ker(\I_q\otimes\bSigma\otimes\bSigma)$. We know ${\Var}(\alpha_j^\top\bV_n^e\,|\,\mathcal{X}_n) =0$ almost surely for any $j \in [r]$. Since $\E^e(\alpha_j^\top\bV_n^e)=0$, then 
\begin{align*}
    \Omega_2 &=\big\{\omega: \pr\{\alpha_j^\top\bV_{n_{k_\ell}}^e=0
        \,|\,\mathcal X_{n_{k_\ell}}=\mathcal X_{n_{k_\ell}}(\omega)\}=1~\text{for all }\ell \ge 1 \text{ and } j \in [r] \big\}
\end{align*}
satisfies $\pr(\Omega_2)=1$.
Hence, for any $\beta \in \ker(\I_q\otimes\bSigma\otimes\bSigma)$ and $\omega \in \Omega_2$, we have $\pr\{\beta^\top\bV_{n_{k_\ell}}^e=0 \,|\,\mathcal X_{n_{k_\ell}}=\mathcal X_{n_{k_\ell}}(\omega)\}=1$ for all $\ell \ge 1$. Hence, \eqref{eq:joint_cond_ge0} also holds with $\Omega_0=\Omega_1\cap\Omega_2$.

Recall
\[
    T_n^e=\frac{1}{\hat\sigma_{n1}}
    \sum_{\tau=1}^q
    \bigg[\frac np\tr\{\bS_\tau^e(\bS_\tau^e)^\top\}
    -M_{n,\tau}\bigg]\,.
\]
Based on \eqref{eq:joint_cond}, we will show in Section \ref{subsec:Tne} that
\begin{align}\label{eq:cond_Tne}
    T_{n_{k_\ell}}^e \,|\,
    \mathcal X_{n_{k_\ell}}=\mathcal X_{n_{k_\ell}}(\omega) \xrightarrow{\mathrm d} T_{\infty}
\end{align}
for every fixed $\omega\in\Omega_0$.
Since the distribution function $F(\cdot)$ of $T_{\infty}$ is continuous, P\'olya's theorem implies that there exists a sequence
$\eta_\ell(\omega)\to0$ such that
\[
    \sup_{x\in\mathbb R}
    \big|\pr\big\{T_{n_{k_\ell}}^e\le x
    \,|\,
    \mathcal X_{n_{k_\ell}}=\mathcal X_{n_{k_\ell}}(\omega)
    \big\}-F(x)\big|<\eta_\ell(\omega)
\]
for all sufficiently large $\ell$.
This implies that $F^{-1}\{1-\alpha-\eta_\ell(\omega)\} \le \tilde{\mathrm{cv}}_{\alpha,n_{k_\ell}}(\omega) \le F^{-1}\{1-\alpha+\eta_\ell(\omega)\}$, where $\tilde{\mathrm{cv}}_{\alpha,n_{k_\ell}}(\omega) = \inf\{t\in \mathbb{R}: \pr\{T_{n_{k_\ell}}^e\le t\,|\,\mathcal X_{n_{k_\ell}}=\mathcal X_{n_{k_\ell}}(\omega)\} \ge 1-\alpha\}$. Since $F(\cdot)$ is continuous and
strictly increasing on the interior of its support and $\eta_\ell(\omega)\to0$ as $\ell \to \infty$, we have $\tilde{\mathrm{cv}}_{\alpha,n_{k_\ell}}(\omega)\to c_\alpha$ as $\ell \to \infty$
for every $\omega\in\Omega_0$.
Since $\pr(\Omega_0)=1$, we have \eqref{eq:as} holds.
We complete the proof of Theorem \ref{tm:boot_size} with fixed $p$. $\hfill\Box$

\subsubsection{Proofs of \eqref{q1_pfix}--\eqref{eq:q4_pfix}}\label{subsubsec:q123_pfix}

Notice that $p$ and $q$ are fixed constants.
For $W_{n1}$, we have
\begin{align*}
    \E(|W_{n1}|) & = \frac{1}{n}\sum_{\tau=1}^q\sum_{t=n-\tau+1}^{n}\sum_{k,l=1}^p \E(x_{k,t}^2)\E(x_{l,t+\tau}^2) \le \frac{K}{n} \sum_{k,l=1}^p  \sigma_{k,k}\sigma_{l,l} = \frac{K}{n}{\tr}^2(\bSigma) \le \frac{K}{n}\,.
\end{align*}
Hence, by Markov's inequality, \eqref{q1_pfix} holds.

For $W_{n2}$, we have
\begin{align*}
     \E(|W_{n2}|) & = \frac{1}{n^2}\sum_{\tau=1}^q \sum_{k,l,k',l'=1}^p\E\bigg[\bigg\{\sum_{j=\tau+1}^n(x_{k,j-\tau} x_{k',j-\tau}x_{l,j}x_{l',j}-\sigma_{k,k'} \sigma_{l,l'})\bigg\}^2\bigg] \\
     & = \frac{1}{n^2}\sum_{\tau=1}^q \sum_{k,l,k',l'=1}^p\sum_{j=\tau+1}^n\E\{(x_{k,j-\tau} x_{k',j-\tau}x_{l,j}x_{l',j}-\sigma_{k,k'} \sigma_{l,l'})^2\}\\
     &~~~~~~+\frac{2}{n^2}\sum_{\tau=1}^q \sum_{k,l,k',l'=1}^p\sum_{j=2\tau+1}^n\E\big\{(x_{k,j-2\tau} x_{k',j-2\tau}x_{l,j-\tau}x_{l',j-\tau}-\sigma_{k,k'} \sigma_{l,l'}) \\
     &~~~~~~~~~~~~~~~~~~~~~~~~~~~~~~~~~~~~~~~~~~~~\times(x_{k,j-\tau} x_{k',j-\tau}x_{l,j}x_{l',j}-\sigma_{k,k'} \sigma_{l,l'})\big\}\\
     & \le \frac{K}{n}\,.
\end{align*}
Hence, by Markov's inequality, \eqref{eq:q2_pfix} holds.

For $W_{n3}(\tau,\tau')$ with any $\tau > \tau'$, we have
\begin{align*}
\E\{|W_{n3}(\tau,\tau')|\}
= &~ \frac{1}{n^2}\sum_{j=1}^n\sum_{k,l,k',l'=1}^p \E(x_{k,j-\tau}^2 x_{l,j}^2x_{k',j-\tau'}^2 x_{l',j}^2)\\
\le &~ \frac{K}{n}\sum_{k,l,k',l'=1}^p \sigma_{k,k}\sigma_{k',k'} \E(x_{l,1}^2 x_{l',1}^2) \le \frac{K}{n}\,.
\end{align*}
Hence, by Markov's inequality, \eqref{eq:q3_pfix} holds.

For $W_{n4}$, we have
\begin{align*}
    \E(|W_{n4}|) & = \frac{1}{n^2} \sum_{\tau=1}^q \sum_{j=\tau+1}^n \sum_{k,l=1}^p \E(x_{k,j-\tau}^4) \E(x_{l,j}^4) \le \frac{K}{n}\,.
\end{align*}
Hence, by Markov's inequality, \eqref{eq:q4_pfix} holds.
$\hfill\Box$

\subsubsection{Proof of \eqref{eq:projection_fixedp}}\label{subsec:projection_pfix}
In this section, we consider the case with $\sigma_{\beta}^2 > 0$. Let $s_{\tau,k,l}^e$ denote the $(k,l)$-th element of $\sqrt{n}\bS_\tau^e$.  
Notice that 
 \begin{align*}
     s_{\tau,k,l}^e
     =&~\frac1{\sqrt{n}}\sum_{t=1}^{n}e_te_{t+\tau} x_{k,t} x_{l,t+\tau}\\
     =&~\frac1{\sqrt{n}}\sum_{t=1}^{n-\tau}e_te_{t+\tau} x_{k,t} x_{l,t+\tau}+\frac1{\sqrt{n}}\sum_{t=n-\tau+1}^{n}e_te_{t+\tau} x_{k,t} x_{l,t+\tau}\,,  
 \end{align*}
which implies 
\begin{align*}
\beta^\top\bV_n^e & = \sum_{\tau=1}^q \sum_{k,l=1}^p c_{\tau,k,l}s_{\tau,k,l}^e \\
&= \frac{1}{\sqrt{n}}\sum_{\tau=1}^q \sum_{k,l=1}^p c_{\tau,k,l}\sum_{t=1}^{n-\tau}e_te_{t+\tau} x_{k,t} x_{l,t+\tau} + \frac{1}{\sqrt{n}}\sum_{\tau=1}^q \sum_{k,l=1}^p c_{\tau,k,l}\sum_{t=n-\tau+1}^{n}e_te_{t+\tau} x_{k,t} x_{l,t+\tau} \\
&= \tilde{J}_1 + \tilde{J}_2\,.
\end{align*}
It is easy to verify that $\E^e(\tilde{J}_2)=0$. By the fact $e_t^2=1$ for $t \in [n]$ and the Cauchy-Schwarz inequality, we have
\begin{align}\label{eq:J2ttt}
    \E^e(\tilde{J}_2^2) & = \frac{1}{n}\sum_{\tau=1}^q \sum_{k,l,k',l'=1}^p c_{\tau,k,l}c_{\tau,k',l'}\sum_{t=n-\tau+1}^{n} x_{k,t} x_{l,t+\tau}x_{k',t} x_{l',t+\tau}\notag\\
    & \le \frac{1}{n}\sum_{\tau=1}^q \sum_{t=n-\tau+1}^{n}\bigg(\sum_{k,l=1}^pc_{\tau,k,l}^2\bigg)\bigg(\sum_{k,l=1}^px_{k,t}^2 x_{l,t+\tau}^2\bigg)\notag\\
    & \le \max_{\tau \in [q]}\sum_{k,l=1}^pc_{\tau,k,l}^2 \cdot\frac{1}{n}\sum_{\tau=1}^q\sum_{t=n-\tau+1}^{n}\sum_{k,l=1}^p x_{k,t}^2 x_{l,t+\tau}^2 =K_1(\beta)W_{n1}\,,
\end{align}
where $W_{n1}$ is defined in \eqref{eq:Wn1}, and $K_1(\beta) = \max_{\tau \in [q]}\sum_{k,l=1}^pc_{\tau,k,l}^2 >0$ is some constant  independent of $n$.
By convention, we set $\x_\ell=0$ for $\ell\leq 0$. Moreover, we have
\begin{align*}
    \tilde{J}_1 = &~ \sum_{j=1}^n \bigg\{\E_j^e \bigg(\frac{1}{\sqrt{n}}\sum_{\tau=1}^q \sum_{k,l=1}^p c_{\tau,k,l}\sum_{t=1}^{n-\tau}e_te_{t+\tau} x_{k,t} x_{l,t+\tau}\bigg)  \\
    &~~~~~~~~~~~~~- \E_{j-1}^e \bigg(\frac{1}{\sqrt{n}}\sum_{\tau=1}^q \sum_{k,l=1}^p c_{\tau,k,l}\sum_{t=1}^{n-\tau}e_te_{t+\tau} x_{k,t} x_{l,t+\tau}\bigg) \bigg\} \\
    = &~ \frac{1}{\sqrt{n}}\sum_{j=1}^n \sum_{\tau=1}^q \sum_{k,l=1}^pc_{\tau,k,l}e_{j-\tau}e_{j} x_{k,j-\tau} x_{l,j} = \sum_{j=1}^n B_{j}^e
\end{align*}
with
\begin{align*}
    B_{j}^e = \frac{1}{\sqrt{n}}\sum_{\tau=1}^q \sum_{k,l=1}^p c_{\tau,k,l}e_{j-\tau}e_{j} x_{k,j-\tau} x_{l,j}\,.
\end{align*}
Given $\mathcal{X}_n$, it can be shown that $\{B_{j}^e\}$ forms a sequence of martingale differences with respect to $\{\mathcal F_j^e\}$. 
Write
\begin{align*}
    \frac{\tilde{J}_1}{\sigma_{\beta}} = \sum_{j=1}^n Y_j^e ~~\mbox{with}~~
    Y_j^e = \frac{B_{j}^e}{\sigma_{\beta}}\,.
\end{align*}
By Theorem 3.9 of \cite{Hall1980} with $\delta=1$, for all $x \in \mathbb{R}$, there exists a universal constant $C>0$ that does not depend on $x$, such that whenever $\tilde{L}_n \le 1$,
\begin{align}\label{eq:BEbound_pfix}
    \bigg|\pr\bigg(\frac{\tilde{J}_1}{\sigma_{\beta}} \le x\,\bigg|\, \mathcal{X}_n\bigg)-\Phi(x)\bigg| \le C\tilde{L}_n^{1/5}(1+|x|^{16/5})^{-1}\,,
\end{align}
 where
\begin{align}\label{eq:Ln_pfix}
    \tilde{L}_n = \E^e\{(\tilde{V}_n^2-1)^2\} + \sum_{j=1}^n \E^e\{(Y_j^e)^4\} ~~\mbox{with}~~ \tilde{V}_n^2 = \sum_{j=1}^n \E^e_{j-1}\{(Y_j^e)^2\}\,.
\end{align}
By the fact $e_\ell^2=1$ for $\ell\in [n]$, we have
\begin{align}\label{eq:Vn2_pfix}
    \tilde{V}_n^2 & = \frac{1}{\sigma_{\beta}^2}\sum_{j=1}^n \E^e_{j-1}\{(B_{j}^e)^2\} \notag\\
    & = \frac{1}{\sigma_{\beta}^2n}\sum_{j=1}^n\sum_{\tau,\tau'=1}^q \sum_{k,l,k',l'=1}^pc_{\tau,k,l}c_{\tau',k',l'}e_{j-\tau}e_{j-\tau'} x_{k,j-\tau} x_{l,j}x_{k',j-\tau'} x_{l',j} \notag\\
    & = \frac{1}{\sigma_{\beta}^2n}\sum_{j=1}^n\sum_{\tau=1}^q \sum_{k,l,k',l'=1}^pc_{\tau,k,l}c_{\tau,k',l'} x_{k,j-\tau} x_{l,j}x_{k',j-\tau} x_{l',j} \\
    &~~~~~+\frac{2}{\sigma_{\beta}^2n}\sum_{j=1}^n\sum_{1 \le \tau' < \tau \le q} \sum_{k,l,k',l'=1}^pc_{\tau,k,l}c_{\tau',k',l'}e_{j-\tau}e_{j-\tau'} x_{k,j-\tau} x_{l,j}x_{k',j-\tau'} x_{l',j} \notag\\
    &:=\tilde{W}_{n2} + \sum_{1 \le \tau' < \tau \le q}\tilde{W}_{n3}^e(\tau,\tau')\,,
\end{align}
where 
\begin{align*}
    \tilde{W}_{n2} & = \frac{1}{\sigma_{\beta}^2n}\sum_{j=1}^n\sum_{\tau=1}^q \sum_{k,l,k',l'=1}^pc_{\tau,k,l}c_{\tau,k',l'} x_{k,j-\tau} x_{l,j}x_{k',j-\tau} x_{l',j}\,,\\
    \tilde{W}_{n3}^e(\tau,\tau') & = \frac{2}{\sigma_{\beta}^2n}\sum_{j=1}^n \sum_{k,l,k',l'=1}^pc_{\tau,k,l}c_{\tau',k',l'}e_{j-\tau}e_{j-\tau'} x_{k,j-\tau} x_{l,j}x_{k',j-\tau'} x_{l',j}\,.
\end{align*}
Notice that $\E^e\{\tilde{W}_{n3}^e(\tau,\tau')\}=0$ for any $1\le \tau' < \tau \le q$, which implies $\E^e(\tilde{V}_n^2)=\tilde{W}_{n2}$.
Hence, it holds that 
\begin{align}\label{eq:EVn2_pfix}
    \E^e\{(\tilde{V}_n^2-1)^2\} & = \E^e[\{\tilde{V}_n^2 -\E^e(\tilde{V}_n^2)\}^2] + \{\E^e(\tilde{V}_n^2) -1\}^2 \notag\\
    &=\E^e\bigg[\bigg\{\sum_{1 \le \tau' < \tau \le q}\tilde{W}_{n3}^e(\tau,\tau')\bigg\}^2\bigg] + (\tilde{W}_{n2}-1)^2 \notag\\
    & \le \frac{q^2(q-1)^2}{4} \max_{1 \le \tau' < \tau \le q}\E^e[\{\tilde{W}_{n3}^e(\tau,\tau')\}^2] + (\tilde{W}_{n2}-1)^2 \,.
\end{align}
Recall $\sigma_{\beta}^2 = \sum_{\tau=1}^q \sum_{u,v,k,l=1}^p c_{\tau,u,v}c_{\tau,k,l}\sigma_{u,k}\sigma_{v,l}$. By convention, we set $\x_t=0$ for $t\leq 0$. Hence, we have
\begin{align*}
    \tilde{W}_{n2}-1 & = \frac{1}{\sigma_{\beta}^2n}\sum_{j=1}^n\sum_{\tau=1}^q \sum_{k,l,k',l'=1}^pc_{\tau,k,l}c_{\tau,k',l'} (x_{k,j-\tau} x_{k',j-\tau}x_{l,j}x_{l',j}-\sigma_{k,k'} \sigma_{l,l'})\\
    & = -\frac{1}{\sigma_{\beta}^2n}\sum_{\tau=1}^q \sum_{j=1}^\tau \sum_{k,l,k',l'=1}^pc_{\tau,k,l}c_{\tau,k',l'} \sigma_{k,k'} \sigma_{l,l'} \\
    &~~~~+\frac{1}{\sigma_{\beta}^2n}\sum_{\tau=1}^q\sum_{j=\tau+1}^n \sum_{k,l,k',l'=1}^pc_{\tau,k,l}c_{\tau,k',l'} (x_{k,j-\tau} x_{k',j-\tau}x_{l,j}x_{l',j}-\sigma_{k,k'} \sigma_{l,l'})\,.
\end{align*}
By the Cauchy-Schwarz inequality, it holds that
\begin{align*}
    (\tilde{W}_{n2}-1)^2 &\le \frac{2}{\sigma_{\beta}^4n^2}\bigg(\sum_{\tau=1}^q \sum_{j=1}^\tau \sum_{k,l,k',l'=1}^pc_{\tau,k,l}c_{\tau,k',l'} \sigma_{k,k'} \sigma_{l,l'}\bigg)^2\\
    &~~~~+\frac{2}{\sigma_{\beta}^4n^2}\bigg\{\sum_{\tau=1}^q\sum_{j=\tau+1}^n \sum_{k,l,k',l'=1}^pc_{\tau,k,l}c_{\tau,k',l'} (x_{k,j-\tau} x_{k',j-\tau}x_{l,j}x_{l',j}-\sigma_{k,k'} \sigma_{l,l'})\bigg\}^2\\
    &\le\frac{q(q+1)p^4}{\sigma_{\beta}^4n^2}\sum_{\tau=1}^q \sum_{j=1}^\tau \sum_{k,l,k',l'=1}^pc_{\tau,k,l}^2c_{\tau,k',l'}^2 \sigma_{k,k'}^2 \sigma_{l,l'}^2\\
    &~~~~+\frac{2}{\sigma_{\beta}^4}\sum_{\tau=1}^q \sum_{k,l,k',l'=1}^pc_{\tau,k,l}^2c_{\tau,k',l'}^2 \\
    &~~~~~~~~~~~~~~~~~~~~\times \frac{1}{n^2}\sum_{\tau=1}^q \sum_{k,l,k',l'=1}^p\bigg\{\sum_{j=\tau+1}^n(x_{k,j-\tau} x_{k',j-\tau}x_{l,j}x_{l',j}-\sigma_{k,k'} \sigma_{l,l'})\bigg\}^2\\
    & \le K_{2,1}(\beta)n^{-2} +K_{2,2}(\beta)W_{n2}\,,
\end{align*}
where $W_{n2}$ is defined in \eqref{eq:Wn2}, $K_{2,1}(\beta)>0$ and $K_{2,2}(\beta)>0$ are some constants independent of $n$.
By the fact $e_\ell^2=1$ for $\ell\in [n]$ and the Cauchy-Schwarz inequality, we have
\begin{align*}
    \E^e[\{\tilde{W}_{n3}^e(\tau,\tau')\}^2] & = \frac{4}{\sigma_{\beta}^4n^2}\sum_{j=1}^n\bigg(\sum_{k,l,k',l'=1}^pc_{\tau,k,l}c_{\tau',k',l'} x_{k,j-\tau} x_{l,j}x_{k',j-\tau'} x_{l',j}\bigg)^2\\
    & \le \frac{4}{\sigma_{\beta}^4}\max_{1 \le \tau' < \tau \le q}\sum_{k,l,k',l'=1}^pc_{\tau,k,l}^2c_{\tau',k',l'}^2 \cdot \frac{1}{n^2}\sum_{j=1}^n\sum_{k,l,k',l'=1}^p x_{k,j-\tau}^2 x_{l,j}^2x_{k',j-\tau'}^2 x_{l',j}^2 \\
    &=K_3(\beta)W_{n3}(\tau,\tau')\,,
\end{align*}
for any $1 \le \tau' < \tau \le q$, where $W_{n3}(\tau,\tau')$ is defined in \eqref{eq:Wn3}, and $$K_3(\beta) = \frac{4}{\sigma_{\beta}^4}\max_{1 \le \tau' < \tau \le q}\sum_{k,l,k',l'=1}^pc_{\tau,k,l}^2c_{\tau',k',l'}^2\ge0$$ is some constant independent of $n$.
Moreover, by the Cauchy-Schwarz inequality and the fact $e_\ell^2=1$ for $\ell\in [n]$, we have 
\begin{align*}
    \sum_{j=1}^n \E^e\{(Y_j^e)^4\} & = \frac{1}{\sigma_{\beta}^4}\sum_{j=1}^n \E^e\{(B_{j}^e)^4\} \notag\\
    & \le \frac{q^2p^4}{\sigma_{\beta}^4n^2}\sum_{j=1}^n \E^e\bigg\{\sum_{\tau=1}^q \sum_{k,l=1}^p c_{\tau,k,l}^2 e_{j-\tau}^2e_{j}^2x_{k,j-\tau}^2 x_{l,j}^2\bigg\}^2 \notag\\
    & = \frac{q^2p^4}{\sigma_{\beta}^4n^2}\sum_{j=1}^n \bigg\{\sum_{\tau=1}^q \sum_{k,l=1}^p c_{\tau,k,l}^2 x_{k,j-\tau}^2 x_{l,j}^2\bigg\}^2 \\
    & \le \frac{q^3p^6}{\sigma_{\beta}^4n^2}\sum_{j=1}^n \sum_{\tau=1}^q \sum_{k,l=1}^p c_{\tau,k,l}^4 x_{k,j-\tau}^4 x_{l,j}^4 \\
    & \le \frac{q^3p^6}{\sigma_{\beta}^4}\max_{\tau\in [q],~k,l \in [p]}c_{\tau,k,l}^4 \cdot \frac{1}{n^2}\sum_{j=1}^n \sum_{\tau=1}^q \sum_{k,l=1}^p x_{k,j-\tau}^4 x_{l,j}^4\\
    &= K_{4}(\beta)W_{n4}\,,
\end{align*}
where $W_{n4}$ is defined in \eqref{eq:Wn4}, and $K_4(\beta) = q^3p^6\sigma_{\beta}^{-4}\max_{\tau\in [q],~k,l \in [p]}c_{\tau,k,l}^4>0$ is some constant  independent of $n$.
Hence, by \eqref{eq:Ln_pfix}, we have 
\begin{align}\label{eq:Ln2_pfix}
    \tilde{L}_n \le \tilde{K}(\beta) \Big\{n^{-2}+\max_{1 \le \tau' < \tau \le q}W_{n3}(\tau,\tau') + W_{n2} + W_{n4}\Big\}\,,
\end{align}
where $\tilde{K}(\beta)>0$ is some constant independent of $n$.
Hence, \eqref{eq:BEbound_pfix} implies that
\begin{align}\label{eq:BE2_pfix}
\sup_{t \in \mathbb{R}}\bigg|\pr(\tilde{J}_{1} \le t\,|\, \mathcal{X}_n)-\Phi\bigg(\frac{t}{\sigma_{\beta}}\bigg)\bigg| &=\sup_{x\in \mathbb{R}}\bigg|\pr\bigg(\frac{\tilde{J}_{1}}{\sigma_{\beta}} \le x\,\bigg|\, \mathcal{X}_n\bigg)-\Phi(x)\bigg| \\
&\le K^\dagger(\beta) \Big\{n^{-2}+\max_{1 \le \tau' < \tau \le q}W_{n3}(\tau,\tau') + W_{n2} + W_{n4}\Big\}^{1/5}\,,\notag
\end{align}
where $K^\dagger(\beta)\geq\max\{C,1\}\tilde{K}(\beta)^{1/5}$ is a constant independent of $n$. Moreover, recall $\beta^\top\bV_n^e - \tilde{J}_1 = \tilde{J}_2$. By the Markov's inequality and \eqref{eq:J2ttt}, for any $\varepsilon > 0$, we have
\begin{align*}
    \pr(|\tilde{J}_2|>\varepsilon\,|\,\mathcal{X}_n) \le \frac{\E^e(\tilde{J}_2^2)}{\varepsilon^2} \le \frac{K_1(\beta)}{\varepsilon^2}W_{n1}\,.
\end{align*}
Hence, we have
\begin{align*}
    \pr\bigg(\frac{\beta^\top\bV_n^e}{\sigma_{\beta}}  \le x\,\bigg|\, \mathcal{X}_n\bigg) & \le \pr\bigg(\frac{\beta^\top\bV_n^e}{\sigma_{\beta}}  \le x,~|\tilde{J}_2|\le\varepsilon \,\bigg|\, \mathcal{X}_n\bigg) + \pr(|\tilde{J}_2|>\varepsilon\,|\,\mathcal{X}_n)\\
    & \le \pr\bigg(\frac{\tilde{J}_1}{\sigma_{\beta}} \le x+\frac{\varepsilon}{\sigma_{\beta}}\,\bigg|\, \mathcal{X}_n\bigg) + \frac{K_1(\beta)}{\varepsilon^2}W_{n1}
\end{align*}
and
\begin{align*}
    \pr\bigg(\frac{\beta^\top\bV_n^e}{\sigma_{\beta}}  \le x\,\bigg|\, \mathcal{X}_n\bigg) & \ge \pr\bigg(\frac{\beta^\top\bV_n^e}{\sigma_{\beta}}  \le x,~|\tilde{J}_2|\le \varepsilon \,\bigg|\, \mathcal{X}_n\bigg) \\
    & \ge \pr\bigg(\frac{\tilde{J}_1}{\sigma_{\beta}} \le x-\frac{\varepsilon}{\sigma_{\beta}}\,\bigg|\, \mathcal{X}_n\bigg) - \pr(|\tilde{J}_2|>\varepsilon\,|\,\mathcal{X}_n)\\
    & \ge \pr\bigg(\frac{\tilde{J}_1}{\sigma_{\beta}} \le x-\frac{\varepsilon}{\sigma_{\beta}}\,\bigg|\, \mathcal{X}_n\bigg) - \frac{K_1(\beta)}{\varepsilon^2}W_{n1}\,,
\end{align*}
which implies that 
\begin{align*}
    \pr\bigg(\frac{\tilde{J}_1}{\sigma_{\beta}}\le x-\frac{\varepsilon}{\sigma_{\beta}}\,\bigg|\, \mathcal{X}_n\bigg) - \frac{K_1(\beta)}{\varepsilon^2}W_{n1} &\le \pr\bigg(\frac{\beta^\top\bV_n^e}{\sigma_{\beta}}  \le x\,\bigg|\, \mathcal{X}_n\bigg) \\
    &\le \pr\bigg(\frac{\tilde{J}_1}{\sigma_{\beta}} \le x+\frac{\varepsilon}{\sigma_{\beta}}\,\bigg|\, \mathcal{X}_n\bigg) + \frac{K_1(\beta)}{\varepsilon^2}W_{n1}\,.
\end{align*}
By \eqref{eq:BE2_pfix}, we have
\begin{align*}
    &\Phi\bigg(x-\frac{\varepsilon}{\sigma_{\beta}}\bigg)-K^\dagger(\beta) \Big\{n^{-2}+\max_{1 \le \tau' < \tau \le q}W_{n3}(\tau,\tau') + W_{n2} + W_{n4}\Big\}^{1/5}- K_1(\beta)\varepsilon^{-2}W_{n1} \\
    &~~~~~~\le \pr\bigg(\frac{\beta^\top\bV_n^e}{\sigma_{\beta}}  \le x\,\bigg|\, \mathcal{X}_n\bigg) \\
    &~~~~~~\le \Phi\bigg(x+\frac{\varepsilon}{\sigma_{\beta}}\bigg)+K^\dagger(\beta) \Big\{n^{-2}+\max_{1 \le \tau' < \tau \le q}W_{n3}(\tau,\tau') + W_{n2} + W_{n4}\Big\}^{1/5}+ K_1(\beta)\varepsilon^{-2}W_{n1}\,.
\end{align*}
Furthermore, since $\Phi(\cdot)$ is Lipschitz continuous, with $\delta_{\varepsilon}=\varepsilon/(\sigma_{\beta}\sqrt{2\pi})\to0$ as $\varepsilon \to0$, we have $|\Phi(x+\sigma_{\beta}^{-1}\varepsilon)-\Phi(x)|\le \delta_{\varepsilon}$. Hence, we have
\begin{align*}
    &\Phi(x)-\delta_{\varepsilon}-K^\dagger(\beta) \Big\{n^{-2}+\max_{1 \le \tau' < \tau \le q}W_{n3}(\tau,\tau') + W_{n2} + W_{n4}\Big\}^{1/5}- K_1(\beta)\varepsilon^{-2}W_{n1} \\
    &~~~~~~\le \pr\bigg(\frac{\beta^\top\bV_n^e}{\sigma_{\beta}}  \le x\,\bigg|\, \mathcal{X}_n\bigg) \\
    &~~~~~~\le \Phi(x)+\delta_{\varepsilon}+K^\dagger(\beta) \Big\{n^{-2}+\max_{1 \le \tau' < \tau \le q}W_{n3}(\tau,\tau') + W_{n2} + W_{n4}\Big\}^{1/5}+ K_1(\beta)\varepsilon^{-2}W_{n1}\,,
\end{align*}
which implies that \eqref{eq:projection_fixedp} holds
for each fixed  $\beta$ with $\sigma_{\beta}^2>0$. $\hfill\Box$

\subsubsection{Proof of \eqref{eq:cond_Tne}}\label{subsec:Tne}

Recall
\begin{align*}
    T_n^e &= \frac{1}{\hat\sigma_{n1}}
    \sum_{\tau=1}^q
    \bigg[\frac np\tr\{\bS_\tau^e(\bS_\tau^e)^\top\}
    -M_{n,\tau}\bigg]\,,\\
    T_\infty
    &=\frac{1}{\sqrt2\,\tr(\bSigma^2)}
    \sum_{\tau=1}^q
    \left\{\tr(\bSigma N_\tau\bSigma N_\tau^\top)
    -{\tr}^2(\bSigma)\right\}\,,
\end{align*}
as given in \eqref{eq:Tne} and \eqref{eq:Tinf_fixedp}, respectively.
Here $N_\tau=(\varepsilon_{i,j,\tau})\in\mathbb R^{p\times p}$, where $\varepsilon_{i,j,\tau}$ are i.i.d.\ $\mathcal{N}(0,1)$ variables for any $i,j\in [p]$ and $\tau\in[q]$. To prove \eqref{eq:cond_Tne}, it suffices to show
\begin{align}\label{eq:all_cf}
    \E\{\exp(\iota tT_{n_{k_\ell}}^e) \,|\,
    \mathcal X_{n_{k_\ell}}=\mathcal X_{n_{k_\ell}}(\omega)\} \to \E\{\exp(\iota tT_\infty)\} ~~\mbox{as}~\ell \to \infty
\end{align}
for every fixed $t \in \mathbb{R}$, where $\iota = \sqrt{-1}$. For simplicity of notation, write $$Y_{n_{k_\ell}} = \frac{n_{k_{\ell}}}{p}\sum_{\tau=1}^q\tr\{\bS_{\tau,n_{k_{\ell}}}^e(\bS_{\tau,n_{k_{\ell}}}^e)^\top\}\,,~~~ Y = \frac{1}{p}\sum_{\tau=1}^q\tr(\bSigma N_\tau\bSigma N_\tau^{\top})\,,$$ 
$c = \sqrt{2}p^{-1}{\tr}(\bSigma^2)$ and $m = qp^{-1}{\tr}^2(\bSigma)$. Here $\{\bS_{\tau,n_{k_{\ell}}}^e\}$ denotes the corresponding subsequence of $\{\bS_{\tau,n}^e\}$, where $\bS_{\tau}^e$ is  written as $\bS_{\tau,n}^e$ to make the dependence of $\bS_\tau^e$ on $n$ explicit.  Hence, 
\begin{align*}
    &\E\{\exp(\iota tT_{n_{k_\ell}}^e) \,|\,
    \mathcal X_{n_{k_\ell}}=\mathcal X_{n_{k_\ell}}(\omega)\} \\
    &~~~~~~~~~~~~~~~~~~= \exp\bigg\{-\frac{\iota t\sum_{\tau=1}^q M_{n_{k_\ell},\tau}(\omega) }{\hat\sigma_{n_{k_\ell}1}(\omega)}\bigg\}\E\bigg\{\exp\bigg(\frac{\iota tY_{n_{k_\ell}}}{\hat\sigma_{n_{k_\ell}1}}\bigg) \,\bigg|\,
    \mathcal X_{n_{k_\ell}}=\mathcal X_{n_{k_\ell}}(\omega)\bigg\}\,, \\
    &\E\{\exp(\iota tT_\infty)\} = \exp\bigg(-\frac{\iota tm }{c}\bigg)\E\bigg\{\exp\bigg(\frac{\iota tY}{c}\bigg)\bigg\}\,.
\end{align*}
Recall $\max_{\tau\in [q]}|M_{n_{k_\ell},\tau}(\omega)-q^{-1}m|
        +|\hat\sigma_{n_{k_\ell}1}(\omega)
        -c|\to0$ for any fixed $\omega \in \Omega_0$, which implies 
\begin{align*}
    \exp\bigg\{-\frac{\iota t\sum_{\tau=1}^q M_{n_{k_\ell},\tau}(\omega) }{\hat\sigma_{n_{k_\ell}1}(\omega)}\bigg\} \to \exp\bigg(-\frac{\iota tm }{c}\bigg) ~~\mbox{as}~\ell \to \infty
\end{align*}
for any fixed $t \in \mathbb{R}$ and $\omega \in \Omega_0$.
Hence, to show \eqref{eq:all_cf}, it suffices to show
\begin{align}\label{eq:all_cf2}
    \E\bigg\{\exp\bigg(\frac{\iota tY_{n_{k_\ell}}}{\hat\sigma_{n_{k_\ell}1}}\bigg) \,\bigg|\,
    \mathcal X_{n_{k_\ell}}=\mathcal X_{n_{k_\ell}}(\omega)\bigg\} \to \E\bigg\{\exp\bigg(\frac{\iota tY}{c}\bigg)\bigg\} ~~\mbox{as}~\ell \to \infty
\end{align}
for any fixed $\omega \in \Omega_0$ and any fixed $t \in \mathbb{R}$.

Analogous to \eqref{eq:SDist}, we can reformulate \eqref{eq:joint_cond} as follows:
\begin{align*}
     \bV_{n_{k_\ell}}^e\,|\,
    \mathcal X_{n_{k_\ell}}=\mathcal X_{n_{k_\ell}}(\omega) \xrightarrow{\mathrm{d}} \begin{pmatrix}
         \bSigma^{1/2}\otimes \bSigma^{1/2} &  &  \\
         & \ddots & \\
         &  &  \bSigma^{1/2}\otimes \bSigma^{1/2}
         \end{pmatrix} \begin{pmatrix}
         \mathrm{vec}(N_1) \\
         \vdots \\
         \mathrm{ vec}( N_q)\end{pmatrix} 
\end{align*}
as $\ell \to \infty$ for any fixed $\omega \in \Omega_0$, which implies
\begin{align*}
    \frac{n_{k_\ell}}{p}\begin{pmatrix}
    \operatorname{tr}\{\bS_{1,n_{k_\ell}}^e(\bS_{1,n_{k_\ell}}^e)^\top\} \\
    \vdots \\
    \operatorname{tr}\{\bS_{q,n_{k_\ell}}^e(\bS_{q,n_{k_\ell}}^e)^\top\} 
    \end{pmatrix} \Bigg|\, \mathcal X_{n_{k_\ell}}=\mathcal X_{n_{k_\ell}}(\omega)  \xrightarrow{\mathrm{d}} \frac{1}{p} \begin{pmatrix}
    \operatorname{tr}(\bSigma N_1\bSigma  N_1^\top ) \notag\\
    \vdots \\
    \operatorname{tr}(\bSigma N_q\bSigma  N_q^\top )
    \end{pmatrix}\,.
\end{align*}
Then we have
\begin{align}\label{eq:Yconverge}
    Y_{n_{k_\ell}} \,|\, \mathcal X_{n_{k_\ell}}=\mathcal X_{n_{k_\ell}}(\omega) \xrightarrow{\mathrm{d}} Y ~~\mbox{as}~\ell \to \infty
\end{align}
for any fixed $\omega \in \Omega_0$,
which implies
\begin{align}\label{eq:cf}
    &\E\bigg\{\exp\bigg(\frac{\iota tY_{n_{k_\ell}}}{c}\bigg)\,\bigg|\,\mathcal X_{n_{k_\ell}}=\mathcal X_{n_{k_\ell}}(\omega)\bigg\}  \to \E\bigg\{\exp\bigg(\frac{\iota tY}{c}\bigg)\bigg\}~~\mbox{as}~\ell \to \infty
\end{align}
for any fixed $\omega \in \Omega_0$ and any $t \in \mathbb{R}$. 
Notice that 
\begin{align}\label{eq:cf2}
    &\bigg| \E\bigg\{\exp\bigg(\frac{\iota tY_{n_{k_\ell}}}{\hat\sigma_{n_{k_\ell}1}}\bigg) \,\bigg|\,
    \mathcal X_{n_{k_\ell}}=\mathcal X_{n_{k_\ell}}(\omega)\bigg\} -  \E\bigg\{\exp\bigg(\frac{\iota tY}{c}\bigg)\bigg\}  \bigg| \notag\\
    &~~~~~\le \bigg| \E\bigg\{\exp\bigg(\frac{\iota tY_{n_{k_\ell}}}{\hat\sigma_{n_{k_\ell}1}}\bigg) \,\bigg|\,
    \mathcal X_{n_{k_\ell}}=\mathcal X_{n_{k_\ell}}(\omega)\bigg\} - \E\bigg\{\exp\bigg(\frac{\iota tY_{n_{k_\ell}}}{c}\bigg) \,\bigg|\,
    \mathcal X_{n_{k_\ell}}=\mathcal X_{n_{k_\ell}}(\omega)\bigg\}\bigg| \notag\\
    &~~~~~~~~~~+\bigg| \E\bigg\{\exp\bigg(\frac{\iota tY_{n_{k_\ell}}}{c}\bigg) \,\bigg|\,
    \mathcal X_{n_{k_\ell}}=\mathcal X_{n_{k_\ell}}(\omega)\bigg\} -\E\bigg\{\exp\bigg(\frac{\iota tY}{c}\bigg)\bigg\} \bigg|\notag\\
    &~~~~~=U_{1,n_{k_\ell}}(\omega,t)+U_{2,n_{k_\ell}}(\omega,t)\,.
\end{align}
By \eqref{eq:cf}, 
$U_{2,n_{k_\ell}}(\omega,t)\to0$ for every fixed
$\omega\in\Omega_0$ and $t\in\mathbb{R}$. In the sequel, we only need to show $U_{1,n_{k_\ell}}(\omega,t)\to0$ for every fixed
$\omega\in\Omega_0$ and $t\in\mathbb{R}$.

For $U_{1,n_{k_\ell}}(\omega,t)$, we have
\begin{align*}
    U_{1,n_{k_\ell}}(\omega,t) & \le  \E\bigg\{\bigg|\exp\bigg(\frac{\iota tY_{n_{k_\ell}}}{\hat\sigma_{n_{k_\ell}1}}\bigg) - \exp\bigg(\frac{\iota tY_{n_{k_\ell}}}{c}\bigg)\bigg| \,\bigg|\,
    \mathcal X_{n_{k_\ell}}=\mathcal X_{n_{k_\ell}}(\omega)\bigg\} \notag\\
    & = \E\bigg\{\bigg|\exp\bigg(\frac{\iota tY_{n_{k_\ell}}}{\hat\sigma_{n_{k_\ell}1}}\bigg) - \exp\bigg(\frac{\iota tY_{n_{k_\ell}}}{c}\bigg)\bigg|I(|Y_{n_{k_\ell}}|\le M) \,\bigg|\,
    \mathcal X_{n_{k_\ell}}=\mathcal X_{n_{k_\ell}}(\omega)\bigg\} \notag\\
    &~~~~~~~+\E\bigg\{\bigg|\exp\bigg(\frac{\iota tY_{n_{k_\ell}}}{\hat\sigma_{n_{k_\ell}1}}\bigg) - \exp\bigg(\frac{\iota tY_{n_{k_\ell}}}{c}\bigg)\bigg|I(|Y_{n_{k_\ell}}|> M) \,\bigg|\,
    \mathcal X_{n_{k_\ell}}=\mathcal X_{n_{k_\ell}}(\omega)\bigg\}\notag\\
    & \le M|t|\bigg|\frac{1}{\hat\sigma_{n_{k_\ell}1}(\omega)}-\frac{1}{c}\bigg| + 2\pr\big\{|Y_{n_{k_\ell}}|> M \,|\, \mathcal X_{n_{k_\ell}}=\mathcal X_{n_{k_\ell}}(\omega)\big\}
\end{align*}
for any $M>0$, where the last inequality is based on the fact $|e^{\iota u}-e^{\iota v}|\leq \min\{2,|u-v|\}$. Hence, for any $\omega \in \Omega_0$ and $t \in \mathbb{R}$, it holds that
\begin{align}\label{eq:U1}
\limsup_{\ell\to\infty}U_{1,n_{k_\ell}}(\omega,t) \le 2 \limsup_{\ell\to\infty} \pr\big\{|Y_{n_{k_\ell}}|> M \,|\, \mathcal X_{n_{k_\ell}}=\mathcal X_{n_{k_\ell}}(\omega)\big\}
\end{align}
for any $M>0$. 
Notice that 
\begin{align*}
    Y \stackrel{\mathrm{d}}{=} \frac{1}{p}\sum_{\tau=1}^q \sum_{i,j=1}^p \lambda_i\lambda_j \varepsilon_{i,j,\tau}^2\,,
\end{align*}
where $\lambda_1\geq\cdots\geq\lambda_p$ are the eigenvalues of $\bSigma$. This implies 
$$\E(|Y|) = \frac{1}{p}\sum_{\tau=1}^q \sum_{i,j=1}^p \lambda_i\lambda_j \E(\varepsilon_{i,j,\tau}^2)=\frac{q}{p}{\tr}^2(\bSigma) < \infty\,.$$ 
Hence, for any $\epsilon>0$, by Markov’s inequality, we have $\pr(|Y|\ge M)<\epsilon/2$ by choosing $M>2\epsilon^{-1}\E(|Y|)$. By \eqref{eq:Yconverge} and the Portmanteau theorem, we have
\begin{align}\label{eq:tail}
    \limsup_{\ell\to\infty}\pr\big\{|Y_{n_{k_\ell}}|>M \,|\,\mathcal X_{n_{k_\ell}}=\mathcal X_{n_{k_\ell}}(\omega)\big\} & \le \limsup_{\ell\to\infty}\pr\big\{|Y_{n_{k_\ell}}|\ge M \,|\, \mathcal X_{n_{k_\ell}}=\mathcal X_{n_{k_\ell}}(\omega)\big\}\notag\\
    & \le \pr(|Y|\ge M)<\frac{\epsilon}{2}\,.
\end{align}
Thus, by \eqref{eq:U1} and \eqref{eq:tail}, we have $\limsup_{\ell\to\infty}U_{1,n_{k_\ell}}(\omega,t) < \epsilon$.
As $\epsilon>0$ is arbitrary and $U_{1,n_{k_\ell}}(\omega,t)\ge0$,
we obtain $U_{1,n_{k_\ell}}(\omega,t)\to0$ for every fixed
$\omega\in\Omega_0$ and $t\in\mathbb{R}$.  Hence, by \eqref{eq:cf2}, we have \eqref{eq:all_cf2} holds. $\hfill\Box$

 \section{Proof of Theorem \ref{boot}}
 
Recall $\mathcal{X}_n = \{\x_1,\ldots,\x_n\}$, $G_{n} = \sum_{\tau=1}^q G_{n,\tau,1}$ and
\begin{align*}
    G_{n}^e = \sum_{\tau=1}^q G_{n,\tau,1}^e \quad \mbox{with} \quad G_{n,\tau,1}^e= \frac 1{pn}\sum_{t\neq s}\y_t^\top \y_s\y_{t+\tau}^\top \y_{s+\tau}\,.
\end{align*} 
For any fixed $\alpha\in (0,1)$, write $\tilde{{\rm cv}}_{G, \alpha} = \inf\{t \in \mathbb{R}: \pr(G_{n}^{e} \le t\,|\, \mathcal{X}_n) \ge 1-\alpha\}$. By Chebyshev's inequality, we have
\begin{align*}
	\pr\bigg\{\bigg|\frac{G_{n}^{e}}{\sqrt{\Var(G_{n}^{e}\,|\,\mathcal{X}_n)}}\bigg|>\frac1{\sqrt{\alpha}}
	\,\bigg|\,\mathcal{X}_n
	\bigg\}\leq \alpha\,,
\end{align*}
which implies 
\begin{align*}
    \tilde{{\rm cv}}_{G, \alpha} \le K_\alpha \sqrt{\Var(G_{n}^{e}\,|\,\mathcal{X}_n)}
\end{align*}
for some constant $K_\alpha >0$. 
Recall $\tilde{{\rm cv}}_{\alpha} = \inf\{t \in \mathbb{R}: \pr(T_{n}^{e} \le t\,|\, \mathcal{X}_n) \ge 1-\alpha\}$. 
As shown in Section \ref{subsec:size_pinf}, $\pr(
T_n > \tilde{{\rm cv}}_{\alpha} ) = \pr(
G_{n}> \tilde{{\rm cv}}_{G,\alpha})$ under both the null and alternative hypotheses. Due to $\mu_n = \E(G_{n})$, then
\begin{align*}	
\pr_{H_1}(
T_n > \tilde{{\rm cv}}_{\alpha} ) 
\geq &~\pr_{H_1}\big\{
G_{n}> K_\alpha \sqrt{\Var(G_{n}^{e}\,|\,\mathcal{X}_n)}
\big\}\\
\geq &~\pr_{H_1}\bigg\{\bigg|\frac{G_{n}-\mu_n}{\sqrt{\Var(G_{n})}}\bigg|
<\frac{\mu_n-K_\alpha \sqrt{\Var(G_{n}^{e}\,|\,\mathcal{X}_n)}}{\sqrt{\Var(G_{n})}}\bigg\}\,.
\end{align*}
Let $\{\eta_n\}$ be a sequence of positive numbers decreasing to zero. By Markov's inequality, it holds that
\begin{align*}
\pr_{H_1}\big[\Var(G_{n}^{e}\,|\,\mathcal{X}_n)>\eta_n^{-2}\E\{\Var(G_{n}^{e}\,|\,\mathcal{X}_n)\}	\big]\leq \eta_n^2\to0\,,
\end{align*}
which implies
\begin{align*}
\pr_{H_1}( T_n > \tilde{{\rm cv}}_{\alpha} ) & \geq \pr_{H_1}\bigg[\bigg|\frac{G_{n}-\mu_n}{\sqrt{\Var(G_{n})}}\bigg|
<\frac{\mu_n-K_\alpha \sqrt{\Var(G_{n}^{e}\,|\,\mathcal{X}_n)}}{\sqrt{\Var(G_{n})}},\,\\
&~~~~~~~~~~~~~~~~~~~~~~~\Var(G_{n}^{e}\,|\,\mathcal{X}_n)\le \eta_n^{-2}\E\{\Var(G_{n}^{e}\,|\,\mathcal{X}_n)\}\bigg]\\
&\geq \mathbb \pr_{H_1}\bigg\{\bigg|\frac{G_{n}-\mu_n}{\sqrt{\Var(G_{n})}}\bigg|
<\frac{\mu_n-\eta_n^{-1}K_\alpha \sqrt{\E\{\Var(G_{n}^{e}\,|\,\mathcal{X}_n)\}}}{\sqrt{\Var(G_{n})}}\bigg\}-o(1)\,.
\end{align*}
Therefore, to prove Theorem \ref{boot}, we only need to show
\begin{align*}
\frac{\mu_n-K_\alpha \eta_n^{-1}\sqrt{\E\{\Var(G_{n}^{e}\,|\,\mathcal{X}_n)\}}}{\sqrt{\Var(G_{n})}}>
\frac1{\eta_n}
\end{align*}
for all large $n$. It suffices  to show
\begin{align}\label{sign-1}
\frac{\mu_n}{\sqrt{\E\{\Var(G_{n}^{e}\,|\,\mathcal{X}_n)\}}+\sqrt{\Var(G_{n})}}\to+\infty\,,
\end{align}
since we can select $\eta_n \to 0$ arbitrarily slow.

We then derive $\E\{\Var(G_{n}^{e}\,|\,\mathcal{X}_n)\}$ and $\Var(G_{n})$. Recall $G_{n} = \sum_{\tau=1}^q G_{n,\tau,1}$ with
\begin{align*}
    G_{n,\tau,1} = \frac {2}{np}\sum_{t < s}\x_t^\top \x_s\x_{t+\tau}^\top \x_{s+\tau}\,,
\end{align*}
where $\sum_{t<s}$ denotes the summation over $t$ and $s$ taking values in $[n]$ with the restriction $t<s$. 
Notice that
\begin{align*}
G_{n,\tau,1} 
=&~ \frac2{np}\sum_{t< s\leq n-\tau}\x_t^\top \x_s\x_{t+\tau}^\top \x_{s+\tau}+ \frac2{np}\sum_{t< s,\, s> n-\tau}\x_t^\top \x_s\x_{t+\tau}^\top \x_{s+\tau}\\
=&~ G_{n,\tau,10}+G_{n,\tau,11}\,,\\
G_{n,\tau,1}^e
=&~ \frac2{np}\sum_{t< s\leq n-\tau}e_te_se_{t+\tau}e_{s+\tau}\x_t^\top \x_s\x_{t+\tau}^\top \x_{s+\tau}
+ \frac2{np}\sum_{t< s,\, s>n-\tau}e_te_se_{t+\tau}e_{s+\tau}\x_t^\top \x_s\x_{t+\tau}^\top \x_{s+\tau}\\
= &~ G_{n,\tau,10}^e+G_{n,\tau,11}^e\,.
\end{align*}

As we will show in Sections \ref{subsubsec:E2} and \ref{subsubsec:E3},
\begin{align}
&\max\big\{\Var(G_{n,\tau,10}), \Var(G_{n,\tau,11})\big\}\leq K\Big(\frac pn+\frac np\Big),\label{varg}\\
&\max\big\{\E\{\Var(G_{n,\tau,10}^{e}\,|\,\mathcal{X}_n)\},~ \E\{\Var(G_{n,\tau,11}^{e}\,|\,\mathcal{X}_n)\}\big\}
\leq K\Big(\frac {p^2}{n}+ 1\Big)\,.\label{ev-hntau}
\end{align}
Notice that $q$ is a fixed constant. Then, by the Cauchy-Schwarz inequality, we have
\begin{align*}
    \Var(G_{n}) & \le q \sum_{\tau=1}^q \Var( G_{n,\tau,1}) \le 2q \sum_{\tau=1}^q \{\Var( G_{n,\tau,10}) + \Var( G_{n,\tau,11})\} \\
    &\le K\Big(\frac pn+\frac np\Big)\,,\\
    \E\{\Var(G_{n}^{e}\,|\,\mathcal{X}_n)\} & \le q \sum_{\tau=1}^q \E\{\Var( G_{n,\tau,1}^e\,|\,\mathcal{X}_n)\}\\
    &\le 2q \sum_{\tau=1}^q \big[\E\{\Var( G_{n,\tau,10}^e\,|\,\mathcal{X}_n)\} + \E\{\Var( G_{n,\tau,11}^e\,|\,\mathcal{X}_n)\}\big] \\
    &\le K\Big(\frac {p^2}{n}+ 1\Big)\,,
\end{align*}
which implies $\sqrt{\E\{\Var(G_{n}^{e}\,|\,\mathcal{X}_n)\}}+\sqrt{\Var(G_{n})} \le K \max\{p/\sqrt{n}, \sqrt{n/p}\}$.
Hence, \eqref{sign-1} holds if $\mu_n/\max\{p/\sqrt{n}, \sqrt{n/p}\}\to\infty$.  We complete the proof of Theorem \ref{boot}. $\hfill\Box$

\subsection{Proof of \eqref{varg}}\label{subsubsec:E2}
 Recall $\x_{t} = \sum_{j=0}^{\infty}\A_j\z_{t-j}$ with $\A_0 = I_p$ and $\A_j \in \mathbb{R}^{p \times p}$. 
Then
\begin{align*}
G_{n,\tau,10} 
=&~ \frac2{np}\sum_{t< s\leq n-\tau}
\sum_{j_1, j_2=0}^\infty \z_{t-j_1}^\top \A_{j_1}^\top \A_{j_2}\z_{s-j_2}
\sum_{j_3, j_4=0}^\infty \z_{t+\tau-j_3}^\top \A_{j_3}^\top \A_{j_4}\z_{s+\tau-j_4} \\
=&~\frac2{np}\sum_{k=1}^{n-\tau-1}\sum_{t=1}^{n-\tau-k}
\sum_{j_1, j_2=0}^\infty \z_{t-j_1}^\top \A_{j_1}^\top \A_{j_2}\z_{t+k-j_2}
\sum_{j_3, j_4=0}^\infty \z_{t+\tau-j_3}^\top \A_{j_3}^\top \A_{j_4}\z_{t+\tau+k-j_4}\\
= &~
 \frac2{np}\sum_{k=1}^{n-\tau-1}\sum_{t=1}^{n-\tau-k}
\big[\{a(t,k)+b(t,k)\}\{a(t+\tau,k)+b(t+\tau,k)\}\big]\,,
\end{align*}
where
\begin{align}\label{eq:abtk}
&~~~~~~~~~~~a(t,k)=\sum_{j=0}^\infty \z_{t-j}^\top \A_{j}^\top \A_{j+k}\z_{t-j}\,,\quad
b(t,k)=\sum_{j_1, j_2=0\atop j_2\neq j_1+k}^\infty \z_{t-j_1}^\top \A_{j_1}^\top \A_{j_2}\z_{t+k-j_2}\,.	
\end{align}
Hence, we have
\begin{align*}
G_{n,\tau,10} = g_{n,\tau,11}+g_{n,\tau,12}+g_{n,\tau,21}+g_{n,\tau,22}\,,	
\end{align*}
where
\begin{align*}
g_{n,\tau,11}&=\frac2{np}\sum_{k=1}^{n-\tau-1}\sum_{t=1}^{n-\tau-k}
a(t,k)a(t+\tau,k)\,,\quad
g_{n,\tau,12}=\frac2{np}\sum_{k=1}^{n-\tau-1}\sum_{t=1}^{n-\tau-k}
a(t,k)b(t+\tau,k)\,,\\
g_{n,\tau,21}&=\frac2{np}\sum_{k=1}^{n-\tau-1}\sum_{t=1}^{n-\tau-k}
b(t,k)a(t+\tau,k)\,,\quad
g_{n,\tau,22}=\frac2{np}\sum_{k=1}^{n-\tau-1}\sum_{t=1}^{n-\tau-k}
b(t,k)b(t+\tau,k)\,.
\end{align*}
Then ${\Var}(G_{n,\tau,10}) \le K\{ {\Var}(g_{n,\tau,11})+{\Var}(g_{n,\tau,12})+{\Var}(g_{n,\tau,21})+{\Var}(g_{n,\tau,22})\}$.
In order to bound ${\Var}(G_{n,\tau,10})$, 
it suffices to derive bounds for ${\Var}(g_{n,\tau,ij})$ for $i,j = 1,2$.  

For $g_{n,\tau,11}$, we have
\begin{align}
\Var(g_{n,\tau,11})
&=\frac4{n^2p^2}\Var\bigg(\sum_{k=1}^{n-\tau-1}\sum_{t=1}^{n-\tau-k}
\sum_{j_1,j_2=0}^\infty 
\z_{t-j_1}^\top \A_{j_1}^\top \A_{j_1+k}\z_{t-j_1} \z_{t+\tau-j_2}^\top\A_{j_2}^\top\A_{j_2+k}\z_{t+\tau-j_2}	
\bigg)\nonumber\\
&=\frac4{n^2p^2}\sum_{k, k'=1}^{n-\tau-1}\sum_{t \in [n-\tau-k] \atop t' \in [n-\tau-k']}
\sum_{j_1,j_2, j_1', j_2'=0}^{\infty} 
{\rm cov}\big(
\z_{t-j_1}^\top \A_{j_1}^\top \A_{j_1+k}\z_{t-j_1} 
\z_{t+\tau-j_2}^\top\A_{j_2}^\top\A_{j_2+k}\z_{t+\tau-j_2},\nonumber\\
&\qquad~~~~~~~~~~~~~~~~~~~~~
\z_{t'-j'_1}^\top \A_{j'_1}^\top \A_{j'_1+k'}\z_{t'-j'_1} \z_{t'+\tau-j'_2}^\top\A_{j'_2}^\top\A_{j'_2+k'}\z_{t'+\tau-j'_2}
 	\big)\,.\label{cov11}
\end{align}
To bound $\Var(g_{n,\tau,11})$, we first present some auxiliary results. More specifically, for any $p\times p$ matrices $\A$ and $\B$, due to $\|\bSigma\|_2 \le K$ as assumed in Assumption \ref{as:A4}, by Lemma \ref{ineq-main}, we have
\begin{align}\label{eq:e4}
    {\E}[\{\z_i^\top\A\z_i-\tr(\A\bSigma)\}^4] & \le K{\tr}^2(\bSigma^{1/2}\A\bSigma{\A}^{\top}\bSigma^{1/2}) + K{\tr}\{(\bSigma^{1/2}\A\bSigma\A^{\top}\bSigma^{1/2})^2\} \notag \\
    & \le K(p^2\|\A\|_2^4+p\|\A\|_2^4) \le Kp^2\|\A\|_2^4\,,
\end{align}
and 
\begin{align}\label{eq:e2}
    {\Var}(\z_i^\top\A\z_i) &= {\E}[\{\z_i^\top\A\z_i-\tr(\A\bSigma)\}^2] \le K{\tr}(\bSigma^{1/2}\A\bSigma\A^{\top}\bSigma^{1/2}) \le Kp\|\A\|_2^2\,,
\end{align}
which implies
\begin{align}\label{eq:var2}
\Var(\z_i^\top\A\z_i \z_j^\top\B\z_j)
\leq &~K\big(
\Var[\{\z_i^\top\A\z_i-\tr(\A\bSigma)\}\{\z_j^\top\B\z_j-\tr(\B\bSigma)\}]\notag\\
&~~~~~~~~+{\tr}^2(\A\bSigma)\Var(\z_j^\top\B\z_j)+{\tr}^2(\B\bSigma)\Var(\z_i^\top\A\z_i)
\big)\notag\\
\leq &~K\big(
{\E}^{1/2}[\{\z_i^\top\A\z_i-\tr(\A\bSigma)\}^4] {\E}^{1/2}[\{\z_j^\top\B\z_j-\tr(\B\bSigma)\}^4]\notag\\
&~~~~~~~~+{\tr}^2(\A\bSigma)\Var(\z_j^\top\B\z_j)+{\tr}^2(\B\bSigma)\Var(\z_i^\top\A\z_i)
\big)\notag\\
\leq &~Kp^3\|\A\|_2^2\|\B\|_2^2\,.
\end{align}
Hence, for any $p\times p$ matrices $\A,\B, \C, \D$, by the Cauchy-Schwarz inequality, we have
\begin{align}\label{var4-1}
|\Cov(\z_i^\top\A\z_i \z_j^\top\B\z_j, \z_k^\top\C\z_k \z_\ell^\top\D\z_\ell)| \leq &~
\Var^{1/2}(\z_i^\top\A\z_i \z_j^\top\B\z_j) \Var^{1/2}(\z_k^\top\C\z_k \z_\ell^\top\D\z_\ell) \nonumber\\
\leq &~ Kp^3\|\A\|_2\|\B\|_2\|\C\|_2\|\D\|_2
\end{align}
for any $\{i, j\}\cap \{k,l\}\neq \emptyset$, and $\Cov(\z_i^\top\A\z_i \z_j^\top\B\z_j, \z_k^\top\C\z_k \z_\ell^\top\D\z_\ell) = 0$ for any $\{i, j\}\cap \{k,l\} = \emptyset$.
Therefore, by Assumption \ref{as:A5}, we have
\begin{align}\label{eq:gn11}
\Var(g_{n,\tau,11})
\leq &~\frac {Kp}{n}\sum_{k, k'=1}^{n-\tau-1}
\sum_{j_1,j_2, j_1',  j_2'=0}^{\infty} 
\|\A_{j_1}^\top \A_{j_1+k}\|_2 \|\A_{j_2}^\top \A_{j_2+k}\|_2
\|\A_{j'_1}^\top \A_{j'_1+k'}\|_2  \|\A_{j'_2}^\top\A_{j'_2+k'}\|_2	\notag\\
\le &~ \frac {Kp}{n}\bigg(\sum_{k=1}^{n-\tau-1}
\sum_{j_1,j_2=0}^{\infty} 
\|\A_{j_1}\|_2 \|\A_{j_1+k}\|_2 \|\A_{j_2}\|_2 \|\A_{j_2+k}\|_2\bigg)^2
 \notag\\
 \leq &~ \frac {Kp}{n}\bigg\{\sum_{k=1}^{n-\tau-1}\bigg(
\sum_{j=0}^{\infty} 
\|\A_{j}\|_2 \|\A_{j+k}\|_2 \bigg)^2\bigg\}^2 \notag\\
\leq &~\frac {Kp}{n}
\bigg\{\sum_{k=1}^{\infty}\bigg(\sum_{j=0}^\infty 
\|\A_{j}\|_2^2\bigg)\bigg(\sum_{j=0}^\infty \|\A_{j+k}\|_2^2\bigg) \bigg\}^2\notag\\
\leq &~\frac {Kp}{n}
\bigg(\sum_{j=0}^\infty 
\|\A_{j}\|_2^2\cdot\sum_{j=0}^\infty j\|\A_{j}\|_2^2 \bigg)^2\leq \frac {Kp}{n}\,,
\end{align}
where the last inequality is based on Assumption \ref{as:A5}.

For $g_{n,\tau,12}$, due to $\E(g_{n,\tau,12})=0$, we have
\begin{align}
&\Var(g_{n,\tau,12}) \nonumber\\
&~~~=\frac 4{n^2p^2}\E\bigg\{\bigg(\sum_{k=1}^{n-\tau-1}\sum_{t=1}^{n-\tau-k}
\sum_{j=0}^\infty 
\sum_{j_3, j_4=0\atop j_4\neq j_3+k}^\infty 
\z_{t-j}^\top \A_{j}^\top \A_{j+k}\z_{t-j}
\z_{t+\tau-j_3}^\top \A_{j_3}^\top \A_{j_4}\z_{t+\tau+k-j_4}	\bigg)^2\bigg\}\nonumber\\
&~~~=\frac 4{n^2p^2}\sum_{k, k'=1}^{n-\tau-1}\sum_{t \in [n-\tau-k] \atop t' \in [n-\tau-k']} 
\sum_{j, j_3, j_4=0 \atop j_4\neq j_3+k}^{\infty} 
\sum_{j', j'_3, j'_4=0 \atop j'_4\neq j'_3+k'}^{\infty}
\E(\z_{t-j}^\top \A_{j}^\top \A_{j+k}\z_{t-j}
\z_{t+\tau-j_3}^\top \A_{j_3}^\top \A_{j_4}\z_{t+\tau+k-j_4}	\nonumber\\
&~~~~~~~~~~~~~~~~~~~~~~~\times
\z_{t'-j'}^\top \A_{j'}^\top \A_{j'+k'}\z_{t'-j'}
\z_{t'+\tau-j'_3}^\top \A_{j'_3}^\top \A_{j'_4}\z_{t'+\tau+k'-j'_4}	) \label{var4-2} \\
&~~~\le \frac K{n^2p^2}\cdot n \max_{t,t' \in [n-\tau-1]}\sum_{k, k'=1}^{n-\tau-1}
\sum_{j, j_3, j_4=0 \atop j_4\neq j_3+k}^{\infty} 
\sum_{j', j'_3, j'_4=0 \atop j'_4\neq j'_3+k'}^{\infty}
\E^{1/2}\{(\z_{t-j}^\top \A_{j}^\top \A_{j+k}\z_{t-j}
\z_{t+\tau-j_3}^\top \A_{j_3}^\top \A_{j_4}\z_{t+\tau+k-j_4})^2\} \nonumber\\
&~~~~~~~~~~~~~~~~~~~~~~~\times
\E^{1/2}\{(\z_{t'-j'}^\top \A_{j'}^\top \A_{j'+k'}\z_{t'-j'}
\z_{t'+\tau-j'_3}^\top \A_{j'_3}^\top \A_{j'_4}\z_{t'+\tau+k'-j'_4}	)^2\}\,, \nonumber
\end{align}
where the last inequality is based on the Cauchy-Schwarz inequality and the fact that the expectation in \eqref{var4-2} is zero if $\{t-j, t+\tau-j_3, t+\tau+k-j_4\}\cap \{t'-j', t'+\tau-j'_3, t'+\tau+k'-j'_4\}= \emptyset$ given $\{\tau, k, k', j, j_3, j_4, j', j_3', j_4'\}$.
Since $t+\tau+k-j_4 \neq t+\tau-j_3$, we can claim that either $t-j \neq t+\tau-j_3$ or $t-j \neq t+\tau+k-j_4$ must be satisfied. If $t-j \neq t+\tau+k-j_4$, by the Cauchy-Schwarz inequality,
we have
\begin{align}\label{three-diff}
&\E\{(\z_{t-j}^\top \A_{j}^\top \A_{j+k}\z_{t-j}
\z_{t+\tau-j_3}^\top \A_{j_3}^\top \A_{j_4}\z_{t+\tau+k-j_4})^2\}\nonumber\\	
&~~~~~=\E\{(\z_{t-j}^\top \A_{j}^\top \A_{j+k}\z_{t-j})^2
(\z_{t+\tau-j_3}^\top \A_{j_3}^\top \A_{j_4}\bSigma \A_{j_4}^\top \A_{j_3}\z_{t+\tau-j_3})\} \nonumber\\
&~~~~~\leq \E^{1/2}\{(\z_{t-j}^\top \A_{j}^\top \A_{j+k}\z_{t-j})^4\}\E^{1/2}\{
(\z_{t+\tau-j_3}^\top \A_{j_3}^\top \A_{j_4}\bSigma \A_{j_4}^\top \A_{j_3}\z_{t+\tau-j_3})^2\} \nonumber\\
&~~~~~\leq Kp^3 \|\A_{j}^\top \A_{j+k}\|_2^2 \| \A_{j_3}^\top \A_{j_4}\|_2^2 \,.
\end{align}
The last inequality is obtained based on the following facts: (i)
\begin{align}\label{eq:e44}
    \E\{(\z_{t-j}^\top \A_{j}^\top \A_{j+k}\z_{t-j})^4\} & \le K\big(\E[\{\z_{t-j}^\top \A_{j}^\top \A_{j+k}\z_{t-j} - {\tr}(\A_{j}^\top \A_{j+k}\bSigma)\}^4] + {\tr}^4(\A_{j}^\top \A_{j+k}\bSigma)\big) \notag\\
    & \le K(p^2\|\A_{j}^\top \A_{j+k}\|_2^4 + p^4\|\A_{j}^\top \A_{j+k}\|_2^4) \le Kp^4\|\A_{j}^\top \A_{j+k}\|_2^4\,,
\end{align}
where the second step is based on \eqref{eq:e4},
and (ii)
\begin{align}\label{eq:e22}
    & \E\{
(\z_{t+\tau-j_3}^\top \A_{j_3}^\top \A_{j_4}\bSigma \A_{j_4}^\top \A_{j_3}\z_{t+\tau-j_3})^2\} \notag\\
&~~~~~~= {\Var}(\z_{t+\tau-j_3}^\top \A_{j_3}^\top \A_{j_4}\bSigma \A_{j_4}^\top \A_{j_3}\z_{t+\tau-j_3}) + {\tr}^2(\A_{j_3}^\top \A_{j_4}\bSigma \A_{j_4}^\top \A_{j_3} \bSigma) \notag\\
&~~~~~~ \le K\big(p \|\A_{j_3}^\top \A_{j_4}\bSigma \A_{j_4}^\top \A_{j_3}\|_2^2 + p^2 \|\A_{j_3}^\top \A_{j_4}\bSigma \A_{j_4}^\top \A_{j_3}\|_2^2\big) \\
&~~~~~~ \le Kp^2\|\A_{j_3}^\top \A_{j_4}\|_2^4\,, \notag
\end{align}
where the second step is based on \eqref{eq:e2}. On the other hand, if $t-j \neq t+\tau-j_3$, we can also show 
\begin{align*}
    \E\{(\z_{t-j}^\top \A_{j}^\top \A_{j+k}\z_{t-j}
\z_{t+\tau-j_3}^\top \A_{j_3}^\top \A_{j_4}\z_{t+\tau+k-j_4})^2\} \le Kp^3 \|\A_{j}^\top \A_{j+k}\|_2^2 \| \A_{j_3}^\top \A_{j_4}\|_2^2 \,.
\end{align*}
Hence, we obtain
\begin{align*}
\Var(g_{n,\tau,12})
\leq &~\frac {Kp}{n}\sum_{k, k'=1}^{n-\tau-1} 
\sum_{j, j_3, j_4=0 \atop j_4\neq j_3+k}^{\infty}
\sum_{j', j'_3, j'_4=0 \atop j'_4\neq j'_3+k'}^{\infty} 
\| \A_{j}^\top \A_{j+k}\|_2 \| \A_{j_3}^\top \A_{j_4}\|_2
\| \A_{j'}^\top \A_{j'+k'}\|_2 \| \A_{j'_3}^\top \A_{j'_4}\|_2\\
\leq &~\frac {Kp}{n}\bigg\{
\sum_{k=1}^\infty \sum_{j=0}^\infty  
(\| \A_{j}\|_2 \|\A_{j+k}\|_2)\cdot \sum_{j_3=0}^\infty \| \A_{j_3}\|_2\cdot\sum_{j_4=0}^\infty\|\A_{j_4}\|_2\bigg\}^2\\
\leq &~\frac {Kp}{n}\bigg( 
\sum_{j=0}^\infty  j \|\A_{j}\|_2\bigg)^2\bigg(\sum_{j=0}^\infty  \|\A_{j}\|_2\bigg)^4 \leq \frac {Kp}{n}\,,
\end{align*}
where the last inequality is based on Assumption \ref{as:A5}. Following the same arguments, we can also show
\begin{align*}
    \Var(g_{n,\tau,21})
\leq \frac {Kp}{n}\,.
\end{align*}

As we will show in Section \ref{subsec:g22}, 
\begin{align}\label{eq:var-g22}
    \Var(g_{n,\tau,22})
\leq K\Big(1+\frac{n}{p}\Big)\,.
\end{align}
Hence, we have 
\begin{align*}
&\Var(G_{n,\tau,10})\leq K\Big(\frac pn+\frac np\Big)\,.
\end{align*}
Using the similar arguments, we can also show
\begin{align*}
&\Var(G_{n,\tau,11})\leq K\Big(\frac pn+\frac np\Big)\,.
\end{align*}
Hence, we have \eqref{varg}. $\hfill\Box$

\subsection{Proof of \eqref{eq:var-g22}}\label{subsec:g22}
For $g_{n,\tau,22}$, we have
\begin{align}
& \Var(g_{n,\tau,22}) \nonumber\\
& =\frac 4{n^2p^2}\Var\bigg(\sum_{k=1}^{n-\tau-1}\sum_{t=1}^{n-\tau-k}
\sum_{j_1, j_2=0\atop j_2\neq j_1+k}^\infty \sum_{j_3, j_4=0\atop j_4\neq j_3+k}^\infty\z_{t-j_1}^\top \A_{j_1}^\top \A_{j_2}\z_{t+k-j_2}  
\z_{t+\tau-j_3}^\top \A_{j_3}^\top \A_{j_4}\z_{t+\tau+k-j_4}	\bigg)\nonumber\\
&=\frac 4{n^2p^2}\sum_{k, k'=1}^{n-\tau-1}
\sum_{{j_1, j_2, j_3,  j_4=0} \atop {j_2\neq j_1+k \atop j_4\neq j_3+k}}^{\infty}
\sum_{{j'_1, j'_2, j'_3,  j'_4=0} \atop {j'_2\neq j'_1+k' \atop j'_4\neq j'_3+k'}}^{\infty}\sum_{t \in [n-\tau-k] \atop t' \in [n-\tau-k']}\Cov\big(
\z_{t-j_1}^\top \A_{j_1}^\top \A_{j_2}\z_{t+k-j_2} 
\z_{t+\tau-j_3}^\top \A_{j_3}^\top \A_{j_4}\z_{t+\tau+k-j_4}	,\nonumber\\
&~~~~~~~~~~~~~~~~~~~~~~~~~~~~~~~~~~~~~~~~~~~~~\z_{t'-j'_1}^\top \A_{j'_1}^\top \A_{j'_2}\z_{t'+k'-j'_2} 
\z_{t'+\tau-j'_3}^\top \A_{j'_3}^\top \A_{j'_4}\z_{t'+\tau+k'-j'_4}	
\big)\nonumber\\
&=\frac 4{n^2p^2}\sum_{k, k'=1}^{n-\tau-1}
\sum_{{j_1, j_2, j_3,  j_4=0} \atop {j_2\neq j_1+k \atop j_4\neq j_3+k}}^{\infty}
\sum_{{j'_1, j'_2, j'_3,  j'_4=0} \atop {j'_2\neq j'_1+k' \atop j'_4\neq j'_3+k'}}^{\infty}\sum_{t \in [n-\tau-k] \atop t' \in [n-\tau-k']}C(\tau,k,k',t,t',j_1, j_2, j_3,  j_4, j'_1, j'_2, j'_3,  j'_4)\,.\label{var4-4}
\end{align}
Write $\mathcal{A} = (a_1,a_2,a_3,a_4) = (t-j_1,t+k-j_2,t+\tau-j_3,t+\tau+k-j_4)$ and $\mathcal{B}=(b_1,b_2,b_3,b_4)=(t'-j_1',t'+k'-j_2',t'+\tau-j_3', t'+\tau+k'-j_4')$. Let $S_{\mathcal{A}}$ and $S_{\mathcal{B}}$ denote the corresponding sets of distinct indices of $\mathcal{A}$ and $\mathcal{B}$, respectively. The restrictions in the summation in \eqref{var4-4} imply that $a_1 \neq a_2$, $a_3 \neq a_4$, $b_1 \neq b_2$ and $b_3 \neq b_4$. Then $|S_{\mathcal{A}}|\ge 2$ and $|S_{\mathcal{B}}|\ge 2$. If an index $i \in S_\mathcal{A}$ (or $S_\mathcal{B}$), then the index $i$ can occur at most twice within $\mathcal{A}$ (or $\mathcal{B}$). Moreover, the following two facts also hold.
\begin{quote}
{\bf Fact 1.} If there exists one of the 8 values $a_1,\ldots,a_4,b_1,\ldots,b_4$ differing from the other 7 values, then $C(\tau,k,k',t,t', j_1, j_2, j_3,  j_4, j'_1, j'_2, j'_3,  j'_4)=0$.

{\bf Fact 2.} If $S_{\mathcal{A}} \cap S_{\mathcal{B}} = \emptyset$, then $C(\tau,k,k',t,t', j_1, j_2, j_3,  j_4, j'_1, j'_2, j'_3,  j'_4)=0$.
\end{quote}
Define 
\begin{align*}
    &\mathcal{C} = \big\{(k,k',t,t',j_1,j_2,j_3,j_4,j'_1, j'_2, j'_3,  j'_4):  k, k' \in [n-\tau-1], \\
    &~~~~~~~~~~~~~~~~~~~~t \in [n-\tau-k],
    t' \in [n-\tau-k'], ~\mbox{and}~j_1,j_2,j_3,j_4, j'_1,j'_2, j'_3, j'_4 \ge 0 \\
    &~~~~~~~~~~~~~~~~~~~~\mbox{such that}~j_2\neq j_1+k, j_4\neq j_3+k, j'_2\neq j'_1+k', j'_4\neq j'_3+k'\big\}\,.
\end{align*} 
Then 
\begin{align*}
\Var(g_{n,\tau,22}) 
= &~\frac 4{n^2p^2}\sum_{(k,k',t,t',j_1,j_2,j_3,j_4,j'_1, j'_2, j'_3,  j'_4) \in \mathcal{C}_1} 
C(\tau,k,k',t,t', j_1, j_2, j_3,  j_4, j'_1, j'_2, j'_3,  j'_4) \\
&~+ \frac 4{n^2p^2}\sum_{(k,k',t,t',j_1,j_2,j_3,j_4,j'_1, j'_2, j'_3,  j'_4) \in \mathcal{C}_2 } 
C(\tau,k,k',t,t', j_1, j_2, j_3,  j_4, j'_1, j'_2, j'_3,  j'_4)\\
&~+\frac 4{n^2p^2}\sum_{(k,k',t,t',j_1,j_2,j_3,j_4,j'_1, j'_2, j'_3,  j'_4) \in \mathcal{C}_3 } 
C(\tau,k,k',t,t', j_1, j_2, j_3,  j_4, j'_1, j'_2, j'_3,  j'_4)\\
=&~ F_1 + F_2 + F_3\,,
\end{align*}
where 
\begin{align}\label{eq:C123}
    \mathcal{C}_1 & = \big\{(k,k',t,t',j_1,j_2,j_3,j_4,j'_1, j'_2, j'_3,  j'_4)\in \mathcal{C}: |S_{\mathcal{A}}|=4\big\}\,,\notag\\
    \mathcal{C}_2 & = \big\{(k,k',t,t',j_1,j_2,j_3,j_4,j'_1, j'_2, j'_3,  j'_4)\in \mathcal{C}: |S_{\mathcal{A}}|=3\big\}\,,\\
    \mathcal{C}_3 & = \big\{(k,k',t,t',j_1,j_2,j_3,j_4,j'_1, j'_2, j'_3,  j'_4)\in \mathcal{C}: |S_{\mathcal{A}}|=2\big\}\,.\notag
\end{align}
As we will show in Sections \ref{subsec:F1}, \ref{subsec:F2} and \ref{subsec:F3}, 
\begin{align}
    &F_1 \le K\,,\label{eq:F1}\\
    &F_2 \le K\Big(1+\frac{n}{p}\Big)\,,\label{eq:F2}\\
    &F_3 \le K\Big(1+\frac{n}{p}\Big)\,.\label{eq:F3}
\end{align}
Hence, \eqref{eq:var-g22} holds. $\hfill\Box$

\subsubsection{Proof of \eqref{eq:F1}}\label{subsec:F1}
In this case, the 4 values in $\mathcal{A}$ are distinct.
By Fact 1, if $S_{\mathcal{B}} \neq S_{\mathcal{A}}$, then $$C(\tau,k,k',t,t', j_1, j_2, j_3,  j_4, j'_1, j'_2, j'_3,  j'_4)=0\,.$$ Hence, for given $(\tau,k,t,j_1,j_2,j_3,j_4,j'_1, j'_2, j'_3,  j'_4)$, if  $C(\tau,k,k',t,t', j_1, j_2, j_3,  j_4, j'_1, j'_2, j'_3,  j'_4) \neq 0$, then $k'$ and $t'$ must be some functions of $(\tau,k,t,j_1,j_2,j_3,j_4,j'_1, j'_2, j'_3,  j'_4)$, which implies
\begin{align}\label{var:CC1}
&\sum_{(k,k',t,t',j_1,j_2,j_3,j_4,j'_1, j'_2, j'_3,  j'_4) \in \mathcal{C}_1} 
C(\tau,k,k',t,t', j_1, j_2, j_3,  j_4, j'_1, j'_2, j'_3,  j'_4) \\
&~~\le O(n^2)\cdot \sum_{j_1, j_2, j_3,  j_4,j'_1, j'_2, j'_3,  j'_4=0}^{\infty}
\max_{(k,k',t,t')~{\rm such~  that}\atop  (k,k',t,t',j_1, j_2, j_3,  j_4, j'_1, j'_2, j'_3,  j'_4)\in \mathcal{C}_1} |C(\tau,k,k',t,t',j_1, j_2, j_3,  j_4, j'_1, j'_2, j'_3,  j'_4)|\,.\notag
\end{align}
Here, by convention, for any given $(j_1,j_2,j_3,j_4,j'_1,j'_2, j'_3,j'_4)$, we set 
\begin{align*}
    \max_{(k,k',t,t')~{\rm such~  that}\atop  (k,k',t,t',j_1, j_2, j_3,  j_4, j'_1, j'_2, j'_3,  j'_4)\in \mathcal{C}_1} |C(\tau,k,k',t,t',j_1, j_2, j_3,  j_4, j'_1, j'_2, j'_3,  j'_4)|=0
\end{align*}
if there does not exist $(k,k',t,t')$ such that $(k,k',t,t',j_1, j_2, j_3,  j_4, j'_1, j'_2, j'_3,  j'_4)\in \mathcal{C}_1$.
To bound 
\begin{align*}
    \max_{(k,k',t,t')~{\rm such~  that}\atop  (k,k',t,t',j_1, j_2, j_3,  j_4, j'_1, j'_2, j'_3,  j'_4)\in \mathcal{C}_1} |C(\tau,k,k',t,t',j_1, j_2, j_3,  j_4, j'_1, j'_2, j'_3,  j'_4)|\,,
\end{align*}
we assume $\mathcal{A} = (1,2,3,4)$ without loss of generality. Then $a_1=1$, $a_2=2$, $a_3=3$ and $a_4=4$. 

\underline{(i) $\{\{b_1,b_2\},\{b_3,b_4\}\}=\{\{1,2\},\{3,4\}\}$.} We first consider $(b_1,b_2,b_3,b_4)=(1,2,3,4)$. It holds that 
\begin{align*}
&|\Cov(
\z_{1}^\top \A_{j_1}^\top \A_{j_2}\z_{2} 
\z_{3}^\top \A_{j_3}^\top \A_{j_4}\z_{4}	, ~\z_{1}^\top \A_{j'_1}^\top \A_{j'_2}\z_{2} 
\z_{3}^\top \A_{j'_3}^\top \A_{j'_4}\z_{4})|\\
&~~~~~~~~~~~~~~~~~~~~~~~~ = 
|\tr(\A_{j_1}^\top \A_{j_2}\bSigma \A_{j'_2}^\top \A_{j'_1}\bSigma)||\tr(\A_{j_3}^\top \A_{j_4}\bSigma \A_{j'_4}^\top \A_{j'_3}\bSigma)|	\\
&~~~~~~~~~~~~~~~~~~~~~~~~ \leq Kp^2
\|\A_{j_1}^\top \A_{j_2}\|_2 \|\A_{j_3}^\top \A_{j_4}\|_2 
\|\A_{j'_1}^\top \A_{j'_2}\|_2 \|\A_{j'_3}^\top \A_{j'_4}\|_2\,.
\end{align*}
This bound also holds for any $(b_1,b_2,b_3,b_4)$ satisfying $\{\{b_1,b_2\},\{b_3,b_4\}\}=\{\{1,2\},\{3,4\}\}$. 

  \underline{(ii) $\{\{b_1,b_2\},\{b_3,b_4\}\} \neq \{\{1,2\},\{3,4\}\}$.} We first consider $(b_1,b_2,b_3,b_4)=(1,3,2,4)$. It holds that \begin{align*}
&|\Cov(
\z_{1}^\top \A_{j_1}^\top \A_{j_2}\z_{2} 
\z_{3}^\top \A_{j_3}^\top \A_{j_4}\z_{4}, ~\z_{1}^\top \A_{j'_1}^\top \A_{j'_2}\z_{3} 
\z_{2}^\top \A_{j'_3}^\top \A_{j'_4}\z_{4})|\\
&~~~~~~~~~~~~~~~~~~~~~~~~ = |\tr(\A_{j_1}^\top \A_{j_2}\bSigma \A_{j'_3}^\top \A_{j'_4}\bSigma \A_{j_4}^\top \A_{j_3}\bSigma \A_{j'_2}^\top \A_{j'_1}\bSigma)|	\\
&~~~~~~~~~~~~~~~~~~~~~~~~ \leq Kp
\|\A_{j_1}^\top \A_{j_2}\|_2 \|\A_{j_3}^\top \A_{j_4}\|_2 
\|\A_{j'_1}^\top \A_{j'_2}\|_2 \|\A_{j'_3}^\top \A_{j'_4}\|_2\,.	
\end{align*}
This bound also holds for any $(b_1,b_2,b_3,b_4)$ satisfying $\{\{b_1,b_2\},\{b_3,b_4\}\}\neq \{\{1,2\},\{3,4\}\}$. 

Hence, we have
\begin{align}\label{eq:boundC1}
&\max_{(k,k',t,t')~{\rm such~  that}\atop  (k,k',t,t',j_1, j_2, j_3,  j_4, j'_1, j'_2, j'_3,  j'_4)\in \mathcal{C}_1} |C(\tau,k,k',t,t',j_1, j_2, j_3,  j_4, j'_1, j'_2, j'_3,  j'_4)|\notag\\
&~~~~~~~~~~~~~~~~\le Kp^2\|\A_{j_1}^\top \A_{j_2}\|_2 \|\A_{j_3}^\top \A_{j_4}\|_2 
\|\A_{j'_1}^\top \A_{j'_2}\|_2 \|\A_{j'_3}^\top \A_{j'_4}\|_2\,.
\end{align}
Together with \eqref{var:CC1}, we have
\begin{align*}
&\sum_{(k,k',t,t',j_1,j_2,j_3,j_4,j'_1, j'_2, j'_3,  j'_4) \in \mathcal{C}_1} 
C(\tau,k,k',t,t', j_1, j_2, j_3,  j_4, j'_1, j'_2, j'_3,  j'_4) \\
&~~~~~~~~~~~~~~~~\le Kn^2p^2\sum_{j_1,j_2,j_3,j_4,j'_1,j'_2,j'_3,j'_4=0}^{\infty}\|\A_{j_1}^\top \A_{j_2}\|_2 \|\A_{j_3}^\top \A_{j_4}\|_2 
\|\A_{j'_1}^\top \A_{j'_2}\|_2 \|\A_{j'_3}^\top \A_{j'_4}\|_2\\
&~~~~~~~~~~~~~~~~\leq Kn^2p^2\bigg(\sum_{j=0}^{\infty}\|\A_{j}\|_2\bigg)^8 \le Kn^2p^2\,,
\end{align*}
where the last step is based on Assumption \ref{as:A5}. Hence, we have \eqref{eq:F1}. $\hfill\Box$

\subsubsection{Proof of \eqref{eq:F2}}\label{subsec:F2}

Recall $|S_{\mathcal{A}}|=3$ for $(k,k',t,t',j_1,j_2,j_3,j_4,j'_1, j'_2, j'_3,  j'_4) \in \mathcal{C}_2$. Hence, there must be a value in $\mathcal{A}$ occurring twice. Since $a_1 \neq a_2$ and $a_3 \neq a_4$, then either $a_2=a_3$, $a_1 = a_4$, $a_1=a_3$, or $a_2 = a_4$ holds. Write 
\begin{align*}
    \mathcal{C}_{2a} & = \big\{(k,k',t,t',j_1,j_2,j_3,j_4,j'_1, j'_2, j'_3,  j'_4)\in \mathcal{C}_2:  a_2=a_3\big\}\,,\\
    \mathcal{C}_{2b} & = \big\{(k,k',t,t',j_1,j_2,j_3,j_4,j'_1, j'_2, j'_3,  j'_4)\in \mathcal{C}_2:  a_1=a_4\big\}\,,\\
    \mathcal{C}_{2c} & = \big\{(k,k',t,t',j_1,j_2,j_3,j_4,j'_1, j'_2, j'_3,  j'_4)\in \mathcal{C}_2:  a_1=a_3\big\}\,,\\
    \mathcal{C}_{2d} & = \big\{(k,k',t,t',j_1,j_2,j_3,j_4,j'_1, j'_2, j'_3,  j'_4)\in \mathcal{C}_2:  a_2=a_4\big\}\,.
\end{align*}
Then
\begin{align}\label{eq:Case2}
&\sum_{(k,k',t,t',j_1,j_2,j_3,j_4,j'_1, j'_2, j'_3,  j'_4) \in \mathcal{C}_2 } 
C(\tau,k,k',t,t', j_1, j_2, j_3,  j_4, j'_1, j'_2, j'_3,  j'_4)\notag\\
&= \sum_{(k,k',t,t',j_1,j_2,j_3,j_4,j'_1, j'_2, j'_3,  j'_4) \in \mathcal{C}_{2a} } 
C(\tau,k,k',t,t', j_1, j_2, j_3,  j_4, j'_1, j'_2, j'_3,  j'_4)\notag\\
&~~~~~~+\sum_{(k,k',t,t',j_1,j_2,j_3,j_4,j'_1, j'_2, j'_3,  j'_4) \in \mathcal{C}_{2b} } 
C(\tau,k,k',t,t', j_1, j_2, j_3,  j_4, j'_1, j'_2, j'_3,  j'_4)\notag\\
&~~~~~~+\sum_{(k,k',t,t',j_1,j_2,j_3,j_4,j'_1, j'_2, j'_3,  j'_4) \in \mathcal{C}_{2c} } 
C(\tau,k,k',t,t', j_1, j_2, j_3,  j_4, j'_1, j'_2, j'_3,  j'_4)\notag\\
&~~~~~~+\sum_{(k,k',t,t',j_1,j_2,j_3,j_4,j'_1, j'_2, j'_3,  j'_4) \in \mathcal{C}_{2d} } 
C(\tau,k,k',t,t', j_1, j_2, j_3,  j_4, j'_1, j'_2, j'_3,  j'_4)\,.
\end{align}
In the sequel, we bound each term on the right-hand side of the above equation.

\underline{\it Case 1: $a_2=a_3$.} Recall $a_2 = t+k-j_2$ and $a_3=t+\tau-j_3$. Hence, $a_2 = a_3$ implies $k = \tau+j_2-j_3$. Since $a_1, a_2$ and  $a_4$ are distinct, by Fact 1, if $a_1 \notin S_{\mathcal{B}}$ or $a_4 \notin S_{\mathcal{B}}$, then $C(\tau,k,k',t,t', j_1, j_2, j_3,  j_4, j'_1, j'_2, j'_3,  j'_4)=0,$ which implies that $a_1,a_4 \in S_{\mathcal{B}}$ provided that
$C(\tau,k,k',t,t', j_1, j_2, j_3,  j_4, j'_1, j'_2, j'_3,  j'_4) \neq 0$. 
Hence, for given $(\tau,k',t,j_1,j_2,j_3,j_4,j'_1, j'_2, j'_3,  j'_4)$ in this case, if  $C(\tau,k,k',t,t', j_1, j_2, j_3,  j_4, j'_1, j'_2, j'_3,  j'_4) \neq 0$,
$t'$ must be some function of $(\tau,k',t,j_1,j_2,j_3,j_4,j'_1, j'_2, j'_3,  j'_4)$, which implies
\begin{align}\label{var:C2a}
&\sum_{(k,k',t,t',j_1,j_2,j_3,j_4,j'_1, j'_2, j'_3,  j'_4) \in \mathcal{C}_{2a}} 
C(\tau,k,k',t,t', j_1, j_2, j_3,  j_4, j'_1, j'_2, j'_3,  j'_4) \\
&~~\le O(n^2)\cdot \sum_{j_1, j_2, j_3,  j_4,j'_1, j'_2, j'_3,  j'_4=0}^{\infty}
\max_{(k,k',t,t')~{\rm such~  that}\atop  (k,k',t,t',j_1, j_2, j_3,  j_4, j'_1, j'_2, j'_3,  j'_4)\in \mathcal{C}_{2a}} |C(\tau,k,k',t,t',j_1, j_2, j_3,  j_4, j'_1, j'_2, j'_3,  j'_4)|\,.\notag
\end{align}
Here, by convention, for any given $(j_1,j_2,j_3,j_4,j'_1,j'_2, j'_3,j'_4)$, we set 
\begin{align*}
    \max_{(k,k',t,t')~{\rm such~  that}\atop  (k,k',t,t',j_1, j_2, j_3,  j_4, j'_1, j'_2, j'_3,  j'_4)\in \mathcal{C}_{2a}} |C(\tau,k,k',t,t',j_1, j_2, j_3,  j_4, j'_1, j'_2, j'_3,  j'_4)|=0
\end{align*}
if there does not exist $(k,k',t,t')$ such that $(k,k',t,t',j_1, j_2, j_3,  j_4, j'_1, j'_2, j'_3,  j'_4)\in \mathcal{C}_{2a}$. Notice that there are 3 distinct values in $\mathcal{A}$ in this case. If $|S_{\mathcal{B}}|=4$, then the 4 values in $\mathcal{B}$ are distinct, which implies there exists a value in $\mathcal{B}$ that is not in $\mathcal{A}$, and then $C(\tau,k,k',t,t', j_1, j_2, j_3,  j_4, j'_1, j'_2, j'_3,  j'_4)=0$. Hence, if  $C(\tau,k,k',t,t', j_1, j_2, j_3,  j_4, j'_1, j'_2, j'_3,  j'_4) \neq 0$ in this case, then $S_{\mathcal{B}}$ must satisfy $|S_{\mathcal{B}}| \le 3$. As we will show in Section \ref{subsec:Case2abound}, 
\begin{align}\label{eq:2abound}
    &\max_{(k,k',t,t')~{\rm such~  that}\atop  (k,k',t,t',j_1, j_2, j_3,  j_4, j'_1, j'_2, j'_3,  j'_4)\in \mathcal{C}_{2a}} |C(\tau,k,k',t,t',j_1, j_2, j_3,  j_4, j'_1, j'_2, j'_3,  j'_4)| \notag\\
    &~~~~~~~~~~~~~~~~~~~~~\le Kp^2\|\A_{j_1}^\top \A_{j_2}\|_2 \|\A_{j_3}^\top \A_{j_4}\|_2 
\|\A_{j'_1}^\top \A_{j'_2}\|_2 \|\A_{j'_3}^\top \A_{j'_4}\|_2\,.
\end{align}
Together with \eqref{var:C2a}, we have
\begin{align*}
&\sum_{(k,k',t,t',j_1,j_2,j_3,j_4,j'_1, j'_2, j'_3,  j'_4) \in \mathcal{C}_{2a}} 
C(\tau,k,k',t,t', j_1, j_2, j_3,  j_4, j'_1, j'_2, j'_3,  j'_4) \\
&~~~~~~~~~~~~~~~~\le Kn^2p^2\sum_{j_1,j_2,j_3,j_4,j'_1,j'_2,j'_3,j'_4=0}^{\infty}\|\A_{j_1}^\top \A_{j_2}\|_2 \|\A_{j_3}^\top \A_{j_4}\|_2 
\|\A_{j'_1}^\top \A_{j'_2}\|_2 \|\A_{j'_3}^\top \A_{j'_4}\|_2\\
&~~~~~~~~~~~~~~~~\leq Kn^2p^2\bigg(\sum_{j=0}^{\infty}\|\A_{j}\|_2\bigg)^8 \le Kn^2p^2\,,
\end{align*}
where the last step is based on Assumption \ref{as:A5}.

\underline{\it Case 2: $a_1=a_4$.} 
Using the similar arguments for Case 1, we can also show
\begin{align*}
&\sum_{(k,k',t,t',j_1,j_2,j_3,j_4,j'_1, j'_2, j'_3,  j'_4) \in \mathcal{C}_{2b}} 
C(\tau,k,k',t,t', j_1, j_2, j_3,  j_4, j'_1, j'_2, j'_3,  j'_4) \le Kn^2p^2\,.
\end{align*}

\underline{\it Case 3: $a_1=a_3$.}
Since $a_1, a_2$ and  $a_4$ are distinct, by Fact 1, if  $a_2 \notin S_{\mathcal{B}}$ or $a_4 \notin S_{\mathcal{B}}$, then $C(\tau,k,k',t,t', j_1, j_2, j_3,  j_4, j'_1, j'_2, j'_3,  j'_4)=0$. Hence, if $C(\tau,k,k',t,t', j_1, j_2, j_3,  j_4, j'_1, j'_2, j'_3,  j'_4) \neq 0$, we must have $a_2, a_4 \in S_{\mathcal{B}}$.
Write 
\begin{align*}
    \mathcal{C}_{2c}^{\rm (i)} & = \big\{(k,k',t,t',j_1,j_2,j_3,j_4,j'_1, j'_2, j'_3,  j'_4)\in \mathcal{C}_{2c}: S_{\mathcal{B}}=S_{\mathcal{A}} ~\mbox{and}~ a_1 ~\mbox{occurs twice in}~ \mathcal{B} \big\}\,,\\
    \mathcal{C}_{2c}^{\rm (ii)} & = \big\{(k,k',t,t',j_1,j_2,j_3,j_4,j'_1, j'_2, j'_3,  j'_4)\in \mathcal{C}_{2c}: S_{\mathcal{B}}=S_{\mathcal{A}} ~\mbox{and}~ a_2 ~\mbox{(or}~ a_4) ~\mbox{occurs twice in}~ \mathcal{B} \big\}\,,\\
    \mathcal{C}_{2c}^{\rm (iii)} & = \big\{(k,k',t,t',j_1,j_2,j_3,j_4,j'_1, j'_2, j'_3,  j'_4)\in \mathcal{C}_{2c}: S_{\mathcal{B}} \neq S_{\mathcal{A}} ~\mbox{and}~ |S_{\mathcal{B}}|=3\big\}\,,\\
    \mathcal{C}_{2c}^{\rm (iv)} & = \big\{(k,k',t,t',j_1,j_2,j_3,j_4,j'_1, j'_2, j'_3,  j'_4)\in \mathcal{C}_{2c}: S_{\mathcal{B}} \neq S_{\mathcal{A}} ~\mbox{and}~ |S_{\mathcal{B}}|=2\big\}\,.
\end{align*}
It then holds that
\begin{align}\label{eq:Case2ccc}
&\sum_{(k,k',t,t',j_1,j_2,j_3,j_4,j'_1, j'_2, j'_3,  j'_4) \in \mathcal{C}_{2c} } 
C(\tau,k,k',t,t', j_1, j_2, j_3,  j_4, j'_1, j'_2, j'_3,  j'_4)\notag\\
&~~~~~= \sum_{(k,k',t,t',j_1,j_2,j_3,j_4,j'_1, j'_2, j'_3,  j'_4) \in \mathcal{C}_{2c}^{\rm (i)} } 
C(\tau,k,k',t,t', j_1, j_2, j_3,  j_4, j'_1, j'_2, j'_3,  j'_4)\notag\\
&~~~~~~~~~~~+\sum_{(k,k',t,t',j_1,j_2,j_3,j_4,j'_1, j'_2, j'_3,  j'_4) \in \mathcal{C}_{2c}^{\rm (ii)} } 
C(\tau,k,k',t,t', j_1, j_2, j_3,  j_4, j'_1, j'_2, j'_3,  j'_4)\notag\\
&~~~~~~~~~~~+\sum_{(k,k',t,t',j_1,j_2,j_3,j_4,j'_1, j'_2, j'_3,  j'_4) \in \mathcal{C}_{2c}^{\rm (iii)} } 
C(\tau,k,k',t,t', j_1, j_2, j_3,  j_4, j'_1, j'_2, j'_3,  j'_4)\notag\\
&~~~~~~~~~~~+\sum_{(k,k',t,t',j_1,j_2,j_3,j_4,j'_1, j'_2, j'_3,  j'_4) \in \mathcal{C}_{2c}^{\rm (iv)} } 
C(\tau,k,k',t,t', j_1, j_2, j_3,  j_4, j'_1, j'_2, j'_3,  j'_4)\,.
\end{align}
If $C(\tau,k,k',t,t', j_1, j_2, j_3,  j_4, j'_1, j'_2, j'_3,  j'_4) \neq 0$ for $(k,k',t,t',j_1,j_2,j_3,j_4,j'_1, j'_2, j'_3,  j'_4) \in \mathcal{C}_{2c}^{\rm (i)} \cup \mathcal{C}_{2c}^{\rm (ii)} \cup \mathcal{C}_{2c}^{\rm (iv)}$, then $k'$ and $t'$ must be some functions of $(\tau,k,t,j_1,j_2,j_3,j_4,j'_1, j'_2, j'_3,  j'_4)$, which implies 
\begin{align*}
&\sum_{(k,k',t,t',j_1,j_2,j_3,j_4,j'_1, j'_2, j'_3,  j'_4) \in \mathcal{C}_{2c}^{\rm (m)}} 
C(\tau,k,k',t,t', j_1, j_2, j_3,  j_4, j'_1, j'_2, j'_3,  j'_4) \\
&~\le O(n^2)\cdot \sum_{j_1, j_2, j_3,  j_4,j'_1, j'_2, j'_3,  j'_4=0}^{\infty}
\max_{(k,k',t,t')~{\rm such~  that}\atop  (k,k',t,t',j_1, j_2, j_3,  j_4, j'_1, j'_2, j'_3,  j'_4)\in \mathcal{C}_{2c}^{\rm (m)}} |C(\tau,k,k',t,t',j_1, j_2, j_3,  j_4, j'_1, j'_2, j'_3,  j'_4) |\notag
\end{align*}
for ${\rm m} \in \{{\rm i}, {\rm ii},{\rm iv}\}$. Here, by convention, for any given $(j_1,j_2,j_3,j_4,j'_1,j'_2, j'_3,j'_4)$, we set 
\begin{align*}
    \max_{(k,k',t,t')~{\rm such~  that}\atop  (k,k',t,t',j_1, j_2, j_3,  j_4, j'_1, j'_2, j'_3,  j'_4)\in \mathcal{C}_{2c}^{\rm (m)}} |C(\tau,k,k',t,t',j_1, j_2, j_3,  j_4, j'_1, j'_2, j'_3,  j'_4)|=0
\end{align*}
if there does not exist $(k,k',t,t')$ such that $(k,k',t,t',j_1, j_2, j_3,  j_4, j'_1, j'_2, j'_3,  j'_4)\in \mathcal{C}_{2c}^{\rm (m)}$ for ${\rm m} \in \{{\rm i}, {\rm ii},{\rm iv}\}$.
If $C(\tau,k,k',t,t', j_1, j_2, j_3,  j_4, j'_1, j'_2, j'_3,  j'_4) \neq 0$ for $(k,k',t,t',j_1,j_2,j_3,j_4,j'_1, j'_2, j'_3,  j'_4)\in \mathcal{C}_{2c}^{\rm (iii)}$, then $t'$ must be some function of $(\tau,k,k',t,j_1,j_2,j_3,j_4,j'_1, j'_2, j'_3,  j'_4)$, which implies 
\begin{align*}
&\sum_{(k,k',t,t',j_1,j_2,j_3,j_4,j'_1, j'_2, j'_3,  j'_4) \in \mathcal{C}_{2c}^{\rm (iii)}} 
C(\tau,k,k',t,t', j_1, j_2, j_3,  j_4, j'_1, j'_2, j'_3,  j'_4) \\
&~\le O(n^3)\cdot \sum_{j_1, j_2, j_3,  j_4,j'_1, j'_2, j'_3,  j'_4=0}^{\infty}
\max_{(k,k',t,t')~{\rm such~  that}\atop  (k,k',t,t',j_1, j_2, j_3,  j_4, j'_1, j'_2, j'_3,  j'_4)\in \mathcal{C}_{2c}^{\rm (iii)}} |C(\tau,k,k',t,t',j_1, j_2, j_3,  j_4, j'_1, j'_2, j'_3,  j'_4)|\,.\notag
\end{align*}
Here, by convention, for any given $(j_1,j_2,j_3,j_4,j'_1,j'_2, j'_3,j'_4)$, we set 
\begin{align*}
    \max_{(k,k',t,t')~{\rm such~  that}\atop  (k,k',t,t',j_1, j_2, j_3,  j_4, j'_1, j'_2, j'_3,  j'_4)\in \mathcal{C}_{2c}^{\rm (iii)}} |C(\tau,k,k',t,t',j_1, j_2, j_3,  j_4, j'_1, j'_2, j'_3,  j'_4)|=0
\end{align*}
if there does not exist $(k,k',t,t')$ such that $(k,k',t,t',j_1, j_2, j_3,  j_4, j'_1, j'_2, j'_3,  j'_4)\in \mathcal{C}_{2c}^{\rm (iii)}$.
As we will show in Section \ref{subsec:Case2c}, 
\begin{align}
&\max_{(k,k',t,t')~{\rm such~  that}\atop  (k,k',t,t',j_1, j_2, j_3,  j_4, j'_1, j'_2, j'_3,  j'_4)\in \mathcal{C}_{2c}^{\rm (i)}} |C(\tau,k,k',t,t',j_1, j_2, j_3,  j_4, j'_1, j'_2, j'_3,  j'_4)|\notag\\
&~~~~~~~~~~~~~~~~~~~~~~~~~~~~\le Kp^2\|\A_{j_1}^\top \A_{j_2}\|_2 \|\A_{j_3}^\top \A_{j_4}\|_2 
\|\A_{j'_1}^\top \A_{j'_2}\|_2 \|\A_{j'_3}^\top \A_{j'_4}\|_2\,,\label{Case2cbound1}\\
&\max_{(k,k',t,t')~{\rm such~  that}\atop  (k,k',t,t',j_1, j_2, j_3,  j_4, j'_1, j'_2, j'_3,  j'_4)\in \mathcal{C}_{2c}^{\rm (ii)}} |C(\tau,k,k',t,t',j_1, j_2, j_3,  j_4, j'_1, j'_2, j'_3,  j'_4)|\notag\\
&~~~~~~~~~~~~~~~~~~~~~~~~~~~~\le Kp^{3/2}\|\A_{j_1}^\top \A_{j_2}\|_2 \|\A_{j_3}^\top \A_{j_4}\|_2 
\|\A_{j'_1}^\top \A_{j'_2}\|_2 \|\A_{j'_3}^\top \A_{j'_4}\|_2\,,\label{Case2cbound2}\\
&\max_{(k,k',t,t')~{\rm such~  that}\atop  (k,k',t,t',j_1, j_2, j_3,  j_4, j'_1, j'_2, j'_3,  j'_4)\in \mathcal{C}_{2c}^{\rm (iii)}} |C(\tau,k,k',t,t',j_1, j_2, j_3,  j_4, j'_1, j'_2, j'_3,  j'_4)|\notag\\
&~~~~~~~~~~~~~~~~~~~~~~~~~~~~\le Kp\|\A_{j_1}^\top \A_{j_2}\|_2 \|\A_{j_3}^\top \A_{j_4}\|_2 
\|\A_{j'_1}^\top \A_{j'_2}\|_2 \|\A_{j'_3}^\top \A_{j'_4}\|_2\,,\label{Case2cbound3}\\
&\max_{(k,k',t,t')~{\rm such~  that}\atop  (k,k',t,t',j_1, j_2, j_3,  j_4, j'_1, j'_2, j'_3,  j'_4)\in \mathcal{C}_{2c}^{\rm (iv)}} |C(\tau,k,k',t,t',j_1, j_2, j_3,  j_4, j'_1, j'_2, j'_3,  j'_4)|\notag\\
&~~~~~~~~~~~~~~~~~~~~~~~~~~~~\le Kp^{3/2}\|\A_{j_1}^\top \A_{j_2}\|_2 \|\A_{j_3}^\top \A_{j_4}\|_2 
\|\A_{j'_1}^\top \A_{j'_2}\|_2 \|\A_{j'_3}^\top \A_{j'_4}\|_2\,.\label{Case2cbound4}
\end{align}
Hence, by \eqref{eq:Case2ccc}, we have
\begin{align*}
&\sum_{(k,k',t,t',j_1,j_2,j_3,j_4,j'_1, j'_2, j'_3,  j'_4) \in \mathcal{C}_{2c}} 
C(\tau,k,k',t,t', j_1, j_2, j_3,  j_4, j'_1, j'_2, j'_3,  j'_4) \\
&~~~~~~~~~~~~~~~~\le K(n^2p^2 + n^3p)\sum_{j_1,j_2,j_3,j_4,j'_1,j'_2,j'_3,j'_4=0}^{\infty}\|\A_{j_1}^\top \A_{j_2}\|_2 \|\A_{j_3}^\top \A_{j_4}\|_2 
\|\A_{j'_1}^\top \A_{j'_2}\|_2 \|\A_{j'_3}^\top \A_{j'_4}\|_2\\
&~~~~~~~~~~~~~~~~\leq K(n^2p^2 + n^3p)\bigg(\sum_{j=0}^{\infty}\|\A_{j}\|_2\bigg)^8 \le K(n^2p^2 + n^3p)\,,
\end{align*}
where the last step is based on Assumption \ref{as:A5}. 

\underline{\it Case 4: $a_2=a_4$.}  Using arguments similar to those for Case 3, we can also show
\begin{align*}
&\sum_{(k,k',t,t',j_1,j_2,j_3,j_4,j'_1, j'_2, j'_3,  j'_4) \in \mathcal{C}_{2d}} 
C(\tau,k,k',t,t', j_1, j_2, j_3,  j_4, j'_1, j'_2, j'_3,  j'_4) \le K(n^2p^2 + n^3p)\,.
\end{align*}

Combining the results of Cases 1--4 together, by \eqref{eq:Case2}, we have \eqref{eq:F2}. $\hfill\Box$

\subsubsection{Proof of \eqref{eq:2abound}}\label{subsec:Case2abound}

Based on the relationship between $S_{\mathcal{A}}$ and $S_{\mathcal{B}}$, we consider five cases: (i)  $S_{\mathcal{A}} = S_{\mathcal{B}}$ and $a_2$ occurs twice in $\mathcal{B}$, (ii) $S_{\mathcal{A}} = S_{\mathcal{B}}$ and $a_1$ occurs twice in $\mathcal{B}$, (iii) $S_{\mathcal{A}} = S_{\mathcal{B}}$ and $a_4$ occurs twice in $\mathcal{B}$, (iv) $S_{\mathcal{A}} \neq S_{\mathcal{B}}$ and $|S_{\mathcal{B}}|=3$, and (v) $S_{\mathcal{A}} \neq S_{\mathcal{B}}$ and $|S_{\mathcal{B}}|=2$. To simplify our presentation, we assume $\mathcal{A} = (1,2,2,3)$ without loss of generality. Then $a_1=1$, $a_2=2=a_3$ and $a_4=3$.

\underline{(i) $S_{\mathcal{A}} = S_{\mathcal{B}}$ and $a_2$ occurs twice in $\mathcal{B}$.} Since $\mathcal{B}=(b_1,b_2,b_3,b_4)$ satisfies $b_1 \neq b_2$ and $b_3 \neq b_4$, then either $b_2=b_3=a_2$,  $b_1 = b_4=a_2$, $b_1=b_3=a_2$, or $b_2 = b_4=a_2$ holds. To prove \eqref{eq:2abound}, we first consider $(b_1,b_2,b_3,b_4) = (1,2,2,3)$ with $b_2=b_3=a_2=2$ without loss of generality.
By the Cauchy-Schwarz inequality, we have
    \begin{align*}
&|\Cov(
\z_{1}^\top \A_{j_1}^\top \A_{j_2}\z_{2} 
\z_{2}^\top \A_{j_3}^\top \A_{j_4}\z_{3}, ~\z_{1}^\top \A_{j'_1}^\top \A_{j'_2}\z_{2} 
\z_{2}^\top \A_{j'_3}^\top \A_{j'_4}\z_{3}	)|\\
&~~~~~~~~~~~~~~~~~~~~~~~~=|\E(
\z_{2}^\top \A_{j_3}^\top \A_{j_4}\bSigma \A_{j'_4}^\top \A_{j'_3}\z_{2}
\z_{2}^\top \A_{j'_2}^\top \A_{j'_1}\bSigma\A_{j_1}^\top \A_{j_2}\z_{2} )|\\
&~~~~~~~~~~~~~~~~~~~~~~~~\le \E^{1/2}\{(
\z_{2}^\top \A_{j_3}^\top \A_{j_4}\bSigma \A_{j'_4}^\top \A_{j'_3}\z_{2})^2\}
\E^{1/2}\{(\z_{2}^\top \A_{j'_2}^\top \A_{j'_1}\bSigma\A_{j_1}^\top \A_{j_2}\z_{2} 
)^2\}\\
&~~~~~~~~~~~~~~~~~~~~~~~~\leq Kp^2
\|\A_{j_1}^\top \A_{j_2}\|_2 \|\A_{j_3}^\top \A_{j_4}\|_2 
\|\A_{j'_1}^\top \A_{j'_2}\|_2 \|\A_{j'_3}^\top \A_{j'_4}\|_2\,,	
\end{align*}
where the last inequality is obtained following the same arguments of  \eqref{eq:e22}. 
This bound also holds in any of the following cases:  $b_1 = b_4=a_2=2$, $b_1=b_3=a_2=2$, or $b_2 = b_4=a_2=2$. Hence, this bound actually holds for any $(b_1,b_2,b_3,b_4)$ taking permutation of $(1,2,2,3)$ satisfying $b_1 \neq b_2$ and $b_3 \neq b_4$. Then we have 
\begin{align*}
    |C(\tau,k,k',t,t',j_1, j_2, j_3,  j_4, j'_1, j'_2, j'_3,  j'_4)| \le Kp^2
\|\A_{j_1}^\top \A_{j_2}\|_2 \|\A_{j_3}^\top \A_{j_4}\|_2 
\|\A_{j'_1}^\top \A_{j'_2}\|_2 \|\A_{j'_3}^\top \A_{j'_4}\|_2
\end{align*}
in this case.

\underline{(ii) $S_{\mathcal{A}} = S_{\mathcal{B}}$ and $a_1$ occurs twice in $\mathcal{B}$.} Since $\mathcal{B}=(b_1,b_2,b_3,b_4)$ satisfies $b_1 \neq b_2$ and $b_3 \neq b_4$, then either $b_1=b_3=a_1$, $b_1 = b_4=a_1$, $b_2=b_3=a_1$, or $b_2 = b_4=a_1$ holds. To prove \eqref{eq:2abound}, we first consider $(b_1,b_2,b_3,b_4) = (1,2,1,3)$ with $b_1=b_3=a_1=1$ without loss of generality. By the Cauchy-Schwarz inequality, we have
\begin{align}\label{eq:Case2a2}
&|\Cov(
\z_{1}^\top \A_{j_1}^\top \A_{j_2}\z_{2} 
\z_{2}^\top \A_{j_3}^\top \A_{j_4}\z_{3}, ~\z_{1}^\top \A_{j'_1}^\top \A_{j'_2}\z_{2} 
\z_{1}^\top \A_{j'_3}^\top \A_{j'_4}\z_{3})|\notag\\
&~~~~~~~~~~~~~~~~~~~~~~~~=|\E(\z_{1}^\top \A_{j_1}^\top \A_{j_2}\z_{2} \z_{1}^\top\A_{j'_1}^\top \A_{j'_2}\z_{2} \z_{2}^\top \A_{j_3}^\top \A_{j_4}\bSigma \A_{j'_4}^\top \A_{j'_3}\z_{1})|\notag\\
&~~~~~~~~~~~~~~~~~~~~~~~~\le \E^{1/4}\{(\z_{1}^\top \A_{j_1}^\top \A_{j_2}\z_{2})^4\}\E^{1/4}\{(\z_{1}^\top\A_{j'_1}^\top \A_{j'_2}\z_{2})^4\} \notag\\
&~~~~~~~~~~~~~~~~~~~~~~~~~~~~~~~~~\times \E^{1/2}\{(\z_{2}^\top \A_{j_3}^\top \A_{j_4}\bSigma \A_{j'_4}^\top \A_{j'_3}\z_{1})^2\}\notag\\
&~~~~~~~~~~~~~~~~~~~~~~~~\leq Kp^{3/2}
\|\A_{j_1}^\top \A_{j_2}\|_2 \|\A_{j_3}^\top \A_{j_4}\|_2 
\|\A_{j'_1}^\top \A_{j'_2}\|_2 \|\A_{j'_3}^\top \A_{j'_4}\|_2\,.	
\end{align}
The last inequality is obtained based on $\z_t = \bSigma^{1/2}W_t$ as assumed in Assumption \ref{as:A5}, and for any $p \times p$ matrix $\A$,
\begin{align}\label{eq:4E2}
    {\E}\{(\z_{1}^\top \A\z_{2})^4\} & = {\E}[\E\{(W_{1}^\top \bSigma^{1/2} \A\bSigma^{1/2}W_{2})^4\,|\,W_1\}] \notag\\
    & \le K \E(| \bSigma^{1/2}\A^\top\bSigma^{1/2}W_{1}|_2^4)  \le K\|\bSigma^{1/2}\A^\top\bSigma^{1/2}\|_2^4 \E(|W_{1}|_2^4) \notag\\
    & \le Kp^2\|\A\|_2^4
\end{align}
by Assumption \ref{as:A4}.
This bound also holds in any of the following cases: $b_1 = b_4=a_1=1$, $b_2=b_3=a_1=1$, or $b_2 = b_4=a_1=1$. Hence, this bound actually holds for any $(b_1,b_2,b_3,b_4)$ taking permutation of $(1,2,1,3)$ satisfying $b_1 \neq b_2$ and $b_3 \neq b_4$. Then we have
\begin{align*}
    |C(\tau,k,k',t,t',j_1, j_2, j_3,  j_4, j'_1, j'_2, j'_3,  j'_4)| \le Kp^{3/2}
\|\A_{j_1}^\top \A_{j_2}\|_2 \|\A_{j_3}^\top \A_{j_4}\|_2 
\|\A_{j'_1}^\top \A_{j'_2}\|_2 \|\A_{j'_3}^\top \A_{j'_4}\|_2
\end{align*}
in this case.

\underline{(iii) $S_{\mathcal{A}} = S_{\mathcal{B}}$ and $a_4$ occurs twice in $\mathcal{B}$.} 
Using the similar arguments of \eqref{eq:Case2a2}, we can also show 
\begin{align*}
    |C(\tau,k,k',t,t',j_1, j_2, j_3,  j_4, j'_1, j'_2, j'_3,  j'_4)| \le Kp^{3/2}
\|\A_{j_1}^\top \A_{j_2}\|_2 \|\A_{j_3}^\top \A_{j_4}\|_2 
\|\A_{j'_1}^\top \A_{j'_2}\|_2 \|\A_{j'_3}^\top \A_{j'_4}\|_2
\end{align*}
in this case.

\underline{(iv) $S_{\mathcal{A}} \neq S_{\mathcal{B}}$ and $|S_{\mathcal{B}}|=3$.} In this case, due to $|S_{\mathcal{A}}|=3$, $\mathcal{B}$ must contain a value not in $S_{\mathcal{A}}$. Without loss of generality, we first assume $b_2 \notin S_{\mathcal{A}}$. By Fact 1, if $b_2$ just occurs once in $\mathcal{B}$, then
$C(\tau,k,k',t,t', j_1, j_2, j_3,  j_4, j'_1, j'_2, j'_3,  j'_4) = 0$. Hence, if $C(\tau,k,k',t,t', j_1, j_2, j_3,  j_4, j'_1, j'_2, j'_3,  j'_4) \neq 0$, $b_2$ must occur twice in $\mathcal{B}$. Since $\mathcal{B}=(b_1,b_2,b_3,b_4)$ satisfies $b_1 \neq b_2$ and $b_3 \neq b_4$, then either $b_2=b_3$ or $b_2 = b_4$ holds. To prove \eqref{eq:2abound}, we first consider $(b_1,b_2,b_3,b_4)=(1,4,4,3)$ with $b_2=b_3=4 \notin S_{\mathcal{A}}$. It holds that
\begin{align*}
&|\Cov(
\z_{1}^\top \A_{j_1}^\top \A_{j_2}\z_{2} 
\z_{2}^\top \A_{j_3}^\top \A_{j_4}\z_{3}, ~\z_{1}^\top \A_{j'_1}^\top \A_{j'_2}\z_{4} 
\z_{4}^\top \A_{j'_3}^\top \A_{j'_4}\z_{3})|\\
&~~~~~~~~~~~~~~~~~~~~~~~~=|\E(
\z_{2}^\top \A_{j_2}^\top \A_{j_1}\bSigma \A_{j'_1}^\top \A_{j'_2}\z_{4}
\z_{2}^\top \A_{j_3}^\top \A_{j_4}\bSigma\A_{j'_4}^\top \A_{j'_3}\z_{4})|\\
&~~~~~~~~~~~~~~~~~~~~~~~~= |\tr(\A_{j_2}^\top \A_{j_1}\bSigma\A_{j'_1}^\top \A_{j'_2}\bSigma
\A_{j'_3}^\top \A_{j'_4}\bSigma \A_{j_4}^\top \A_{j_3}\bSigma)| \\
&~~~~~~~~~~~~~~~~~~~~~~~~\leq Kp
\|\A_{j_1}^\top \A_{j_2}\|_2 \|\A_{j_3}^\top \A_{j_4}\|_2 
\|\A_{j'_1}^\top \A_{j'_2}\|_2 \|\A_{j'_3}^\top \A_{j'_4}\|_2\,,
\end{align*}
where the last inequality is based on Assumption \ref{as:A4}.
This bound also holds when $b_2 = b_4=4$.
Then we consider the case when $b_1 =4 \notin S_{\mathcal{A}}$ occurs twice in $\mathcal{B}$. Since $\mathcal{B}=(b_1,b_2,b_3,b_4)$ satisfies $b_1 \neq b_2$ and $b_3 \neq b_4$, then either $b_1 = b_4=4$ or $b_1=b_3=4$ holds.
Using the similar arguments, we can show that the above bound also holds when $b_1 = b_4=4$ or $b_1=b_3=4$.
Hence, this bound actually holds for any $(b_1,b_2,b_3,b_4)$ taking permutation of $(1,4,4,3)$ satisfying $b_1 \neq b_2$ and $b_3 \neq b_4$. Then we have 
\begin{align*}
    |C(\tau,k,k',t,t',j_1, j_2, j_3,  j_4, j'_1, j'_2, j'_3,  j'_4)| \le Kp
\|\A_{j_1}^\top \A_{j_2}\|_2 \|\A_{j_3}^\top \A_{j_4}\|_2 
\|\A_{j'_1}^\top \A_{j'_2}\|_2 \|\A_{j'_3}^\top \A_{j'_4}\|_2
\end{align*}
in this case.

\underline{(v)  $S_{\mathcal{A}} \neq S_{\mathcal{B}}$ and $|S_{\mathcal{B}}|=2$.} Recall $a_1 \neq a_4$. As shown in Case 1 of Section \ref{subsec:F2},
$a_1,a_4 \in S_{\mathcal{B}}$ if
$C(\tau,k,k',t,t', j_1, j_2, j_3,  j_4, j'_1, j'_2, j'_3,  j'_4) \neq 0$, 
which implies that $S_{\mathcal{B}}=\{a_1,a_4\}$ in this case. Notice that $\mathcal{B}=(b_1,b_2,b_3,b_4)$ satisfies $b_1\neq b_2$ and $b_3 \neq b_4$. Hence, $a_1$ and $a_4$ both occur twice in $\mathcal{B}$.
To prove \eqref{eq:2abound}, we first consider $(b_1,b_2,b_3,b_4)=(1,3,1,3)$. By the Cauchy-Schwarz inequality, we have
\begin{align}\label{eq:Case2a5}
&|\Cov(
\z_{1}^\top \A_{j_1}^\top \A_{j_2}\z_{2} \z_{2}^\top \A_{j_3}^\top \A_{j_4}\z_{3}, 
~\z_{1}^\top \A_{j'_1}^\top \A_{j'_2}\z_{3} \z_{1}^\top \A_{j'_3}^\top \A_{j'_4}\z_{3})|\notag\\
&~~~~~~~~~~=|\E(
\z_{1}^\top \A_{j_1}^\top \A_{j_2}\bSigma \A_{j_3}^\top \A_{j_4}\z_{3}	 
\z_{1}^\top \A_{j'_1}^\top \A_{j'_2}\z_{3}\z_{1}^\top \A_{j'_3}^\top \A_{j'_4}\z_{3})|\notag\\
&~~~~~~~~~~\leq {\E}^{1/2}\{(
\z_{1}^\top \A_{j_1}^\top \A_{j_2}\bSigma \A_{j_3}^\top \A_{j_4}\z_{3})^2\}	 
{\E}^{1/4}\{(\z_{1}^\top \A_{j'_1}^\top \A_{j'_2}\z_{3})^4\}{\E}^{1/4}\{(
\z_{1}^\top \A_{j'_3}^\top \A_{j'_4}\z_{3})^4\}\notag\\
&~~~~~~~~~~\leq K p^{3/2}
\|\A_{j_1}^\top \A_{j_2}\|_2 \|\A_{j_3}^\top \A_{j_4}\|_2 
\|\A_{j'_1}^\top \A_{j'_2}\|_2 \|\A_{j'_3}^\top \A_{j'_4}\|_2\,,	
\end{align} 
where the last inequality is based on \eqref{eq:4E2}.
This bound also holds for any $(b_1,b_2,b_3,b_4)$ taking permutation of $(1,3,1,3)$ such that $b_1 \neq b_2$ and $b_3 \neq b_4$. Hence,
\begin{align*}
    |C(\tau,k,k',t,t',j_1, j_2, j_3,  j_4, j'_1, j'_2, j'_3,  j'_4)| \le Kp^{3/2}
\|\A_{j_1}^\top \A_{j_2}\|_2 \|\A_{j_3}^\top \A_{j_4}\|_2 
\|\A_{j'_1}^\top \A_{j'_2}\|_2 \|\A_{j'_3}^\top \A_{j'_4}\|_2
\end{align*}
in this case. $\hfill\Box$

\subsubsection{Proofs of \eqref{Case2cbound1}--\eqref{Case2cbound4}}\label{subsec:Case2c}

To simplify our presentation for the proofs of \eqref{Case2cbound1}--\eqref{Case2cbound4}, we assume $\mathcal{A} = (1,2,1,3)$ without loss of generality. Then $a_1=1=a_3$, $a_2=2$ and $a_4=3$.

Recall $S_{\mathcal{B}}=S_{\mathcal{A}}$ and $a_1$ occurs twice in $\mathcal{B}$ for $(k,k',t,t', j_1, j_2, j_3,  j_4, j'_1, j'_2, j'_3,  j'_4)\in \mathcal{C}_{2c}^{\rm (i)}$. Since $\mathcal{B}=(b_1,b_2,b_3,b_4)$ satisfies $b_1 \neq b_2$ and $b_3 \neq b_4$, then either  $b_1=b_3=a_1$, $b_1 = b_4=a_1$, $b_2=b_3=a_1$, or $b_2 = b_4=a_1$ holds. To prove \eqref{Case2cbound1}, we first consider $(b_1,b_2,b_3,b_4)=(1,2,1,3)$ with $b_1=b_3=a_1=1$. By the Cauchy-Schwarz inequality, we have
\begin{align*}
&|\Cov(
\z_{1}^\top \A_{j_1}^\top \A_{j_2}\z_{2} \z_{1}^\top \A_{j_3}^\top \A_{j_4}\z_{3}, 
~\z_{1}^\top \A_{j'_1}^\top \A_{j'_2}\z_{2} \z_{1}^\top \A_{j'_3}^\top \A_{j'_4}\z_{3})|\\
&~~~~~~~~~~=|\E(
\z_{1}^\top \A_{j_1}^\top \A_{j_2}\bSigma \A_{j'_2}^\top \A_{j'_1}\z_{1}\z_{1}^{\top}\A_{j_3}^\top \A_{j_4}
\bSigma \A_{j'_4}^\top \A_{j'_3}\z_{1})|\\
&~~~~~~~~~~\leq {\E}^{1/2}\{(
\z_{1}^\top \A_{j_1}^\top \A_{j_2}\bSigma \A_{j'_2}^\top \A_{j'_1}\z_{1})^2\}	 
{\E}^{1/2}\{(\z_{1}^{\top}\A_{j_3}^\top \A_{j_4}
\bSigma \A_{j'_4}^\top \A_{j'_3}\z_{1})^2\}\\
&~~~~~~~~~~\leq K p^{2}
\|\A_{j_1}^\top \A_{j_2}\|_2 \|\A_{j_3}^\top \A_{j_4}\|_2 
\|\A_{j'_1}^\top \A_{j'_2}\|_2 \|\A_{j'_3}^\top \A_{j'_4}\|_2\,,	
\end{align*} 
where the last inequality is obtained following the same arguments of  \eqref{eq:e22}. This bound also holds in any of the following cases: $b_1 = b_4=a_1=1$, $b_2=b_3=a_1=1$, or $b_2 = b_4=a_1=1$. Hence, 
this bound actually holds for any $(b_1,b_2,b_3,b_4)$ taking permutation of $(1, 2, 1, 3)$ satisfying $b_1 \neq b_2$ and $b_3 \neq b_4$. Then \eqref{Case2cbound1} holds.

Recall $S_{\mathcal{B}}=S_{\mathcal{A}}$ and $a_2$ (or $a_4$) occurs twice in $\mathcal{B}$ for $(k,k',t,t', j_1, j_2, j_3,  j_4, j'_1, j'_2, j'_3,  j'_4)\in \mathcal{C}_{2c}^{\rm (ii)}$. We first consider the case when $a_2$ occurs twice in $\mathcal{B}$. Since $\mathcal{B}=(b_1,b_2,b_3,b_4)$ satisfies $b_1 \neq b_2$ and $b_3 \neq b_4$, then either $b_2=b_3=a_2$, $b_1 = b_4=a_2$, $b_1=b_3=a_2$, or $b_2 = b_4=a_2$ holds. 
To prove \eqref{Case2cbound2}, we first consider $(b_1,b_2,b_3,b_4)= (1,2,2,3)$ with $b_2=b_3=a_2=2$.  
Following the same arguments of \eqref{eq:Case2a2}, we have
\begin{align*}
&|\Cov(
\z_{1}^\top \A_{j_1}^\top \A_{j_2}\z_{2} \z_{1}^\top \A_{j_3}^\top \A_{j_4}\z_{3}, 
~\z_{1}^\top \A_{j'_1}^\top \A_{j'_2}\z_{2} \z_{2}^\top \A_{j'_3}^\top \A_{j'_4}\z_{3})|\\
&~~~~~~~~~~\leq K p^{3/2}
\|\A_{j_1}^\top \A_{j_2}\|_2 \|\A_{j_3}^\top \A_{j_4}\|_2 
\|\A_{j'_1}^\top \A_{j'_2}\|_2 \|\A_{j'_3}^\top \A_{j'_4}\|_2\,.	
\end{align*} 
This bound also holds in any of the following cases: $b_1 = b_4=a_2=2$, $b_1=b_3=a_2=2$, or $b_2 = b_4=a_2=2$. Hence,
this bound actually holds for any $(b_1,b_2,b_3,b_4)$ taking permutation of $(1, 2, 2, 3)$ satisfying $b_1 \neq b_2$ and $b_3 \neq b_4$. Moreover, this bound also holds when $a_4$ occurs twice in $(b_1,b_2,b_3,b_4)$, which satisfies $b_1 \neq b_2$ and $b_3 \neq b_4$. Then  
\eqref{Case2cbound2} holds.

Due to $|S_{\mathcal{A}}|=3$, for $(k,k',t,t', j_1, j_2, j_3,  j_4, j'_1, j'_2, j'_3,  j'_4)\in \mathcal{C}_{2c}^{\rm (iii)}$, 
$\mathcal{B}$ must contain a value not in $S_{\mathcal{A}}$. Without loss of generality, we first assume $b_1 \notin S_{\mathcal{A}}$. Since $|S_{\mathcal{B}}|=3$ in this case, each value in $\mathcal{B}$ occurs at most twice. By Fact 1, if $b_1$ just occurs once in $\mathcal{B}$, then
$C(\tau,k,k',t,t', j_1, j_2, j_3,  j_4, j'_1, j'_2, j'_3,  j'_4) = 0$.  Hence, if $C(\tau,k,k',t,t', j_1, j_2, j_3,  j_4, j'_1, j'_2, j'_3,  j'_4) \neq 0$, then $b_1$ must occur twice in $\mathcal{B}$.  Since $\mathcal{B}=(b_1,b_2,b_3,b_4)$ satisfies $b_1 \neq b_2$ and $b_3 \neq b_4$, then either $b_1=b_3$ or $b_1 = b_4$ holds.
To prove \eqref{Case2cbound3}, we first consider $(b_1,b_2,b_3,b_4)= (4,2,4,3)$ with $b_1=b_3=4 \notin S_{\mathcal{A}}$. Then
\begin{align*}
&|\Cov(
\z_{1}^\top \A_{j_1}^\top \A_{j_2}\z_{2} \z_{1}^\top \A_{j_3}^\top \A_{j_4}\z_{3}, 
~\z_{4}^\top \A_{j'_1}^\top \A_{j'_2}\z_{2} \z_{4}^\top \A_{j'_3}^\top \A_{j'_4}\z_{3})|\\
&~~~~~~~~~~=|\tr(
 \A_{j_1}^\top \A_{j_2}\bSigma
 \A_{j'_2}^\top \A_{j'_1}\bSigma\A_{j'_3}^\top \A_{j'_4}\bSigma\A_{j_4}^\top \A_{j_3}\bSigma)|\\
&~~~~~~~~~~\leq K p
\|\A_{j_1}^\top \A_{j_2}\|_2 \|\A_{j_3}^\top \A_{j_4}\|_2 
\|\A_{j'_1}^\top \A_{j'_2}\|_2 \|\A_{j'_3}^\top \A_{j'_4}\|_2\,.	
\end{align*} 
This bound also holds when $b_1 = b_4=4$.
For the case when $b_2 =4 \notin S_{\mathcal{A}}$ occurs twice in $\mathcal{B}$, since $\mathcal{B}=(b_1,b_2,b_3,b_4)$ satisfies $b_1 \neq b_2$ and $b_3 \neq b_4$, then either $b_2 = b_4=4$ or $b_2=b_3=4$ holds.
Using the similar arguments, we can show that the above bound also holds when $b_2 = b_4=4$ or $b_2=b_3=4$. Hence,
this bound actually holds for any $(b_1,b_2,b_3,b_4)$ taking permutation of $(4, 2, 4, 3)$ satisfying $b_1 \neq b_2$ and $b_3 \neq b_4$. Then \eqref{Case2cbound3} holds.

Recall that, for $(k,k',t,t', j_1, j_2, j_3,  j_4, j'_1, j'_2, j'_3,  j'_4)\in \mathcal{C}_{2c}^{\rm (iv)}$, $a_2 \neq a_4$ and $a_2,a_4 \in S_{\mathcal{B}}$ if
$C(\tau,k,k',t,t', j_1, j_2, j_3,  j_4, j'_1, j'_2, j'_3,  j'_4) \neq 0$. 
If $|S_{\mathcal{B}}|=2$, then $S_{\mathcal{B}}=\{a_2,a_4\}$. 
Notice that $\mathcal{B}=(b_1,b_2,b_3,b_4)$ satisfies $b_1\neq b_2$ and $b_3 \neq b_4$. Hence, if $|S_{\mathcal{B}}|=2$, then $a_2$ and $a_4$ both occur twice in $\mathcal{B}$. To prove \eqref{Case2cbound4}, we first consider $(b_1,b_2,b_3,b_4)=(2,3,2,3)$. Using the similar arguments for \eqref{eq:Case2a5}, we have
\begin{align*}
&|\Cov(
\z_{1}^\top \A_{j_1}^\top \A_{j_2}\z_{2} \z_{1}^\top \A_{j_3}^\top \A_{j_4}\z_{3}, 
~\z_{2}^\top \A_{j'_1}^\top \A_{j'_2}\z_{3} \z_{2}^\top \A_{j'_3}^\top \A_{j'_4}\z_{3})|\notag\\
&~~~~~~~~~~\leq K p^{3/2}
\|\A_{j_1}^\top \A_{j_2}\|_2 \|\A_{j_3}^\top \A_{j_4}\|_2 
\|\A_{j'_1}^\top \A_{j'_2}\|_2 \|\A_{j'_3}^\top \A_{j'_4}\|_2\,.	
\end{align*} 
This bound also holds for any $(b_1,b_2,b_3,b_4)$ taking permutation of $(2, 3, 2, 3)$ satisfying $b_1 \neq b_2$ and $b_3 \neq b_4$. Hence, \eqref{Case2cbound4} holds. $\hfill\Box$

\subsubsection{Proof of \eqref{eq:F3}}\label{subsec:F3}

Recall $|S_{\mathcal{B}}|\ge 2$. Notice that $|S_{\mathcal{A}}|=2$ for $(k,k',t,t', j_1, j_2, j_3,  j_4, j'_1, j'_2, j'_3,  j'_4) \in \mathcal{C}_3$.  If $|S_{\mathcal{B}}|=4$, then the 4 values in $\mathcal{B}$ are distinct, which implies there exists a value in $\mathcal{B}$ that is not in $\mathcal{A}$, and then $C(\tau,k,k',t,t', j_1, j_2, j_3,  j_4, j'_1, j'_2, j'_3,  j'_4)=0$ by Fact 1. Hence, we have
\begin{align}\label{eq:Case3}
&\sum_{(k,k',t,t',j_1,j_2,j_3,j_4,j'_1, j'_2, j'_3,  j'_4) \in \mathcal{C}_{3} } 
C(\tau,k,k',t,t', j_1, j_2, j_3,  j_4, j'_1, j'_2, j'_3,  j'_4)\notag\\
&= \sum_{(k,k',t,t',j_1,j_2,j_3,j_4,j'_1, j'_2, j'_3,  j'_4) \in \mathcal{C}_{3}^{\rm (i)} } 
C(\tau,k,k',t,t', j_1, j_2, j_3,  j_4, j'_1, j'_2, j'_3,  j'_4)\notag\\
&~~~~~~+\sum_{(k,k',t,t',j_1,j_2,j_3,j_4,j'_1, j'_2, j'_3,  j'_4) \in \mathcal{C}_{3}^{\rm (ii)} } 
C(\tau,k,k',t,t', j_1, j_2, j_3,  j_4, j'_1, j'_2, j'_3,  j'_4)\notag\\
&~~~~~~+\sum_{(k,k',t,t',j_1,j_2,j_3,j_4,j'_1, j'_2, j'_3,  j'_4) \in \mathcal{C}_{3}^{\rm (iii)} } 
C(\tau,k,k',t,t', j_1, j_2, j_3,  j_4, j'_1, j'_2, j'_3,  j'_4)\,,
\end{align}
where
\begin{align*}
    \mathcal{C}_{3}^{\rm (i)} & = \big\{(k,k',t,t',j_1,j_2,j_3,j_4,j'_1, j'_2, j'_3,  j'_4)\in \mathcal{C}_{3}: S_{\mathcal{B}}=S_{\mathcal{A}}\big\}\,,\\
    \mathcal{C}_{3}^{\rm (ii)} & = \big\{(k,k',t,t',j_1,j_2,j_3,j_4,j'_1, j'_2, j'_3,  j'_4)\in \mathcal{C}_{3}: S_{\mathcal{B}}\neq S_{\mathcal{A}} ~\mbox{and}~ |S_{\mathcal{B}}|=2 \big\}\,,\\
    \mathcal{C}_{3}^{\rm (iii)} & = \big\{(k,k',t,t',j_1,j_2,j_3,j_4,j'_1, j'_2, j'_3,  j'_4)\in \mathcal{C}_{3}: S_{\mathcal{B}} \neq S_{\mathcal{A}} ~\mbox{and}~ |S_{\mathcal{B}}|=3\big\}\,.
\end{align*}

Notice that $\mathcal{B}=(b_1,b_2,b_3,b_4)$ satisfies $b_1\neq b_2$ and $b_3 \neq b_4$. Then each value in $\mathcal{B}$ occurs at most twice.  Hence, for any $(k,k',t,t',j_1,j_2,j_3,j_4,j'_1, j'_2, j'_3,  j'_4) \in \mathcal{C}_{3}^{\rm (i)}$, the two values in $S_{\mathcal{B}}$ must both occur twice in $\mathcal{B}$, which implies that $k'$ and $t'$ must be some functions of $(\tau,k,t,j_1,j_2,j_3,j_4,j'_1, j'_2, j'_3,  j'_4)$ and then 
\begin{align*}
&\sum_{(k,k',t,t',j_1,j_2,j_3,j_4,j'_1, j'_2, j'_3,  j'_4) \in \mathcal{C}_{3}^{\rm (i)}} 
C(\tau,k,k',t,t', j_1, j_2, j_3,  j_4, j'_1, j'_2, j'_3,  j'_4) \\
&~\le O(n^2)\cdot \sum_{j_1, j_2, j_3,  j_4,j'_1, j'_2, j'_3,  j'_4=0}^{\infty}
\max_{(k,k',t,t')~{\rm such~  that}\atop  (k,k',t,t',j_1, j_2, j_3,  j_4, j'_1, j'_2, j'_3,  j'_4)\in \mathcal{C}_{3}^{\rm (i)}} |C(\tau,k,k',t,t',j_1, j_2, j_3,  j_4, j'_1, j'_2, j'_3,  j'_4)|\,. \notag
\end{align*}
Here, by convention, for any given $(j_1,j_2,j_3,j_4,j'_1,j'_2, j'_3,j'_4)$, we set 
\begin{align*}
    \max_{(k,k',t,t')~{\rm such~  that}\atop  (k,k',t,t',j_1, j_2, j_3,  j_4, j'_1, j'_2, j'_3,  j'_4)\in \mathcal{C}_{3}^{\rm (i)}} |C(\tau,k,k',t,t',j_1, j_2, j_3,  j_4, j'_1, j'_2, j'_3,  j'_4)|=0
\end{align*}
if there does not exist $(k,k',t,t')$ such that $(k,k',t,t',j_1, j_2, j_3,  j_4, j'_1, j'_2, j'_3,  j'_4)\in \mathcal{C}_3^{\rm (i)}$.

Recall $|S_{\mathcal{A}}|=2$. Notice that $|S_{\mathcal{B}}|=2$ for any $(k,k',t,t',j_1,j_2,j_3,j_4,j'_1, j'_2, j'_3,  j'_4) \in \mathcal{C}_{3}^{\rm (ii)}$. By Fact 2, we must have $S_{\mathcal{B}} \cap S_{\mathcal{A}} \neq \emptyset$ for  any $(k,k',t,t',j_1,j_2,j_3,j_4,j'_1, j'_2, j'_3,  j'_4) \in \mathcal{C}_{3}^{\rm (ii)}$ such that  $C(\tau,k,k',t,t', j_1, j_2, j_3,  j_4, j'_1, j'_2, j'_3,  j'_4) \neq 0$, which implies $\mathcal{B}$ must contain a value not in $S_{\mathcal{A}}$.  Since each value in $\mathcal{B}$ occurs at most twice, the two values in $S_{\mathcal{B}}$ must both occur twice in $\mathcal{B}$.
Hence, for any $(k,k',t,t',j_1,j_2,j_3,j_4,j'_1, j'_2, j'_3,  j'_4) \in \mathcal{C}_{3}^{\rm (ii)}$, if $C(\tau,k,k',t,t', j_1, j_2, j_3,  j_4, j'_1, j'_2, j'_3,  j'_4) \neq 0$, then $t'$ must be some function of $(\tau,k,k',t,j_1,j_2,j_3,$ $j_4,j'_1, j'_2, j'_3,  j'_4)$, which implies  
\begin{align*}
&\sum_{(k,k',t,t',j_1,j_2,j_3,j_4,j'_1, j'_2, j'_3,  j'_4) \in \mathcal{C}_{3}^{\rm (ii)}} 
C(\tau,k,k',t,t', j_1, j_2, j_3,  j_4, j'_1, j'_2, j'_3,  j'_4) \\
&~\le O(n^3)\cdot \sum_{j_1, j_2, j_3,  j_4,j'_1, j'_2, j'_3,  j'_4=0}^{\infty}
\max_{(k,k',t,t')~{\rm such~  that}\atop  (k,k',t,t',j_1, j_2, j_3,  j_4, j'_1, j'_2, j'_3,  j'_4)\in \mathcal{C}_{3}^{\rm (ii)}} |C(\tau,k,k',t,t',j_1, j_2, j_3,  j_4, j'_1, j'_2, j'_3,  j'_4)|\,. \notag
\end{align*}
Here, by convention, for any given $(j_1,j_2,j_3,j_4,j'_1,j'_2, j'_3,j'_4)$, we set 
\begin{align*}
    \max_{(k,k',t,t')~{\rm such~  that}\atop  (k,k',t,t',j_1, j_2, j_3,  j_4, j'_1, j'_2, j'_3,  j'_4)\in \mathcal{C}_{3}^{\rm (ii)}} |C(\tau,k,k',t,t',j_1, j_2, j_3,  j_4, j'_1, j'_2, j'_3,  j'_4)|=0
\end{align*}
if there does not exist $(k,k',t,t')$ such that $(k,k',t,t',j_1, j_2, j_3,  j_4, j'_1, j'_2, j'_3,  j'_4)\in \mathcal{C}_3^{\rm (ii)}$.

For any  $(k,k',t,t',j_1,j_2,j_3,j_4,j'_1, j'_2, j'_3,  j'_4) \in \mathcal{C}_{3}^{\rm (iii)}$, due to $|S_{\mathcal{A}}|=2$ and $|S_{\mathcal{B}}|=3$, $S_{\mathcal{B}}$ must contain a value that is not in $S_\mathcal{A}$. Without loss of generality, we assume $b_1 \notin S_{\mathcal{A}}$. Recall each value in $\mathcal{B}$ occurs at most twice. By Fact 1, if $b_1$ just occurs once in $\mathcal{B}$, then
$C(\tau,k,k',t,t', j_1, j_2, j_3,  j_4, j'_1, j'_2, j'_3,  j'_4) = 0$.  Hence, if $C(\tau,k,k',t,t', j_1, j_2, j_3,  j_4, j'_1, j'_2, j'_3,  j'_4) \neq 0$, then $b_1$ must occur twice in $\mathcal{B}$. Notice that $|S_{\mathcal{B}} \cap S_{\mathcal{A}}| \le 2$. 
By Fact 2, if $S_{\mathcal{B}} \cap S_{\mathcal{A}} = \emptyset$, then $C(\tau,k,k',t,t', j_1, j_2, j_3,  j_4, j'_1, j'_2, j'_3,  j'_4) = 0$. If $|S_{\mathcal{B}} \cap S_{\mathcal{A}}| = 1$, there must exist a value in $\mathcal{B}$ which does not belong to $S_{\mathcal{A}}$ and is also different from $b_1$. Without loss of generality, we assume such value is $b_2$. If $b_2$ just occurs once in $\mathcal{B}$, by Fact 1, we know $C(\tau,k,k',t,t', j_1, j_2, j_3,  j_4, j'_1, j'_2, j'_3,  j'_4) = 0$. Hence, if $C(\tau,k,k',t,t', j_1, j_2, j_3,  j_4, j'_1, j'_2, j'_3,  j'_4) \neq 0$, then $b_2$ must also occur twice in $\mathcal{B}$ when $|S_{\mathcal{B}}\cap S_{\mathcal{A}}|=1$. Hence, in order to make $C(\tau,k,k',t,t', j_1, j_2, j_3,  j_4, j'_1, j'_2, j'_3,  j'_4) \neq 0$ under the scenario $|S_{\mathcal{B}}\cap S_{\mathcal{A}}|=1$, we need to require both $b_1$ and $b_2$ occur twice in $\mathcal{B}$, which implies $|S_{\mathcal{B}}|=2$. This is a contradiction of $|S_{\mathcal{B}}|=3$ required in the definition of $\mathcal{C}_{3}^{\rm (iii)}$. Then    $C(\tau,k,k',t,t', j_1, j_2, j_3,  j_4, j'_1, j'_2, j'_3,  j'_4) = 0$ if $|S_{\mathcal{B}} \cap S_{\mathcal{A}}| = 1$. We can claim that $S_{\mathcal{A}}\subset S_{\mathcal{B}}$ for any $(k,k',t,t',j_1,j_2,j_3,j_4,j'_1, j'_2, j'_3,  j'_4) \in \mathcal{C}_{3}^{\rm (iii)}$ such that $C(\tau,k,k',t,t', j_1, j_2, j_3,  j_4, j'_1, j'_2, j'_3,  j'_4) \neq 0$.
Hence, for any $(k,k',t,t',j_1,j_2,j_3,j_4,j'_1, j'_2, j'_3,$ $j'_4) \in \mathcal{C}_{3}^{\rm (iii)}$ such that $C(\tau,k,k',t,t', j_1, j_2, j_3,  j_4, j'_1, j'_2, j'_3,  j'_4) \neq 0$, $k$ and $t$ must be some functions of $(\tau,k',t',j_1,j_2,j_3,j_4,j'_1, j'_2, j'_3,  j'_4)$,  which implies
\begin{align*}
&\sum_{(k,k',t,t',j_1,j_2,j_3,j_4,j'_1, j'_2, j'_3,  j'_4) \in \mathcal{C}_{3}^{\rm (iii)}} 
C(\tau,k,k',t,t', j_1, j_2, j_3,  j_4, j'_1, j'_2, j'_3,  j'_4) \\
&~\le O(n^2)\cdot \sum_{j_1, j_2, j_3,  j_4,j'_1, j'_2, j'_3,  j'_4=0}^{\infty}
\max_{(k,k',t,t')~{\rm such~  that}\atop  (k,k',t,t',j_1, j_2, j_3,  j_4, j'_1, j'_2, j'_3,  j'_4)\in \mathcal{C}_{3}^{\rm (iii)}} |C(\tau,k,k',t,t',j_1, j_2, j_3,  j_4, j'_1, j'_2, j'_3,  j'_4)|\,. \notag
\end{align*}
Here, by convention, for any given $(j_1,j_2,j_3,j_4,j'_1,j'_2, j'_3,j'_4)$, we set 
\begin{align*}
    \max_{(k,k',t,t')~{\rm such~  that}\atop  (k,k',t,t',j_1, j_2, j_3,  j_4, j'_1, j'_2, j'_3,  j'_4)\in \mathcal{C}_{3}^{\rm (iii)}} |C(\tau,k,k',t,t',j_1, j_2, j_3,  j_4, j'_1, j'_2, j'_3,  j'_4)|=0
\end{align*}
if there does not exist $(k,k',t,t')$ such that $(k,k',t,t',j_1, j_2, j_3,  j_4, j'_1, j'_2, j'_3,  j'_4)\in \mathcal{C}_3^{\rm (iii)}$.

As we will show in Section \ref{subsec:Case3},
\begin{align}
&\max_{(k,k',t,t')~{\rm such~  that}\atop  (k,k',t,t',j_1, j_2, j_3,  j_4, j'_1, j'_2, j'_3,  j'_4)\in \mathcal{C}_{3}^{\rm (i)}} |C(\tau,k,k',t,t',j_1, j_2, j_3,  j_4, j'_1, j'_2, j'_3,  j'_4)|\notag\\
&~~~~~~~~~~~~~~~~\le Kp^2\|\A_{j_1}^\top \A_{j_2}\|_2 \|\A_{j_3}^\top \A_{j_4}\|_2 
\|\A_{j'_1}^\top \A_{j'_2}\|_2 \|\A_{j'_3}^\top \A_{j'_4}\|_2\,,\label{eq:Case3bound1}\\
&\max_{(k,k',t,t')~{\rm such~  that}\atop  (k,k',t,t',j_1, j_2, j_3,  j_4, j'_1, j'_2, j'_3,  j'_4)\in \mathcal{C}_{3}^{\rm (ii)}} |C(\tau,k,k',t,t',j_1, j_2, j_3,  j_4, j'_1, j'_2, j'_3,  j'_4)|\notag\\
&~~~~~~~~~~~~~~~~\le Kp\|\A_{j_1}^\top \A_{j_2}\|_2 \|\A_{j_3}^\top \A_{j_4}\|_2 
\|\A_{j'_1}^\top \A_{j'_2}\|_2 \|\A_{j'_3}^\top \A_{j'_4}\|_2\,,\label{eq:Case3bound2}\\
&\max_{(k,k',t,t')~{\rm such~  that}\atop  (k,k',t,t',j_1, j_2, j_3,  j_4, j'_1, j'_2, j'_3,  j'_4)\in \mathcal{C}_{3}^{\rm (iii)}} |C(\tau,k,k',t,t',j_1, j_2, j_3,  j_4, j'_1, j'_2, j'_3,  j'_4)|\notag\\
&~~~~~~~~~~~~~~~~\le Kp^{3/2}\|\A_{j_1}^\top \A_{j_2}\|_2 \|\A_{j_3}^\top \A_{j_4}\|_2 
\|\A_{j'_1}^\top \A_{j'_2}\|_2 \|\A_{j'_3}^\top \A_{j'_4}\|_2\,.\label{eq:Case3bound3}
\end{align}
Hence, by \eqref{eq:Case3}, we have
\begin{align*}
&\sum_{(k,k',t,t',j_1,j_2,j_3,j_4,j'_1, j'_2, j'_3,  j'_4) \in \mathcal{C}_{3}} 
C(\tau,k,k',t,t', j_1, j_2, j_3,  j_4, j'_1, j'_2, j'_3,  j'_4)\notag\\
&~~~~~~~~~~~~~~~~\le K(n^2p^2+n^3p)\sum_{j_1,j_2,j_3,j_4,j'_1,j'_2,j'_3,j'_4=0}^{\infty}\|\A_{j_1}^\top \A_{j_2}\|_2 \|\A_{j_3}^\top \A_{j_4}\|_2 
\|\A_{j'_1}^\top \A_{j'_2}\|_2 \|\A_{j'_3}^\top \A_{j'_4}\|_2\\
&~~~~~~~~~~~~~~~~\leq K(n^2p^2+n^3p)\bigg(\sum_{j=0}^{\infty}\|\A_{j}\|_2\bigg)^8 \le K(n^2p^2+n^3p)\,,
\end{align*}
where the last inequality is based on Assumption \ref{as:A5}.
Hence, \eqref{eq:F3} holds. $\hfill\Box$

\subsubsection{Proofs of \eqref{eq:Case3bound1}--\eqref{eq:Case3bound3}}\label{subsec:Case3}

To simplify our presentation for the proofs of \eqref{eq:Case3bound1}--\eqref{eq:Case3bound3}, we assume $\mathcal{A} = (1,2,1,2)$ without loss of generality. Then $a_1=1=a_3$ and $a_2=2=a_4$.

Recall that $\mathcal{B}=(b_1,b_2,b_3,b_4)$ satisfies $b_1\neq b_2$ and $b_3 \neq b_4$. Then each value in $\mathcal{B}$ occurs at most twice. For any $(k,k',t,t',j_1,j_2,j_3,j_4,j'_1, j'_2, j'_3,  j'_4) \in \mathcal{C}_{3}^{\rm (i)}$,  $S_{\mathcal{B}} = S_{\mathcal{A}} = \{1,2\}$, which implies the two values in $S_{\mathcal{B}}$ must occur twice in $\mathcal{B}$.
To prove \eqref{eq:Case3bound1}, we first consider $(b_1,b_2,b_3,b_4)=(1,2,1,2)$. By the Cauchy-Schwarz inequality, we have
\begin{align*}
&|\Cov(
\z_{1}^\top \A_{j_1}^\top \A_{j_2}\z_{2} \z_{1}^\top \A_{j_3}^\top \A_{j_4}\z_{2}	, 
~\z_{1}^\top \A_{j'_1}^\top \A_{j'_2}\z_{2} \z_{1}^\top \A_{j'_3}^\top \A_{j'_4}\z_{2})|\\
&~~~~~~~~\leq 
\Var^{1/2}(\z_{1}^\top \A_{j_1}^\top \A_{j_2}\z_{2} \z_{1}^\top \A_{j_3}^\top \A_{j_4}\z_{2})	 
\Var^{1/2}(\z_{1}^\top \A_{j'_1}^\top \A_{j'_2}\z_{2} \z_{1}^\top \A_{j'_3}^\top \A_{j'_4}\z_{2}	
) \\
&~~~~~~~~\leq \E^{1/4}\{(\z_{1}^\top \A_{j_1}^\top \A_{j_2}\z_{2})^4\}\E^{1/4}\{(\z_{1}^\top \A_{j_3}^\top \A_{j_4}\z_{2})^4\}\E^{1/4}\{(\z_{1}^\top \A_{j'_1}^\top \A_{j'_2}\z_{2})^4\} \\
&~~~~~~~~~~~~~~~\times \E^{1/4}\{(\z_{1}^\top \A_{j'_3}^\top \A_{j'_4}\z_{2})^4\}\\
&~~~~~~~~ \leq Kp^{2}
\|\A_{j_1}^\top \A_{j_2}\|_2 \|\A_{j_3}^\top \A_{j_4}\|_2 
\|\A_{j'_1}^\top \A_{j'_2}\|_2 \|\A_{j'_3}^\top \A_{j'_4}\|_2\,,	
\end{align*}
where the last inequality is based on \eqref{eq:4E2}. This bound also holds for any $(b_1,b_2,b_3,b_4)$ taking permutation of $(1, 2, 1, 2)$ such that $b_1 \neq b_2$ and $b_3 \neq b_4$. Hence, \eqref{eq:Case3bound1} holds.

Recall that, for any $(k,k',t,t',j_1,j_2,j_3,j_4,j'_1, j'_2, j'_3,  j'_4) \in \mathcal{C}_{3}^{\rm (ii)}$, if $C(\tau,k,k',t,t', j_1, j_2, j_3,  j_4,$ $j'_1, j'_2, j'_3,  j'_4) \neq 0$, the two values in $S_{\mathcal{B}}$ must occur twice in $\mathcal{B}$ and $S_{\mathcal{B}} \cap S_{\mathcal{A}} \neq \emptyset$. Recall $S_{\mathcal{A}}=\{1,2\}$, $S_{\mathcal{B}}\neq S_{\mathcal{A}}$ and $|S_{\mathcal{B}}|=2$. Therefore, if $(k,k',t,t',j_1,j_2,j_3,j_4,j'_1, j'_2, j'_3,  j'_4) \in \mathcal{C}_{3}^{\rm (ii)}$ such that $C(\tau,k,k',t,t', j_1, j_2, j_3,  j_4,$ $j'_1, j'_2, j'_3,  j'_4) \neq 0$, we have $S_{\mathcal{B}}\cap S_{\mathcal{A}}=\{1\}$ or $\{2\}$. 
 When  $S_{\mathcal{B}}\cap S_{\mathcal{A}}=\{1\}$, we assume $S_{\mathcal{B}}=\{1,3\}$ without loss of generality. To prove \eqref{eq:Case3bound2},  we first consider $(b_1,b_2,b_3,b_4)=(1,3,1,3)$. By the Cauchy-Schwarz inequality, we have
\begin{align}\label{eq:1313}
        &|\Cov(\z_{1}^\top \A_{j_1}^\top \A_{j_2}\z_{2} \z_{1}^\top \A_{j_3}^\top \A_{j_4}\z_{2}, 
        ~\z_{1}^\top \A_{j'_1}^\top \A_{j'_2}\z_{3} \z_{1}^\top \A_{j'_3}^\top \A_{j'_4}\z_{3})|\notag\\
        &~~~~~~~~=|\Cov(\z_{1}^\top \A_{j_1}^\top \A_{j_2}\bSigma \A_{j_4}^\top \A_{j_3}\z_{1}	, 
 ~\z_{1}^\top \A_{j'_1}^\top \A_{j'_2}\bSigma \A_{j'_4}^\top \A_{j'_3}\z_{1})|\notag\\
        &~~~~~~~~\leq  \Var^{1/2}(\z_{1}^\top \A_{j_1}^\top \A_{j_2}\bSigma \A_{j_4}^\top \A_{j_3}\z_{1})  \Var^{1/2}(\z_{1}^\top \A_{j'_1}^\top \A_{j'_2}\bSigma \A_{j'_4}^\top \A_{j'_3}\z_{1})\notag\\
        &~~~~~~~~ \leq Kp
        \|\A_{j_1}^\top \A_{j_2}\|_2 \|\A_{j_3}^\top \A_{j_4}\|_2 
        \|\A_{j'_1}^\top \A_{j'_2}\|_2 \|\A_{j'_3}^\top \A_{j'_4}\|_2\,,	
    \end{align}
where the last inequality is based on \eqref{eq:e2}.
This bound also holds for any $(b_1,b_2,b_3,b_4)$ taking permutation of $(1, 3, 1, 3)$ such that $b_1 \neq b_2$ and $b_3 \neq b_4$. When $S_{\mathcal{B}}\cap S_{\mathcal{A}}=\{2\}$, we assume $S_{\mathcal{B}}=\{2,3\}$ without loss of generality. Using the similar arguments, we can show that the upper bound in \eqref{eq:1313} also holds when $(b_1,b_2,b_3,b_4)$ taking permutation of $(2, 3, 2, 3)$ such that $b_1 \neq b_2$ and $b_3 \neq b_4$. Hence, \eqref{eq:Case3bound2} holds.

Recall that, for any $(k,k',t,t',j_1,j_2,j_3,j_4,j'_1, j'_2, j'_3,  j'_4) \in \mathcal{C}_{3}^{\rm (iii)}$, if $C(\tau,k,k',t,t', j_1, j_2, j_3,  j_4,$ $j'_1, j'_2, j'_3,  j'_4) \neq 0$, we must have $S_{\mathcal{A}} \subset S_{\mathcal{B}}$ and $\mathcal{B}$ contains a value that is not in $S_{\mathcal{A}}$ and occurs twice in $\mathcal{B}$. Notice that $S_{\mathcal{A}}=\{1,2\}$. We assume $S_{\mathcal{B}}= \{1,2,3\}$ without loss of generality. Based on the above discussion, we know the value 3 must occur twice in $\mathcal{B}$.
To prove \eqref{eq:Case3bound3}, we first consider $(b_1,b_2,b_3,b_4)=(1,3,3,2)$. By the Cauchy-Schwarz inequality, we have
    \begin{align*}
&|\Cov(
\z_{1}^\top \A_{j_1}^\top \A_{j_2}\z_{2} 
\z_{1}^\top \A_{j_3}^\top \A_{j_4}\z_{2}, ~\z_{1}^\top \A_{j'_1}^\top \A_{j'_2}\z_{3} 
\z_{3}^\top \A_{j'_3}^\top \A_{j'_4}\z_{2}	
)|\\
&~~~~~~~~~~~~~~~~~~~~~~~~=|\E(\z_{1}^\top \A_{j_1}^\top \A_{j_2}\z_{2} 
\z_{1}^\top \A_{j_3}^\top \A_{j_4}\z_{2} \z_{1}^\top \A_{j'_1}^\top \A_{j'_2}\bSigma \A_{j'_3}^\top \A_{j'_4}\z_{2})| \\
&~~~~~~~~~~~~~~~~~~~~~~~~\le \E^{1/4}\{(\z_{1}^\top \A_{j_1}^\top \A_{j_2}\z_{2})^4\}\E^{1/4}\{(\z_{1}^\top \A_{j_3}^\top \A_{j_4}\z_{2})^4\} \\
&~~~~~~~~~~~~~~~~~~~~~~~~~~~~~~~~~\times \E^{1/2}\{(\z_{1}^\top \A_{j'_1}^\top \A_{j'_2}\bSigma \A_{j'_3}^\top \A_{j'_4}\z_{2})^2\}\\
&~~~~~~~~~~~~~~~~~~~~~~~~\leq Kp^{3/2}
\|\A_{j_1}^\top \A_{j_2}\|_2 \|\A_{j_3}^\top \A_{j_4}\|_2 
\|\A_{j'_1}^\top \A_{j'_2}\|_2 \|\A_{j'_3}^\top \A_{j'_4}\|_2\,,
\end{align*}
where the last inequality is based on \eqref{eq:4E2}. This bound also holds for any $(b_1,b_2,b_3,b_4)$ taking permutation of $(1, 3, 3, 2)$ such that $b_1 \neq b_2$ and $b_3 \neq b_4$. Hence, \eqref{eq:Case3bound3} holds. $\hfill\Box$

\subsection{Proof of \eqref{ev-hntau}}\label{subsubsec:E3}
Notice that $\E(G_{n,\tau,10}^e\,|\,\mathcal{X}_n) = 0$. Hence, we have
\begin{align}\label{eq:Cond_var}
    \Var(G_{n,\tau,10}^e\,|\,\mathcal{X}_n) = &~ \E\{(G_{n,\tau,10}^e)^2 \,|\,\mathcal{X}_n\} \notag\\
    = &~ \E\bigg\{\bigg( \frac2{np}\sum_{k=1}^{n-\tau-1}\sum_{t=1}^{n-\tau-k}
e_te_{t+k}e_{t+\tau}e_{t+k+\tau}\x_t^\top \x_{t+k}\x_{t+\tau}^\top \x_{t+\tau+k}\bigg)^2 \,\bigg|\, \mathcal{X}_n \bigg\} \notag\\
    = &~ \frac{4}{n^2p^2}\sum_{k=1}^{n-\tau-1}\sum_{t=1}^{n-\tau-k}(\x_t^\top \x_{t+k}\x_{t+\tau}^\top \x_{t+\tau+k})^2  \,.
\end{align}
We decompose $X_t^\top X_{t+k}$ ($t+k\leq n$) as
\begin{align*}
\x_t^\top \x_{t+k}
=&~
\sum_{j=0}^\infty \z_{t-j}^\top  \A_j^\top \A_{j+k}\z_{t-j}
+\sum_{j, \ell=0\atop \ell \neq k+j}^\infty \z_{t-j}^\top  \A_j^\top \A_\ell\z_{t+k-\ell} = a(t,k)+b(t,k)\,,
\end{align*}
where $a(t,k)$ and $b(t,k)$ are defined in \eqref{eq:abtk}. 
By \eqref{eq:Cond_var}, we have 
\begin{align*}
    \Var(G_{n,\tau,10}^e\,|\,\mathcal{X}_n) = &~ \frac{4}{n^2p^2}\sum_{k=1}^{n-\tau-1}\sum_{t=1}^{n-\tau-k} \big[\{a(t,k)+b(t,k)\}\{a(t+\tau,k)+b(t+\tau,k)\} \big]^2 \\
    \leq &~ \frac{K}{n^2p^2} \sum_{k=1}^{n-\tau-1}\sum_{t=1}^{n-\tau-k} \big\{a^2(t,k)a^2(t+\tau,k) + a^2(t,k)b^2(t+\tau,k) \\
    &~~~~~~~~~~~~~~~~~~~~~~~~~+ a^2(t+\tau,k)b^2(t,k) + b^2(t,k)b^2(t+\tau,k)\big\}\,.
\end{align*}
Hence, it holds that
\begin{align}\label{eq:Gn10aa}
\E\{\Var(G_{n,\tau,10}^e \,|\,\mathcal{X}_n)\}
\leq &~
\frac{K}{n^2p^2} \sum_{k=1}^{n-\tau-1}\sum_{t=1}^{n-\tau-k} \big[\E\{a^2(t,k)a^2(t+\tau,k)\} + \E\{a^2(t,k)b^2(t+\tau,k)\} \notag\\
    &~~~~~~~~~~~~~~~~~~~~~~~~+ \E\{a^2(t+\tau,k)b^2(t,k)\} + \E\{b^2(t,k)b^2(t+\tau,k)\}\big]\notag\\
= &~ K \sum_{i,j=1}^2 V_{n,\tau,ij}^{(1)}+ K \sum_{i,j=1}^2 V_{n,\tau,ij}^{(2)}\,,
\end{align}
where 
\begin{align*}
V_{n,\tau,11}^{(1)}&=\frac1{n^2p^2}\sum_{k=1}^{n-\tau-1}\sum_{t=1}^{n-\tau-k}
{\Var}\{a(t,k)a(t+\tau,k)\}\,,\\
V_{n,\tau,12}^{(1)}&=\frac1{n^2p^2}\sum_{k=1}^{n-\tau-1}\sum_{t=1}^{n-\tau-k}
{\Var}\{a(t,k)b(t+\tau,k)\}\,,\\
V_{n,\tau,21}^{(1)}&=\frac1{n^2p^2}\sum_{k=1}^{n-\tau-1}\sum_{t=1}^{n-\tau-k}
{\Var}\{a(t+\tau,k)b(t,k)\}\,,\\
V_{n,\tau,22}^{(1)}&=\frac1{n^2p^2}\sum_{k=1}^{n-\tau-1}\sum_{t=1}^{n-\tau-k}
{\Var}\{b(t,k)b(t+\tau,k)\}\,,\\
V_{n,\tau,11}^{(2)}&=\frac1{n^2p^2}\sum_{k=1}^{n-\tau-1}\sum_{t=1}^{n-\tau-k}
\E^2\{a(t,k)a(t+\tau,k)\}\,,\\
V_{n,\tau,12}^{(2)}&=\frac1{n^2p^2}\sum_{k=1}^{n-\tau-1}\sum_{t=1}^{n-\tau-k}
\E^2\{a(t,k)b(t+\tau,k)\}\,,\\
V_{n,\tau,21}^{(2)}&=\frac1{n^2p^2}\sum_{k=1}^{n-\tau-1}\sum_{t=1}^{n-\tau-k}
\E^2\{a(t+\tau,k)b(t,k)\}\,,\\
V_{n,\tau,22}^{(2)}&=\frac1{n^2p^2}\sum_{k=1}^{n-\tau-1}\sum_{t=1}^{n-\tau-k}
\E^2\{b(t,k)b(t+\tau,k)\}\,.
\end{align*}
As we will show in Sections \ref{subsec:V1} and \ref{subsec:V2},
\begin{align}
    &V_{n,\tau,11}^{(1)}, V_{n,\tau,12}^{(1)}, V_{n,\tau,21}^{(1)}  \le \frac{Kp}{n}\,, \label{eq:V11to13-1}\\
    &~~~~~~~~~~V_{n,\tau,22}^{(1)}  \le K\,, \label{eq:V22-1}\\
    &~~~~~~~~~~V_{n,\tau,11}^{(2)}  \le \frac{Kp^2}{n}\,, \label{eq:V11-2}\\
    &~~~V_{n,\tau,12}^{(2)} =0\,, ~ V_{n,\tau,21}^{(2)}=0\,, \label{eq:V1221-2}\\
    &~~~~~~~~~~V_{n,\tau,22}^{(2)} \le K\,.\label{eq:V22-2}
\end{align}
By \eqref{eq:Gn10aa}, we have 
\begin{align*}
    \E\{\Var(G_{n,\tau,10}^{e}\,|\,\mathcal{X}_n)\}\leq K\Big(\frac {p^2}{n}+ 1\Big)\,.
\end{align*}
Using the similar arguments, we can also show 
\begin{align*}
    \E\{\Var(G_{n,\tau,11}^{e}\,|\,\mathcal{X}_n)\}\leq K\Big(\frac {p^2}{n}+ 1\Big)\,.
\end{align*}
Hence, we have \eqref{ev-hntau}.
$\hfill\Box$

\subsubsection{Proofs of \eqref{eq:V11to13-1} and \eqref{eq:V22-1}}\label{subsec:V1}

Recall $\x_{t} = \sum_{j=0}^{\infty}\A_j\z_{t-j}$ with $\A_0 = I_p$ and $\A_j \in \mathbb{R}^{p \times p}$.
For $V_{n,\tau,11}^{(1)}$, we have
\begin{align*}
V_{n,\tau,11}^{(1)}&=\frac1{n^2p^2}\sum_{k=1}^{n-\tau-1}\sum_{t=1}^{n-\tau-k}
{\rm var}\{a(t,k)a(t+\tau,k)\} \\
&=\frac1{n^2p^2}\sum_{k=1}^{n-\tau-1}\sum_{t=1}^{n-\tau-k}
\sum_{j_1,j_2, j_1', j_2'=0}^{\infty} 
\Cov\big(
\z_{t-j_1}^\top \A_{j_1}^\top \A_{j_1+k}\z_{t-j_1} 
\z_{t+\tau-j_2}^\top\A_{j_2}^\top\A_{j_2+k}\z_{t+\tau-j_2},\nonumber\\
&\qquad~~~~~~~~~~~~~~~~~~~~~
\z_{t-j'_1}^\top \A_{j'_1}^\top \A_{j'_1+k}\z_{t-j'_1} \z_{t+\tau-j'_2}^\top\A_{j'_2}^\top\A_{j'_2+k}\z_{t+\tau-j'_2}
 	\big)\,.
\end{align*}
By \eqref{var4-1}, we have 
\begin{align*}
  &  |\Cov\big(
\z_{t-j_1}^\top \A_{j_1}^\top \A_{j_1+k}\z_{t-j_1} 
\z_{t+\tau-j_2}^\top\A_{j_2}^\top\A_{j_2+k}\z_{t+\tau-j_2},
\z_{t-j'_1}^\top \A_{j'_1}^\top \A_{j'_1+k}\z_{t-j'_1} \z_{t+\tau-j'_2}^\top\A_{j'_2}^\top\A_{j'_2+k}\z_{t+\tau-j'_2}
 	\big)|\\
    &~~~~~~\leq \begin{cases} Kp^3\|\A_{j_1}^\top \A_{j_1+k}\|_2\|\A_{j_2}^\top\A_{j_2+k}\|_2\|\A_{j'_1}^\top \A_{j'_1+k}\|_2\|\A_{j'_2}^\top\A_{j'_2+k}\|_2\,,\\
    ~~~~~~~~~~~~~~~~~~~~~~~~~~~~~~~~~~~~~~~~~~\mbox{if}~ \{t-j_1, t+\tau-j_2\}\cap \{t-j'_1,t+\tau-j'_2\}\neq \emptyset\,,\\
    ~~~0\,,~~~~~~~~~~~~~~~~~~~~~~~~~~~~~~~~~~~\,\mbox{if}~ \{t-j_1, t+\tau-j_2\}\cap \{t-j'_1,t+\tau-j'_2\}= \emptyset\,,
    \end{cases}
\end{align*}
which implies, by Assumption \ref{as:A5}, that
\begin{align*}
V_{n,\tau,11}^{(1)}
\leq &~\frac {Kp}{n}\sum_{k=1}^{n-\tau-1}
\sum_{j_1,j_2, j_1', j_2'=0}^{\infty} 
\|\A_{j_1}^\top \A_{j_1+k}\|_2 \|\A_{j_2}^\top \A_{j_2+k}\|_2
\|\A_{j'_1}^\top \A_{j'_1+k}\|_2  \|\A_{j'_2}^\top\A_{j'_2+k}\|_2	\notag\\
\le &~ \frac {Kp}{n}\sum_{k=1}^{n-\tau-1}\bigg(
\sum_{j=0}^{\infty} 
\|\A_{j}\|_2 \|\A_{j+k}\|_2 \bigg)^4
 \notag\\
 \leq &~ \frac {Kp}{n}\sum_{k=1}^{n-\tau-1}\bigg(
\sum_{j=0}^{\infty} 
\|\A_{j}\|_2^2\bigg)^2 \bigg(
\sum_{j=0}^{\infty} \|\A_{j+k}\|_2^2 \bigg)^2 \notag\\
\leq &~\frac {Kp}{n}
\sum_{k=1}^{\infty}\bigg(\sum_{j=0}^\infty 
\|\A_{j}\|_2^2\bigg)^3\bigg(\sum_{j=0}^\infty \|\A_{j+k}\|_2^2\bigg) \notag\\
\leq &~\frac {Kp}{n}
\sum_{j=0}^\infty j\|\A_{j}\|_2^2 \leq \frac {Kp}{n}\,.
\end{align*}

For $V_{n,\tau,12}^{(1)}$, due to $\E\{a(t,k)b(t+\tau,k)\}=0$ for any $k \in [n-\tau-1]$ and $t \in [n-\tau-k]$, we have
\begin{align*}
    V_{n,\tau,12}^{(1)}&=\frac1{n^2p^2}\sum_{k=1}^{n-\tau-1}\sum_{t=1}^{n-\tau-k}
{\rm var}\{a(t,k)b(t+\tau,k)\} \\
&=\frac1{n^2p^2}\sum_{k=1}^{n-\tau-1}\sum_{t=1}^{n-\tau-k} E\bigg\{\bigg(
\sum_{j=0}^\infty 
\sum_{j_3, j_4=0\atop j_4\neq j_3+k}^\infty 
Z_{t-j}^\top A_{j}^\top A_{j+k}\z_{t-j}
Z_{t+\tau-j_3}^\top A_{j_3}^\top A_{j_4}Z_{t+\tau+k-j_4}	\bigg)^2\bigg\} \\
&=\frac {1}{n^2p^2}\sum_{k=1}^{n-\tau-1}\sum_{t=1}^{n-\tau-k} 
\sum_{j, j_3, j_4=0 \atop j_4\neq j_3+k}^{\infty} 
\sum_{j', j'_3, j'_4=0 \atop j'_4\neq j'_3+k}^{\infty}
E(\z_{t-j}^\top A_{j}^\top A_{j+k}\z_{t-j}
Z_{t+\tau-j_3}^\top A_{j_3}^\top A_{j_4}Z_{t+\tau+k-j_4}	\\
&~~~~~~~~~~~~~~~~~~~~~~~\times
Z_{t-j'}^\top A_{j'}^\top A_{j'+k}\z_{t-j'}
Z_{t+\tau-j'_3}^\top A_{j'_3}^\top A_{j'_4}Z_{t+\tau+k-j'_4}	) \\
&\le \frac K{np^2} \max_{t \in [n-\tau-1]}\sum_{k=1}^{n-\tau-1}
\sum_{j, j_3, j_4=0 \atop j_4\neq j_3+k}^{\infty} 
\sum_{j', j'_3, j'_4=0 \atop j'_4\neq j'_3+k}^{\infty}
E^{1/2}\{(Z_{t-j}^\top A_{j}^\top A_{j+k}\z_{t-j}
Z_{t+\tau-j_3}^\top A_{j_3}^\top A_{j_4}Z_{t+\tau+k-j_4})^2\} \\
&~~~~~~~~~~~~~~~~~~~~~~~\times
E^{1/2}\{(Z_{t-j'}^\top A_{j'}^\top A_{j'+k}Z_{t-j'}
Z_{t+\tau-j'_3}^\top A_{j'_3}^\top A_{j'_4}Z_{t+\tau+k-j'_4})^2\}\,,
\end{align*}
where the last inequality is based on the Cauchy-Schwarz inequality. Since $j_4\neq j_3+k$, we have $t+\tau+k-j_4 \neq t+\tau -j_3$, which implies that at least one of the conditions $t-j \neq t+\tau+k-j_4$ and $t-j \neq t+\tau -j_3$ must hold. By \eqref{three-diff}, the following bound holds if either $t-j \neq t+\tau+k-j_4$ or $t-j \neq t+\tau -j_3$:
\begin{align*}
    E\{(Z_{t-j}^\top A_{j}^\top A_{j+k}Z_{t-j}
Z_{t+\tau-j_3}^\top A_{j_3}^\top A_{j_4}Z_{t+\tau+k-j_4})^2\} &\le Kp^3 \|A_{j}^\top A_{j+k}\|_2^2 \| A_{j_3}^\top A_{j_4}\|_2^2\,.
\end{align*}
Hence, by Assumption \ref{as:A5}, we have
\begin{align*}
    V_{n,\tau,12}^{(1)}
\leq &~\frac {Kp}{n}\sum_{k=1}^{n-\tau-1} 
\sum_{j, j_3, j_4=0 \atop j_4\neq j_3+k}^{\infty}
\sum_{j', j'_3, j'_4=0 \atop j'_4\neq j'_3+k}^{\infty} 
\| A_{j}^\top A_{j+k}\|_2 \| A_{j_3}^\top A_{j_4}\|_2
\| A_{j'}^\top A_{j'+k}\|_2 \| A_{j'_3}^\top A_{j'_4}\|_2\\
\leq &~\frac {Kp}{n}
\sum_{k=1}^\infty \bigg\{\sum_{j=0}^\infty  
(\| A_{j}\|_2 \|A_{j+k}\|_2) \sum_{j_3=0}^\infty \| A_{j_3}\|_2\sum_{j_4=0}^\infty\|A_{j_4}\|_2\bigg\}^2\\
\leq &~ \frac {Kp}{n}
\sum_{k=1}^\infty \bigg(\sum_{j=0}^\infty 
\|A_{j}\|_2^2\bigg)\bigg(\sum_{j=0}^\infty \|A_{j+k}\|_2^2\bigg) \bigg(\sum_{j=0}^\infty  \|A_{j}\|_2\bigg)^4\\
\leq &~\frac {Kp}{n} 
\sum_{j=0}^\infty  j \|A_{j}\|_2^2 \leq \frac {Kp}{n}\,.
\end{align*}
Using the similar arguments, we can also show $V_{n,\tau,21}^{(1)} \le Kpn^{-1}$. Hence, \eqref{eq:V11to13-1} holds.

For $V_{n,\tau,22}^{(1)}$, we have
\begin{align*}
    V_{n,\tau,22}^{(1)}&=\frac1{n^2p^2}\sum_{k=1}^{n-\tau-1}\sum_{t=1}^{n-\tau-k}
\sum_{{j_1, j_2, j_3,  j_4=0} \atop {j_2\neq j_1+k \atop j_4\neq j_3+k}}^{\infty}\sum_{{j'_1, j'_2, j'_3,  j'_4=0} \atop {j'_2\neq j'_1+k \atop j'_4\neq j'_3+k}}^{\infty} \Cov\big(
Z_{t-j_1}^\top A_{j_1}^\top A_{j_2}Z_{t+k-j_2} 
Z_{t+\tau-j_3}^\top A_{j_3}^\top A_{j_4}Z_{t+\tau+k-j_4}	,\nonumber\\
&~~~~~~~~~~~~~~~~~~~~~~~~~~~~~~~~~~~~~~~~~~~~~~~~~~~~~~~~Z_{t-j'_1}^\top A_{j'_1}^\top A_{j'_2}Z_{t+k-j'_2} 
Z_{t+\tau-j'_3}^\top A_{j'_3}^\top A_{j'_4}Z_{t+\tau+k-j'_4}	
\big)\nonumber\\
&=\frac{1}{n^2p^2}\sum_{k=1}^{n-\tau-1}\sum_{t=1}^{n-\tau-k}
\sum_{{j_1, j_2, j_3,  j_4=0} \atop {j_2\neq j_1+k \atop j_4\neq j_3+k}}^{\infty}\sum_{{j'_1, j'_2, j'_3,  j'_4=0} \atop {j'_2\neq j'_1+k \atop j'_4\neq j'_3+k}}^{\infty}C(\tau,k,k,t,t,j_1, j_2, j_3,  j_4, j'_1, j'_2, j'_3,  j'_4)\,,
\end{align*}
where $C(\tau,k,k',t,t',j_1, j_2, j_3,  j_4, j'_1, j'_2, j'_3,  j'_4)$ is defined in \eqref{var4-4} in Section \ref{subsec:g22}. Define 
\begin{align*}
    &\tilde{\mathcal{C}} = \big\{(k,k',t,t',j_1,j_2,j_3,j_4,j'_1, j'_2, j'_3,  j'_4):  k=k' \in [n-\tau-1], \\
    &~~~~~~~~~~~~~~~~~~~~t=t' \in [n-\tau-k] ~\mbox{and}~j_1,j_2,j_3,j_4, j'_1,j'_2, j'_3, j'_4 \ge 0 \\
    &~~~~~~~~~~~~~~~~~~~~\mbox{such that}~j_2\neq j_1+k, j_4\neq j_3+k, j'_2\neq j'_1+k, j'_4\neq j'_3+k\big\}\,.
\end{align*} 
Notice that, for any fixed $(j_1,j_2,j_3,j_4,j'_1, j'_2, j'_3,  j'_4)$, the number of tuples $(k,k',t,t')$ such that $(k,k',t,t',j_1,j_2,j_3,j_4,j'_1, j'_2, j'_3,  j'_4) \in \tilde{\mathcal{C}}$ is $O(n^2)$.
Hence, $V_{n,\tau,22}^{(1)}$ can be reformulated and bounded as follows:
\begin{align}\label{eq:222}
& V_{n,\tau,22}^{(1)} \notag\\
&~~= \frac1{n^2p^2}\sum_{(k,k',t,t',j_1,j_2,j_3,j_4,j'_1, j'_2, j'_3,  j'_4) \in \tilde{\mathcal{C}}} 
C(\tau,k,k',t,t', j_1, j_2, j_3,  j_4, j'_1, j'_2, j'_3,  j'_4) \\
&~~\le \frac{K}{p^2} \sum_{j_1, j_2, j_3,  j_4,j'_1, j'_2, j'_3,  j'_4=0}^{\infty}
\max_{(k,k',t,t')~{\rm such~  that}\atop  (k,k',t,t',j_1, j_2, j_3,  j_4, j'_1, j'_2, j'_3,  j'_4)\in \tilde{\mathcal{C}}} |C(\tau,k,k',t,t',j_1, j_2, j_3,  j_4, j'_1, j'_2, j'_3,  j'_4)|\,.\notag
\end{align}
Here, by convention, for any given $(j_1,j_2,j_3,j_4,j'_1,j'_2, j'_3,j'_4)$, we set 
\begin{align*}
    \max_{(k,k',t,t')~{\rm such~  that}\atop  (k,k',t,t',j_1, j_2, j_3,  j_4, j'_1, j'_2, j'_3,  j'_4)\in \tilde{\mathcal{C}}} |C(\tau,k,k',t,t',j_1, j_2, j_3,  j_4, j'_1, j'_2, j'_3,  j'_4)|=0
\end{align*}
if there does not exist $(k,k',t,t')$ such that $(k,k',t,t',j_1, j_2, j_3,  j_4, j'_1, j'_2, j'_3,  j'_4)\in \tilde{\mathcal{C}}$. 
Recall 
\begin{align*}
    \mathcal{C} & = \big\{(k,k',t,t',j_1,j_2,j_3,j_4,j'_1, j'_2, j'_3,  j'_4):  k, k' \in [n-\tau-1], \\
    &~~~~~~~~~~~~~~~~t \in [n-\tau-k],
    t' \in [n-\tau-k'], ~\mbox{and}~j_1,j_2,j_3,j_4, j'_1,j'_2, j'_3, j'_4 \ge 0 \\
    &~~~~~~~~~~~~~~~~\mbox{such that}~j_2\neq j_1+k, j_4\neq j_3+k, j'_2\neq j'_1+k', j'_4\neq j'_3+k'\big\}\\
    & = \mathcal{C}_1 \cup \mathcal{C}_2 \cup \mathcal{C}_3 \,,
\end{align*} 
where $\mathcal{C}_1$, $\mathcal{C}_2$ and $\mathcal{C}_3$ are defined in \eqref{eq:C123} in Section \ref{subsec:g22}. In the sequel, by convention, for any $\mathcal{G} \subset \mathcal{C}$ and any given $(j_1,j_2,j_3,j_4,j'_1,j'_2, j'_3,j'_4)$, we set 
\begin{align*}
    \max_{(k,k',t,t')~{\rm such~  that}\atop  (k,k',t,t',j_1, j_2, j_3,  j_4, j'_1, j'_2, j'_3,  j'_4)\in \mathcal{G}} |C(\tau,k,k',t,t',j_1, j_2, j_3,  j_4, j'_1, j'_2, j'_3,  j'_4)|=0
\end{align*}
if there does not exist $(k,k',t,t')$ such that $(k,k',t,t',j_1, j_2, j_3,  j_4, j'_1, j'_2, j'_3,  j'_4)\in \mathcal{G}$. 
Notice that $\tilde{\mathcal{C}} \subset \mathcal{C}$. Then
\begin{align}\label{eq:max}
    &\max_{(k,k',t,t')~{\rm such~  that}\atop  (k,k',t,t',j_1, j_2, j_3,  j_4, j'_1, j'_2, j'_3,  j'_4)\in \tilde{\mathcal{C}}} |C(\tau,k,k',t,t',j_1, j_2, j_3,  j_4, j'_1, j'_2, j'_3,  j'_4)|\\
     &~~~~~~~~~~~~~~~~~~~~~\le \max_{(k,k',t,t')~{\rm such~  that}\atop  (k,k',t,t',j_1, j_2, j_3,  j_4, j'_1, j'_2, j'_3,  j'_4)\in \mathcal{C}} |C(\tau,k,k',t,t',j_1, j_2, j_3,  j_4, j'_1, j'_2, j'_3,  j'_4)|\notag\\
     &~~~~~~~~~~~~~~~~~~~~~= \max_{(k,k',t,t')~{\rm such~  that}\atop  (k,k',t,t',j_1, j_2, j_3,  j_4, j'_1, j'_2, j'_3,  j'_4)\in \mathcal{C}_1} |C(\tau,k,k',t,t',j_1, j_2, j_3,  j_4, j'_1, j'_2, j'_3,  j'_4)|\notag\\
     &~~~~~~~~~~~~~~~~~~~~~~~~~~~\vee \max_{(k,k',t,t')~{\rm such~  that}\atop  (k,k',t,t',j_1, j_2, j_3,  j_4, j'_1, j'_2, j'_3,  j'_4)\in \mathcal{C}_{2}} |C(\tau,k,k',t,t',j_1, j_2, j_3,  j_4, j'_1, j'_2, j'_3,  j'_4)|\notag\\
     &~~~~~~~~~~~~~~~~~~~~~~~~~~~\vee \max_{(k,k',t,t')~{\rm such~  that}\atop  (k,k',t,t',j_1, j_2, j_3,  j_4, j'_1, j'_2, j'_3,  j'_4)\in \mathcal{C}_{3}} |C(\tau,k,k',t,t',j_1, j_2, j_3,  j_4, j'_1, j'_2, j'_3,  j'_4)|\,.\notag
\end{align}
As shown in \eqref{eq:boundC1}, we have 
\begin{align*}
    &\max_{(k,k',t,t')~{\rm such~  that}\atop  (k,k',t,t',j_1, j_2, j_3,  j_4, j'_1, j'_2, j'_3,  j'_4)\in \mathcal{C}_1} |C(\tau,k,k',t,t',j_1, j_2, j_3,  j_4, j'_1, j'_2, j'_3,  j'_4)|\notag\\
&~~~~~~~~~~~~~~~~\le Kp^2\|A_{j_1}^\top A_{j_2}\|_2 \|A_{j_3}^\top A_{j_4}\|_2 
\|A_{j'_1}^\top A_{j'_2}\|_2 \|A_{j'_3}^\top A_{j'_4}\|_2\,.
\end{align*}
Notice that $\mathcal{C}_2 = \mathcal{C}_{2a} \cup \mathcal{C}_{2b}\cup \mathcal{C}_{2c}\cup \mathcal{C}_{2d}$, where $\mathcal{C}_{2a}$, $\mathcal{C}_{2b}$, $\mathcal{C}_{2c}$ and $\mathcal{C}_{2d}$ are defined at the beginning of Section \ref{subsec:F2}. Combining the bounds established in \eqref{eq:2abound} and 
\eqref{Case2cbound1}--\eqref{Case2cbound4} in Section \ref{subsec:F2}, we have 
\begin{align*}
    &\max_{(k,k',t,t')~{\rm such~  that}\atop  (k,k',t,t',j_1, j_2, j_3,  j_4, j'_1, j'_2, j'_3,  j'_4)\in \mathcal{C}_2} |C(\tau,k,k',t,t',j_1, j_2, j_3,  j_4, j'_1, j'_2, j'_3,  j'_4)|\notag\\
    &~~~~~~~~~~~~~~~~ = \max_{(k,k',t,t')~{\rm such~  that}\atop  (k,k',t,t',j_1, j_2, j_3,  j_4, j'_1, j'_2, j'_3,  j'_4)\in \mathcal{C}_{2a}} |C(\tau,k,k',t,t',j_1, j_2, j_3,  j_4, j'_1, j'_2, j'_3,  j'_4)|\notag\\
    &~~~~~~~~~~~~~~~~~~~~~~~~~ \vee \max_{(k,k',t,t')~{\rm such~  that}\atop  (k,k',t,t',j_1, j_2, j_3,  j_4, j'_1, j'_2, j'_3,  j'_4)\in \mathcal{C}_{2b}} |C(\tau,k,k',t,t',j_1, j_2, j_3,  j_4, j'_1, j'_2, j'_3,  j'_4)|\notag\\
    &~~~~~~~~~~~~~~~~~~~~~~~~~ \vee \max_{(k,k',t,t')~{\rm such~  that}\atop  (k,k',t,t',j_1, j_2, j_3,  j_4, j'_1, j'_2, j'_3,  j'_4)\in \mathcal{C}_{2c}} |C(\tau,k,k',t,t',j_1, j_2, j_3,  j_4, j'_1, j'_2, j'_3,  j'_4)|\notag\\
    &~~~~~~~~~~~~~~~~~~~~~~~~~ \vee \max_{(k,k',t,t')~{\rm such~  that}\atop  (k,k',t,t',j_1, j_2, j_3,  j_4, j'_1, j'_2, j'_3,  j'_4)\in \mathcal{C}_{2d}} |C(\tau,k,k',t,t',j_1, j_2, j_3,  j_4, j'_1, j'_2, j'_3,  j'_4)|\notag\\
&~~~~~~~~~~~~~~~~\le Kp^2\|A_{j_1}^\top A_{j_2}\|_2 \|A_{j_3}^\top A_{j_4}\|_2 
\|A_{j'_1}^\top A_{j'_2}\|_2 \|A_{j'_3}^\top A_{j'_4}\|_2\,.
\end{align*}
Moreover, $\mathcal{C}_3 = \mathcal{C}_{3}^{\rm (i)} \cup \mathcal{C}_{3}^{\rm (ii)} \cup \mathcal{C}_{3}^{\rm (iii)}$, where $\mathcal{C}_{3}^{\rm (i)}$, $\mathcal{C}_{3}^{\rm (ii)}$ and $\mathcal{C}_{3}^{\rm (iii)}$ are defined in \eqref{eq:Case3} in Section \ref{subsec:F3}.  Combining the bounds established in \eqref{eq:Case3bound1}--\eqref{eq:Case3bound3} in Section \ref{subsec:F3}, we have 
\begin{align*}
    &\max_{(k,k',t,t')~{\rm such~  that}\atop  (k,k',t,t',j_1, j_2, j_3,  j_4, j'_1, j'_2, j'_3,  j'_4)\in \mathcal{C}_3} |C(\tau,k,k',t,t',j_1, j_2, j_3,  j_4, j'_1, j'_2, j'_3,  j'_4)|\notag\\
    &~~~~~~~~~~~~~~~~ = \max_{(k,k',t,t')~{\rm such~  that}\atop  (k,k',t,t',j_1, j_2, j_3,  j_4, j'_1, j'_2, j'_3,  j'_4)\in \mathcal{C}_{3}^{\rm (i)}} |C(\tau,k,k',t,t',j_1, j_2, j_3,  j_4, j'_1, j'_2, j'_3,  j'_4)|\notag\\
    &~~~~~~~~~~~~~~~~~~~~~~~~~ \vee \max_{(k,k',t,t')~{\rm such~  that}\atop  (k,k',t,t',j_1, j_2, j_3,  j_4, j'_1, j'_2, j'_3,  j'_4)\in \mathcal{C}_{3}^{\rm (ii)}} |C(\tau,k,k',t,t',j_1, j_2, j_3,  j_4, j'_1, j'_2, j'_3,  j'_4)|\notag\\
    &~~~~~~~~~~~~~~~~~~~~~~~~~ \vee \max_{(k,k',t,t')~{\rm such~  that}\atop  (k,k',t,t',j_1, j_2, j_3,  j_4, j'_1, j'_2, j'_3,  j'_4)\in \mathcal{C}_{3}^{\rm (iii)}} |C(\tau,k,k',t,t',j_1, j_2, j_3,  j_4, j'_1, j'_2, j'_3,  j'_4)|\notag\\
&~~~~~~~~~~~~~~~~\le Kp^2\|A_{j_1}^\top A_{j_2}\|_2 \|A_{j_3}^\top A_{j_4}\|_2 
\|A_{j'_1}^\top A_{j'_2}\|_2 \|A_{j'_3}^\top A_{j'_4}\|_2\,.
\end{align*}
Together with \eqref{eq:max}, we have
\begin{align*}
    &\max_{(k,k',t,t')~{\rm such~  that}\atop  (k,k',t,t',j_1, j_2, j_3,  j_4, j'_1, j'_2, j'_3,  j'_4)\in \tilde{\mathcal{C}}} |C(\tau,k,k',t,t',j_1, j_2, j_3,  j_4, j'_1, j'_2, j'_3,  j'_4)|\\
    &~~~~~~~~~~~~~~~~~~~~~\le Kp^2\|A_{j_1}^\top A_{j_2}\|_2 \|A_{j_3}^\top A_{j_4}\|_2 
\|A_{j'_1}^\top A_{j'_2}\|_2 \|A_{j'_3}^\top A_{j'_4}\|_2\,.
\end{align*}
Hence, by \eqref{eq:222} and Assumption \ref{as:A5}, we have
\begin{align*}
 V_{n,\tau,22}^{(1)} &\le  K\sum_{j_1,j_2,j_3,j_4,j'_1,j'_2,j'_3,j'_4=0}^{\infty}\|A_{j_1}^\top A_{j_2}\|_2 \|A_{j_3}^\top A_{j_4}\|_2 
\|A_{j'_1}^\top A_{j'_2}\|_2 \|A_{j'_3}^\top A_{j'_4}\|_2\\
&\leq K\bigg(\sum_{j=0}^{\infty}\|A_{j}\|_2\bigg)^8 \le K\,.
\end{align*}
Then \eqref{eq:V22-1} holds. $\hfill\Box$

\subsubsection{Proofs of \eqref{eq:V11-2}--\eqref{eq:V22-2}}\label{subsec:V2}
To bound $V_{n,\tau,11}^{(2)}$, we have
\begin{align*}
|E \{a(t,k)a(t+\tau,k)\}|
=&~ 
\bigg|\sum_{j, \ell=0}^\infty 
E(Z_{t-j}^\top  A_j^\top A_{j+k}Z_{t-j}
Z_{t+\tau-\ell}^\top  A_\ell^\top A_{\ell+k}Z_{t+\tau-\ell})\bigg|\\
\le&~
\sum_{j=0}^\infty \big|\Cov(
Z_{t-j}^\top  A_j^\top A_{j+k}Z_{t-j},~
Z_{t-j}^\top  A_{j+\tau}^\top A_{j+\tau+k}Z_{t-j})\big|\\
&~~~+\sum_{j, \ell=0}^\infty |\tr(  A_j^\top A_{j+k}\Sigma)|
|\tr( A_\ell^\top A_{\ell+k}\Sigma)|\\
\le &~ \sum_{j=0}^\infty{\rm var}^{1/2}(
Z_{t-j}^\top  A_j^\top A_{j+k}Z_{t-j}){\rm var}^{1/2}(
Z_{t-j}^\top  A_{j+\tau}^\top A_{j+\tau+k}Z_{t-j})\\
&~~~+\sum_{j, \ell=0}^\infty |\tr(A_j^\top A_{j+k}\Sigma)|
|\tr(A_\ell^\top A_{\ell+k}\Sigma)|\\
\le &~ Kp\sum_{j=0}^\infty\|A_j\|_2\|A_{j+k}\|_2\|A_{j+\tau}\|_2\|A_{j+\tau+k}\|_2 \\
&~~~+Kp^2\sum_{j, \ell=0}^\infty \|A_j\|_2\|A_{j+k}\|_2\|A_{\ell}\|_2\|A_{\ell+k}\|_2\,,
\end{align*}
where the second inequality is based on the Cauchy-Schwarz inequality, and the last inequality is based on \eqref{eq:e2}, Assumption \ref{as:A4}, and the inequality $|\tr(A\Sigma)| \le p\|A\Sigma\|_2 \le p\|A\|_2\|\Sigma\|_2$. Hence, under Assumption \ref{as:A5}, we have 
\begin{align*}
    V_{n,\tau,11}^{(2)} &\le \frac{Kp^2}{n}\sum_{k=1}^{n-\tau-1}\bigg(\sum_{j, \ell=0}^\infty \|A_j\|_2\|A_{j+k}\|_2\|A_{\ell}\|_2\|A_{\ell+k}\|_2\bigg)^2 \\
    &\le \frac{Kp^2}{n}\sum_{k=1}^{\infty} \bigg(\sum_{j=0}^{\infty}\|A_{j}\|_2 \|A_{j+k}\|_2 \bigg)^4\\
     & \le \frac {Kp^2}{n}\sum_{k=1}^{\infty}\bigg(\sum_{j=0}^{\infty} \|A_{j}\|_2^2\bigg)^2 \bigg(
\sum_{j=0}^{\infty} \|A_{j+k}\|_2^2 \bigg)^2 \notag\\
  & \le \frac {Kp^2}{n}
\sum_{k=1}^{\infty}\bigg(\sum_{j=0}^\infty 
\|A_{j}\|_2^2\bigg)^3\bigg(\sum_{j=0}^\infty \|A_{j+k}\|_2^2\bigg) \notag\\
& \leq \frac {Kp^2}{n}\sum_{j=0}^\infty j\|A_{j}\|_2^2 \le \frac{Kp^2}{n}\,,
\end{align*}
where the third inequality is based on the Cauchy-Schwarz inequality. Then \eqref{eq:V11-2} holds.

To bound $V_{n,\tau,12}^{(2)}$, $V_{n,\tau,21}^{(2)}$ and $V_{n,\tau,22}^{(2)}$, we decompose $b(t,k)$ as 
\begin{align*}
    b(t,k)  =Z_{t}^\top Z_{t+k} + \sum_{j, \ell=0\atop \ell \neq k+j,~ j+\ell> 0}^\infty Z_{t-j}^\top  A_j^\top A_\ell Z_{t+k-\ell} = b_1(t,k)+b_2(t,k)\,.
\end{align*}
Hence, we have
\begin{align*}
    V_{n,\tau,12}^{(2)}&\le \frac{K}{n^2p^2}\sum_{k=1}^{n-\tau-1}\sum_{t=1}^{n-\tau-k}
E^2\{a(t,k)b_1(t+\tau,k)\} + \frac{K}{n^2p^2}\sum_{k=1}^{n-\tau-1}\sum_{t=1}^{n-\tau-k}
E^2\{a(t,k)b_2(t+\tau,k)\}\,,\\
V_{n,\tau,21}^{(2)}&\le \frac{K}{n^2p^2}\sum_{k=1}^{n-\tau-1}\sum_{t=1}^{n-\tau-k}
E^2\{a(t+\tau,k)b_1(t,k)\}+\frac{K}{n^2p^2}\sum_{k=1}^{n-\tau-1}\sum_{t=1}^{n-\tau-k}
E^2\{a(t+\tau,k)b_2(t,k)\}\,,\\
V_{n,\tau,22}^{(2)}&\le \frac{K}{n^2p^2}\sum_{k=1}^{n-\tau-1}\sum_{t=1}^{n-\tau-k}
E^2\{b_1(t,k)b_1(t+\tau,k)\} + \frac{K}{n^2p^2}\sum_{k=1}^{n-\tau-1}\sum_{t=1}^{n-\tau-k}
E^2\{b_1(t,k)b_2(t+\tau,k)\}\\
&~~+\frac{K}{n^2p^2}\sum_{k=1}^{n-\tau-1}\sum_{t=1}^{n-\tau-k}
E^2\{b_1(t+\tau,k)b_2(t,k)\} + \frac{K}{n^2p^2}\sum_{k=1}^{n-\tau-1}\sum_{t=1}^{n-\tau-k}
E^2\{b_2(t,k)b_2(t+\tau,k)\}\,.
\end{align*}
Notice that $E\{a(t,k)b_1(t+\tau,k)\}=0$, $E\{a(t,k)b_2(t+\tau,k)\}=0$, $E\{a(t+\tau,k)b_1(t,k)\}=0$, and $E\{a(t+\tau,k)b_2(t,k)\}=0$. Hence, we have $V_{n,\tau,12}^{(2)}=0$ and $V_{n,\tau,21}^{(2)}=0$. Then \eqref{eq:V1221-2} holds.

For $V_{n,\tau,22}^{(2)}$, we have $E\{b_1(t,k)b_1(t+\tau,k)\}=0$ and $E\{b_1(t+\tau,k)b_2(t,k)\}=0$. In addition, it holds that
\begin{align*}
|E \{b_1(t,k)b_2(t+\tau,k)\}|
=&~|E (Z_{t}^\top Z_{t+k} Z_{t}^\top  A_\tau^\top A_\tau Z_{t+k})+
E (Z_{t}^\top Z_{t+k} Z_{t+k}^\top  A_{\tau-k}^\top A_{\tau+k}Z_{t})|\\
=&~|\tr(A_\tau^\top A_\tau\Sigma^2)+
\tr (A_{\tau-k}^\top A_{\tau+k}\Sigma^2)|\\
\le &~ Kp\|A_\tau\|_2^2 + Kp\|A_{\tau-k}\|_2 \|A_{\tau+k}\|_2
\end{align*}
for $k\le \tau$,  and
\begin{align*}
|E \{b_1(t,k)b_2(t+\tau,k)\}|
=&~|E (Z_{t}^\top Z_{t+k} Z_{t}^\top  A_\tau^\top A_\tau Z_{t+k})|=\tr(A_\tau^\top A_\tau\Sigma^2)\le  Kp\|A_\tau\|_2^2
\end{align*}
for $k > \tau$. Furthermore, by convention, we set $A_{j}=0$ for $j < 0$.
Then 
\begin{align*}
|E \{b_2(t,k)b_2(t+\tau,k)\}|
=&~ 
\bigg|\sum_{j, \ell=0\atop \ell \neq k+j,~ j+\ell> 0}^\infty 
\big\{ E(Z_{t-j}^\top  A_j^\top A_\ell Z_{t+k-\ell}
Z_{t-j}^\top  A_{\tau+j}^\top A_{\tau+\ell}Z_{t+k-\ell})\\
&~~~~~~~~~~~~~~~~~~~
+E(Z_{t-j}^\top  A_j^\top A_\ell Z_{t+k-\ell}
Z_{t+k-\ell}^\top  A_{\ell+\tau-k}^\top A_{\tau+k+j}Z_{t-j}) \big\}\bigg|\\
=&~ 
\bigg|\sum_{j, \ell=0\atop \ell \neq k+j,~ j+\ell> 0}^\infty 
\big\{ \tr( A_j^\top A_\ell\Sigma  A_{\tau+\ell}^\top A_{\tau+j}\Sigma)
+\tr(  A_j^\top A_\ell\Sigma A_{\ell+\tau-k}^\top A_{\tau+k+j}\Sigma ) \big\}\bigg|\\
\le &~ Kp \sum_{j, \ell=0\atop \ell \neq k+j,~ j+\ell> 0}^\infty\|A_j\|_2\|A_{\ell}\|_2\|A_{\tau+\ell}\|_2\|A_{\tau+j}\|_2 \\
&~~~~~+Kp \sum_{j, \ell=0\atop \ell \neq k+j,~ j+\ell> 0}^\infty\|A_j\|_2\|A_{\ell}\|_2\|A_{\ell+\tau-k}\|_2\|A_{\tau+k+j}\|_2\,.
\end{align*}
Hence, we have
\begin{align*}
    V_{n,\tau,22}^{(2)} &\le \frac{K}{n}\sum_{k=1}^{\tau} (\|A_\tau\|_2^2 +
\|A_{\tau-k}\|_2 \|A_{\tau+k}\|_2)^2 + \frac{K}{n}\sum_{k=\tau+1}^{n-\tau-1}\|A_\tau\|_2^4\\
    &~~~~~+\frac{K}{n}\sum_{k=1}^{n-\tau-1}\bigg(\sum_{j, \ell=0\atop \ell \neq k+j,~ j+\ell> 0}^\infty\|A_j\|_2\|A_{\ell}\|_2\|A_{\tau+\ell}\|_2\|A_{\tau+j}\|_2 \\
    &~~~~~~~~~~~~~~~~~~~~~~~~~~~~~+\sum_{j, \ell=0\atop \ell \neq k+j,~ j+\ell> 0}^\infty\|A_j\|_2\|A_{\ell}\|_2\|A_{\ell+\tau-k}\|_2\|A_{\tau+k+j}\|_2\bigg)^2 \\
    & \le K + \frac{K}{n}\sum_{k=1}^{n-\tau-1}\bigg(\sum_{j, \ell=0\atop \ell \neq k+j,~ j+\ell> 0}^\infty\|A_j\|_2\|A_{\ell}\|_2\|A_{\tau+\ell}\|_2\|A_{\tau+j}\|_2\bigg)^2\\
    &~~~~~ + \frac{K}{n}\sum_{k=1}^{n-\tau-1}\bigg(\sum_{j, \ell=0\atop \ell \neq k+j,~ j+\ell> 0}^\infty\|A_j\|_2\|A_{\ell}\|_2\|A_{\ell+\tau-k}\|_2\|A_{\tau+k+j}\|_2\bigg)^2 \\
    & \le K + K \bigg(\sum_{j=0}^\infty\|A_j\|_2\bigg)^4 \le K\,,
\end{align*}
where the last inequality is based on Assumption \ref{as:A5}. Then \eqref{eq:V22-2} holds. $\hfill\Box$

	\end{document}